\documentclass[paper=A4,fontsize=11pt,bibliography=totoc,twoside]{scrbook}
\usepackage{amsmath}
\usepackage{amsthm}
\usepackage{amssymb}
\usepackage{amsfonts}
\usepackage{amsxtra}
\usepackage{color}
\usepackage{graphicx}
\usepackage{ngerman}
\usepackage{caption}
\usepackage[german,intoc,noprefix]{nomencl}

\makenomenclature
\newtheorem{theorem}{Theorem}[chapter]
\newtheorem{lemma}[theorem]{Lemma}

\newtheorem{proposition}[theorem]{Proposition}

\newtheorem{korollar}[theorem]{Korollar}
\newtheorem{definition}[theorem]{Definition}

\newtheorem{behauptung}[theorem]{Behauptung}
\newtheorem{bemerkung}[theorem]{Bemerkung}

\newcommand{\norm}[1]{\lVert #1\rVert}
\newcommand{\Bignorm}[1]{\big\lVert #1\big\rVert}
\newcommand{\setR}{\mathbb{R}}
\newcommand{\setN}{\mathbb{N}}
\newcommand{\embedding}{\hookrightarrow}
\newcommand{\compactembedding}{\hookrightarrow\hookrightarrow}
\newcommand{\bu}{{\bf u}}

\newcommand{\be}{{\bf e}}
\newcommand{\bg}{{\bf g}}
\newcommand{\bv}{{\bf v}}
\newcommand{\bw}{{\bf w}}
\newcommand{\bX}{{\bf X}}
\newcommand{\bY}{{\bf Y}}
\newcommand{\bW}{{\bf W}}
\newcommand{\bA}{{\bf A}}
\newcommand{\bB}{{\bf B}}
\newcommand{\ff}{{\bf f}}
\newcommand{\R}{\mathcal{R}}
\newcommand{\F}{\mathcal{F}}
\newcommand{\M}{\mathcal{M}}
\newcommand{\T}{\mathcal{T}}
\newcommand{\bnu}{\boldsymbol{\nu}}
\newcommand{\bnue}{\boldsymbol{\nu}_{\eta(t)}}
\newcommand{\bphi}{\boldsymbol{\phi}}
\newcommand{\bpsi}{\boldsymbol{\psi}}
\newcommand{\balpha}{\boldsymbol{\alpha}}
\newcommand{\pa}{\partial}
\newcommand{\iot}{\int_{\Omega_{\eta(t)}}}

\newcommand{\iorte}{\int_{\Omega_{\mathcal{R}_\epsilon\eta(t)}}}
\newcommand{\ios}{\int_{\Omega_{\eta(s)}}}

\newcommand{\iortd}{\int_{\Omega_{\mathcal{R}\delta(t)}}}
\newcommand{\iortdn}{\int_{\Omega_{\mathcal{R}\delta_n(t)}}}
\newcommand{\idot}{\int_{\pa\Omega_{\eta(t)}\setminus\Gamma}}
\newcommand{\im}{\int_M}

\newcommand{\onorm}[1]{\norm{#1}_{L^2(\Omega_{\eta(t)})}}
\newcommand{\mnorm}[1]{\norm{#1}_{L^2(M)}}
\newcommand{\oet}{\Omega_{\eta(t)}}

\newcommand{\beweis}{\noindent{\bf Beweis:}\ }

\newcommand{\loc}{\text{loc}}
\renewcommand{\phi}{\varphi}

\DeclareMathOperator{\dv}{div}

\DeclareMathOperator{\tr}{tr_\eta}
\DeclareMathOperator{\trt}{tr_{\tilde\eta}}
\DeclareMathOperator{\trent}{tr_{\tilde\eta_n}}
\DeclareMathOperator{\trace}{tr}
\DeclareMathOperator{\trnormal}{tr^n_\eta}
\DeclareMathOperator{\trnormaln}{tr^n_{\eta_0}}
\DeclareMathOperator{\trnormale}{tr^n_{\R_\epsilon\eta_0}}
\DeclareMathOperator{\trnormaleo}{tr^n_{\R\eta_0}}
\DeclareMathOperator{\trnormald}{tr^n_{\R\delta}}
\DeclareMathOperator{\tren}{tr_{\eta_n}}

\DeclareMathOperator{\trre}{tr_{\mathcal{R}_\epsilon\eta}}

\DeclareMathOperator{\trrd}{tr_{\mathcal{R}\delta}}
\DeclareMathOperator{\trrdt}{tr_{\mathcal{R}\delta(t)}}
\DeclareMathOperator{\trd}{tr_{\delta}}
\DeclareMathOperator{\trrdn}{tr_{\mathcal{R}\delta_n}}

\DeclareMathOperator{\grad}{grad}
\DeclareMathOperator{\id}{id}

\DeclareMathOperator{\spann}{span}
\DeclareMathOperator{\inn}{int}
\DeclareMathOperator{\supp}{supp}
\begin{document}
%%%%%%%%%%%%%%%%%%%%%%%%%%% titlepage %%%%%%%%%%%%%%%%%%%%%%%%%%%
\pagenumbering{roman}

\begin{titlepage}\begin{addmargin}[1.5cm]{0cm}
\vspace*{0.2cm}
\begin{center}
\begin{minipage}{0.5\textwidth}
\centering
\includegraphics[width=0.8\textwidth]{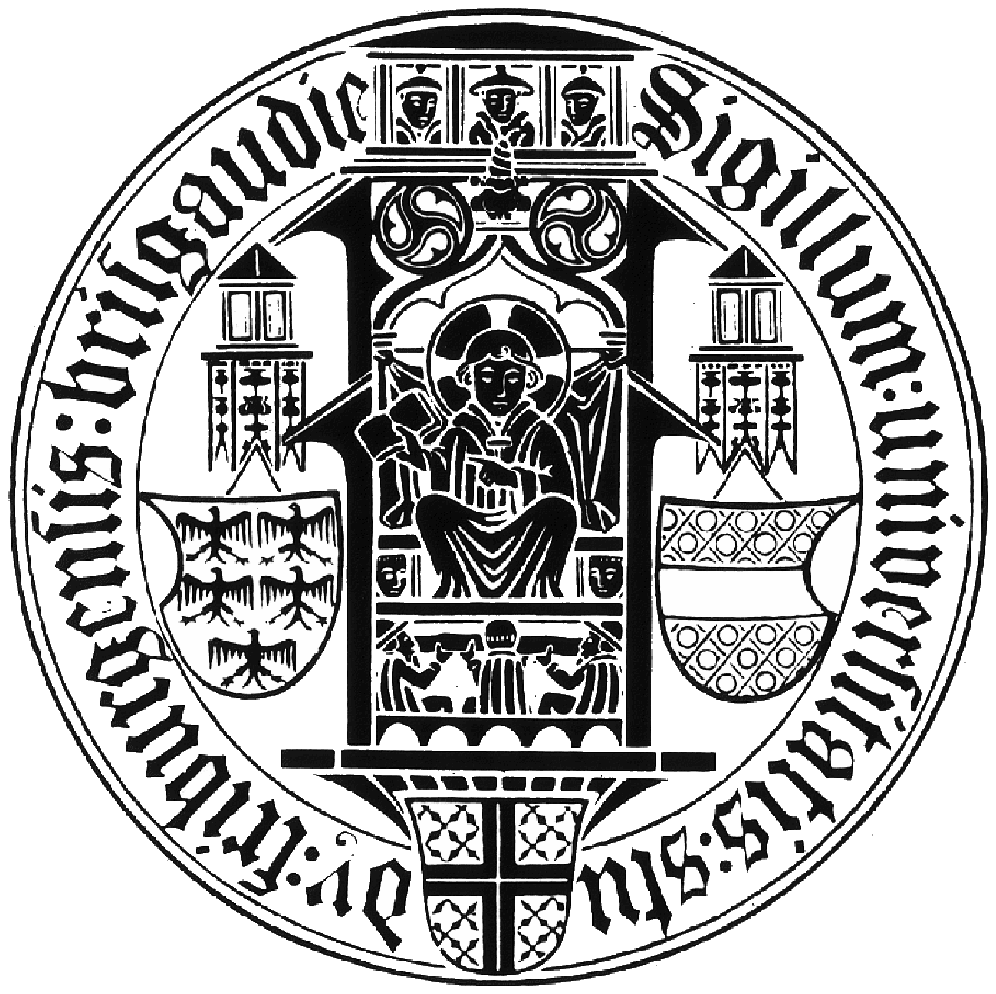}\\[0.2cm]
Albert-Ludwigs-Universit"at Freiburg\\
Fakult"at f"ur Mathematik und Physik
\end{minipage}\\[2.5cm]
{\LARGE\textbf{Globale Existenz f"ur die Interaktion eines Navier-Stokes-Fluids mit einer
linear elastischen Schale}}\\[1cm]
Dissertation zur Erlangung des  Doktorgrades\\
 der Fakult"at f"ur Mathematik und Physik \\
der Albert-Ludwigs-Universit"at Freiburg im Breisgau\\[1cm]
vorgelegt von Dipl.-Math. Daniel Lengeler\\[2cm]
  betreut von Prof. Dr. Michael R{\r{u}}{\v{z}}i{\v{c}}ka
\end{center}\end{addmargin}
\newpage
\thispagestyle{empty}
\vspace*{19cm}
\hspace*{-1.5cm} \begin{tabular}{ll}
    Dekan: & Prof. Dr. Kay K"onigsmann
    \\
    Referent: & Prof. Dr. Michael R{\r{u}}{\v{z}}i{\v{c}}ka
    \\
    Koreferent: & Prof. Dr. Helmut Abels
    \\
    Datum der Promotion: & 15. Dezember 2011
 \end{tabular} 
\end{titlepage}
\chapter*{Erkl"arung}
% Hiermit erkl"are ich, dass ich diese Arbeit selbst"andig verfasst habe, keine anderen als
% die angegebenen Quellen/Hilfsmittel verwendet habe und alle Stellen, die w"ortlich oder
% sinngem"a"s
% aus ver"offentlichten Schriften entnommen wurden, als solche kenntlich gemacht habe. Dar"uber
% hinaus
% erkl"are ich, dass diese Arbeit nicht, auch nicht auszugsweise, bereits f"ur eine andere
% Pr"ufung angefertigt wurde. 
Ich erkl"are hiermit, dass ich die vorliegende Arbeit ohne unzul"assige Hilfe Dritter und ohne
Benutzung anderer als der angegebenen Hilfsmittel angefertigt habe. Die aus anderen Quellen direkt
oder indirekt "ubernommenen Daten und Konzepte sind unter Angabe der Quelle gekennzeichnet.
Insbesondere habe ich hierf"ur nicht die entgeltliche Hilfe von Vermittlungs- bzw. Beratungsdiensten
(Promotionsberater oder anderer Personen) in Anspruch genommen. Die Arbeit wurde bisher weder im
In- noch im Ausland in gleicher oder "ahnlicher Form einer anderen Pr"ufungsbeh"orde vorgelegt.
\vspace{3cm}

\begin{center}
\begin{tabular}{p{6.5cm}p{4cm}}
\rule{3cm}{0.4pt} & \rule{4cm}{0.4pt} \\
{\small Ort, Datum} & {\small Unterschrift}\\
\end{tabular}
\end{center}
\chapter*{Danksagung}
Zun"achst m"ochte ich meiner Frau Wiebke und meinen Eltern danken. Alle drei haben mich immer
"au"serst wohlwollend unterst"utzt. Zudem danke ich meinem Betreuer Prof. Dr.
Michael {R{\r u}{\v z}i{\v c}ka}\ f"ur die fruchtbare Aufgabenstellung, zahlreiche hilfreiche
Gespr"ache und die freundschaftliche Zusammenarbeit. Ebenso danke ich Prof. Dr. Lars Diening und
meinem Mitdoktoranden Dipl.-Math. Philipp N"agele f"ur regelm"a"sige Diskussionen, die ebenfalls
sehr zum Gelingen dieser Arbeit beigetragen haben. 
% 
% Zun"achst m"ochte ich mich bei meinen Betreuern Prof. Dr. Michael {R{\r u}{\v z}i{\v c}ka}\ und
% Prof. Dr. Lars Diening bedanken, die immer ein offenes Ohr f"ur die zahlreichen Fragen und
% Probleme hatten und ohne deren Hilfe diese Arbeit sicherlich nicht zustande gekommen
% w"are. Desweiteren m"ochte ich meinem Mitdoktoranden Philipp N"agele f"ur die zahlreichen
% anregenden und hilfreichen Diskussionen danken. Nicht zuletzt gilt mein Dank meinen Eltern, die
% mich
% immer "au"serst gro"sz"ugig und wohlwollend unterst"utzt haben, und insbesondere meiner Frau
% Wiebke
% f"ur die moralische Unterst"utzung und das Ertragen mancher Laune.

%%%%%%%%%%%%%%%%%%%%%%%% tableofcontents %%%%%%%%%%%%%%%%%%%%%%%%

\tableofcontents
\newpage~
\thispagestyle{empty}
%%%%%%%%%%%%%%%%%%%%%%%%%%%%% input %%%%%%%%%%%%%%%%%%%%%%%%%%%%%
\pagenumbering{arabic}
\setcounter{page}{0}

\nomenclature[00]{$B'$}{Dualraum eines Banachraums $B$}
\nomenclature[01]{$B_0\embedding B_1$, $B_0\compactembedding B_1$}{Stetige bzw. kompakte Einbettung
von
Banach-R\"aumen}
\nomenclature[02]{$(B_0,B_1)_{\theta,p}$, $[B_0,B_1]_{\theta}$}{Reeller bzw. komplexer
Interpolationsraum}
\nomenclature[03]{$L^p(\Omega)$, $W^{s,p}(\Omega)$}{Lebesgue- bzw. Sobolev-Slobodetskii-Raum. Im
Falle $\Omega\subset\setR^3$ k\"onnen diese Bezeichner auch f\"ur $L^p(\Omega,\setR^3)$,
$W^{s,p}(\Omega,\setR^3)$ stehen.}
\nomenclature[03]{$W^{1,p}_{\dv}(\Omega)$, $W^{1,p}_{0,\dv}(\Omega)$}{Raum der divergenzfreien
Elemente von $W^{1,p}(\Omega,\setR^3)$ (mit verschwindenden Randwerten), wobei
$\Omega\subset\setR^3$ (mit Lipschitz-Rand)}
\nomenclature[04]{$H^s(\Omega)$}{$L^2$-Skala $W^{s,2}(\Omega)$}
\nomenclature[05]{$\nabla$}{Levi-Civita-Zusammenhang. Dieser kann je nach Kontext (euklidischer Raum
oder Fl\"ache) variieren. Ist $\bu$ ein Vektorfeld auf $\setR^3$, so kann $\nabla\bu$ mit der
Funktionalmatrix identifiziert werden.}
\nomenclature[06]{$\Delta$}{Dem Levi-Civita-Zusammenhang zugeordneter Laplace-Operator}
\nomenclature[07]{$\dv$, $\nabla^*$}{Dem Levi-Civita-Zusammenhang zugeordnete Divergenz}
\nomenclature[08]{$D\bu$}{Symmetrischer Anteil der Funktionalmatrix eines Vektorfeldes $\bu$
auf $\setR^3$, d.h. $D\bu:=1/2(\nabla\bu+(\nabla\bu)^T)$}
\nomenclature[09]{$d\Phi$}{Differential einer Abbildung $\Phi$ zwischen Mannigfaltigkeiten
(Teilmengen des euklidischen Raums eingeschlossen). Dieses kann im euklidischen Fall mit der
Funktionalmatrix identifiziert werden.}
\nomenclature[10]{$\bu\cdot\bv$}{$\setR^3$-Skalarpodukt}
\nomenclature[11]{$A:B$}{Matrixskalarpodukt}
\nomenclature[12]{$\pa M$}{Rand einer Mannigfaltigkeit $M$ im Sinne von Mannigfaltigkeiten}
\nomenclature[13]{$\inn M$}{Inneres einer Mannigfaltigkeit $M$ im Sinne von Mannigfaltigkeiten, d.h.
$\inn M:=M\setminus\pa M$}
\nomenclature[14]{$p'$}{Definiert durch
$\frac{1}{p}+\frac{1}{p'}=1$, wobei $1\le p\le\infty$}
\nomenclature[15]{$\supp f$}{Tr\"ager einer Funktion $f$ (per Definition abgeschlossen)}
\nomenclature[16]{$C^k_0(U)$}{Menge der $C^k$-Funktionen $f$ mit $\supp
f\subset U$, wobei $U$ Teilmenge einer Mannigfaltigkeit (euklidischer Raum
eingeschlossen)}
\nomenclature[17]{$\chi_A$}{Charakteristische Funktion einer Menge $A$}
\nomenclature[18]{$c$}{Generische Konstante}
\nomenclature[19]{$\kappa$, $\Lambda$, $B_\alpha$, $\Omega_\eta$}{Definition auf S.5}
\nomenclature[20]{$\Psi_\eta$, $\Phi_\eta$}{Definition auf S.6}
\nomenclature[21]{$\tau(\eta)$}{Definition auf S.7}
\nomenclature[22]{$\tr$}{Definition auf S.8}
\nomenclature[23]{$\gamma(\eta)$}{Definition auf S.10}
\nomenclature[24]{$\trnormal$}{Definition auf S.11}
\nomenclature[25]{$\F_\eta$, $H(\Omega_\eta)$}{Definition auf S.12}
\nomenclature[26]{$\Omega_\eta^I$, $L^p(I,L^r(\oet))$}{Definition auf S.16}
\nomenclature[27]{$K$}{Definition auf S.20}
\nomenclature[28]{$Y^I$, $X_\eta^I$, $T^I_\eta$}{Definition auf S.27}
\nomenclature[29]{$c_0$}{Definition auf S.30}
\nomenclature[30]{$\R_\epsilon$}{Definition auf S.40}
\nomenclature[31]{$\bu_0^\epsilon$, $\eta_1^\epsilon$}{Definition auf S.41}
\nomenclature[32]{$\T_\delta$}{Definition auf S.76}
\nomenclature[33]{$\M_\eta$}{Definition auf S.77}
%Symbolverzeichniss in Notationen umbenennen
\renewcommand{\nomname}{Konventionen und Notation}
% Notationen einleitender Text
\renewcommand{\nompreamble}{
In dieser Arbeit werden Vektorfelder durch Fettdruck
gekennzeichnet, w"ahrend alle anderen Tensorfelder normal gedruckt werden. Eine Ausnahme von dieser
Regel bildet Anhang A.1, wo alle Tensorfelder normal gedruckt werden. Wir verwenden zudem die
Summenkonvention, nach der "uber identische in einem Monom auftretende Indizes summiert wird.
Schlie"slich setzen wir auf Teilmengen von Mannigfaltigkeiten (euklidischer Raum eingeschlossen)
definierte Funktionen stets stillschweigend durch $0$ auf den gesamten Raum fort, sofern nicht
explizit Anderes gesagt wird.}
%Kopf anpassen
\cleardoublepage% or \clearpage
\markboth{\nomname}{\nomname}
% Notationen (Verzeichniss einfügen) [2cm] <-- Spaltenbreite
\printnomenclature[3.8cm]

\chapter{Einleitung}

Die mathematische Analysis der Interaktion von Fluiden mit Festk"orpern bei zeitlich variablen
Phasengrenzen ist seit den $90$er Jahren des letzten Jahrhunderts Gegenstand intensiver Forschung.
Die ersten Arbeiten untersuchten die Bewegung starrer K"orper in viskosen Fluiden; siehe zum
Beispiel \cite{b35}, \cite{b36}. Die Einbeziehung elastischer K"orper ist wesentlich schwieriger,
bedingt durch potentielle Regularit"atsinkompatibilit"aten des zumeist hyperbolischen
Festk"orperanteils
und des parabolischen Fluidanteils der Gleichungen. So verwendeten die ersten Existenzresultate in
diese Richtung Regularisierungen des Festk"orperanteils in Form von D"ampfungstermen, siehe
\cite{b37}, \cite{b27}, \cite{b38}, oder einer endlichen Anzahl von Moden, siehe \cite{b39}.
Speziell
in \cite{b27} wird die Interaktion eines Navier-Stokes-Fluids mit einer ged"ampften elastischen
Platte untersucht und die Langzeitexistenz schwacher L"osungen gezeigt. Ein Durchbruch wurde in
\cite{b40}, \cite{b41} erzielt. Dort wird f"ur die Interaktion eines
Navier-Stokes-Fluids mit einem dreidimensionalen elastischen K"orper, beschrieben durch lineare bzw.
quasilineare hyperbolische Gleichungen ohne Regularisierung, die Kurzzeitexistenz
f"ur sehr regul"are Daten gezeigt. Eine der zentralen Ideen dabei ist die Verwendung eines
funktionalen
Rahmens, der sich am hyperbolischen Anteil des Systems orientiert. Eine weitere wichtige Arbeit ist
\cite{b28}. Dort werden Absch"atzungen gezeigt, die es erm"oglichen auf die D"ampfung
in \cite{b27} zu verzichten. In \cite{b42} schlie"slich wird die Interaktion eines
Navier-Stokes-Fluids mit einer
ann"ahernd zweidimensionalen, im Ruhezustand gekr"ummten elastischen Struktur, einer sogenannten
Schale, untersucht. Die elastische Energie der Schale wird dabei durch die Koiter-Energie
modelliert; siehe \cite{b54}, \cite{b20}, \cite{b45}, \cite{b23} und die dortigen Referenzen. Im
Vergleich mit
\cite{b41} tritt dabei die zus"atzliche Schwierigkeit
auf, dass der (quasilineare) Gradient der Koiter-Energie in den Richtungen tangential zur
ausgelenkten Schale elliptisch degeneriert ist. Allerdings konnte die Kurzzeitexistenz nur unter
Vernachl"assigung der Massentr"agheit der Schale gezeigt werden. In diesem Fall wird die Auslenkung
der Schale durch eine elliptische Gleichung mit zeitabh"angiger rechter Seite beschrieben.

In der vorliegenden Arbeit zeigen wir die Langzeitexistenz schwacher L"osungen f"ur die Interaktion
eines Navier-Stokes-Fluids mit einer elastischen Koiter-Schale unter Ber"ucksichtigung der
Massentr"agheit, wobei wir die Gleichungen der elastischen Struktur linearisieren und
ihre Bewegung auf transversale Auslenkungen einschr"anken. Die L"osungen existieren,
solange wir sicherstellen k"onnen, dass sich verschiedene Teile der Schale 
nicht ber"uhren. Die Gleichung der elastischen Struktur, die wir erhalten, ist eine
Verallgemeinerung der instation"aren, linearen Kirchhoff-Love-Plattengleichung f"ur transversale
Auslenkungen. Mithin ist unser Resultat eine direkte Verallgemeinerung von \cite{b27}, \cite{b28}.
Unsere Konstruktion schwacher L"osungen folgt in den grundlegenden Z"ugen der in
\cite{b27}. Die zus"atzliche Schwierigkeit unserer Aufgabe im Vergleich zu dieser
Arbeit besteht nicht in der (wenig) komplizierteren Gleichung der
elastischen Struktur, sondern
vielmehr in der allgemeinen Geometrie der Schale. Diese erzwingt die Entwicklung neuer Techniken,
die "uber die vorliegende Arbeit hinaus von Interesse sein d"urften. Als Nebenprodukt verk"urzt
unser Vorgehen den Beweis in \cite{b27}, \cite{b28}.\footnote{Angesichts der Seitenzahl der
vorliegenden Arbeit mag man sich "uber diese Aussage wundern. W"urde man sich jedoch auf den Fall
einer Platte beschr"anken, so w"urden einige der folgenden Aussagen und Konstruktionen trivial
oder gar hinf"allig, und die Seitenzahl w"urde sich erheblich reduzieren.} Im scharfen
Gegensatz zu diesen Arbeiten tritt
in unserem Beweis nicht eine einzige l"angere Absch"atzung auf. Zudem ben"otigen wir an keiner
Stelle eine Regularisierung der Schalengleichung. Schlie"slich ist unser Beweis so gestaltet, dass
eine "Ubertragung
auf verallgemeinerte Newton'sche Fluide in greifbare N"ahe r"uckt.

Eine der Schwierigkeiten dieser Arbeit besteht in der geringen Regularit"at der Auslenkungen der
Schale, die den mathematischen Rand des Fluids darstellt. In Kapitel 2 werden wir deshalb zun"achst
einige Vor"uberlegungen bez"uglich Gebieten mit geringer Randregularit"at anstellen. In Kapitel 3
stellen wir die Koiter-Energie vor, und in Kapitel 4 folgt eine Spezifizierung der Problemstellung.
Anschlie"send geben wir in Kapitel 5 das zentrale Existenzresultat dieser Arbeit an und f"uhren den
zugeh"origen Beweis. In Kapitel 6 skizzieren wir die ersten Schritte einer "Ubertragung des Beweises
auf verallgemeinerte Newton'sche Fluide, und Kapitel 7 gibt schlie"slich einen Ausblick auf
m"ogliche zuk"unftige Forschung.
\chapter{Variable Gebiete}

Es sei $\Omega\subset\setR^3$ ein beschr"anktes, nichtleeres Gebiet mit $C^4$-Rand und "au"serer
Einheitsnormale $\bnu\in C^3(\pa\Omega,\setR^3)$. Offenbar ist $\pa\Omega$ eine
geschlossene\footnote{D.h. kompakt
und nicht berandet.}, nicht notwendig zusammenh"angende Fl"ache. Wir bezeichnen mit $dA$ das
Fl"achenma"s von $\pa\Omega$ und f"ur $\alpha>0$ mit $S_{\alpha}$ den offenen
$\alpha$-Schlauch um $\pa\Omega$. Es existiert ein maximales $\kappa>0$ derart, dass die Abbildung 
\begin{equation*}
 \begin{aligned}
\Lambda: \partial\Omega\times (-\kappa,\kappa)\rightarrow
S_{\kappa},\
(q,s)\mapsto q + s\,\bnu(q)  
 \end{aligned}
\end{equation*}
ein $C^3$-Diffeomorphismus ist; siehe zum Beispiel Theorem
10.19 in \cite{b1}. F"ur die Inverse $\Lambda^{-1}$ schreiben wir auch $x\mapsto(q(x),s(x))$; vgl.
Abbildung \ref{Bild1}. Man beachte, dass $\kappa$
nicht notwendig klein ist. Ist $\Omega$ der Ball vom
Radius $R$, so ist $\kappa=R$. $\Lambda$ wird zum Rand seines Definitionsbereiches hin
singul"ar.
\vspace{0.5cm}
\begin{figure}[h]
\centering
\input{Bild1b.pstex_t}
\caption{} 
\label{Bild1}
\end{figure}
\vspace{0.5cm}
\hspace{-0.2cm}Wir setzen $B_\alpha:=\Omega\cup S_{\alpha}$ f"ur $0<\alpha<\kappa$. Die Abbildung
$\Lambda(\,\cdot\, ,\alpha): \pa\Omega\rightarrow\pa B_\alpha$ ist ebenfalls ein
$C^3$-Diffeomorphismus, sodass $B_\alpha$ ein beschr"anktes Gebiet mit $C^3$-Rand ist. F"ur
stetiges $\eta:\pa\Omega\rightarrow (-\kappa,\kappa)$ setzen wir
\begin{equation*}
 \begin{aligned}
  \Omega_\eta:=\Omega\setminus S_\kappa\ \cup \{x\in S_\kappa\ |\ s(x)<\eta(q(x))\};
 \end{aligned}
\end{equation*}
vgl. Abbildung \ref{Bild1}. Offenbar ist $\Omega_\eta$ offen. Ist $\eta\in C^k(\pa\Omega)$,
$k\in\{1,2,3\}$, so ist $\pa\Omega_\eta$ eine
$C^k$-Fl"ache, wie wir in K"urze sehen werden. In
diesem Falle seien $\bnu_\eta$ die "au"sere
Einheitsnormale und $dA_\eta$ das Fl"achenma"s von $\pa\Omega_\eta$. Wir wollen nun einen
Hom"oomorphismus von $\overline\Omega$ auf
$\overline{\Omega_\eta}$ konstruieren, die sogenannte Hanzawa-Transformation. Sei dazu
$\beta\in C^\infty(\setR)$ mit $\beta=0$ in
einer
Umgebung von $-1$ und $\beta=1$ in einer Umgebung von $0$. F"ur stetiges $\eta:\pa\Omega\rightarrow
(-\kappa,\kappa)$ definieren wir 
\[\Psi_\eta: \overline\Omega\rightarrow\overline{\Omega_\eta}\]
in $S_\kappa\cap\overline\Omega$ durch
\begin{equation}\label{eqn:defpsi}
 \begin{aligned}
x&\mapsto x +
\bnu(q(x))\,
\eta(q(x))\,\beta(s(x)/\kappa)=q(x)+\bnu(q(x))\Big(s(x)+\eta(q(x))\,\beta(s(x)/\kappa)\Big). 
 \end{aligned}
\end{equation}
In $\Omega\setminus S_\kappa$ sei $\Psi_\eta$ die Identit"at. Punkte
$x\in S_\kappa\cap\overline\Omega$ werden durch $\Psi_\eta$ im Wesentlichen um die L"ange
$\eta(q(x))$ in Richtung $\bnu(q(x))$ verschoben. Allerdings m"ussen wir die L"ange, um die
translatiert wird, zum "`inneren"' Rand von $S_\kappa\cap\overline\Omega$ hin geeignet gegen $0$
gehen lassen, um eine Bijektion zu erhalten. Das erreichen wir durch die Multiplikation mit
der Abschneidefunktion $\beta$, die noch n"aher zu spezifizieren ist. Punkte $x\in
S_\kappa\cap\overline\Omega$ mit festem Fu"spunkt $q(x)=q\in\pa\Omega$ werden durch
\[\beta_q: s\mapsto s+\eta(q)\,\beta(s/\kappa),\ s=s(x)\in [-\kappa,0]\]
abgebildet. Insbesondere gilt $\beta_q(-\kappa)=-\kappa$ und $\beta_q(0)=\eta(q)$. Damit
$\beta_q$ eine (stetig differenzierbare) Inverse besitzt, muss
\[0<\beta_q'(s)=1+\eta(q)/\kappa\ \beta'(s/\kappa)\]
gelten. Wir fordern dazu $|\beta'(s)|<\kappa/|\eta(q)|$ f"ur alle $s\in [-1,0]$ und alle
$q\in\pa\Omega$, was wegen
$\norm{\eta}_{L^\infty(\pa\Omega)}<\kappa$ m"oglich ist. $\Psi_\eta$ ist dann bijektiv und die
Inverse ist in $S_\kappa\cap\overline\Omega$ durch
\[\Psi_\eta^{-1}:x\mapsto q(x)+\bnu(q(x))\,\beta_{q(x)}^{-1}(s(x))\]
gegeben. Offenbar h"angt $\beta_q^{-1}$ stetig und im Falle $\eta\in
C^k(\pa\Omega)$, $k\in\setN$, sogar $k$-mal stetig
differenzierbar von $q$ ab; bezeichnen wir mit $d$ das Differential bez"uglich der Variablen $q$
und mit $'$ die Ableitung nach der skalaren Variable, so gilt wegen
\[  0= ds=d(\beta_q^{-1}(\beta_q(s)))=(d\beta_q^{-1})(\beta_q(s)) +
(\beta_q^{-1})'(\beta_q(s))\,(d\beta_q)(s)\]
die Identit"at 
\begin{equation}\label{eqn:ableitung}
 \begin{aligned}
(d\beta_q^{-1})(\beta_q(s))=-\frac{\beta(s/\kappa)}{1+\eta(q)/\kappa\ \beta'(s/\kappa)}d\eta(q).  
 \end{aligned}
\end{equation}
$\Psi_\eta$ ist somit ein
Hom"oomorphismus und im Falle $\eta\in C^k(\pa\Omega)$, $k\in\{1,2,3\}$, sogar ein
$C^k$-Diffeomorphismus. Auch der
Hom"oomorphismus \[\Phi_\eta:=\Psi_\eta|_{\pa\Omega}:\pa\Omega\rightarrow\pa\Omega_\eta,\
q\mapsto q+\eta(q)\,\bnu(q)\]
mit der Inversen $x\mapsto q(x)$ ist f"ur $\eta\in C^k(\pa\Omega)$, $k\in\{1,2,3\}$, ein
$C^k$-Diffeomorphismus. Die Kettenregel und \eqref{eqn:ableitung} zeigen, dass die Eintr"age der
Funktionalmatrizen von
$\Psi_\eta$, $\Phi_\eta$ und ihren Inversen von der Form 
\begin{equation}\label{eqn:funk}
 \begin{aligned}
b_0+\sum_i b_i\ (\pa_i\eta)\circ q  
 \end{aligned}
\end{equation}
sind mit stetigen und beschr"ankten Funktionen $b_0$, $b_i$, die f"ur $\tau(\eta)\rightarrow\infty$
mit
\begin{equation}\label{eqn:taueta}
 \begin{aligned}
\tau(\eta):=\left\{\begin{array}{cl} (1-\norm{\eta}_{L^\infty(\pa\Omega)}/\kappa)^{-1} &
\text{, falls }\norm{\eta}_{L^\infty(\pa\Omega)}<\kappa\\
\infty & \text{, sonst}\end{array}\right.
 \end{aligned}
\end{equation}
teilweise singul"ar werden, weil dann die Funktionen $\beta_q$
singul"ar werden und der Diffeomorphismus $\Lambda$ in der N"ahe seiner Singularit"aten ausgewertet
wird. Zudem sind die Tr"ager der $b_i$ in $S_\kappa$ enthalten, wobei der Abstand zum Rand von
$S_\kappa$ f"ur $\tau(\eta)\rightarrow\infty$ gegen $0$ geht. F"ur $\tau(\eta)\rightarrow\infty$
wird somit die Abbildung $q$ in \eqref{eqn:funk} in der N"ahe ihrer Singularit"aten
ausgewertet. Neben m"oglichen Irregularit"aten der Auslenkung $\eta$ k"onnen die Abbildungen
$\Psi_\eta$ und $\Phi_\eta$ also auch dadurch singul"ar werden, dass die maximale Auslenkung an
$\kappa$ heranr"uckt. Aus diesem Grund werden die Stetigkeitskonstanten der im Folgenden
konstuierten linearen Abbildungen zwischen Funktionenr"aumen stets von $\tau(\eta)$ abh"angen.

Die Abbildung $\Psi_\eta$ h"angt von der Abschneidefunktion $\beta$ ab, die wiederum in
Abh"angigkeit von $\tau(\eta)$ gew"ahlt werden kann. Wann immer wir mit Folgen $(\eta_n)$ von
Auslenkungen mit $\sup_n\tau(\eta_n)<\infty$ zu tun haben werden, wollen und k"onnen wir $\beta$
unabh"angig vom Folgenindex w"ahlen. Dadurch stellen wir sicher, dass die Folge $(\Psi_{\eta_n})$
konvergiert, falls $(\eta_n)$ konvergiert.

Eine Bi-Lipschitz-Abbildung von Definitionsbereichen induziert Isomorphismen der jeweiligen
$L^p$- und $W^{1,p}$-R"aume. F"ur $\eta\in H^2(\pa\Omega)$ ist wegen der Einbettung
$H^2(\pa\Omega)\embedding C^{0,\theta}(\pa\Omega)$ f\"ur $\theta<1$ die Abbildung
$\Psi_\eta$ gerade nicht bi-Lipschitz-stetig, sodass ein kleiner Verlust unter dieser
Transformation entsteht.
\begin{lemma}\label{lemma:psi}
Es seien $1<p\le\infty$ und $\eta\in H^2(\pa\Omega)$ mit $\norm{\eta}_{L^\infty(\pa\Omega)}<\kappa$.
Dann ist die lineare Abbildung $v\mapsto v\circ\Psi_\eta$ stetig von $L^p(\Omega_\eta)$ nach
$L^r(\Omega)$ und von $W^{1,p}(\Omega_\eta)$ nach $W^{1,r}(\Omega)$ f"ur alle $1\le r<p$. Eine
analoge Aussage gilt mit $\Psi_\eta^{-1}$ anstelle von $\Psi_\eta$. Die
Stetigkeitskonstanten h"angen nur von $\Omega$, $p$, $\norm{\eta}_{H^2(\pa\Omega)}$, $\tau(\eta)$
und
$r$ ab; sie bleiben beschr"ankt, falls $\norm{\eta}_{H^2(\pa\Omega)}$ und $\tau(\eta)$
beschr"ankt bleiben.
% Konvergiert die Folge
% $(\eta_n)$ schwach gegen $\eta$ in $H^2(\pa\Omega)$, so gilt f"ur $v\in W^{1,p}(\Omega_\eta)$
% \[v\circ\Psi_{\eta_n}\rightarrow v\circ\Psi\text{ in } W^{1,r}(\Omega).\]
\end{lemma}
\beweis Wir k"onnen ohne Einschr"ankung $p<\infty$ annehmen. Wir approximieren $\eta$ durch
Funktionen $(\eta_n)\subset
C^2(\pa\Omega)$ in
$H^2(\pa\Omega)$, insbesondere gleichm"a"sig. Wegen \eqref{eqn:funk} und der
Einbettung \[H^2(\pa\Omega)\embedding W^{1,s}(\pa\Omega) \text{ f"ur alle }1\le s<\infty\]
sind die Eintr"age der Funktionalmatrix von $\Psi_{\eta_n}^{-1}$
und somit auch die
Funktionaldeterminante in $L^s(\Omega_{\eta_n})$ f"ur jedes $s<\infty$ in Abh"angigkeit von
$\tau(\eta_n)$ beschr"ankt. Wir erhalten also f"ur $v\in C_0^\infty(\setR^3)$, $r<p$ und
$1/r=1/p+1/s$ unter Verwendung des Transformationssatzes und der H"older-Ungleichung 
\begin{equation*}
 \begin{aligned}
\norm{v\circ\Psi_{\eta_n}}_{L^r(\Omega)}=\norm{v\,
(\det d\Psi_{\eta_n}^{-1})^{1/r}}_{L^r(\Omega_{\eta_n})}&\le \norm{(\det
d\Psi_{\eta_n}^{-1})^{1/r}}_{L^s(\Omega_{\eta_n})}\,\norm{v}_{L^p(\Omega_{
\eta_n})}.
 \end{aligned}
\end{equation*}
Aus der Konvergenz von $(\eta_n)$ in $H^2(\pa\Omega)$ folgt
\begin{equation*}
 \begin{aligned}
\norm{v\circ\Psi_{\eta}}_{L^r(\Omega)}&\le \norm{(\det
d\Psi_{\eta}^{-1})^{1/r}}_{L^s(\Omega_{\eta})}\,\norm{v}_{L^p(\Omega_{
\eta})}\le c(\Omega,p,\norm{\eta}_{H^2(\pa\Omega)},\tau(\eta),r)\,\norm{v}_{L^p(\Omega_{
\eta})}.
 \end{aligned}
\end{equation*}
Mit Hilfe der Dichtheit glatter Funktionen in $L^p(\Omega_\eta)$ folgern wir die Stetigkeit
bez"uglich der Lebesgue-R"aume.                                

Aufgrund der Kettenregel und \eqref{eqn:funk} gilt
\[\norm{\nabla(v\circ\Psi_{\eta_n})}_{L^r(\Omega)}\le
c(\Omega,\norm{\eta_n}_{H^2(\pa\Omega)},\tau(\eta_n),r)\,\norm{(\nabla 
v)\circ\Psi_{\eta_n}}_{L^{(r+p)/2}(\Omega)}.\]
Zudem folgt aus der gleichm"a"sigen Konvergenz von $((\nabla v)\circ\Psi_{\eta_n})$ und der
Konvergenz der Eintr"age der Funktionalmatrix von $\Psi_{\eta_n}$ in $L^s(\Omega)$ f"ur $s<\infty$
die Konvergenz von $(\nabla(v\circ\Psi_{\eta_n}))$ gegen $\nabla(v\circ\Psi_{\eta})$
in $L^r(\Omega)$. Wir erhalten somit die Stetigkeit bez"uglich der Sobolev-R"aume aus
dem bereits Gezeigten sowie der Dichtheit glatter Funktionen in $W^{1,p}(\Omega_\eta)$; siehe
Proposition \ref{theorem:dicht}.

Der Beweis der analogen Aussage mit $\Psi_\eta^{-1}$ anstelle von $\Psi_\eta$ geht genauso.
\qed\\

\begin{bemerkung}\label{bem:konv}
Die Folge $(\eta_n)$ konvergiere gegen ein $\eta$ schwach in $H^2(\pa\Omega)$ und aufgrund der aus
dem Satz von Arzela-Ascoli folgenden kompakten Einbettung $(\theta<1)$ 
\[H^2(\pa\Omega)\embedding C^{0,\theta}(\pa\Omega)\compactembedding C(\pa\Omega)\]
insbesondere gleichm"a"sig. Zudem gelte $\sup_n\tau(\eta_n)<\infty$. Setzen wir $v\in
L^p(\Omega_\eta)$ durch $0$ auf $\setR^3$ fort, so konvergiert die Folge
$(v\circ\Psi_{\eta_n})$ gegen $v\circ\Psi_{\eta}$ in $L^r(\Omega)$, $r<p$. Das folgt mit Lemma
\ref{lemma:psi} und der Approximierbarkeit von $v$
durch $C_0^\infty(\setR^3)$-Funktionen $\tilde v$ aus der Absch"atzung (in der $L^r(\Omega)$-Norm)
\begin{equation*}
 \begin{aligned}
  \norm{v\circ\Psi_\eta-v\circ\Psi_{\eta_n}}\le \norm{\tilde v\circ\Psi_\eta-\tilde
v\circ\Psi_{\eta_n}} + \norm{(v-\tilde v)\circ\Psi_{\eta}} +
\norm{(v-\tilde v)\circ\Psi_{\eta_n}}.
 \end{aligned}
\end{equation*}
Ist $v\in L^p(\Omega)$ und setzen wir die Funktionen $v\circ\Psi^{-1}_{\eta_n}$ und
$v\circ\Psi^{-1}_{\eta}$ durch $0$ auf $\setR^3$ fort, so l"asst sich auf vollkommen analoge Weise
die Konvergenz der Folge $(v\circ\Psi^{-1}_{\eta_n})$ gegen $v\circ\Psi^{-1}_{\eta}$ in
$L^r(\setR^3)$, $r<p$, einsehen.
\end{bemerkung}
\vspace{0.2cm}

Wir konstruieren nun einen Spuroperator f"ur Auslenkungen $\eta\in H^2(\pa\Omega)$. Man beachte
dabei, dass $\pa\Omega_\eta$ wegen der Einbettung $H^2(\pa\Omega)\embedding C^{0,\theta}(\pa\Omega)$
f\"ur $\theta<1$ gerade nicht Lipschitz-stetig ist. Die Einschr"ankung $\cdot\,|_{\pa\Omega}$ ist
fortan im Sinne des "ublichen Spuroperators f"ur regul"are
R"ander zu verstehen.
\begin{korollar}\label{lemma:spur}
Es seien $1<p\le\infty$ und $\eta\in H^2(\pa\Omega)$ mit $\norm{\eta}_{L^\infty(\pa\Omega)}<\kappa$.
Dann ist die lineare Abbildung $\tr: v\mapsto (v\circ\Psi_\eta)|_{\pa\Omega}$ wohldefiniert und
stetig von $W^{1,p}(\Omega_\eta)$ nach $W^{1-1/r,r}(\pa\Omega)$ f"ur alle $1<r<p$. Die
Stetigkeitskonstante h"angt nur von $\Omega$, $\norm{\eta}_{H^2(\pa\Omega)}$, $\tau(\eta)$ und $r$
ab; sie bleibt beschr"ankt, falls $\norm{\eta}_{H^2(\pa\Omega)}$ und $\tau(\eta)$
beschr"ankt bleiben.
\end{korollar}
\beweis Die Behauptung folgt sofort aus Lemma \ref{lemma:psi} und den Stetigkeitseigenschaften
des
"ublichen Spuroperators; siehe Theorem \ref{theorem:spur}.
\qed\\

Der Operator $\tr$ ist nichts anderes als die stetige Fortsetzung der "`zur"uckgeholten"' Spur
\[v\mapsto (q\mapsto v(q+\eta(q)\,\bnu(q)), q\in\pa\Omega),\]
die f"ur glatte Funktionen $v$ wohldefiniert ist. Dabei bleibt die lokale Fl"achenverzerrung der
Abbildung $\Phi_\eta: q\mapsto q+\eta(q)\,\bnu(q)$ au"sen vor. Wir messen also die Spuren mit dem
Fl"achenma"s von $\pa\Omega$ und nicht mit dem von $\pa\Omega_\eta$. 
% Das ist schon deshalb sinnvoll,
% weil das zweidimensionale Hausdorff-Ma"s auf $\pa\Omega_\eta$ als Graph einer nicht
% Lipschitz-stetigen Funktion im Allgemeinen nicht (lokal) endlich ist. 
Dieser Operator stellt genau die richtige Konstruktion f"ur den Vergleich $\tr\bu=\pa_t\eta\,\bnu$
der
Geschwindigkeiten von Fluid
und Schale am Rand dar. Man beachte, dass die Gleichungen der Elastizit"at in Lagrange-Koordinaten
formuliert werden.

Aus Lemma \ref{lemma:psi} und den Sobolev-Einbettungen f"ur Gebiete mit regul"arem Rand folgen
Sobolev-Einbettungen f"ur unsere speziellen Gebiete.
\begin{korollar}\label{lemma:sobolev}
Es seien $1<p<3$ und $\eta\in H^2(\pa\Omega)$ mit $\norm{\eta}_{L^\infty(\pa\Omega)}<\kappa$. Dann
gilt die kompakte Einbettung \[W^{1,p}(\Omega_\eta)\compactembedding L^s(\Omega_\eta)\]
f"ur $1\le s < p^*=3p/(3-p)$. Die Stetigkeitskonstante h"angt nur von
$\Omega$, $\norm{\eta}_{H^2(\pa\Omega)}$, $\tau(\eta)$, $p$ und $s$ ab; sie bleibt beschr"ankt,
falls $\norm{\eta}_{H^2(\pa\Omega)}$ und $\tau(\eta)$
beschr"ankt bleiben.
\end{korollar}
\vspace{0.2cm}
% \beweis Wie im Beweis von Lemma \ref{lemma:spur} zeigen wir f"ur $v\in W^{1,p}(\Omega_\eta)$ und
% $r<p$ die Absch"atzung
% \[\norm{v\circ\Psi_\eta}_{W^{1,r}(\Omega)}\le c
% \norm{v}_{W^{1,p}(\Omega_\eta)}.\]
% Die "ubliche Soboleveinbettung, siehe ANHANG, liefert f"ur $\tilde s<r$
% \[\norm{v\circ\Psi_\eta}_{L^{\tilde s}(\Omega)}\le c
% \norm{v\circ\Psi_\eta}_{W^{1,r}(\Omega)}.\]
% Desweiteren erhalten wir unter Verwendung des Transformationssatz und mit Hilfe der Argumente aus
% dem Beweis von Lemma \ref{lemma:spur} f"ur $s<\tilde s$
% \[\norm{v}_{L^s(\Omega)}\le c \norm{v\circ\Psi_\eta}_{L^{\tilde s}(\Omega)}.\]
% \qed 
% Wir wollen nun den Kern des Spuroperators im Falle
% divergenzfreier Vektorfelder charakterisieren.
% \begin{lemma}
% Es seien $1<p<\infty$ und $\eta\in H^2(\pa\Omega)$ mit $|\eta| <\kappa$. Dann liegt die Menge
% $C_{0,\dv}^\infty(\Omega_\eta):=\{\bphi\in C_0^\infty(\Omega_\eta)|\dv\bphi=0\}$ dicht im
% Kern des Operators $\tr: W^{1,p}_{\dv}(\Omega_\eta):=\{\bphi\in
% W^{1,p}(\Omega_\eta)|\dv\bphi=0\}\rightarrow W^{1-1/r,r}(\pa\Omega)$, wobei $r<p$. 
% \end{lemma}
% \beweis \marginpar{NOCH UNKLAR! EVT. GEN"UGT WENIGER!}\qed\\
Es seien im Folgenden stets $H$ die mittlere Kr"ummung und $G$ die Gau"s'sche K"ummung von
$\pa\Omega$.

\begin{proposition}\label{lemma:partInt}
Es seien $1<p\le\infty$ und $\eta\in H^2(\pa\Omega)$ mit $\norm{\eta}_{L^\infty(\pa\Omega)}<\kappa$.
Dann gilt f"ur
$\bphi\in W^{1,p}(\Omega_\eta)$ mit $\tr \bphi=b\,\bnu$, $b$ eine skalare Funktion, und $\psi\in
C^1(\overline{\Omega_\eta})$
\[\int_{\Omega_\eta} \bphi\cdot\nabla\psi\ dx=-\int_{\Omega_\eta} \dv\bphi\ \psi\ dx+
\int_{\pa\Omega}b\, (1-2H\eta+G\,\eta^2)\, \tr \psi\ dA.\]
\end{proposition}
\beweis Wir k"onnen ohne Einschr"ankung $p<\infty$ annehmen. Wir approximieren $\bphi$ durch
$(\bphi_k)\subset C^\infty_0(\setR^3)$ in $W^{1,p}(\Omega_\eta)$ und $\eta$ durch
$(\eta_n)\subset C^2(\pa\Omega)$ in $H^2(\pa\Omega)$. Dann erhalten wir
mittels partieller Integration
\[\int_{\Omega_{\eta_n}} \bphi_k\cdot\nabla\psi\ dx=-\int_{\Omega_{\eta_n}} \dv\bphi_k\
 \psi\ dx + \int_{\pa\Omega_{\eta_n}} \bphi_k\cdot\bnu_{\eta_n}\,
\psi\
dA_{\eta_n}.\] 
Anwenden des Transformationssatzes, siehe \eqref{eqn:trafo}, auf das Randintegral ergibt
\[\int_{\pa\Omega}\tren\bphi_k\cdot(\bnu_{\eta_n}\circ \Phi_{\eta_n})\
\tren\psi\ |\det d\Phi_{\eta_n}|\ dA.\]
Mit Hilfe des Gram-Schmidt-Algorithmus lassen sich in einer Umgebung eines jeden Punktes von
$\pa\Omega$ orthonormale, tangentiale $C^3$-Vektorfelder $\be_1,\be_2$
konstruieren; siehe zum Beispiel \cite{b1}. Ohne Einschr"ankung gelte $\be_1\times\be_2=\bnu$. Wir
setzen
$\bv_i^n:=d\Phi_{\eta_n}\be_i$. Die Funktionen $\bv_i^n$ lassen sich als nichttangentiale lokale
Vektorfelder auf
$\pa\Omega$ auffassen. Die
lokale Fl"achenverzerrung von $\Phi_{\eta_n}$ ist identisch der Fl"ache des durch die $\bv_i^n$
aufgespannten Parallelepipeds, also mit $|\bv_1^n\times
\bv_2^n|$, w"ahrend die Normale $\bnu_{\eta_n}\circ\Phi_{\eta_n}$
durch $(\bv_1^n\times
\bv_2^n)/|\bv_1^n\times \bv_2^n|$ gegeben ist. Die Felder $\bv_1^n\times
\bv_2^n$ sind offenbar unabh"angig von der konkreten Wahl der $\be_i$, sodass wir durch die lokale
Definition tats"achlich globale, d.h. auf  $\pa\Omega$ definierte, Vektorfelder $\bv^n$ erhalten.
Das
Randintegral l"a"st sich also schreiben als
\[\int_{\pa\Omega}\tren\bphi_k\cdot\bv^n\,
\tren\psi\ dA.\]
Ist $q\in\pa\Omega$ und $c$ eine Kurve in $\pa\Omega$ mit $c(0)=q$ und $\frac{d}{dt}\big|_{t=0}\
c(t) = \be_i(q)$, so ist
\begin{equation}\label{eqn:v}
 \begin{aligned}
\bv_i^n(q)&=\frac{d}{dt}\Big|_{t=0}\Phi_{\eta_n}(c(t))=\frac{d}{dt}\big|_{t=0}
(c(t)+\eta_n(c(t))\,\bnu(c(t)))\\
&=\be_i(q)+d\eta_n\be_i(q)\, \bnu(q)+\eta_n(q)\,\frac{d}{dt}\Big|_{t=0}
\bnu(c(t))\\
&=\be_i(q)+d\eta_n\be_i(q)\, \bnu(q)-\eta_n(q)\, h_i^j(q)\, \be_j(q),
 \end{aligned}
\end{equation}
wobei $h^j_i$ die Komponenten der Weingartenabbildung bez"uglich der
Orthonormalbasis $\be_1,\be_2$ sind. Da $\eta_n$ in $H^2(\pa\Omega)$ konvergiert, konvergiert
$\bv_i^n$ in $L^r$ auf seinem Definitionsbereich f"ur alle $r<\infty$. Somit konvergiert $\bv^n$
gegen
ein $\bv$ in $L^r(\pa\Omega)$ f"ur alle $r<\infty$. Lassen wir nun zun"achst $n$ und anschlie"send
$k$ gegen unendlich gehen, so erhalten wir
\[\int_{\Omega_\eta} \bphi\cdot\nabla\psi\ dx=-\int_{\Omega_\eta} \dv\bphi\
 \psi\ dx + \int_{\pa\Omega} b\,\bnu\cdot\bv\,
\tr\psi\ dA.\]
Es bleibt, $\bnu\cdot\bv$ zu bestimmen. Wir lesen von Gleichung \eqref{eqn:v} ab, dass 
\[\bnu\cdot(\bv_1^n\times\bv_2^n)=1-(h_1^1+h^2_2)\,\eta_n+(h_1^1h_2^2-
h_1^2h_2^1)\,\eta_n^2=1-2H\eta_n+G\,\eta_n^2\]
gilt, woraus wir durch Grenz"ubergang die Behauptung folgern.\qed\\

Das im Beweis konstruierte Feld $\bv\in L^r(\pa\Omega)$, $r<\infty$, wird auch im Folgenden noch
Verwendung finden. Um die
Abh"angigkeit von $\eta$ auszudr"ucken, setzen wir $\bv=\bv_\eta$. Da die Felder $\bv^n$ unabh"angig
von der Wahl der $\be_i$ sind, gilt dies ebenso f"ur $\bv_\eta$. Dies deckt
sich mit der Koordinateninvarianz der Interpretation von $\bv_\eta$ als "`"au"sere Normale"' an
$\pa\Omega_\eta$, deren
L"ange identisch der "`lokalen Fl"achenverzerrung"' von $\Phi_\eta$ ist.

\begin{bemerkung}\label{bem:trafo}
 Der Beweis von Proposition \ref{lemma:partInt} zeigt f"ur $\psi\in L^1(\pa\Omega)$ und $\eta$
hinreichend glatt die Identit"at
\[\int_{\pa\Omega_{\eta}} \big((\psi\,\bnu)\circ\Phi_{\eta}^{-1}\big)\cdot\bnu_\eta\
dA_\eta=\int_{\pa\Omega} \psi\,\bnu\cdot\bv_\eta\ dA=\int_{\pa\Omega} \psi\,(1-2H\eta+G\,\eta^2)\
dA.\]
\end{bemerkung}
\vspace{0.2cm}

\begin{bemerkung}\label{bem:groessernull}
Die Gr"o"se $\gamma(\eta):=1-2H\eta+G\,\eta^2$ ist positiv, solange $|\eta|<\kappa$. Die Nullstellen
dieses Polynoms liegen f"ur $G\not=0$ n"amlich bei
\[\frac{H}{G}\pm\frac{\sqrt{H^2-G}}{|G|}=\frac{1}{2}\bigg(\frac{h_1+h_2}{h_1h_2}\pm\frac{|h_1-h_2|}{
|h_1h_2|}\bigg),\]
wobei $h_1$, $h_2$ die Hauptkr"ummungen sind. Die Betr"age der Nullstellen sind somit identisch
\[\frac{|h_1+h_2\pm(h_1-h_2)|}{2|h_1h_2|}=|h_1|^{-1},|h_2|^{-1}.\]
F"ur $G\not=0$ ist die Behauptung somit bewiesen, wenn wir zeigen k"onnen,
dass $\kappa\le\min(|h_1(q)|^{-1},|h_2(q)|^{-1})$ f"ur alle $q\in\pa\Omega$. F"ur $G=0$ ist der
Beweis dann offensichtlich. Um die Absch"atzung von $\kappa$ einzusehen, betrachten wir den
Diffeomorphismus
\[\Lambda_\alpha:=\Lambda(\,\cdot\, ,\alpha):\pa\Omega\rightarrow
\Lambda_\alpha(\pa\Omega)\subset\setR^3\]
mit $\alpha\in(-\kappa,\kappa)$. Die Absch"atzung von $\kappa$ folgt nun aus der Tatsache, dass das
Differential von $\Lambda_\alpha$ an der Stelle $q\in \pa\Omega$ singul"ar wird f"ur
$|\alpha|\nearrow\min(|h_1(q)|^{-1},|h_2(q)|^{-1})$. Ist n"amlich $\be_1,\be_2$ eine Eigenbasis
der
Weingarten-Abbildung in $q$, so gilt $d\Lambda_\alpha\,\be_i=(1-\alpha\, h_i(q))\,\be_i$.
\end{bemerkung}
\vspace{0.2cm}

Wir betrachten nun den kanonisch normierten Raum 
\[E^p(\Omega_\eta):=\{\bphi\in
L^p(\Omega_\eta)\ |\ \dv \bphi\in
L^p(\Omega_\eta)\}\]
mit $1\le p\le\infty$. F"ur solche Vektorfelder l"asst sich immerhin noch eine Spur in
Normalenrichtung konstruieren.
\begin{proposition}\label{lemma:nspur}
Es seien $1<p<\infty$ und $\eta\in H^2(\pa\Omega)$ mit $\norm{\eta}_{L^\infty(\pa\Omega)}<\kappa$.
Dann existiert ein stetiger, linearer Operator 
\[\trnormal: E^p(\Omega_\eta)\rightarrow (W^{1,p'}(\pa\Omega))'\]
derart, dass f"ur
$\bphi\in E^p(\Omega_\eta)$ und $\psi\in C^1(\overline{\Omega_\eta})$ gilt
\[\int_{\Omega_\eta} \bphi\cdot\nabla\psi\ dx=-\int_{\Omega_\eta} \dv\bphi\
 \psi\ dx + \langle \trnormal\bphi,\tr\psi\rangle_{W^{1,p'}(\pa\Omega)}.\]
Die Stetigkeitskonstante h"angt nur von $\Omega$, $\tau(\eta)$ und $p$ ab; sie bleibt beschr"ankt,
falls $\tau(\eta)$ beschr"ankt bleibt.
\end{proposition}
\beweis Es gen"ugt $\bphi\in C^1(\overline{\Omega_\eta})$ zu betrachten, da diese Funktionen dicht
in $E^p(\Omega_\eta)$ liegen; siehe Proposition \ref{theorem:dicht}. Wir erhalten analog zum
Vorgehen im Beweis
von Proposition \ref{lemma:partInt} f"ur $\psi\in C^1(\overline\Omega_\eta)$ die Identit"at
\begin{equation}\label{eqn:approxnspur}
 \begin{aligned}
\int_{\Omega_\eta} \bphi\cdot\nabla\psi\ dx+\int_{\Omega_\eta} \dv\bphi\
 \psi\ dx = \int_{\pa\Omega} \tr\bphi\cdot\bv_\eta\,\tr\psi\ dA.  
 \end{aligned}
\end{equation}
Die linke Seite ist durch $2\norm{\bphi}_{E^p
(\Omega_\eta)}\norm{\psi}_{W^{1,p'}(\Omega_\eta)}$ abgesch"atzt. Offenbar gilt die Gleichung sogar
f"ur $\psi\in W^{1,p'}(\Omega_\eta)$. 

Die Zuordnung $b\mapsto b\circ q$ definiert einen stetigen
Fortsetzungsoperator von $W^{1,p'}(\pa\Omega)$ nach $W^{1,p'}(S_\alpha\cap\Omega_\eta)$ f"ur ein
festes $\alpha$ mit $\norm{\eta}_{L^\infty(\pa\Omega)}<\alpha<\kappa$. Der Transformationssatz,
Gleichung \eqref{eqn:trafo}, gibt uns n"amlich
\[\int_{S_\alpha\cap\Omega_\eta}|b\circ q|^{p'}\ dx=\int_{\pa\Omega}|b|^{p'}\int_{-\alpha}^\eta|\det
d\Lambda|\ dsdA,\]
sodass die $L^{p'}$-Norm von $b\circ q$ in Abh"angigkeit von $\Omega$ und $\alpha$ durch die
$L^{p'}$-Norm von $b$ abgesch"atzt ist. Durch Approximation von $b$ durch hinreichend glatte
Funktionen
erhalten
wir die Kettenregel
\[\nabla(b\circ q)=(\nabla b)\circ q\ dq,\]
sodass auch die $L^{p'}$-Norm von $\nabla(b\circ q)$ in Abh"angigkeit von $\Omega$, $\alpha$
und
$p'$
durch die $L^{p'}$-Norm von $\nabla b$ abgesch"atzt ist. Ebenfalls durch Approximation von $b$
erhalten wir die Identit"at $\tr (b\circ q)=b$. Durch Multiplikation der konstruierten Fortsetzungen
mit der Abschneidefunktion $x\mapsto\beta(|s(x)|)$, wobei
$\beta\in C^\infty(\setR)$, $\beta=1$ in einer Umgebung von $[0,\norm{\eta}_{L^\infty(\pa\Omega)}]$
und $\beta=0$ in einer Umgebung von $\alpha$, erhalten wir schlie"slich
% Wenn
% wir anschlie"send noch die Spur $b\circ \Phi_\delta$, $\delta\equiv\alpha$, dieser Fortsetzung auf
% $\pa(\Omega\setminus
% \overline{S_\alpha})$ mit Hilfe eines Fortsetzungsoperators f"ur regul"are R"ander, siehe Theorem
% \ref{theorem:spur}, zu
% einer $W^{1,p'}(\Omega\setminus\overline{S_\alpha})$-Funktion fortsetzen, so haben wir insgesamt
einen stetigen, linearen Fortsetzungsoperator von $W^{1,p'}(\pa\Omega)$ nach
$W^{1,p'}(\Omega_\eta)$. Damit ist
klar, dass die rechte Seite in obiger Identit"at einen Spuroperator mit den behaupteten
Eigenschaften definiert.
\qed\\

\begin{proposition}\label{lemma:divdicht}
Sei $\eta\in H^2(\pa\Omega)$ mit $\norm{\eta}_{L^\infty(\pa\Omega)}<\kappa$. Dann liegt der
Teilraum der Funktionen mit in $\Omega_\eta$ enthaltenen Tr"agern dicht im kanonisch normierten
Raum
\[H(\Omega_\eta):=\{\bphi\in L^2(\Omega_\eta)\ |\ \dv\bphi=0,\ \trnormal\bphi =0 \}.\]
\end{proposition}
\beweis F"ur $\bphi\in H(\Omega_\eta)$ gelte 
\begin{equation*}
 \begin{aligned}
  \int_{\Omega_\eta} \bphi\cdot\bpsi\ dx =0
 \end{aligned}
\end{equation*}
f"ur alle $\bpsi\in L^2(\Omega_\eta)$ mit $\dv\bpsi=0$ und $\supp \bpsi\subset\Omega_\eta$. Gem"a"s
Theorem \ref{theorem:derham} existiert eine Funktion $p\in L^2_{\loc}(\Omega_\eta)$ mit
$\bphi=\nabla p$. Aus Proposition \ref{lemma:nspur} folgt zudem
\begin{equation*}
 \begin{aligned}
  \int_{\Omega_\eta} \bphi\cdot\nabla \psi\ dx =0 
 \end{aligned}
\end{equation*}
f"ur alle $\psi\in C^1(\overline{\Omega_\eta})$. Der Beweis von Proposition \ref{theorem:dicht}
zeigt, dass $\nabla p$ durch Funktionen $\nabla\psi$, $\psi\in C^1(\overline{\Omega_\eta})$, in
$L^2(\Omega_\eta)$ approximierbar ist, woraus wir
\begin{equation*}
 \begin{aligned}
  \int_{\Omega_\eta} |\bphi|^2\ dx = \int_{\Omega_\eta} \bphi\cdot\nabla p\ dx = 0 
 \end{aligned}
\end{equation*} 
folgern.
\qed\\

Nun konstruieren wir einen Operator, der geeignete Randwerte $b\,\bnu$ zu einem divergenzfreien
Vektorfeld fortsetzt. Dazu skalieren wir das Feld $(b\,\bnu)\circ q$ derart, dass es divergenzfrei
wird. Ein entprechender Ansatz f"uhrt auf eine gew"ohnliche Differentialgleichung und deren L"osung
auf Definition \eqref{eqn:fort} unten.
\begin{proposition}\label{lemma:FortVonRand}
Es seien $1<p<\infty$, $\eta\in H^2(\pa\Omega)$ mit $\norm{\eta}_{L^\infty(\pa\Omega)}<\kappa$ und
$\alpha$ eine Zahl mit $\norm{\eta}_{L^\infty(\pa\Omega)}<\alpha<\kappa$. Dann existiert ein
stetiger, linearer Fortsetzungsoperator
\[\F_\eta: \Big\{b\in W^{1,p}(\pa\Omega)\ \big|\ \int_{\pa\Omega}b\,\gamma(\eta)\
dA=0\Big\}\rightarrow W^{1,p}_{\dv}(B_\alpha);\]
insbesondere gilt $\tr\F_\eta b=b\,\bnu.$ Die Stetigkeitskonstante h"angt nur von $\Omega$,
$\norm{\eta}_{H^2(\pa\Omega)}$, $\alpha$ und $p$ ab; sie bleibt beschr"ankt, falls
$\norm{\eta}_{H^2(\pa\Omega)}$ und $\tau(\alpha)$ beschr"ankt bleiben.
\end{proposition}
\beweis Sei $b\in W^{1,p}(\pa\Omega)$ mit 
\begin{equation}\label{eqn:integralnull}
 \begin{aligned}
\int_{\pa\Omega}b\,\gamma(\eta)\ dA=0.  
 \end{aligned}
\end{equation}
F"ur $x\in S_\alpha$ setzen wir 
\begin{equation}\label{eqn:fort}
 \begin{aligned}
(\F_\eta b)(x):=\exp\Big(\int_{\eta(q(x))}^{s(x)}
\beta(q(x)+\tau\,\bnu(q(x)))\ d\tau\Big)\,
(b\,\bnu)(q(x)),  
 \end{aligned}
\end{equation}
wobei $\beta:=-\dv(\bnu\circ q)$. Offenbar ist $\beta$ eine
$C^2$-Funktion. Mit Hilfe der Kettenregel sehen
wir, dass $\F_\eta b$ schwach differenzierbar ist und dass 
\begin{equation}\label{eqn:fortnabla}
 \begin{aligned}
  \pa_i\, \F_\eta b = \Big[\pa_i((b\,\bnu)\circ q) +
(b\,\bnu)\circ q\ \Big(&\int_{\eta\circ q}^s \pa_i(\beta(q+\tau\,\bnu\circ q))\ d\tau +
\beta(q+s\,\bnu\circ q)\ \pa_i s\\
 &- \beta(q+(\eta\,\bnu)\circ q)\ \pa_i (\eta\circ q)\Big)\Big]\,e^{\int_{\eta\circ q}^s
\beta(q+\tau\bnu\circ q) d\tau}
 \end{aligned}
\end{equation}
gilt. Da $\pa_i((b\,\bnu)\circ q)\in L^p(S_\alpha)$, $(b\,\bnu)\circ q\in L^{r}(S_\alpha)$ f"ur ein
$r>p$, $\pa_i (\eta\circ q)\in L^r(S_\alpha)$ f"ur alle $r<\infty$ und alle anderen Terme
beschr"ankt sind, ist $\F_\eta b$ in $W^{1,p}(S_\alpha)$ abgesch"atzt. Desweiteren
gilt\footnote{F"ur ein skalares
Feld $f$ und ein Vektorfeld $\bX$: $\dv (f\,\bX)=f\dv\bX +
\nabla f\cdot\bX$.}  
\begin{equation*}
 \begin{aligned}
  \dv\F_\eta b&=\dv (\bnu\circ q)\, e^{\int_{\eta\circ q}^s
\beta(q+\tau\bnu\circ q) d\tau}\,b\circ q + (\bnu\circ q)\cdot \nabla (e^{\int_{\eta\circ q}^s
\beta(q+\tau\bnu\circ q) d\tau}\,b\circ q)\\
&= (-\beta\, e^{\int_{\eta\circ q}^s\beta(q+\tau\bnu\circ q) d\tau} + \pa_s e^{\int_{\eta\circ q}^s
\beta(q+\tau\bnu\circ q) d\tau})\, b\circ q
=0.
\end{aligned}
\end{equation*}
F"ur das zweite Gleichheitszeichen haben wir die Definition von $\beta$ und die Tatsache,
dass f"ur $x\in S_\alpha$ 
\[dq\, \bnu(q(x))=\frac{d}{dt}\Big|_{t=0}q(x+t\,\bnu(q(x)))=\frac{d}{dt}\Big|_{t=0}q(x)=0\]
und 
\[ds\, \bnu(q(x))=\frac{d}{dt}\Big|_{t=0}s(x+t\,\bnu(q(x)))=\frac{d}{dt}\Big|_{t=0}s(x)+t=1\]
gilt, verwendet.

Approximieren wir $\eta$ und $b$ durch $C^2$-Funktionen $(\eta_n)$ in
$H^2(\pa\Omega)$ bzw. durch $C^1$-Funktionen $(b_n)$ in $W^{1,p}(\pa\Omega)$, so konvergieren die
$C^1(\overline{S_\alpha})$-Funktionen
\[\bphi_n:=\exp\Big(\int_{\eta_n\circ q}^{s}
\beta(q+\tau\,\bnu\circ q)\
d\tau\Big)\, (b_n\,\bnu)\circ q\]
gegen $\F_\eta b$ in $W^{1,p}(S_\alpha)$. Zudem konvergieren die Spuren
\begin{equation}\label{eqn:approxspur}
 \begin{aligned}
\tr\bphi_n=\exp\Big(\int_{\eta_n}^{\eta}
\beta(q+\tau\,\bnu\circ q)\
d\tau\Big)\,b_n\,\bnu  
 \end{aligned}
\end{equation}
gegen $b\,\bnu$, woraus wir $\tr\F_\eta b=b\,\bnu$ folgern.

Es bleibt noch, $\F_\eta b$ auf $\Omega\setminus\overline{S_\alpha}$ fortzusetzen. Unter Verwendung
von
Proposition \ref{lemma:partInt} folgern wir, dass
\begin{equation*}
 \begin{aligned}
\int_{\pa(\Omega\setminus\overline{S_\alpha})}(\F_\eta b)\cdot\bnu\
dA=-\int_{\pa\Omega}b\, \gamma(\eta)\ dA+\int_{S_\alpha\cap\Omega_\eta} \dv\F_\eta b\ dx=0,
 \end{aligned}
\end{equation*}
wobei $\bnu$ hier die innere Einheitsnormale an
$\pa(\Omega\setminus\overline{S_\alpha})$ bezeichnet. Da
$\pa(\Omega\setminus\overline{S_\alpha})$ ein $C^3$-Rand ist, k"onnen wir aufgrund dieser
Identit"at 
das Stokes-System in $\Omega\setminus \overline{S_\alpha}$ mit Randwerten
$(\F_\eta b)|_{\pa(\Omega\setminus\overline{S_\alpha})}$ l"osen;
siehe Theorem \ref{theorem:stokes}. Dies liefert uns die gew"unschte
Fortsetzung.\qed\\

\begin{bemerkung}\label{bem:nspur}
Die Funktion $b\in L^p(\pa\Omega)$, $1<p<\infty$, erf"ulle die Identit"at
\eqref{eqn:integralnull}. Dann k"onnen wir eine Fortsetzung $\F_\eta b\in L^p(B_\alpha)$ mit
verschwindender (distributioneller) Divergenz wie im Beweis der obigen Proposition konstruieren. Die
lineare Abbildung 
\[\F_\eta: \Big\{b\in L^{p}(\pa\Omega)\ \big|\ \int_{\pa\Omega}b\, \gamma(\eta)\
dA=0\Big\}\rightarrow
\{\bphi\in L^p(B_\alpha)\ |\ \dv\bphi=0\}\]
definiert wiederum einen stetigen, linearen Operator, dessen Stetigkeitskonstante wie in
Proposition \ref{lemma:FortVonRand} von den Daten abh"angt. 

Bei der Konstruktion der Fortsetzung wenden wir den L"osungsoperator des Stokes-Systems auf die
formale Spur 
\[\bphi|_{\pa(\Omega\setminus\overline{S_\alpha})}=\exp\Big(\int_{\eta\circ
q}^{-\alpha}
\beta(q+\tau\,\bnu\circ q)\ d\tau\Big)\, (b\,\bnu)\circ q\]
an. Diese ist offenbar in $L^p(\pa(\Omega\setminus\overline{S_\alpha}))$ abgesch"atzt und erf"ullt
die
Identit"at
\begin{equation}\label{eqn:formalmitt}
 \begin{aligned}
 \int_{\pa(\Omega\setminus\overline{S_\alpha})}\bphi\cdot\bnu\
dA=0.
 \end{aligned}
\end{equation}
Zum anschlie"senden Nachweis der Divergenzfreiheit des resultierenden Vektorfelds auf $B_\alpha$
verwenden wir die Tatsache, dass f"ur $\psi\in C_0^\infty(\Omega)$ die Gleichung
\begin{equation}\label{eqn:horst}
 \begin{aligned}
 \int_{S_\alpha}\bphi\cdot\nabla\psi\
dx=\int_{\pa(\Omega\setminus\overline{S_\alpha})}\bphi\cdot\bnu\, \psi\
dA=-\int_{\Omega\setminus\overline{S_\alpha}}\bphi\cdot\nabla\psi\
dx
 \end{aligned}
\end{equation}
gilt. Die Identit"aten \eqref{eqn:formalmitt} und \eqref{eqn:horst} lassen sich leicht durch
Approximation von $b$ durch hinreichend glatte Funktionen einsehen. Ebenfalls durch Approximation
von $b$ erhalten wir aus den Gleichungen
\eqref{eqn:approxspur} und
\eqref{eqn:approxnspur} die
Identit"at 
\begin{equation}\label{eqn:nspur}
 \begin{aligned}
\trnormal \bphi=b\, \gamma(\eta),
 \end{aligned}
\end{equation}
die die Bezeichnung \emph{Fortsetzungsoperator} rechtfertigt.
\end{bemerkung}
% \begin{lemma}
% Es seien $1<p<\infty$, $\eta\in H^2(\pa\Omega)$ mit $\norm{\eta}_{L^\infty(\pa\Omega)}<\kappa$
% sowie
% $\Gamma\subset\pa\Omega$
% offen und nichtleer. Dann existiert eine Konstante
% $c=c(\Omega,\Gamma,\norm{\eta}_{H^2(\pa\Omega)},p)$ derart, dass f"ur $\bphi\in
% W^{1,p}(\Omega_\eta)$ mit $\tr\bphi=0$ auf $\Gamma$ gilt
% \[\norm{\nabla\bphi}_{L^p(\Omega_\eta)}\le c\norm{D\bphi}_{L^p(\Omega_\eta)}.\]
% Dabei ist $D\bphi=\frac{1}{2}(\nabla\bphi+(\nabla\bphi)^T)$ der symmetrische Anteil des
% Gradienten.
% \end{lemma}
% \beweis \marginpar{NOCH UNKLAR! GEN"UGT F"UR DIVFREIE! LARS FRAGEN! F"UR P=2 WOHL NICHT
% N"OTIG!}\qed\\
\vspace{0.2cm}

\begin{proposition}\label{lemma:nullfort}
Es seien $1<p<3$ und $\eta\in H^2(\pa\Omega)$ mit
$\norm{\eta}_{L^\infty(\pa\Omega)}<\kappa$.
Dann ist die Fortsetzung durch $0$ eine stetige,
lineare Abbildung von $W^{1,p}(\Omega_\eta)$ nach $W^{1/4,p}(\setR^3)$. Die Stetigkeitskonstante
h"angt nur von $\Omega$, $p$, $\norm{\eta}_{H^2(\pa\Omega)}$ und $\tau(\eta)$ ab; sie bleibt
beschr"ankt, falls $\norm{\eta}_{H^2(\pa\Omega)}$ und $\tau(\eta)$ beschr"ankt bleiben.
% Konvergiert die Folge
% $(\eta_n)$ schwach gegen $\eta$ in $H^2(\pa\Omega)$, so gilt f"ur $v\in W^{1,p}(\Omega_\eta)$
% \[v\circ\Psi_{\eta_n}\rightarrow v\circ\Psi\text{ in } W^{1,r}(\Omega).\]
\end{proposition}
\beweis Wir zeigen zun"achst, dass die Fortsetzung durch $0$ eine stetige, lineare Abbildung von
$W^{1,r}(\Omega)$ nach $W^{s,r}(\setR^3)$ mit $1\le r<3$ und $s<1/3$ ist. Dazu gen"ugt es f"ur $v\in
W^{1,r}(\Omega)$ das Integral
\begin{equation*}
 \begin{aligned}
  \int_{\setR^3}\int_{\setR^3} \frac{|v(x)-v(y)|^r}{|x-y|^{3+sr}}\ dydx =
\int_{\Omega}\int_{\Omega} \frac{|v(x)-v(y)|^r}{|x-y|^{3+sr}}\ dydx +2
\int_{\Omega}\int_{\setR^3\setminus\Omega} \frac{|v(x)|^r}{|x-y|^{3+sr}}\ dydx\\
=\int_{\Omega}\int_{\Omega} \frac{|v(x)-v(y)|^r}{|x-y|^{3+sr}}\ dydx+ 2\int_{\Omega} |v(x)|^r
\int_{\setR^3\setminus\Omega} \frac{1}{|x-y|^{3+sr}}\ dydx
 \end{aligned}
\end{equation*}
abzusch"atzen. W"ahrend der erste Term auf der rechten Seite durch
$c\,\norm{v}_{W^{1,r}(\Omega)}$
majorisiert ist, l"asst sich das innere Integral des zweiten Terms durch
\[\int_{|z|>d(x)}\frac{1}{|z|^{3+sr}}\ dz=\frac{c(s,r)}{d(x)^{sr}}\]
absch"atzen, wobei $d(x)$ den Abstand von $x$ zum Rand $\pa\Omega$ bezeichnet. Eine Anwendung der
H"older-Ungleichung mit den Exponenten $r^*/r$ und $(r^*/r)'$  auf den zweiten Summanden zeigt
somit, dass dieser durch $c(s,r)\,\norm{v}_{r^*}^r\,\norm{d(\cdot)^{-s}}_{L^3(\Omega)}^r$ dominiert
ist. Die Identit"at
\[\int_{S_{\kappa/2}\cap\Omega}|d(x)|^{-3s}\
dx=\int_{\pa\Omega}\int_{-\kappa/2}^0|\det
d\Lambda|\, \alpha^{-3s} \ d\alpha dA,\]
eine Konsequenz des Transformationssatzes, zeigt zusammen mit der Ungleichung $3s<1$, dass der
Faktor $\norm{d(\cdot)^{-s}}_{L^3(\Omega)}^r$ endlich ist. Damit ist die
Behauptung bewiesen.

Verketten der Abbildung aus Lemma \ref{lemma:psi} mit obiger Fortsetzung ergibt nun eine
stetige, lineare Abbildung
\[W^{1,p}(\Omega_\eta)\rightarrow W^{1,r}(\Omega)\rightarrow W^{s,r}(\setR^3)\]
f"ur $r<p$ und $s<1/3$. Sei nun $\delta\in C^4(\pa\Omega)$ mit $\eta<\delta<\kappa$
in $\pa\Omega$. Setzen wir die Abbildung $\Psi_\delta$ auf $\overline{B_\alpha}$, $\alpha>0$
hinreichend klein, fort, indem wir die Definition \eqref{eqn:defpsi} $x\in S_\kappa\cap
\overline{B_\alpha}$ verwenden, so erhalten wir einen $C^3$-Diffeomorphismus
\[\widetilde\Psi_\delta: \overline{B_\alpha}\rightarrow \overline{\Omega_{\delta+\alpha}}.\] 
Da die gebrochenen Sobolev-R"aume Interpolationsr"aume sind, ist die lineare
Abbildung $v\mapsto v\circ\widetilde\Psi_\delta^{-1}$ stetig von $W^{s,r}(B_\alpha)$
nach $W^{s,r}(\Omega_{\delta+\alpha})$. Aus demselben Grund ist die
Einschr"ankung von Funktionen stetig von $W^{s,r}(\setR^3)$ nach
$W^{s,r}(B_\alpha)$. F"ur geeignete $s<1/3$ und $r<p$ gilt zudem die Einbettung 
\[W^{s,r}(\Omega_{\delta+\alpha})\embedding W^{1/4,p}(\Omega_{\delta+\alpha}),\]
siehe Theorem \ref{theorem:emb}. Verkn"upfen der obigen Abbildungen zeigt, dass die
Fortsetzung durch $0$ eine stetige, lineare Abbildung von $W^{1,p}(\Omega_\eta)$ nach
$W^{1/4,p}(\Omega_{\delta+\alpha})$ ist. Da $\Omega_\eta$ positiven
Abstand von $\setR^3\setminus\Omega_{\delta+\alpha}$ hat, folgt die Behauptung der Proposition
nun aus der Absch"atzung
\begin{equation*}
 \begin{aligned}
|v|_{1/4,p;\setR^3}^p&=|v|_{1/4,p;\Omega_{\delta+\alpha}}^p+\int_{\Omega_\eta}\int_{
\setR^3\setminus\Omega_{\delta+\alpha}}
\frac{|v(x)|^p}{|x-y|^{3+p/4}}\ dydx\\
&\le |v|_{1/4,p;\Omega_{\delta+\alpha}}^p + c\,\norm{v}_{L^p(\Omega_\eta)}^p.
 \end{aligned}
\end{equation*}
\qed\\

F"ur Funktionen, die auf zeitlich variablen Gebieten definiert sind, sind die "ublichen
Bochner-R"aume nicht die richtigen Objekte. Deshalb geben wir nun eine naheliegende, auf unsere
Bed"urfnisse zugeschnittene Konstruktion an. F"ur $I:=(0,T)$, $T>0$, und $\eta\in C(\bar
I\times\pa\Omega)$ mit
$\norm{\eta}_{L^\infty(I\times\pa\Omega)}<\kappa$ setzen wir $\Omega_\eta^I:=\bigcup_{t\in
I}\, \{t\}\times \Omega_{\eta(t)}$.
$\Omega_\eta^I$ ist ein Gebiet des $\setR^4$. Wir definieren nun f"ur $1\le p,r\le\infty$
\begin{equation*}
 \begin{aligned}
 L^p(I,L^r(\oet))&:=\{v\in L^1(\Omega_\eta^I)\ |\
v(t,\cdot)\in
L^r(\Omega_{\eta(t)})\text{ f"ur fast alle $t$ und}\\
&\hspace{4.9cm}\norm{v(t,\cdot)}_{L^r(\Omega_{\eta(t)})}\in
L^p(I)\},\\
L^p(I,W^{1,r}(\Omega_{\eta(t)}))&:=\{v\in L^p(I,L^r(\oet))\ |\ \nabla v\in
L^p(I,L^r(\oet))\},\\
L^p(I,W_{\dv}^{1,r}(\Omega_{\eta(t)}))&:=\{\bv\in L^p(I,W^{1,r}(\Omega_{\eta(t)}))\ |\
\dv \bv=0\},\\
W^{1,p}(I,W^{1,r}(\Omega_{\eta(t)}))&:=\{\bv\in L^p(I,W^{1,r}(\Omega_{\eta(t)}))\ |\
\pa_t\bv\in L^p(I,W^{1,r}(\Omega_{\eta(t)}))\}.
 \end{aligned}
\end{equation*}
Dabei wirken $\nabla$ und $\dv$ nur auf die r"aumlichen Koordinaten. Wir setzen zudem
\[\Psi_\eta: \bar I\times \overline\Omega\rightarrow \overline{\Omega_\eta^I},\ (t,x)\mapsto
(t,\Psi_{\eta(t)}(x))\]
und
\[\Phi_\eta: \bar I\times \pa\Omega\rightarrow\bigcup_{t\in
\bar I}\, \{t\}\times \pa\Omega_{\eta(t)},\ (t,x)\mapsto
(t,\Phi_{\eta(t)}(x)).\]
Gilt zus"atzlich $\eta\in L^\infty(I,H^2(\pa\Omega))$, so erhalten wir "`instation"are"' Versionen
der obigen Aussagen, wenn wir diese zu (fast) jedem Zeitpunkt verwenden. Zum Beispiel folgt
aus Korollar \ref{lemma:sobolev} die Einbettung  
\[L^2(I,H^{1}(\Omega_{\eta(t)}))\embedding L^2(I,L^s(\Omega_{\eta(t)}))\]
f"ur $1\le s<2^*$. 

Man beachte, dass obige Konstruktion nicht ohne Weiteres einen Ersatz f"ur
Bochner-R"aume von
Funktionen mit Werten in einem Dualraum bietet. Solche R"aume spielen in der Theorie schwacher
L"osungen von Evolutionsgleichungen eine wichtige Rolle, weil den Zeitableitungen solcher L"osungen
typischerweise nur als Elemente von R"aumen dieser Art Sinn gegeben werden kann. Funktionale aus
$(W^{1,r}_0)'$ lassen sich in der Form $f+\pa_i f^i$ mit $f,f^i\in L^{r'}$ schreiben;
siehe zum Beispiel \cite{b4}. Liegt eine derartige Charakterisierung des fraglichen Dualraums vor,
so kann man versuchen, darauf aufbauend einen Ersatz nach obigem Schema zu
entwickeln. Dennoch dr"angt sich die Frage auf, ob nicht ein nat"urlichere Konstruktion denkbar
ist. Man k"onnte auf die Idee kommen, R"aume von
Funktionen, die jeden
Zeitpunkt in
einen anderen Banach-Raum abbilden, zu betrachten. Objekte dieser Art treten in der
Differentialgeometrie auf, und zwar in Form von Vektorb"undeln. Ein Vektorb"undel ordnet
jedem Punkt einer Mannigfaltigkeit einen anderen Vektorraum zu. Um differenzierbaren Abbildungen
von der Mannigfaltigkeit in das B"undel Sinn zu geben, ben"otigt letzteres nat"urlich eine
differenzierbare Struktur. In unserem Falle m"usste die differenzierbare
Struktur\footnote{bzw. die lokalen Trivialisierungen des B"undels, die "ublicherweise die
differenzierbare Struktur induzieren} so gew"ahlt werden, dass sie den Sinn
der Zeitableitung, n"amlich in
einem inneren Punkt des sich bewegenden Gebiets die zeitliche "Anderung der Funktion anzugeben, 
erfasst. Die Konstruktion einer derartigen Struktur ist bislang nicht gelungen. Eine verwandte
M"oglichkeit, die Problematik anzugehen, besteht auf den ersten Blick darin, das
Gebiet auf einen Raumzeitzylinder zu transformieren und auf diesem die
"ublichen Bochner-R"aume zu betrachten. Das scheitert in unserem Fall an der geringen Regularit"at,
die unser variabler Rand haben wird. Zudem verf"alscht eine solche Transformation den Sinn der
Zeitableitung. Prinzipiell scheint die Aussage zu gelten, dass durch eine Transformation auf einen
Raumzeitzylinder nichts gewonnen wird, sondern dass im Gegenteil die Sache komplizierter wird.
Die obige Konstruktion scheint somit das einzig sinnvolle Substitut f"ur die "ublichen
Bochner-R"aume
zu sein. Die Frage nach einer M"oglichkeit, den Verlust des "ublichen Begriffs der schwachen
Zeitableitung zu kompensieren, wird uns deshalb noch besch"aftigen.
% Diese wird "ublicherweise durch lokale Trivialisierungen induziert. Das
% sind Abbildungen, die das B"undel lokal auf das Kreuzprodukt einer offenen Teilmenge der
% Mannigfaltigkeit mit der sogenannten typischen Faser, einem festen Vektorraum, abbilden. 
% ergibt sich durch
% Anwenden des obigen Spuroperators zu jedem Zeitpunkt ein stetiger Operator $\tr$ von
% $L^p(I,W^{1,p}(\Omega_\eta(t)))$ nach $L^p(I,W^{1-1/r,r}(\pa\Omega))$ f"ur alle $r<p$.
% \marginpar{PR"UFEN!}

F"ur alle $1/2<\theta<1$ gilt 
\begin{equation}\label{eqn:hoeldereinb}
 \begin{aligned}
  W^{1,\infty}(I,L^2(\pa\Omega))\cap L^\infty(I,H^2(\pa\Omega)) \embedding C^{0,1-\theta}(\bar
I,H^{2\theta}(\pa\Omega))\embedding C^{0,1-\theta}(\bar I, C^{0,2\theta -1}(\pa\Omega)).
 \end{aligned}
\end{equation}
W"ahrend wir f"ur die zweite Einbettung Theorem \ref{theorem:emb} verwendet haben, folgt die erste
aus der elementaren Absch"atzung
\begin{equation*}
 \begin{aligned}
\norm{u(t)-u(s)}_{(L^2(\pa\Omega),H^2(\pa\Omega))_{\theta,2}}&\le\norm{u(t)-u(s)}
^{\theta}_{H^2(\pa\Omega)}\,\norm{u(t)-u(s)}_{L^2(\pa\Omega)}^{1-\theta}\\ &\le c
\norm{u}_{L^\infty(I,H^2(\pa\Omega))}^{\theta}\,\norm{u}_{W^{1,\infty}(I,L^2(\pa\Omega))}^{1-\theta}
\,|t-s|^{1-\theta}.
 \end{aligned}
\end{equation*}
\vspace{0.2cm}

\begin{proposition}\label{lemma:FortVonRandZeit}
Es seien $\eta\in W^{1,\infty}(I,L^2(\pa\Omega))\cap L^\infty(I,H^2(\pa\Omega))$ mit
$\norm{\eta}_{L^\infty(I\times\pa\Omega)}<\kappa$ und
$\alpha$ eine Zahl mit $\norm{\eta}_{L^\infty(I\times\pa\Omega)}<\alpha<\kappa$. Dann definiert die
Anwendung des Fortsetzungsoperators aus Proposition \ref{lemma:FortVonRand} zu (fast) allen Zeiten
einen stetigen, linearen Fortsetzungsoperator $\F_\eta$ von
\[\Big\{b\in H^{1}(I,L^2(\pa\Omega))\cap L^2(I,H^2(\pa\Omega))\ |\
\int_{\pa\Omega}b(t,\cdot)\, \gamma(\eta(t,\cdot))\
dA=0\text{ f"ur fast alle }t\in I\Big\}\]
nach
\[\{\bphi\in H^{1}(I,L^2(B_\alpha))\cap C(\bar I,H^1(B_\alpha))\ |\ \dv\bphi=0\}.\]
Die Stetigkeitskonstante h"angt nur von $\Omega$,
$\norm{\eta}_{ W^{1,\infty}(I,L^2(\pa\Omega))\cap L^\infty(I,H^2(\Omega))}$ und $\alpha$ ab;
sie bleibt beschr"ankt, falls
$\norm{\eta}_{ W^{1,\infty}(I,L^2(\pa\Omega))\cap L^\infty(I,H^2(\Omega))}$ und $\tau(\alpha)$
beschr"ankt bleiben.
\end{proposition}
\beweis Die Stetigkeit nach $\{\bphi\in L^\infty(I,H^1(B_\alpha))\ |\ \dv\bphi=0\}$ folgt unter
Beachtung der Einbettung
\begin{equation}\label{eqn:einbschnurz}
 \begin{aligned}
  H^1(I,L^2(\pa\Omega))\cap L^2(I,H^2_0(\pa\Omega))\embedding C(\bar I,H^1(\pa\Omega)),
 \end{aligned}
\end{equation}
eine Konsequenz von Proposition \ref{theorem:bochneremb},\footnote{Die
Identit"at
$[L^2(\pa\Omega),H^2(\pa\Omega)]_\frac{1}{2}=H^1(\pa\Omega)$ folgt dabei unter Verwendung eines
endlichen
Atlas mit
untergeordneter Zerlegung der Eins und der in \cite{b4} konstruierten
Fortsetzungsoperatoren aus Theorem 6.4.5 in \cite{b48}.} sofort aus Proposition
\ref{lemma:FortVonRand}. Aus \eqref{eqn:einbschnurz}, \eqref{eqn:hoeldereinb} und
Theorem \ref{theorem:emb} erhalten wir zudem die Einbettungen
\begin{equation*}
 \begin{aligned}
 H^1(I,L^2(\pa\Omega))\cap L^2(I,H^2_0(\pa\Omega))\embedding C(\bar I,L^4(\pa\Omega)),\\
 W^{1,\infty}(I,L^2(\pa\Omega))\cap L^\infty(I,H^2(\pa\Omega)) \embedding C(\bar
I,W^{1,4}(\pa\Omega)).
 \end{aligned}
\end{equation*}
Diese zeigen zusammen mit \eqref{eqn:fort}, \eqref{eqn:fortnabla} und den Stetigkeitseigenschaften
des L"osungsoperators des Stokes-Systems, dass $\F_\eta b$ in $C(\bar I,H^1(B_\alpha))$ liegt. Es
bleibt also lediglich $\F_\eta b$ in $H^1(I,L^2(B_\alpha))$ abzusch"atzen. In
$I\times S_\alpha$ gilt
\begin{equation}\label{eqn:klopps}
 \begin{aligned}
\pa_t\, (\F_\eta b)(t,\cdot):= \Big[(\pa_tb)(t,q)\ \bnu\circ q - \beta(q+\eta(t,q)\, \bnu\circ q)\ 
\pa_t\eta(t,q) \  b(t,q)\  \bnu\circ q\Big]\\
e^{\int_{\eta\circ q}^s
\beta(q+\tau\bnu\circ q) d\tau}.
\end{aligned}
\end{equation}
Wegen $\pa_tb\in L^2(I,L^2(\pa\Omega))$ ist der erste Summand in $L^2(I,L^2(S_\alpha))$
abgesch"atzt; die Absch"atzung des zweiten Summanden in $L^2(I,L^2(S_\alpha))$ folgt aus
$\pa_t\eta\in
L^\infty(I,L^2(\pa\Omega))$ und $b\in L^2(I,L^\infty(\pa\Omega))$. Die Zeitableitung der Spur von
$\F_\eta b$ auf $I\times\pa(\Omega\setminus\overline{S_\alpha})$ ist identisch
\[\Big[(\pa_tb)(t,q)\ \bnu\circ q - \beta(q+\eta(t,q)\,\bnu\circ q)\
\pa_t\eta(t,q) \ b(t,q)\ \bnu\circ q\Big]\,e^{\int_{\eta\circ q}^{-\alpha}
\beta(q+\tau\bnu\circ q) d\tau}\]
und somit in $L^2(I,L^2(\pa(\Omega\setminus\overline{S_\alpha})))$ abgesch"atzt. Aus den
Stetigkeitseigenschaften des L"osungsoperators des Stokes-Systems folgt nun die Behauptung.
\qed\\

\begin{bemerkung}
Die Argumente in obigem Beweis zeigen, dass die Anwendung des Fortsetzungsoperators aus
Bemerkung \ref{bem:nspur} zu (fast) allen Zeiten unter den Voraussetzungen von Proposition
\ref{lemma:FortVonRandZeit} einen stetigen, linearen Fortsetzungsoperator $\F_\eta$ von
\[\Big\{b\in C(\bar I,L^2(\pa\Omega))\ |\
\int_{\pa\Omega}b(t,\cdot)\, \gamma(\eta(t,\cdot))\
dA=0\text{ f"ur fast alle }t\in I\Big\}\]
nach
\[\{\bphi\in C(\bar I,L^2(B_\alpha))\ |\ \dv\bphi=0\}\]
definiert. Die Stetigkeitskonstante h"angt wie in Proposition \ref{lemma:FortVonRandZeit} von den
Daten ab.
\end{bemerkung}

\chapter{Koiter-Energie}

Es sei $M$ eine kompakte, orientierte und berandete $C^4$-Fl"ache im
$\setR^3$ mit durch den umgebenden euklidischen Raum induzierter erster und zweiter
Fundamentalform $g$ bzw. $h$ und Fl"achenma"s $dA$. $M$ stelle die
Mittelfl"ache einer elastischen Schale der Dicke $2\epsilon_0>0$ in ihrer Ruhelage dar, wobei
$\epsilon_0$ klein sei gegen die Kehrwerte der Hauptkr"ummungen von $M$. Die
elastische Schale bestehe aus einem homogenen und isotropen Material, dessen lineares elastisches
Verhalten durch die
Lam\'{e}-Konstanten $\lambda$ und $\mu$ charakterisiert sei. Deformationen der Mittelfl"{a}che
und damit der Schale beschreiben wir durch ein (hinreichend glattes) Vektorfeld $\boldsymbol{\eta}:\
M \rightarrow
\setR^3$. Wir bezeichnen mit $g(\boldsymbol{\eta})$ und $h(\boldsymbol{\eta})$ die auf $M$
zur"uckgeholte Metrik bzw. zweite Fundamentalform der gem"a"s $\boldsymbol{\eta}$ deformierten
Fl"ache. Die elastische Energie einer solchen Verschiebung kann durch die \emph{Koiter-Energie f"ur
eine nichtlinear elastische Schale}
\[K_N(\boldsymbol{\eta})=\frac{1}{2}\int_M \epsilon_0\,\langle C,
\Sigma(\boldsymbol{\eta})\otimes\Sigma(\boldsymbol{\eta})\rangle + \frac{\epsilon_0^3}{3}\,\langle
C,
\Xi(\boldsymbol{\eta})\otimes\Xi(\boldsymbol{\eta})\rangle\ dA \]
modelliert werden; siehe \cite{b54}, \cite{b20}, \cite{b45}, \cite{b23} und die dortigen Referenzen.
Dabei ist
\begin{equation*}
 \begin{aligned}
  C_{\alpha\beta\gamma\delta}:=\frac{4\lambda\mu}{\lambda+2\mu}\,g_{\alpha\beta}\,g_{\gamma\delta} +
2\mu\,(g_{\alpha\gamma}\,g_{\beta\delta} + g_{\alpha\delta}\,g_{\beta\gamma}) 
 \end{aligned}
\end{equation*}
das Elastizit"atstensorfeld der Schale, und $\Sigma(\boldsymbol{\eta}):=1/2\,
(g(\boldsymbol{\eta})-g)$ sowie $\Xi(\boldsymbol{\eta}):=h(\boldsymbol{\eta})-h$ bezeichnen die
Differenzen der ersten bzw. zweiten Fundamentalformen. Diese Energie wird in \cite{b20} auf Basis
der dreidimensionalen
Elastizit"at unter den zus"atzlichen Annahmen kleiner
Dehnungen und ebener Spannungszust"ande parallel zur Mittelfl"ache hergeleitet.  
% Die erste ist die
% Love-Kirchhoff-Annahme, die besagt,
% dass die (zun"achst dreidimensionale) Deformation die Dicke der Schale unver"andert
% l"asst und zudem die Normalen der urspr"unglichen
% Mittelfl"ache auf die jeweiligen Normalen der deformierten Mittelfl"ache abbildet. Nach der
% zweiten
% Annahme ist der Spannungszustand in der Schale eben und parallel zur
% Mittelfl"ache. 
Der mit $\epsilon_0$ skalierende Anteil der
Koiter-Energie, die L"angungs- oder Membranenergie, erfasst ausschlie"slich nichtisometrische
Deformationen.
Der mit $\epsilon_0^3$ skalierende Anteil kann als Biegeenergie interpretiert werden. F"ur
die rigorose Rechtfertigung derartiger "`zweidimensionaler Theorien"' auf
Basis der dreidimensionalen Elastizit"at verweisen wir auf \cite{b45}, \cite{b21}, \cite{b22}.
Setzen wir
\[\widetilde C:=\frac{4\lambda\mu}{\lambda+2\mu}\, g\otimes g,\]
so ist offenbar der Anteil 
\begin{equation*}
 \begin{aligned}
\frac{1}{2}\int_M \frac{\epsilon_0^3}{3}\,\langle \widetilde C,
h(\boldsymbol{\eta})\otimes
h(\boldsymbol{\eta})\rangle\ dA&=\frac{\epsilon_0^3}{3}\frac{4\lambda\mu}{\lambda+2\mu}\
\frac{1}{2}\int_M \langle g,
h(\boldsymbol{\eta})\rangle^2\ dA
 \end{aligned}
\end{equation*}
der Biegeenergie verwandt mit der Willmore-Energie 
\[W(\boldsymbol{\eta})=\frac{1}{2}\int_M H^2(\boldsymbol{\eta})\ dA_{\boldsymbol{\eta}}\]
der deformierten Mittelfl"ache. Hierbei bezeichnen $H(\boldsymbol{\eta})=\frac{1}{2}\langle
g(\boldsymbol{\eta}),h(\boldsymbol{\eta})\rangle$ die auf $M$
zur"uckgeholte mittlere Kr"ummung der deformierten Mittelfl"ache und $dA_{\boldsymbol{\eta}}$ das
Fl"achenma"s bez"uglich der Metrik $g(\boldsymbol{\eta})$. 

In der Koiter-Energie besteht "uber das Hooke'sche Gesetz zwar ein linearer Zusammenhang zwischen
Dehnungen und Spannungen, aber der Zusammenhang zwischen Deformationen und Dehnungen ist
nichtlinear. 
% Wenn wir bei der Herleitung der Bewegungsgleichungen der Schale zus"atzlich zur
% elastischen Energie der Fl"ache "uber die
% Massentr"agheit noch eine kinetische Energie in das Wirkungsfunktional aufnehmen, erhalten wir 
% Euler-Lagrange-Gleichungen in Form eines quasilinearen, dispersiven Systems; vgl. die Rechnungen
% unten. Da wir an zeitlich globalen L"osungen (f"ur beliebige, nicht notwendig kleine
% Anfangsdaten) interessiert sind, ist es geboten Vereinfachungen vornehmen. 
Wir linearisieren nun die Abh"angigkeit der Tensorfelder
$\Sigma(\boldsymbol{\eta})$ und $\Xi(\boldsymbol{\eta})$ von $\boldsymbol{\eta}$ an der Stelle
$\boldsymbol{\eta}=\boldsymbol{0}$
und erhalten die Tensorfelder $\sigma(\boldsymbol{\eta})$ und $\xi(\boldsymbol{\eta})$. Entsprechend
definieren wir die \emph{Koiter-Energie (f"ur eine linear elastische Schale)}
\begin{equation*}
 \begin{aligned}
  K(\boldsymbol{\eta})=K(\boldsymbol{\eta},\boldsymbol{\eta})=\frac{1}{2}\int_M \epsilon_0\,\langle
C,
\sigma(\boldsymbol{\eta})\otimes\sigma(\boldsymbol{\eta})\rangle + \frac{\epsilon_0^3}{3}\,\langle
C,
\xi(\boldsymbol{\eta})\otimes\xi(\boldsymbol{\eta})\rangle\ dA. 
 \end{aligned}
\end{equation*}
$K$ ist eine quadratische Form in $\boldsymbol{\eta}$, und es gilt $K\ge 0$; siehe Theorem 4.4-1 in
\cite{b23}. Sie ist unter Umst"anden sogar koerziv; siehe Theorem 4.4-2 in \cite{b23}. Der
Gradient dieser Energie ist ein elliptischer Operator vierter Ordnung, der in an $M$ tangentiale
Richtung degeneriert ist.\footnote{Diese Beobachtung spielt in \cite{b42} eine prominente Rolle.
Da dort die Koiter-Energie f"ur nichtlinear elastische Schalen betrachtet wird, h"angen die
Entartungsrichtungen von der L"osung ab, sind also a-priori unbekannt.} Speziell ist er von
lediglich dritter Ordnung auf dem tangentialen Anteil. Der
normale Anteil der Auslenkung ist deshalb im Allgemeinen regul"arer als der tangentiale; vgl.
Theorem 4.4-2 in \cite{b23}. Motiviert durch das Vorgehen in \cite{b27} schr"anken
% Da wir die Kopplung der Schale mit einem Fluid analysieren wollen,
% ben"otigen wir eine gewisse Regularit"at der Deformation $\boldsymbol{\eta}$. 
wir im Rahmen dieser Arbeit die Verschiebungen auf die Richtung der Einheitsnormale $\bnu$
von $M$ ein. Neben dem
Elliptizit"ats- und Regularit"atsgewinn hat die Einschr"ankung den Vorteil, dass wir
ein einfaches Kriterium daf"ur angeben k"onnen, dass sich verschiedene Teile der
Schale nicht ber"uhren. Wir k"onnen die Auslenkung $\boldsymbol{\eta}=:\eta\,\bnu$ durch eine
skalare Funktion $\eta$ beschreiben. Dadurch integrieren wir die vorliegende Zwangsbedingung in den
Phasenraum, sodass kein Lagrange-Multiplikator in den Gleichungen auftreten wird.  Wir
setzen zudem $K(\eta):=K(\eta\,\bnu)$. Die
Tensorfelder $\sigma(\eta):=\sigma(\eta\,\bnu)$ und $\xi(\eta):=\xi(\eta\,\bnu)$ nehmen mit der
Definition $k_{\alpha\beta}=h_\alpha^\sigma\, h_{\sigma\beta}$ die Form
\begin{equation*}
 \begin{aligned}
  \sigma(\eta)=-h\, \eta,\ \xi(\eta)=\nabla^2 \eta - k\, \eta
 \end{aligned}
\end{equation*}
an; siehe Theorem 4.2-1 und Theorem 4.2-2 in \cite{b23}.
Damit gilt f"ur (hinreichend glatte) skalare Felder $\eta,\zeta$
\begin{equation*}
 \begin{aligned}
  K(\eta,\zeta)&=\frac{1}{2}\int_M \epsilon_0\,\langle C,
\sigma(\eta)\otimes\sigma(\zeta)\rangle + \frac{\epsilon_0^3}{3}\,\langle C,
\xi(\eta)\otimes\xi(\zeta)\rangle\ dA \\
&= \frac{1}{2}\int_M \epsilon_0 \, \eta\,\zeta\, \Big(\frac{16\lambda\mu}{\lambda+2\mu} H^2 + 4\mu
|h|^2\Big)\\
&\hspace{1.3cm}+ \frac{\epsilon_0^3}{3} \frac{4\lambda\mu}{\lambda+2\mu}(\Delta \eta\,\Delta\zeta -
|h|^2\ (\zeta\, \Delta \eta + \eta\,\Delta\zeta) -\eta\,\zeta\, |h|^4)\\
&\hspace{1.3cm} + \frac{\epsilon_0^3}{3}\,4\mu\,(\langle \nabla^2\eta,\nabla^2\zeta\rangle - \zeta\,
\langle \nabla^2\eta,k\rangle - \eta\,
\langle \nabla^2\zeta,k\rangle - \eta\,\zeta\, |k|^2)\ dA,
 \end{aligned}
\end{equation*}
wobei $H$ die mittlere Kr"ummung von $M$ bezeichnet. F"ur die zweite Gleichung haben wir
beispielsweise
\[\langle g\otimes g,\nabla^2\eta\otimes \nabla^2\zeta\rangle=\langle g,\nabla^2\eta\rangle\,\langle
g,\nabla^2\zeta\rangle=\Delta\eta\,\Delta\zeta\]
und
\[g^{\alpha\bar\alpha}\,g^{\gamma\bar\gamma}\,g^{\beta\bar\beta}\,g^{\delta\bar\delta}\,
g_{\alpha\gamma}\,g_{\beta\delta}\, \nabla^2_{\bar\alpha\bar\beta}\eta\,
\nabla^2_{\bar\gamma\bar\delta
}\zeta=g^{\alpha\gamma}\,g^{\beta\delta}\,\nabla^2_{\alpha\beta}\eta\, \nabla^2_{\gamma\delta}
\eta=\langle\nabla^2\eta,\nabla^2\zeta\rangle\]
verwendet. Wenn wir annehmen, dass $\zeta$ und $\nabla\zeta$ auf dem Rand vom $M$
verschwinden, erhalten wir durch
partielle
Integration, Gleichung \eqref{eqn:pisatz}, den $L^2$-Gradienten
der
Energie
\begin{equation}\label{eqn:grad}
 \begin{aligned}
dK(\eta)\,\zeta&=2\,K(\eta,\zeta)\\
&= \int_M \Big(\epsilon_0\,\Big(\frac{16\lambda\mu}{\lambda+2\mu} H^2 + 4\mu\, |h|^2\Big)\,\eta\\
&\hspace{1.3cm}+ \frac{\epsilon_0^3}{3}\, \frac{4\lambda\mu}{\lambda+2\mu}(\Delta^2 \eta -
|h|^2\ \Delta \eta - \Delta (|h|^2\eta) -\eta\, |h|^4)\\
&\hspace{1.3cm} + \frac{\epsilon_0^3}{3}\,4\mu\,( \Delta^2\eta + \nabla^*[\Delta,\nabla]\eta -
\langle
\nabla^2\eta,k\rangle -
(\nabla^*)^2 (k\,\eta) - \eta\,|k|^2)\Big)\, \zeta\ dA\\
&=:(\grad_{L^2}K(\eta),\zeta)_{L^2}.
 \end{aligned}
\end{equation}
Zum Beispiel erhalten wir durch jeweils zweimalige Verwendung von
\eqref{eqn:pisatz}\footnote{$\trace_g$ bezeichnet hier die Spurbildung im Tensorsinne, nicht im
Sinne von Randwerten.}
\begin{equation*}
 \begin{aligned}
\int_M \Delta\eta\,\Delta\zeta\ dA&=\int_M \Delta\eta\, \trace_g\nabla^2\zeta\
dA=-\int_M\langle\nabla\Delta\eta,\nabla\zeta\rangle\ dA\\
&=\int_M\trace_g\nabla^2\Delta\eta\, \zeta\ dA=\int_M\Delta^2\eta\, \zeta\ dA
 \end{aligned}
\end{equation*}
und
\begin{equation*}
 \begin{aligned}
\int_M \langle\nabla^2\eta,\nabla^2\zeta\rangle\ dA=-\int_M
\langle\Delta\nabla\eta,\nabla\zeta\rangle\ dA&=-\int_M
\langle \nabla\Delta\eta + [\Delta,\nabla]\eta,\nabla\zeta\rangle\ dA\\
&=\int_M(\Delta^2\eta + \nabla^*[\Delta,\nabla]\eta)\, \zeta\ dA.
 \end{aligned}
\end{equation*}
Der Kommutator l"asst sich in lokalen Koordinaten in der Form 
\[([\Delta,\nabla]\eta)_\alpha = (2 H\, h_{\alpha}^\beta - k_\alpha^\beta)\, \pa_\beta\eta\]
schreiben; siehe \eqref{eqn:kommutator}, \eqref{eqn:ricci}. Insbesondere ist er ein
Differentialoperator erster Ordnung, der im Falle $h=0$ verschwindet. 

Wir nehmen nun an, dass die Massenverteilung der Schale durch eine konstante Fl"achenmassendichte
$\epsilon_0\rho_S$ der Mittelfl"ache $M$ beschrieben werden kann. Das Hamilton'sche Prinzip besagt,
dass die
Bewegung der Schale, hier gegeben durch ein zeitabh"angiges skalares Feld $\eta$, ein station"arer
Punkt des Wirkungsintegrals
\[\mathcal{A}(\eta)=\int_I \epsilon_0\rho_S\int_M \frac{(\pa_t\eta(t,\cdot))^2}{2}\
dA-K(\eta(t,\cdot))\ dt\]
ist, wobei $I:=(0,T)$, $T>0$. Der zeitliche Integrand ist die Differenz von
kinetischer und potentieller Energie der Schale. Die Ableitung von $\mathcal{A}$ an der Stelle
$\eta$ in
Richtung eines skalaren Feldes $\zeta$, das zusammen
mit $\nabla\zeta$ auf dem Rand von $I\times M$ identisch $0$ sei, ist durch
\begin{equation*}
 \begin{aligned}
d\mathcal{A}(\eta)\,\zeta=\int_I\epsilon_0\rho_S\int_M \pa_t\eta(t,\cdot)\, \pa_t\zeta(t,\cdot)\ dA
-dK(\eta(t,\cdot))\,\zeta(t,\cdot) \ dt
 \end{aligned}
\end{equation*}
gegeben. Da die erste Variation verschwinden soll, erhalten wir mittels partieller Integration in
der Zeit und \eqref{eqn:grad} die Gleichung
\begin{equation*}
 \begin{aligned}
0&=\epsilon_0\rho_S\,\pa^2_t\eta +
\grad_{L^2}K(\eta)=\epsilon_0\rho_S\,\pa^2_t\eta+\epsilon_0^3\,\frac{8\mu(\lambda+\mu)}{
3(\lambda+2\mu)}\,
\Delta^2\eta+B\eta\text { in
} I\times M.
 \end{aligned}
\end{equation*}
Hierbei ist $B$ ein Operator zweiter Ordnung, der im Falle $h=0$ verschwindet. Es liegt also eine
Verallgemeinerung der linearen Kirchhoff-Love-Plattengleichung f"ur transversale Auslenkungen vor;
vgl. zum Beispiel \cite{b44}.
% \begin{equation}\label{eqn:K}
%  \begin{aligned}
%   K(\eta,\zeta)&= \frac{1}{2}\int_M \eta\zeta \big(H^2 + |h|^2\big)+\frac{1}{2}\Delta
% \eta\Delta\zeta -
% \trace_g k\ (\zeta \Delta \eta + \eta\Delta\zeta) -\eta\zeta (\trace_g k)^2\\
% &\hspace{1.3cm} + \frac{1}{2}\langle \nabla^2\eta,\nabla^2\zeta\rangle - \zeta
% \langle \nabla^2\eta,k\rangle - \eta
% \langle \nabla^2\zeta,k\rangle - \eta\zeta |k|^2\ dA.
%  \end{aligned}
% \end{equation}
% Die Bewegungsgleichung nimmt dann die Form
% \begin{equation*}
%  \begin{aligned}
%   \pa^2_t\eta + \Delta&^2\eta + B\eta = 0
%  \end{aligned}
% \end{equation*}
% an, wobei $B\eta := - 2\trace_g k\ \Delta\eta - \langle \nabla^2\eta,k\rangle + \eta(H^2 + |h|^2 -
% (\trace_g k)^2 - (\nabla^*)^2 k - |k|^2)$.
Man sollte sich klar machen, dass diese Gleichung dispersiv, aber nicht hyperbolisch ist. Der
Hauptteil faktorisiert in zwei Schr"odinger-Operatoren
\[\pa^2_t +
\pa^4_x=(i\pa_t+\pa^2_x)(-i\pa_t+\pa^2_x);\] im Gegensatz zum d'Alembert-Operator, der in
zwei Transportoperatoren faktorisiert
\[\pa^2_t - \pa^2_x=(\pa_t+\pa_x)(\pa_t-\pa_x).\] Die
Gleichung l"asst mithin eine unendliche Ausbreitungsgeschwindigkeit zu. Der Einfachheit halber
setzen wir fortan $\epsilon_0\rho_S=1$.

\chapter{Problemstellung}

Wir "ubernehmen die Bezeichnungen aus den vorangegangenen Abschnitten. Insbesondere
sei $\Omega\subset\setR^3$ ein beschr"anktes, nichtleeres Gebiet mit $C^4$-Rand. Der Rand
$\pa\Omega$ stelle die Mittelfl"ache einer elastischen Schale in ihrer Ruhelage dar, deren
Bewegung auf Auslenkungen
l"angs der "au"seren Einheitsnormale $\bnu$ von $\pa\Omega$ eingeschr"ankt sei. $\Gamma$ sei eine
Vereinigung von
Gebieten in $\pa\Omega$ mit $C^{1,1}$-Rand, die nichttrivialen Schnitt mit allen
Zusammenhangskomponenten von $\pa\Omega$ habe. Diesen Teil der Mittelfl"ache
setzen wir als fixiert
voraus. Mit $M$ wollen wir den beweglichen Teil bezeichnen, d.h.
$M:=\pa\Omega\setminus\Gamma$. Die Mannigfaltigkeit $M$ ist kompakt und berandet. Desweiteren sei
$I:=(0,T)$, $T>0$. Die zeitabh"angige Auslenkung des Gebiets
beschreiben wir durch eine Funktion
$\eta: \bar I\times M\rightarrow (-\kappa,\kappa)$, die wir (wie immer) durch
$0$ auf $\bar I\times\pa\Omega$ fortsetzen. Desweiteren nehmen wir an,
dass das Innere des variablen Gebiets $\Omega_{\eta}^I$ durch ein inkompressibles, viskoses
Fluid bef"ullt sei, dessen Geschwindigkeitsfeld $\bu$ und Druckfeld $\pi$ durch die
\emph{Navier-Stokes-Gleichungen} beschrieben werden kann, d.h. es gelte
\begin{equation}\label{eqn:fluid}
 \begin{aligned}
  \rho_F\big(\pa_t \bu + (\bu\cdot\nabla)\bu\big) &= \dv (2\sigma D\bu - \pi\id) + \ff &&\mbox{ in
}
\Omega_{\eta}^I, \\
 \dv \bu &= 0 &&\mbox{ in }  \Omega_{\eta}^I,\\
 \bu(\,\cdot\,,\,\cdot\, + \eta\,\bnu) &= \pa_t\eta\,\bnu &&\mbox{ auf } I\times M,\\
  \bu &= 0 &&\mbox{ auf } I\times\Gamma.
\end{aligned}
\end{equation}
Die Bezeichnung $\id$ steht dabei f"ur die $3\times 3$-Einheitsmatrix, und $\ff$
ist eine gegebene "au"sere Kraftdichte. Zudem ist $(\bu\cdot\nabla)\bu:=u^i\,\pa_i\bu.$
%=\nabla_\bu\bu
Die konstante Dichte $\rho_F$
und die konstante Viskosit"at $\sigma$ setzen wir der Einfachheit halber fortan identisch $1$.
$2D\bu - \pi\id$ ist dann der Spannungstensor, und $2D\bu$ ist der viskose Spannungstensor. Aufgrund
der Divergenzfreiheit von $\bu$ gilt $\dv
2D\bu=\Delta\bu$. \eqref{eqn:fluid}$_{3,4}$
ist die no-slip-Bedingung im Falle eines sich bewegenden Randes,
d.h.
die Geschwindigkeit des Fluids stimmt am Rand mit der Geschwindigkeit des Randes "uberein. Die
vom Fluid auf die Schale ausge"ubte Kraft ist gegeben durch die
Auswertung des Spannungstensors am Rand in Richtung der inneren Normale
$-\bnu_{\eta(t)}$, also durch
\begin{equation}\label{eqn:kraft}
 \begin{aligned}
-2D\bu(t,\cdot)\,\bnu_{\eta(t)} + \pi(t,\cdot)\,\bnu_{\eta(t)}.  
 \end{aligned}
\end{equation}
In die Gleichungen f"ur die Auslenkung der Schale ($M:=\pa
\Omega\setminus\Gamma$)
\begin{equation}\label{eqn:shell}
\begin{aligned}
  \pa^2_{t} \eta + \grad_{L^2}K(\eta)  &= g + {\bf F}\cdot\bnu &&\text{ in } I\times M, \\
 \eta=0,\ \nabla \eta &= 0 &&\text{ auf } I\times\pa M
\end{aligned}
\end{equation}
geht diese Kraft zus"atzlich zu einer gegebenen "au"seren Kraftdichte $g$ in der Form
\[{\bf F}(t,\cdot) = \big(-2D\bu(t,\cdot)\,\bnu_{\eta(t)} +
\pi(t,\cdot)\,\bnu_{\eta(t)}\big)\circ\Phi_{\eta(t)}\ |\det d\Phi_{\eta(t)}|\]
ein. Da die Gleichungen der Elastizit"at in Lagrangekoordinaten, d.h. im Referenzgebiet $M$,
formuliert werden, muss die auf $\pa\Omega_\eta$ registrierte Fl"achenkraftdichte \eqref{eqn:kraft}
mit der Abbildung $\Phi_{\eta(t)}$ verkn"upft und mit der lokalen Fl"achenverzerrung skaliert
werden. Schlie"slich geben wir noch Anfangsdaten 
\begin{equation}\label{eqn:data}
 \begin{aligned}
\eta(0,\cdot)=\eta_0,\ \pa_t\eta(0,\cdot)=\eta_1 \text{ und } \bu(0,\cdot)=\bu_0
 \end{aligned}
\end{equation}
vor. Man beachte, dass in dieses System nicht nur der Druckgradient eingeht, wie es bei einem festen
Rand der Fall ist, sondern auch der Druck selber. So wie der Druckgradient als
Lagrange-Multiplikator bez"uglich der Zwangsbedingung $\dv\bu=0$ interpretiert werden kann, l"asst
sich sein Mittelwert als der mit der Zwangsbedingung
\[\im \pa_t\eta\, \gamma(\eta)\ dA=0\]
assoziierte Lagrange-Multiplikator verstehen. Diese Zwangsbedingung folgt mittels Lemma
\ref{lemma:partInt} aus der Divergenzfreiheit von $\bu$ und der Kopplung \eqref{eqn:fluid}$_3$.

Wir wollen nun rein formal, d.h. unter Vernachl"assigung von Regularit"atsfragen,
Energieabsch"atzungen f"ur dieses parabolisch-dispersive System herleiten. Die Konstruktion
schwacher L"osungen beruht wesentlich auf diesen Absch"atzungen. Nach dem Noether-Theorem
ist der Energiesatz
eng verkn"upft mit der zeitlichen Translationsinvarianz des physikalischen Systems. Wir erhalten
einen lokalen Energiesatz, wenn wir die Gleichung mit dem Vektorfeld multiplizieren, das sich
ergibt, wenn wir den infinitesimalen Erzeuger der zeitlichen Translationen, also $\pa_t$, auf die
Auslenkung anwenden. Wir multiplizieren also Gleichung \eqref{eqn:fluid}$_1$ mit $\bu$ (und sp"ater
Gleichung \eqref{eqn:shell}$_1$ mit $\pa_t\eta$). Da wir an einem globalen
Energiesatz interessiert sind, integrieren wir die resultierende Identit"at "uber $\Omega_{\eta(t)}$
und erhalten nach partieller Integration der Terme des Spannungstensors\footnote{Der
"Ubersichtlichkeit halber unterdr"ucken wir hier und im Folgenden in der Notation die unabh"angigen
Variablen, d.h.
wir setzen $\bu=\bu(t,\cdot)$, etc.}
\begin{equation}\label{eqn:mult}
 \begin{aligned}
  \iot \pa_t\bu\cdot\bu\ dx + \iot (\bu\cdot\nabla)\bu\cdot\bu\ &dx = -\iot |\nabla\bu|^2\
dx+ \iot \ff\cdot\bu\ dx\\
&+\idot
(2D\bu\, \bnue - \pi\,\bnue)\cdot \bu\ dA_{\eta(t)}.
 \end{aligned}
\end{equation}
Aufgrund der Divergenzfreiheit von $\bu$ verschwindet das r"aumliche Integral des Druckterms.
Bez"uglich des r"aumlichen Integrals des viskosen Spannungstensors haben wir die Korn'sche Gleichung
\[2\int_{\Omega_{\eta(t)}}D\bu:D\bu\
dx=2\int_{\Omega_{\eta(t)}}D\bu:\nabla\bu\ dx=\int_{\Omega_{\eta(t)}}\nabla\bu:\nabla\bu\ dx\]
verwendet. F"ur das erste
Gleichheitszeichen in dieser Identit"at ist zu beachten, dass die Kontraktion eines symmetrischen
Tensors mit einem antisymmetrischen Tensor verschwindet, sodass wir $\nabla\bu$ durch seinen
symmetrischen Anteil ersetzen k"onnen. F"ur das zweite Gleichheitszeichen sei auf Bemerkung
\ref{bem:korn} im Anhang verwiesen. Die ersten beiden Integrale in \eqref{eqn:mult}
k"onnen wir mit Hilfe des Reynolds'schen Transporttheorems, Proposition \ref{theorem:reynolds},
zusammenfassen, wobei hier
$\bv=\bu$ gilt und wir $\xi=\frac{|\bu|^2}{2}$ w"ahlen. Partielle Integration
des zweiten Integrals ergibt n"amlich aufgrund
der Divergenzfreiheit von $\bu$
\begin{equation}\label{eqn:wirbelid}
 \begin{aligned}
\iot (\bu\cdot\nabla)\bu\cdot\bu\ dx=-\iot (\bu\cdot\nabla)\bu\cdot\bu\ dx +
\int_{\pa\Omega_{\eta(t)}}
\bu\cdot\bnu_{\eta(t)}|\bu|^2\ dA_{\eta(t)}.  
 \end{aligned}
\end{equation}
Wir erhalten dann
\begin{equation}\label{eqn:energiefluid}
 \begin{aligned}
  \frac{1}{2}\frac{d}{dt} \iot |\bu|^2\ dx = &-\iot
|\nabla\bu|^2\ dx+ \iot \ff\cdot\bu\ dx\\
&+\idot
(2D\bu\, \bnue - \pi\,\bnue)\cdot \bu\ dA_{\eta(t)}.
 \end{aligned}
\end{equation}
% und verwenden
% wir die Koerzitit"at von $S$, so erhalten wir 
% \begin{equation}\label{ap:fluid}
%  \begin{aligned}
%   \frac{1}{2}\frac{d}{dt} &\iot |u|^2\ dx + c_1\iot
% |D\bu|^p\ dx - \iot\kappa_2\ dx\\ &-\idot
% (S(D\bu)\cdot \bnue + \pi\bnue)\cdot \bu\ dA_{\eta(t)} \le \iot \ff\cdot\bu\ dx.
%  \end{aligned}
% \end{equation}
Multiplikation der Gleichung \eqref{eqn:shell}$_1$ mit $\pa_t\eta$, Integration "uber $M$ und
partielle Integration ergibt wegen $(\grad_{L^2}K(\eta),\pa_t\eta)_{L^2}=2K(\eta,\pa_t\eta)$
\begin{equation}\label{eqn:energieshell}
 \begin{aligned}
  \frac{1}{2}\frac{d}{dt} \im |\pa_t\eta|^2\ dA + \frac{d}{dt}K(\eta) = \im g\, \pa_t\eta\ dA + \im
{\bf F}\cdot\bnu\, \pa_t\eta\ dA.
 \end{aligned}
\end{equation}
% \begin{equation}\label{ap:shell}
%  \begin{aligned}
%   \frac{1}{2}\frac{d}{dt} \im |\pa_t\eta|^2\ dA + \frac{1}{2}\frac{d}{dt}\im &|\Delta\eta|^2\ dA +
% \im B\eta\ \pa_t\eta\ dA\\
%  &= \im g\ \pa_t\eta\ dA + \im {\bf F}\cdot\bnu \pa_t\eta\ dA.
%  \end{aligned}
% \end{equation}
% Offenbar gilt 
% \begin{equation*}
%  \begin{aligned}
%   \im B\eta\ \pa_t\eta\ dA \le c\big(\norm{\eta}_{H^2(M)}^2 + \norm{\pa_t \eta}_{L^2(M)}^2\big).
%  \end{aligned}
% \end{equation*}
Durch Addition von \eqref{eqn:energiefluid} und \eqref{eqn:energieshell} erhalten wir unter
Verwendung der
Definition von ${\bf F}$, der Randbedingung \eqref{eqn:fluid}$_3$ und des Transformationssatzes den
Energiesatz
\begin{equation}\label{eqn:energiesatz}
 \begin{aligned}
  \frac{1}{2}\frac{d}{dt} \iot |\bu|^2\ dx + \frac{1}{2}\frac{d}{dt} &\im |\pa_t\eta|^2\ dA +
\frac{d}{dt}K(\eta)\\
 &= - \iot |\nabla\bu|^2\ dx + \iot \ff\cdot\bu\ dx + \im g\, \pa_t\eta\ dA.
 \end{aligned}
\end{equation}
Die zeitliche "Anderung der Gesamtenergie des Systems, zusammengesetzt aus den kinetischen Energien
von Fluid und Schale sowie der potentiellen Energie der Schale, ist identisch der negativen
Energiedissipation durch das Fluid und der Leistung der "au"seren Kr"afte. Wir definieren den Term
$R$ durch
\[\frac{d}{dt}K(\eta)=2\,K(\eta,\pa_t\eta)=:\epsilon_0^3\,\frac{8\mu(\lambda+\mu)}{6(\lambda+2\mu)}
\,\frac{d}{dt}\int_M
|\Delta\eta|^2\ dA + R;\]
man beachte die Identit"at \eqref{eqn:grad}. Offenbar gilt
\[|R|\le c\,(\norm{\eta}_{H^2(M)}^2 + \norm{\pa_t\eta}_{L^2(M)}^2).\]
% Dabei haben wir zum Beispiel die Absch"atzung
% \[\big|\int_M \eta\langle\nabla^2\pa_t\eta,k\rangle\ dA\big| = \big|\int_M
% \pa_t\eta\ (\nabla^*)^2 (k\eta)\ dA\big|\le c(\norm{\eta}_{H^2(M)}^2 +
% \norm{\pa_t\eta}_{L^2(M)}^2)\]
% verwendet.
Wir erhalten somit aus \eqref{eqn:energiesatz}
\begin{equation*}
 \begin{aligned}
 \frac{d}{dt} \iot |u|^2\ dx &+ \iot |\nabla\bu|^2\ dx + \frac{d}{dt} \im
|\pa_t\eta|^2\ dA + \frac{d}{dt}\im |\Delta\eta|^2\ dA \\ &\le
c\,\big(\norm{\ff}_{L^2(\Omega_{\eta(t)})}^2 +
\norm{\bu}_{L^2(\Omega_{\eta(t)})}^2 + \norm{g}_{L^2(M)}^2+ \norm{\pa_t
\eta}_{L^2(M)}^2 + \norm{\eta}_{H^2(M)}^2 \big).
 \end{aligned}
\end{equation*}
$\mnorm{\Delta\,\cdot\,}$ definiert eine
"aquivalente Norm auf $H^2_0(M)$; siehe Anhang A.2. Mit Hilfe des Gronwall'schen Lemmas, Proposition
\ref{lemma:gronwall}, erhalten wir
also 
\begin{equation}\label{ab:apriori}
 \begin{aligned}
&\onorm{\bu(t,\cdot)}^2 + \int_0^t\norm{\nabla\bu(s,\cdot)}^2_{L^2(\Omega_{\eta(s)})}\ ds +
\mnorm{\pa_t\eta(t,\cdot)}^2 + \norm{\eta(t,\cdot)}_{H^2(M)}^2 \\ 
&\hspace{4cm}\le \big(\norm{\bu_0}_{L^2(\Omega_{\eta_0})}^2 + \mnorm{\eta_1}^2 +
\norm{\eta_0}_{H^2(M)}^2\big)\,e^{ct}\\ 
&\hspace{4.5cm} + \int_0^t
\big(\norm{\ff(s,\cdot)}_{L^2(\Omega_{\eta(s)})}^2 + \mnorm{g(s,\cdot)}^2\big)\,e^{c(t-s)}\ ds.
 \end{aligned}
\end{equation}
Es gilt somit
\begin{equation*}
 \begin{aligned}
 \norm{\eta}_{W^{1,\infty}(I,L^2(M))\cap L^\infty(I,H_0^2(M))} +
\norm{\bu}_{L^\infty(I,L^2(\Omega_{\eta(t)}))\cap L^2(I,H^1(\oet))}\le c(T,\text{Daten}).
 \end{aligned}
\end{equation*}
Wir werden schwache L"osungen in dieser Regularit"atsklasse konstruieren. Die Einbettung
\eqref{eqn:hoeldereinb} zeigt, dass der Rand unseres zeitlich
variablen Gebietes Graph einer
H"older-stetigen, aber nicht Lipschitz-stetigen Funktion sein wird.
\chapter{Existenz}

Wir setzen
\begin{equation*}
 \begin{aligned}
 Y^I&:=W^{1,\infty}(I,L^2(M))\cap L^\infty(I,H_0^2(M)),
 \end{aligned}
\end{equation*}
und f"ur $\eta\in Y^I$ mit $\norm{\eta}_{L^\infty(I\times M)}<\kappa$
\begin{equation*}
 \begin{aligned}
X_\eta^I&:=L^\infty(I,L^2(\Omega_{\eta(t)}))\cap L^2(I,H^1_{\dv}(\oet)).
%  T_\eta&:=\{(b,\bphi)\in H^1(I,L^2(M))\cap L^2(I,H^2_0(M))\times H^1(\Omega_{\eta}^I)\cap
% L^\infty(I,L^4(\oet))\\
% &\hspace{1.0cm}| \bphi-Fb\in (H^1(\Omega_{\eta}^I)\cap
% L^\infty(I,L^4(\oet)))_{0,w^*}, b(T,\cdot)=0,
% \bphi(T,\cdot)=\boldsymbol{0},\\ 
% &\hspace{1.1cm}\dv\bphi=0 \}.
%  X^p_\eta&:=L^\infty(I,L^2(\Omega_{\eta(t)}))\cap L^p(I,W^{1,p}_{\dv}(\oet)),\\
%  T^p_\eta&:=\{(b,\bphi)\in (H^1(I,L^2(M))\cap
% L^p(I,H^2_0(M)))\times (H^1(\Omega_{\eta}^I)\\
% &\hspace{0.7cm}\cap L^p(I,W^{1,p}_{\dv} (\Omega_{\eta(t)})))| \bphi(T,\cdot)=\boldsymbol{0}, \tr
% \bphi = b\bnu \}.
 \end{aligned}
\end{equation*}
Wir wollen nun einen geeigneten Raum von Testfunktionen $b$ auf $M$ und $\bphi$ auf
$\Omega_\eta^I$ definieren. Nat"urlich m"ussen die Funktionen $\bphi$ divergenzfrei sein, und es
soll $\tr\bphi=b\,\bnu$ f"ur eine Testfunktion $b$ gelten. "Ublicherweise w"urde man $b$ und $\bphi$
hinreichend oft stetig differenzierbar w"ahlen. Aufgrund der geringen Regularit"at von $\eta$ ist
aber nicht klar, ob "uberhaupt ein ausreichender Satz solcher Testfunktionen existiert. Insbesondere
ist die Fortsetzung in Proposition \ref{lemma:FortVonRand} auch f"ur $C^4$-Randwerte $b$ im
Allgemeinen nicht stetig differenzierbar. Wir bezeichnen deshalb mit $T_\eta^I$ den
kanonisch normierten Raum aller Tupel
\[(b,\bphi)\in \big(H^1(I,L^2(M))\cap L^2(I,H^2_0(M))\big)\times
\big(H^1(\Omega_{\eta}^I)\cap
L^\infty(I,L^4(\oet))\big)\]
mit $b(T,\cdot)=0$, $\bphi(T,\cdot)=0$\footnote{Die Auswertung von $\bphi$ bei $t=T$ ist
sinnvoll, weil diese Funktion in $H^1((t,T)\times Q)$ f"ur jeden
offenen Ball $Q\subset\subset\Omega_{\eta(T)}$ und hinreichend gro"ses $t<T$
liegt.}, $\dv\bphi=0$\footnote{Der Divergenzoperator wirkt wie
immer nur auf die r"aumlichen Koordinaten.} und 
$\bphi-\F_\eta b\in H_0$. Dabei ist $H_0$ der Abschluss in $H^1(\Omega_{\eta}^I)\cap
L^\infty(I,L^4(\oet))$ der bei $t=T$ verschwindenden, divergenzfreien Elemente dieses Raums mit in
$\Omega_{\eta}^{\bar I}$ enthaltenen Tr"agern. Aus
der letzten Forderung folgt offenbar $\tr\bphi=\tr F_\eta b=b\,\bnu$. Die Umkehrung gilt
zumindest dann, wenn $\eta$ hinreichend glatt ist und wir den Testfunktionen den kleineren Raum
\[H^1(I,H^2_0(M))\times H^1(I,H^1(\Omega_{\eta(t)}))\]
zugrunde legen.\footnote{Dass dieser Raum kleiner ist, ist nicht offensichtlich. Wir zeigen aber in
Bemerkung \ref{bem:ausw}, dass f"ur $(b,\bphi)\in T^I_\eta$ die Fortsetzung von $\bphi$ durch
$(b\,\bnu)\circ q$
in $H^1(I,L^2(B_\alpha))$ (f"ur geeignetes $\alpha$) liegt. Dasselbe Argument zeigt, dass $\bphi$
in $H^1(I,L^4(B_\alpha))$ liegt, wenn wir den Testfunktionen obigen Raum zugrunde legen.} Unter
Verwendung der Abbildung $\T_{\eta}$ aus Bemerkung \ref{bem:tdelta} sehen wir n"amlich, dass es bei
hinreichend glattem $\eta$ gen"ugt, die analoge
Situation im Raumzeitzylinder $I\times\Omega$ zu betrachten. Dort k"onnen wir standardm"a"sig unter
Verwendung des L"osungsoperators der Divergenzgleichung geeignete Approximationen konstruieren;
siehe zum Beispiel III.4.1 in \cite{b18}. F"ur jedes $t\in\bar I$ liegen diese in
$C_0^\infty(\Omega)$ und konvergieren in $H^1(\Omega)$. Da die dortige zeitunabh"angige Konstruktion
mit der
Zeitableitung kommutiert, folgt die Konvergenz in $H^1(I,H^1(\Omega_{\eta(t)}))$ in trivialer
Weise.\footnote{Diese einfache Schlussfolgerung ist f"ur den Raum $H^1(\Omega_{\eta}^I)\cap
L^\infty(I,L^4(\oet))$ nicht m"oglich.} Gilt lediglich $\eta\in Y^I$, so ist die
Abbildung $\T_\eta$ nicht mehr anwendbar. Die Verwendung des
L"osungsoperators der Divergenzgleichung in $\Omega_\eta^I$ ist dann aber auch problematisch, weil
dessen "ubliche Stetigkeitseigenschaften in Gebieten, deren Rand nicht Lipschitz-stetig ist, im
Allgemeinen nicht gelten; vgl. zum Beispiel \cite{b26}. M"oglicherweise ist die Umkehrung in diesem
Fall falsch; vgl. die Diskussion in Abschnitt III.4 in \cite{b18}. Wie immer setzen wir
auf $M$ definierte Funktionen stillschweigend durch $0$ auf $\pa\Omega$ fort. Die Bedingung $\tr
\bphi = b\,\bnu$ besagt dann insbesondere, dass $\bphi$ auf $\Gamma$ verschwindet. 

Wir haben durch ein lokales Argument gezeigt, dass die Auswertung von $\bphi$ bei $t=T$ sinnvoll
ist. Dieses Argument zeigt allgemeiner, dass die Auswertung bei $t\in\bar I$ in
$L^2_{\loc}(\Omega_{\eta(t)})$ liegt. F"ur die nachfolgende
Definition ben"otigen wir jedoch $\bphi(0,\cdot)\in L^2(\Omega_{\eta_0})$. Die G"ultigkeit dieser
Aussage sehen wir ein, indem wir $\bphi$ durch $(b\,\bnu)\circ q$ auf
$I\times B_\alpha$, $\norm{\eta}_{L^\infty(I\times M)}<\alpha<\kappa$, fortsetzen und zeigen, dass
diese
Fortsetzung in $H^1(I,L^2(B_\alpha))$ liegt. Insbesondere ist die Auswertung von $\bphi$ bei
fester Zeit $t\in \bar I$ sinnvoll, und es gilt $\bphi(t,\cdot)\in L^2(\Omega_{\eta(t)})$. Die
Details finden sich in Bemerkung \ref{bem:ausw} im Anhang. 

Felder $\ff\in L_{\text{loc}}^{2}([0,\infty)\times \setR^3)$,
% $\ff\in L_{\text{loc}}^{p'}([0,\infty),L^{p'}(\setR^3))$,
$g\in L_{\text{loc}}^2([0,\infty)\times M)$, $\eta_0\in H^2_0(M)$ mit
$\norm{\eta_0}_{L^\infty(M)}<\kappa$, $\eta_1\in L^2(M)$ und
$\bu_0\in L^2(\Omega_{\eta_0})$
mit $\dv \bu_0=0$ sowie $\trnormaln\bu_0=\eta_1\,\gamma(\eta_0)$ nennen wir \emph{zul"assige Daten}.

\begin{definition} Ein Tupel $(\eta,\bu)$ hei"st schwache
L"osung von \eqref{eqn:fluid}, \eqref{eqn:shell} und \eqref{eqn:data}
mit den zul"assigen Daten $(\ff,g,\bu_0,\eta_0,\eta_1)$ auf dem Intervall $I$, falls
$\eta\in Y^I$ mit $\norm{\eta}_{L^\infty(I\times M)}<\kappa$ und $\eta(0,\cdot)=\eta_0$, $\bu\in
X_\eta^I$ mit $\tr \bu =
\pa_t\eta\,\bnu$ und
\begin{equation}\label{eqn:schwach}
\begin{aligned}
% &\iot \bu(t,\cdot)\cdot\bphi(t,\cdot)\ dx 
&- \int_I\iot \bu\cdot\pa_t\bphi\ dxdt - \int_I\im(\pa_t\eta)^2\, b\, \gamma(\eta)\
dAdt+\int_I\iot(\bu\cdot\nabla)\bu\cdot\bphi\ dxdt\\ &+ \int_I\iot \nabla\bu:\nabla\bphi\ dxdt
% +\im\pa_t\eta(t,\cdot)b(t,\cdot)\ dA
-\int_I\im\pa_t\eta\, \pa_tb\ dAdt + 2\int_I K(\eta,b)\ dt \\
&\hspace{1cm}=\int_I\iot \ff\cdot\bphi\ dxdt + \int_I\im g\, b\
dAdt+\int_{\Omega_{\eta_0}}\bu_0\cdot\bphi(0,\cdot)\ dx + \im\eta_1\, b(0,\cdot)\ dA
 \end{aligned}
\end{equation}
f"ur alle Testfunktionen $(b,\bphi)\in T_\eta^I$.
% \begin{definition} Sei $p\ge 2$. Ein Tupel $(\eta,\bu)$ hei"st schwache
% L"osung (von \eqref{eqn:shell}, \eqref{eqn:fluid} und \eqref{eqn:data}) auf dem Intervall $I$,
% falls
% $\eta\in Y$ mit $|\eta|<\kappa$ und $\eta(0,\cdot)=\eta_0$, $\bu\in X^p_\eta$ mit $\tr \bu =
% \pa_t\eta\bnu$ und
% \begin{equation}\label{eqn:schwach}
% \begin{aligned}
% % &\iot \bu(t,\cdot)\cdot\bphi(t,\cdot)\ dx 
% &- \int_I\iot \bu\cdot\pa_t\bphi\ dxdt - \int_I\im(\pa_t\eta)^2b(1-2H\eta+G\eta^2)\ dAdt\\
% &+\int_I\iot(\bu\cdot\nabla)\bu\cdot\bphi\ dxdt + \int_I\iot S(D\bu):D\bphi\ dxdt\\ 
% % +\im\pa_t\eta(t,\cdot)b(t,\cdot)\ dA
% &-\int_I\im\pa_t\eta\pa_tb\ dAdt + \int_I\im K(\eta,b)\
% dAdt \\
% &\hspace{3cm}=\int_I\iot \ff\cdot\bphi\ dxdt + \int_I\im gb\
% dAdt\\
% &\hspace{3.5cm} +\int_{\Omega_{\eta_0}}\bu_0\cdot\bphi(0,\cdot)\ dx + \im\eta_1b(0,\cdot)\ dA
%  \end{aligned}
% \end{equation}
% f"ur alle Testfunktionen $(b,\bphi)\in T^p_\eta$.
\end{definition}
\vspace{0.2cm}

% \marginpar{WIRBELTERM! DEN KOENNTE MAN HIER NOCH UMFORMEN
% ABER FUER DIE ANALOGE UMFORMULIERUNG AUF DER NAECHSTEN SEITE DENNOCH NOETIG!}
Gleichung \eqref{eqn:schwach} ergibt sich formal
durch Multiplikation von
\eqref{eqn:fluid} mit
einer Testfunktion $\bphi$, Integration "uber Ort und Zeit, partielle Integration und Verwenden von
\eqref{eqn:shell}. Genauer wird der Spannungstensor im Ort partiell integriert, Bemerkung
\ref{bem:korn} verwendet, das auftretende
Randintegral mittels \eqref{eqn:shell}$_1$ ersetzt und die Identit"at 
$(\grad_{L^2}K(\eta),b)_{L^2}=2K(\eta,b)$ ausgenutzt. Desweiteren werden die Terme mit
den ersten Zeitableitungen von $\bu$ und den zweiten Zeitableitungen von $\eta$ zeitlich partiell
integriert. Bei dem $\bu$-Term tritt dabei ein Randintegral auf, das sich mit Hilfe des
Reynolds'schen Transporttheorems berechnen l"asst. Mit $\bv=\bu=(\pa_t\eta\bnu)\circ
\Phi_{\eta(t)}^{-1}$ und $\xi=\bu\cdot\bphi$ erhalten wir 
\begin{equation*}
 \begin{aligned}
  \frac{d}{dt}\iot\bu\cdot\bphi\ dx &= \iot \pa_t\bu\cdot\bphi\ dx + \iot \bu\cdot\pa_t\bphi\
dx\\
&\hspace{0.5cm}+\int_{\pa\Omega_{\eta(t)}}\bu\cdot\bphi\
\big((\pa_t\eta\,\bnu)\circ\Phi_{\eta(t)}^{-1}\big)\cdot\bnu_{\eta(t)}\ dA_{\eta(t)}.
% &\hspace{0.5cm}+\int_{\pa\Omega_{\eta(t)}\setminus\Gamma}
% ((\pa_t\eta)^2\ b\bnu)\circ(\Phi_{\eta(t)})^{-1}\cdot\bnu_{\eta(t)}\ dA_{\eta(t)}dt
\end{aligned}
\end{equation*}
Gem"a"s Bemerkung \ref{bem:trafo} ist das Randintegral identisch
\[\im(\pa_t\eta)^2\, b\, \gamma(\eta)\ dA,\]
sodass sich \eqref{eqn:schwach} ergibt. Wir wollen uns noch davon "uberzeugen, dass der dritte Term
in \eqref{eqn:schwach} wohldefiniert ist. Es ist nicht schwer zu sehen, dass die Fortsetzung von
$\bphi$ durch $(b\,\bnu)\circ q$ in
$H^1(I\times B_\alpha)$, $\norm{\eta}_{L^\infty(I\times M)}<\alpha<\kappa$, liegt, sodass aus
Proposition \ref{theorem:bochneremb} folgt, dass $\bphi\in
L^\infty(I,L^3(\oet))$ gilt. W"are $\bu\in L^2(I,L^6(\oet))$, so w"urde dies gen"ugen, damit
$(\bu\cdot\nabla)\bu\cdot\bphi$ in $L^1(\Omega_\eta^I)$ liegt. Aufgrund der geringen
Randregularit"at gilt aber lediglich $\bu\in L^2(I,L^r(\oet))$ f"ur
alle $r<6$; siehe Korollar \ref{lemma:sobolev}. Deshalb ist die zus"atzliche Bedingung $\bphi\in
L^\infty(I,L^4(\oet))$ n"otig.
\vspace{0.2cm}

\begin{bemerkung}\label{bem:grenzen}
F"ur jede schwache L"osung $(\eta,\bu)$, alle Testfunktionen $(b,\bphi)\in T^I_\eta$ und fast alle
$t\in I$ gilt
\begin{equation}\label{eqn:schwacht}
\begin{aligned}
% &\iot \bu(t,\cdot)\cdot\bphi(t,\cdot)\ dx 
&- \int_{0}^t\ios \bu\cdot\pa_t\bphi\ dxds - \int_{0}^t\im(\pa_t\eta)^2\, b\, \gamma(\eta)\
dAds + \int_{0}^t\ios(\bu\cdot\nabla)\bu\cdot\bphi\ dxds\\
& + \int_{0}^t\ios \nabla\bu:\nabla\bphi\ dxds 
% +\im\pa_t\eta(t,\cdot)b(t,\cdot)\ dA
-\int_{0}^t\im\pa_t\eta\, \pa_tb\ dAds + 2\int_{0}^t K(\eta,b)\ ds \\
&\hspace{1.5cm}=\int_{0}^t\ios \ff\cdot\bphi\ dxds + \int_{0}^t\im g\, b\
dAds +\int_{\Omega_{\eta_0}}\bu_0\cdot\bphi(0,\cdot)\ dx\\
&\hspace{2.0cm}+ \im\eta_1\, b(0,\cdot)\ dA
-\int_{\Omega_{\eta(t)}}\bu(t,\cdot)\cdot\bphi(t,\cdot)\
dx -
\im\pa_t\eta(t,\cdot)\, b(t,\cdot)\ dA.
 \end{aligned}
\end{equation}
Tats"achlich k"onnen wir auf die Forderung $b(T,\cdot)=0$, $\bphi(T,\cdot)=0$ verzichten. Um
die Identit"at \eqref{eqn:schwacht} einzusehen, verwenden wir n"amlich in \eqref{eqn:schwach} die
Testfunktion $(b\,\rho^t_\epsilon,\bphi\,\rho^t_\epsilon)$, wobei $\rho\in C^\infty(\setR)$,
$\rho(s)=1$ f"ur $s\le0$ und $\rho(s)=0$ f"ur $s\ge 1$ sowie
$\rho^t_\epsilon(s)=\rho(\epsilon^{-1}(s-t))$. Dann gilt 
\begin{equation*}
 \begin{aligned}
  -\int_I\im\pa_t\eta\, &\pa_s(b\rho^t_\epsilon)\ dAds\\
&= -\int_I\rho^t_\epsilon\im\pa_t\eta\, \pa_tb\
dAds-\int_I\epsilon^{-1}\rho'(\epsilon^{-1}(s-t))\im\pa_t\eta\, b\ dAds\\
&\rightarrow -\int_0^t\im\pa_t\eta\, \pa_tb\ dAds+\im\pa_t\eta(t,\cdot)\, b(t,\cdot)\ dA
 \end{aligned}
\end{equation*}
f"ur $\epsilon \rightarrow 0$ und fast alle $t\in I$. Die Konvergenz des ersten Terms folgt aus dem
Satz "uber die dominierte Konvergenz, w"ahrend wir f"ur die Konvergenz des zweiten Terms den
Lebesgue'schen Differentiationssatz, der in Proposition \ref{theorem:faltung} enthalten ist, und die
Identit"at $\int_\setR \rho'\ ds=-1$ verwendet
haben. Analog l"asst sich die Konvergenz der anderen Terme in
\eqref{eqn:schwach} einsehen, sodass wir \eqref{eqn:schwacht} erhalten.
\end{bemerkung}
\vspace{0.2cm}

Wir kommen nun zum zentralen Existenzresultat dieser Arbeit.
\begin{theorem}\label{theorem:hs}
F"ur beliebige zul"assige Daten $(\ff,g,\bu_0,\eta_0,\eta_1)$ existieren ein $T^*\in (0,\infty]$
und ein Tupel $(\eta,\bu)$ derart, dass $(\eta,\bu)$ f"ur
alle $T<T^*$ eine schwache L"osung von \eqref{eqn:fluid}, \eqref{eqn:shell} und \eqref{eqn:data} auf
dem Intervall $I=(0,T)$ ist und die
Absch"atzung\footnote{Genau genommen verwenden wir auf $X_\eta^I$ die "aquivalente Norm
$\norm{\,\cdot\,}_{L^\infty(I,L^2(\Omega_{\eta(t)}))} +
\norm{\nabla\,\cdot\,}_{L^2(\Omega_{\eta}^I)}$.} 
\begin{equation}\label{ab:hs}
 \begin{aligned}
&\norm{\eta}_{Y^I}^2 + \norm{\bu}_{X_\eta^I}^2\le \big(\norm{\bu_0}_{L^2(\Omega_{\eta_0})}^2 +
\mnorm{\eta_1}^2 +
\norm{\eta_0}_{H^2(M)}^2\big)\,e^{cT}\\
&\hspace{4.5cm} + \int_0^T
\big(\norm{\ff(t,\cdot)}_{L^2(\Omega_{\eta(t)})}^2 + \mnorm{g(t,\cdot)}^2\big)\,e^{c(T-t)}\ dt\\
 \end{aligned}
\end{equation}
erf"ullt. Entweder gilt
$T^*=\infty$ oder $\lim_{t\rightarrow
T^*}\norm{\eta(t,\cdot)}_{L^\infty(M)}=\kappa$.
\end{theorem}
\vspace{0.2cm}

Die rechte Seite von \eqref{ab:hs} als Funktion von $T$, $\Omega_\eta^I$ und den Daten wird im
Folgenden wiederholt auftreten und deshalb kurz
mit $c_0(T,\Omega_\eta^I,\ff,g,\bu_0,\eta_0,\eta_1)$ bezeichnet. 

Ein wichtiger Schritt im Beweis von Theorem \ref{theorem:hs} ist der Nachweis der relativen
Kompaktheit beschr"ankter Folgen approximativer L"osungen in $L^2$. Im Falle der
Navier-Stokes-Gleichungen auf einem Raumzeitzylinder folgt diese sofort aus dem Satz von
Aubin-Lions, Proposition \ref{theorem:aubinlions}. In diesem Fall l"asst sich n"amlich direkt aus
der Gleichung ablesen,
dass die L"osungen beschr"ankte Zeitableitungen in $L^{4/3}(I,(H^{1}_{0,\dv})')$
besitzen. Auf den ersten Blick scheint eine "ahnliche Aussage auch f"ur unser gekoppeltes System zu
gelten. Betrachten wir \eqref{eqn:schwach} mit $b(0,\cdot)=0$ und $\bphi(0,\cdot)=0$, so erhalten
wir
\begin{equation*}
\begin{aligned}
% &\iot \bu(t,\cdot)\cdot\bphi(t,\cdot)\ dx 
- \int_I\iot &\bu\cdot\pa_t\bphi\ dxdt - \int_I\im(\pa_t\eta)^2\, b\, \gamma(\eta)\
dAdt -\int_I\im\pa_t\eta\, \pa_tb\ dAdt\\
& = -\int_I\iot(\bu\cdot\nabla)\bu\cdot\bphi\ dxdt -
\int_I\iot \nabla\bu:\nabla\bphi\ dxdt
% +\im\pa_t\eta(t,\cdot)b(t,\cdot)\ dA
 - 2\int_I K(\eta,b)\ dt\\
&\hspace{0.5cm} +\int_I\iot \ff\cdot\bphi\ dxdt + \int_I\im g\, b\
dAdt.
 \end{aligned}
\end{equation*}
Es ist nun naheliegend, die rechte Seite als Definition der Zeitableitung des Tupels
$(\pa_t\eta,\bu)$ zu interpretieren und zu folgern, dass diese durch die r"aumlichen Regularit"aten
der L"osung beschr"ankt in einem geeigneten Dualraum liegt. Dieser Idee eine rigorose Bedeutung zu
geben, ist allerdings schwierig. Erstens haben wir bislang kein Substitut f"ur Bochner-R"aume von
Funktionen mit Werten in Dualr"aumen konstruiert. Zweitens ist unklar, in welchem
konkreten Sinne man "uberhaupt von einer Zeitableitung des Tupels $(\pa_t\eta,\bu)$ sprechen
k"onnte. Vollends un"ubersichtlich wird die Situation schlie"slich durch den Umstand, dass die
involvierten Funktionenr"aume im Allgemeinen mit der L"osung variieren und dadurch Folgen von
(Dual-)R"aumen
auftreten. 

Angesichts dieser Schwierigkeiten wird in \cite{b27}, \cite{b28} ein anderer Zugang gew"ahlt; das
System wird mit zeitlichen Differenzen(quotienten) getestet. Die Konstruktion der richtigen
Testfunktionen ist allerdings subtil, da die L"osung zum Zeitpunkt $t_0$ im Allgemeinen keine
zul"assige Testfunktion zum Zeitpunkt $t_1\not=t_0$ ist. Auch die Analyse der resultierenden
Identit"at ist nicht einfach und ben"otigt viele Seiten langwieriger Absch"atzungen. Im ersten
Anlauf \cite{b27} erforderte dieser Zugang die Einf"uhrung eines D"ampfungsterms. In unserem Fall
einer allgemeinen Schalengeometrie ist die Konstruktion der richtigen Testfunktionen nochmals
deutlich schwieriger und bislang daran gescheitert, dass, wie bereits angemerkt, die "ublichen
Stetigkeitseigenschaften des L"osungsoperators der Divergenzgleichung in Gebieten, deren Rand nicht
Lipschitz-stetig ist, nicht gelten.

Wir wollen uns deshalb erneut dem nat"urlicher erscheinenden Zugang "uber die Zeitableitung
zuwenden. Man k"onnte auf die Idee kommen, die Kompaktheit der Fluidgeschwindigkeit mit Hilfe des
Satzes von Aubin-Lions zun"achst lokal, d.h. in kleinen Raumzeitzylindern, zu zeigen, und dann
"uber Interpolation und Spurbildung auf die globale Kompaktheit der Fluidgeschwindigkeit und die
Kompaktheit der Schalengeschwindigkeit zu schlie"sen. Durch die Inkompressibilit"at und die
resultierende unendliche Schallgeschwindigkeit ist unser System jedoch hochgradig nichtlokal, sodass
lokale Argumente an dieser Stelle vermutlich nicht zum Erfolg f"uhren k"onnen. Wir m"ussen das
Argument "uber die Zeitableitung also global durchf"uhren. Das ist tats"achlich m"oglich, und der im
Folgenden dargelegte Zugang ist au"serdem direkter und k"urzer als der in \cite{b27}, \cite{b28}.
Wir wollen den Beweis des Satzes von Aubin-Lions auf unsere Situation "uber"-tragen. Aus diesem
Grund
skizzieren wir nun die Beweisvariante, die in \cite{b58} zu finden ist. Diese basiert auf dem Satz
von Arzela-Ascoli und einem Interpolationsresultat. Die "Ubertragung anderer Beweisvarianten, siehe
zum Beispiel \cite{b57}, ist zwar auch m"oglich, jedoch scheint die Variante in \cite{b58} am besten
geeignet zu sein.

\begin{proposition}\label{theorem:aubinlions}
Es seien $I\subset\setR$ ein offenes, beschr"anktes Intervall und $1\le p\le\infty$, $1<r\le\infty$.
F"ur die Banach-R"aume $B_0$, $B$ und $B_1$ gelte
\[B_0\compactembedding B\embedding B_1.\]
Dann ist die Einbettung
\[W:=\big\{v\in L^p(I,B_0)\ \big|\ v'\in L^r(I,B_1)\big\}\embedding L^p(I,B)\]
kompakt.
\end{proposition}
\noindent {\bf Beweisskizze:}\ Die Folge $(v_n)\subset W$ sei beschr"ankt. Es gen"ugt zu zeigen,
dass eine Teilfolge in $C(\bar I,B_1)$ konvergiert, denn aus dem
Ehrling-Lemma erhalten wir f"ur jedes $\epsilon>0$ und $n,m\in\setN$ die Absch"atzung
\[\norm{v_n-v_m}_{L^p(I,B)}\le \epsilon\,\norm{v_n-v_m}_{L^p(I,B_0)} + c(\epsilon)\,
\norm{v_n-v_m}_{L^p(I,B_1)}.\]
Offenbar liegt $(v_n)$ beschr"ankt in $C^{0,1-1/r}(\bar I,B_1)$. Der Satz von Arzela-Ascoli liefert
somit die Konvergenz einer Teilfolge in $C(\bar I,B_1)$, wenn wir zeigen k"onnen, dass f"ur
alle $t$ aus einer dichten Teilmenge von $I$ die Folge $(v_n(t))_n$ relativ kompakt in $B_1$ liegt.
Das folgt aber aus den Einbettungen
\begin{equation*}
 \begin{aligned}
W\embedding C(\bar I,(B_0,B_1)_{\theta,1/\theta})\ \text{ und }\
(B_0,B_1)_{\theta,1/\theta}\compactembedding B_1
 \end{aligned}
\end{equation*}
f"ur ein geeignetes $0<\theta<1$; siehe Theorem $33$ in \cite{b59} und Theorem 3.8.1 in
\cite{b48}.
\qed\\ 

Weder wollen wir versuchen, ein Substitut f"ur den Raum $L^r(I,B_1)$ zu konstruieren, noch werden
wir den Begriff der Zeitableitung des Tupels $(\pa_t\eta,\bu)$ konkretisieren. Stattdessen werden
wir direkt mit der schwachen Formulierung, d.h. mit der getesteten formalen Zeitableitung, arbeiten
und dabei stets zeigen, dass die jeweiligen Aussagen unabh"angig von den Testfunktionen gelten,
solange diese in geeigneten Normen beschr"ankt bleiben. Die bereits angesprochene Problematik des
Auftretens von
Folgen formaler Dualr"aume bekommen wir durch Verwenden zweier spezieller Familien von
Testfunktionen in den Griff. Nat"urlich m"ussen wir auch den Beweis des Satzes von Arzela-Ascoli auf
die vorliegende Situation "ubertragen. Das Analogon der relativen Kompaktheit der Folgen
$(v_n(t))_n$ in $B_1$ erhalten wir hingegen sehr einfach aus der Schranke
f"ur die kinetischen Energien.

Die Aussage der folgenden Proposition wird im Beweis von Theorem \ref{theorem:hs} nicht verwendet.
Jedoch ist der Beweis der Proposition fast w"ortlich auf die im Beweis von Theorem \ref{theorem:hs}
auftretenden Situationen "ubertragbar.

\begin{proposition}\label{lemma:komp}
Es seien $(\ff,g,\bu_0^n,\eta_0^n,\eta_1^n)$ eine Folge zul"assiger Daten mit
\begin{equation}\label{eqn:abeschraenkt}
 \begin{aligned}
\sup_n\big(\tau(\eta_0^n)+\norm{\eta_0^n}_{H^2_0(M)}+\norm{\eta_1^n}_{L^2(M)}+\norm{
\bu^n_0}_{
L^2(\Omega_ {\eta_0^n})}\big)<\infty
 \end{aligned}
\end{equation} 
und $(\eta_n,\bu_n)$ eine Folge schwacher L"osungen von \eqref{eqn:fluid},
\eqref{eqn:shell} und \eqref{eqn:data} mit den Daten $(\ff,g,\bu_0^n,\eta_0^n,\eta_1^n)$
auf dem Intervall $I=(0,T)$ mit\,\footnote{F"ur zeitabh"angige Funktionen $\eta$ ersetzen wir  in
der Definition von $\tau(\eta)$ die Norm $\norm{\eta}_{L^\infty(\pa\Omega)}$ durch
$\norm{\eta}_{L^\infty(I\times \pa\Omega)}$.}
\begin{equation}\label{ab:beschraenkt}
 \begin{aligned}
\sup_n\big(\tau(\eta_n) + \norm{\eta_n}_{Y^I} +
\norm{\bu_n}_{X_{\eta_n}^I}\big)<\infty.
 \end{aligned}
\end{equation}
Dann liegt die Folge $(\pa_t\eta_n,\bu_n)$ relativ kompakt in $L^2(I\times M)\times L^2(I\times
\setR^3)$.\footnote{Wie immer setzen wir die Felder $\bu_n$ durch
$\boldsymbol{0}$ auf $I\times\setR^3$ fort.}
\end{proposition}
\beweis
Wegen \eqref{ab:beschraenkt} erhalten wir f"ur eine Teilfolge\footnote{Hier und im Rest der
Arbeit werden wir bei der Auswahl von Teilfolgen stets stillschweigend bez"uglich aller beteiligten
Folgen zu derselben Teilfolge "ubergehen und diese wieder mit $n$ indizieren.} die Konvergenzen
\begin{equation}\label{eqn:schwkonv}
 \begin{aligned}
  \eta_n&\rightarrow\eta &&\text{ schwach$^*$ in }L^\infty(I,H^2_0(M))\text{ und gleichm"a"sig},\\
  \pa_t\eta_n&\rightarrow\pa_t\eta &&\text{ schwach$^*$ in }L^\infty(I,L^2(M)),\\
  \bu_n&\rightarrow\bu &&\text{ schwach$^*$ in } L^\infty(I,L^2(\setR^3)),\\ 
  \nabla\bu_n&\rightarrow\xi &&\text{ schwach in } L^2(I\times \setR^3),
\end{aligned}
\end{equation}
wobei wir die zun"achst auf $\Omega_{\eta_n}^I$ definierten Felder
$\nabla\bu_n$ durch $0$ auf $I\times \setR^3$ fortsetzen. Wie man leicht "uberpr"uft ist der
Grenzwert $\xi$ nichts anderes als das Feld $\nabla\bu$, wenn wir dieses ebenfalls durch $0$
fortsetzen. Die gleichm"a"sige Konvergenz der Folge $(\eta_n)$ folgt aus der kompakten Einbettung $
(1/2<\theta<1)$
\[Y^I\embedding  C^{0,1-\theta}(\bar I, C^{0,2\theta -1}(\pa\Omega)) \compactembedding C(\bar
I\times\pa\Omega).\]
K"onnen wir die Konvergenz
\begin{equation}\label{eqn:finalkonv}
 \begin{aligned}
\int_I\int_{\Omega_{\eta_n(t)}}|\bu_n|^2\ dxdt + \int_I\im|\pa_t\eta_n|^2\ dAdt\rightarrow
\int_I\int_{\Omega_{\eta(t)}}|\bu|^2\ dxdt + \int_I\im|\pa_t\eta|^2\
dAdt
 \end{aligned}
\end{equation}
zeigen, so folgt die Behauptung der Proposition mit Hilfe
der Konvergenzen \eqref{eqn:schwkonv} und der trivialen Identit"at
\begin{equation*}
 \begin{aligned}
&\int_I\int_{\setR^3} |\bu_n-\bu|^2\ dxdt + \int_I\im |\pa_t\eta_n-\pa_t\eta|^2\ dAdt \\
& = \int_I\int_{\Omega_{\eta_n(t)}} |\bu_n|^2\ dxdt + \int_I\im |\pa_t\eta_n|^2\ dAdt +
\int_I\int_{\Omega_{\eta(t)}}
|\bu|^2\ dxdt + \int_I\im |\pa_t\eta|^2\ dAdt\\
&\hspace{5mm} - 2\int_I\int_{\setR^3} \bu_n\cdot\bu\ dxdt  -
2\int_I\im \pa_t\eta_n\ \pa_t\eta\ dAdt.
 \end{aligned}
\end{equation*}
\eqref{eqn:finalkonv} wiederum ist eine Konsequenz der Konvergenzen
\begin{equation}\label{eqn:l2konv1}
 \begin{aligned}
 \int_I\int_{\Omega_{\eta_n(t)}}\bu_n\cdot \F_{\eta_n}\pa_t\eta_n\ dxdt
+\int_I\im|\pa_t\eta_n|^2\ dAdt\\
\rightarrow \int_I\int_{\Omega_{\eta(t)}}\bu\cdot \F_{\eta}\pa_t\eta\ dxdt 
+ \int_I\im|\pa_t\eta|^2\ dAdt
 \end{aligned}
\end{equation}
und
\begin{equation}\label{eqn:l2konv2}
 \begin{aligned}
\int_I\int_{\Omega_{\eta_n(t)}} \bu_n\cdot(\bu_n -
\F_{\eta_n}\pa_t\eta_n)\ dxdt \rightarrow\int_I \int_{\Omega_{\eta(t)}} \bu\cdot(\bu -
\F_{\eta}\pa_t\eta)\ dxdt.
 \end{aligned}
\end{equation}
Die Zahl $\alpha$ in der Definition der Fortsetzungsoperatoren $\F$, siehe Proposition
\ref{lemma:FortVonRandZeit}, gen"uge dabei der Ungleichung
$\sup_n\norm{\eta_n}_{L^\infty(I\times M)}<\alpha<\kappa$. 

Die Beweise von \eqref{eqn:l2konv1} und \eqref{eqn:l2konv2} verlaufen sehr "ahnlich. Wir beginnen
mit \eqref{eqn:l2konv1}. Eine beliebige Funktion $b\in H^2_0(M)$ verletzt im Allgemeinen die
Mittelwertbedingung \eqref{eqn:integralnull} bez"uglich $\eta_n(t,\cdot)$ und ist somit nicht
divergenzfrei auf $\Omega_{\eta_n(t)}$ fortsetzbar. Aus diesem Grund ben"otigen wir die Operatoren
$\M_{\eta_n}$ aus Lemma \ref{lemma:mittelwert}. Aus diesem Lemma und Proposition
\ref{lemma:FortVonRandZeit} folgt unter Beachtung von
\eqref{ab:beschraenkt} die Absch"atzung
\begin{equation}\label{ab:b}
 \begin{aligned}
  \norm{\M_{\eta_n}b}_{H^1(I,L^2(M))\cap
L^2(I,H^2_0(M))}+\norm{\F_{\eta_n}\M_{\eta_n}b}_{H^1(I,L^2(B_\alpha))\cap
C(\bar I,H^1(B_\alpha))}\le c\,\norm{b}_{H^2_0(M)}.
 \end{aligned}
\end{equation}
%  der Einbettung
% \begin{equation*}
% H^1(I,L^2(M))\cap L^2(I,H^2_0(M)) \embedding C(\bar I,H^1(M)),
% \end{equation*}
% siehe Proposition \ref{theorem:bochneremb},\footnote{Die Identit"at
% $[L^2(M),H^2(M)]_\frac{1}{2}=H^1(M)$ folgt unter Verwendung eines endlichen Atlas mit
% untergeordneter Zerlegung der Eins und des in \cite{b55} konstruierten
% Fortsetzungsoperators aus Theorem 6.4.5 in \cite{b48}.}
Betrachten wir die Identit"at \eqref{eqn:schwacht} mit den L"osungen
$(\eta_n,\bu_n)$ und den Testfunktionen $(\M_{\eta_n}b,\F_{\eta_n}\M_{\eta_n}b)\in T_{\eta_n}^I$, so
folgern wir aus \eqref{ab:beschraenkt} und \eqref{ab:b}, dass die Integranden der auftretenden
zeitlichen Integrale unabh"angig von $b$ und $n$ in $L^{12/11}(I)$
beschr"ankt sind, sofern $\norm{b}_{H^2_0(M)}\le 1$ gilt.\footnote{Anstelle von "`unabh"angig von
$b$, sofern $\norm{b}_{H^2_0(M)}\le 1$"' schreiben wir fortan abk"urzend "`unabh"angig von
$\norm{b}_{H^2_0(M)}\le 1$"'.} Nur f"ur den Wirbelterm ist diese Behauptung nicht
offensichtlich. F"ur diesen folgt sie aus der Absch"atzung
\begin{equation*}
 \begin{aligned}
\Bignorm{\int_{\Omega_{\eta_n(t)}}(\bu_n\cdot\nabla)\bu_n&\cdot\F_{\eta_n}\M_{\eta_n}b\
dx}_{L^{12/11}(I)}\\
&\le \norm{\bu_n}_{L^{12/5}(I,L^4(\Omega_{\eta_n(t)}))}\,
\norm{\nabla\bu_n}_{L^{2}(\Omega_{\eta_n}^I)}\,
\norm{\F_{\eta_n}\M_{\eta_n}b}_{L^{\infty}(I,L^4(\Omega_{\eta_n(t)}))}  
 \end{aligned}
\end{equation*}
unter Beachtung der Sobolev-Einbettung
\[H^1(\Omega_{\eta_n(t)})\embedding L^5(\Omega_{\eta_n(t)}),\]
die gleichm"a"sig in $n$ und $t$ gilt, siehe Korollar \ref{lemma:sobolev}, und der
Interpolationseinbettung ($\theta=5/6$)
\[L^\infty(I,L^2(\Omega_{\eta_n(t)}))\cap L^2(I,L^5(\Omega_{\eta_n(t)}))\embedding
L^{12/5}(I,L^4(\Omega_{\eta_n(t)})).\]
Letztere ist eine einfache Folgerung aus der H"older'schen Ungleichung; vgl. Proposition
\ref{theorem:bochnerinterpol}. Die ersten acht Summanden der Identit"at \eqref{eqn:schwacht} 
sind somit unabh"angig von $\norm{b}_{H^2_0(M)}\le 1$ und $n$ in $C^{0,1/12}(\bar I)$ beschr"ankt.
Dieselbe Aussage gilt mithin f"ur die letzten beiden Summanden dieser Identit"at, d.h. f"ur die
Funktionen
\[c_{b,n}(t):=\int_{\Omega_{\eta_n(t)}}\bu_n(t,\cdot)\cdot(\F_{\eta_n}\M_{\eta_n}b)(t,\cdot)\ dx +
\im\pa_t\eta_n(t,\cdot)\ (\M_{\eta_n}b)(t,\cdot)\ dA,\]
da der neunte und zehnte Summand wegen \eqref{eqn:abeschraenkt} unabh"angig von
$\norm{b}_{H^2_0(M)}\le 1$ und $n$ beschr"ankte reelle Zahlen sind. Aufgrund
der Konvergenzen \eqref{eqn:schwkonv} und Lemma \ref{lemma:konvergenzen} $(1.a)$, $(2.a)$
konvergiert die Folge $(c_{b,n})_n$ f"ur festes $b\in H^2_0(M)$ im Distributionssinne gegen die
Funktion
\[c_{b}(t):= \int_{\Omega_{\eta(t)}}\bu(t,\cdot)\cdot(\F_\eta\M_\eta b)(t,\cdot)\ dx +
\im\pa_t\eta(t,\cdot)\, (\M_\eta b)(t,\cdot)\ dA.\]
Eine Anwendung des Satzes von Arzela-Ascoli zeigt, dass $(c_{b,n})_n$ sogar gleichm"a"sig in $\bar
I$ gegen $c_b$ konvergiert. 

Wir wollen nun zeigen, dass diese gleichm"a"sige Konvergenz unabh"angig von
$\norm{b}_{H^2_0(M)}\le 1$ ist, d.h. dass die Funktionen
\begin{equation*}
\begin{aligned}
&h_n(t):=&\sup_{\norm{b}_{H^2_0(M)}\le
1} \big(c_{b,n}(t)-c_{b}(t)\big)
\end{aligned}
\end{equation*} 
gleichm"a"sig in $\bar I$ gegen $0$ konvergieren. Wegen \eqref{ab:beschraenkt} gilt f"ur fast alle
$t\in I$
\[\sup_n \big(\norm{\bu_n(t)}_{L^2(\setR^3)}+\norm{\pa_t\eta_n(t)}_{L^2(M)}\big)<\infty.\] 
Mit Hilfe eines Diagonalfolgenarguments schlie"sen wir, dass eine abz"ahlbare, dichte Teilmenge
$I_0$ von $I$ und eine Teilfolge von $(\eta_n,\bu_n)$ existieren derart, dass $(\bu_n(t,\cdot))_n$
und $(\pa_t\eta_n(t,\cdot))_n$ f"ur alle $t\in I_0$ schwach in $L^2(\setR^3)$ bzw. $L^2(M)$
konvergieren. Wir zeigen nun, dass die Folge $(c_{b,n}(t))_n$ f"ur festes, aber beliebiges $t\in
I_0$ unabh"angig von $\norm{b}_{H^2_0(M)}\le 1$ konvergiert. Bezeichnet $\eta^*$ den schwachen
Grenzwert von $(\pa_t\eta_n(t,\cdot))_n$ in $L^2(M)$, so folgt
diese Behauptung f"ur den zweiten Summanden in $c_{b,n}(t)$ aus der Absch"atzung
\begin{equation*}
 \begin{aligned}
 &\big|\im \pa_t\eta_n(t,\cdot)\, (\M_{\eta_n} b)(t,\cdot)-\eta^*\,(\M_\eta b)(t,\cdot)\
dA\big|\\
&\hspace{5.0cm}\le \big|\im\big(\pa_t\eta_n(t,\cdot)-\eta^*\big)\, (\M_{\eta_n} b)(t,\cdot)\
dA\big|\\
&\hspace{5.5cm}+ \big|\im\eta^*\, \big((\M_{\eta_n}b)(t,\cdot)-(\M_{\eta}
b)(t,\cdot)\big)\ dA\big|\\
&\hspace{5.0cm}\le \norm{\pa_t\eta_n(t,\cdot)-\eta^*}_{(H^1(M))'}\,
\norm{(\M_{\eta_n}b)(t,\cdot)}_{H^1(M)}\\
&\hspace{5.5cm}+\norm{\eta^*}_{L^2(M)}\norm{(\M_{\eta_n}b)(t,\cdot)-(\M_{\eta}b)(t,\cdot)}_{L^2(M)}
 \end{aligned}
\end{equation*}
unter Beachtung von \eqref{ab:b}, Einbettung \eqref{eqn:einbschnurz}, Lemma
\ref{lemma:konvergenzen} $(1.a)$ und der kompakten Einbettung
\[L^2(M)\compactembedding (H^1(M))'.\]
Letztere folgt unter Verwendung des Satzes von Schauder "uber die Kompaktheit adjungierter
Operatoren aus den "ublichen Sobolev-Einbettungen. Bezeichnet $\bu^*$ den schwachen
Grenzwert von $(\bu_n(t,\cdot))_n$ in $L^2(\setR^3)$, so folgt
die Behauptung f"ur den zweiten Summanden in $c_{b,n}(t)$ vollkommen analog aus der Absch"atzung
\begin{equation*}
 \begin{aligned}
 &\big|\int_{\Omega_{\eta_n(t)}} \bu_n(t,\cdot)\cdot(\F_{\eta_n}\M_{\eta_n}
b)(t,\cdot)\ dx-\int_{\Omega_{\eta(t)}} \bu^*\cdot(\F_{\eta}\M_{\eta}
b)(t,\cdot)\ dx\big|\\
&\hspace{5.0cm}\le \big|\int_{B_\alpha}\big(\bu_n(t,\cdot)-\bu^*\big)\cdot(\F_{\eta_n}\M_{\eta_n}
b)(t,\cdot)\ dx\big|\\
&\hspace{5.5cm}+
\big|\int_{B_\alpha}\bu^*\cdot\big((\F_{\eta_n}\M_{\eta_n}b)(t,\cdot)-(\F_{\eta}\M_{\eta}
b)(t,\cdot)\big)\ dx\big|\\
&\hspace{5.0cm}\le \norm{\bu_n(t,\cdot)-\bu^*}_{(H^1(B_\alpha))'}\,
\norm{(\F_{\eta_n}\M_{\eta_n}b)(t,\cdot)}_{H^1(B_\alpha)}\\
&\hspace{5.5cm}+\norm{\bu^*}_{L^2(B_\alpha)}\norm{(\F_{\eta_n}\M_{\eta_n}b)(t,\cdot)-(\F_{\eta}\M_
{\eta}b)(t, \cdot)}_{L^2(B_\alpha)}
 \end{aligned}
\end{equation*}
unter Beachtung von \eqref{ab:b}, Lemma \ref{lemma:konvergenzen} $(2.a)$ und der
kompakten Einbettung
\begin{equation}\label{eqn:keinb}
 \begin{aligned}
 L^2(B_\alpha)\compactembedding (H^1(B_\alpha))'.
 \end{aligned}
\end{equation}
Wir zeigen schlie"slich, dass die von $\norm{b}_{H^2_0(M)}\le 1$ unabh"angige Konvergenz der Folge
$(c_{b,n}(t))_n$ nicht nur f"ur $t\in I_0$, sondern sogar gleichm"a"sig f"ur alle $t\in \bar I$
gilt. Aufgrund der gleichm"a"sigen Beschr"anktheit von $(c_{b,n})$ in
$C^{0,1/12}(\bar I)$ gilt f"ur alle $t,t'\in \bar I$ und $n,m\in\setN$
\begin{equation*}
 \begin{aligned}
  |c_{b,n}(t)-c_{b,m}(t)|&\le |c_{b,n}(t)-c_{b,n}(t')| + |c_{b,n}(t')-c_{b,m}(t')| +
|c_{b,m}(t')-c_{b,m}(t)|\\
&\le c\,|t-t'|^{1/12} + |c_{b,n}(t')-c_{b,m}(t')|. 
 \end{aligned}
\end{equation*}
Zu einem gegebenen $\epsilon>0$ k"onnen wir nun eine endliche Menge $I_0^\epsilon\subset I_0$
finden derart, dass zu jedem $t\in \bar I$ ein $t'\in I_0^\epsilon$ existiert mit $c\,|t-t'|^{1/12}
<
\epsilon/2.$ Zudem haben wir gerade gezeigt, dass $|c_{b,n}(t')-c_{b,m}(t')|<\epsilon/2$ gilt, falls
$t'\in I_0^\epsilon$ und $n,m\ge N$, wobei $N$ von $\epsilon$, aber nicht von $t'$ und
$\norm{b}_{H^2_0(M)}\le 1$ abh"angt. Somit konvergiert die Folge $(c_{b,n})_n$ gegen $c_b$
gleichm"a"sig in $\bar I$ unabh"angig von $\norm{b}_{H^2_0(M)}\le 1$, d.h. die Folge $(h_n)$
konvergiert gleichm"a"sig in $\bar I$ gegen $0$.

Mit der Definition
\begin{equation*}
\begin{aligned}
&g_n(t):=\sup_{\norm{b}_{L^2(M)}\le
1}\big(c_{b,n}(t)-c_{b}(t)\big)
 \end{aligned}
\end{equation*}
existiert wegen Lemma \ref{lemma:ehrling}, einer Aussage vom Typ des Ehrling-Lemmas, f"ur jedes
$\epsilon>0$ eine Konstante $c(\epsilon)$ derart, dass
\begin{equation*}
 \begin{aligned}
  \int_I g_n(t)\ dt \le \epsilon\,
c\,\big(\norm{\bu_n}_{L^2(I,H^1(\Omega_{\eta_n(t)}))}+\norm{\bu}_{L^2(I,H^1(\Omega_{\eta(t)}))}
\big) + c(\epsilon)\int_I h_n(t)\ dt
 \end{aligned}
\end{equation*}
gilt. Aufgrund der gleichm"a"sigen Konvergenz von $(h_n)$ gegen $0$ folgt
\begin{equation}\label{eqn:gn}
 \begin{aligned}
  \lim_n \int_I g_n(t)\ dt = 0.
 \end{aligned}
\end{equation}
Durch Nulladdition erhalten wir die Identit"at
\begin{equation}\label{eqn:nulladd}
 \begin{aligned}
&  \int_{\Omega_{\eta_n(t)}}\bu_n(t,\cdot)\cdot (\F_{\eta_n}\pa_t\eta_n)(t,\cdot)\ dx
+\im|\pa_t\eta_n(t,\cdot)|^2\ dA\\
& -\int_{\Omega_{\eta(t)}}\bu(t,\cdot)\cdot (\F_{\eta}\pa_t\eta)(t,\cdot)\ dx-
\im|\pa_t\eta(t,\cdot)|^2\ dA\\
&\hspace{1cm}= \int_{\Omega_{\eta_n(t)}}\bu_n(t,\cdot)\cdot (\F_{\eta_n}\pa_t\eta_n)(t,\cdot)\ dx +
\im\pa_t\eta_n(t,\cdot)\, \pa_t\eta_n(t,\cdot)\ dA\\
&\hspace{1.5cm}-\int_{\Omega_{\eta(t)}}\bu(t,\cdot)\cdot (\F_{\eta}\M_{\eta}\pa_t\eta_n)(t,\cdot)\
dx -\im\pa_t\eta(t,\cdot)\,
(\M_{\eta}\pa_t\eta_n)(t,\cdot)\ dA\\
&\hspace{1.5cm}+ \int_{\Omega_{\eta(t)}}\bu(t,\cdot)\cdot
\big((\F_{\eta}\M_{\eta}\pa_t\eta_n)(t,\cdot)-\F_{\eta}\pa_t\eta)(t,\cdot)\big)\ dx\\
&\hspace{1.5cm} + \im\pa_t\eta(t,\cdot)\,
\big((\M_\eta\pa_t\eta_n)(t,\cdot)-\pa_t\eta(t,\cdot)\big)\ dA.
 \end{aligned}
\end{equation}
Aufgrund der Gleichung $\M_{\eta_n}\pa_t\eta_n=\pa_t\eta_n$ sind die ersten beiden Zeilen der
rechten Seite identisch $c_{b,n}(t)-c_{b}(t)$ mit $b=\pa_t\eta_n(t,\cdot)$. Wegen
\eqref{ab:beschraenkt} ist ihr Betrag somit f"ur fast
alle $t$ durch $c\,g_n(t)$ abgesch"atzt. Integrieren wir \eqref{eqn:nulladd}
"uber $I$ und verwenden \eqref{eqn:gn} sowie die aus der schwachen Konvergenz von $(\pa_t\eta_n)$ in
$L^2(I\times M)$
resultierenden\footnote{Ein linearer Operator zwischen Banach-R"aumen ist genau dann stetig
bez"uglich der Normtopologien, wenn er
stetig ist bez"uglich der schwachen Topologien.} schwachen Konvergenzen von $(\M_\eta
\pa_t\eta_n)$
und $(\F_\eta\M_\eta
\pa_t\eta_n)$ in $L^2(I\times M)$ bzw. $L^2(I\times B_\alpha)$, so erhalten wir
\eqref{eqn:l2konv1}; man beachte dabei $\M_{\eta}\pa_t\eta=\pa_t\eta$.

Wir kommen nun zum Beweis von \eqref{eqn:l2konv2}. Es sei $\sigma>0$ hinreichend klein und
$\delta_\sigma\in C^4(\bar I\times \pa\Omega)$ mit $\norm{\delta_\sigma-\eta}_{L^\infty(I\times
\pa\Omega)}<\sigma$ und $\delta_\sigma<\eta$ in $\bar I\times\pa\Omega$. F"ur $\bphi\in
H(\Omega)$ setzen wir
\[c_{\bphi,n}^\sigma(t):=\int_{\Omega_{\eta_n(t)}}\bu_n(t,\cdot)\cdot \T_{\delta_\sigma(t)}\bphi\
dx,\quad
c_{\bphi}^\sigma(t):=\int_{\Omega_{\eta_n(t)}}\bu(t,\cdot)\cdot \T_{\delta_\sigma(t)}\bphi\ dx;\]
siehe Bemerkung \ref{bem:tdelta}. Wir k"onnen nun wie gehabt vorgehen, um zu zeigen, dass die
Funktionen
% \footnote{Bei der Konstruktion von $\Omega_{\eta(t_0)-\sigma}$ wird
% \emph{zuerst} $\eta(t_0,\cdot)$ durch $0$ auf $\pa\Omega$ fortgesetzt und \emph{anschlie"send}
% $\sigma$ subtrahiert.}
\[h_n^\sigma(t):=\sup_{\norm{\bphi}_{H^1_{0,\dv}(\Omega)}\le
1} \big(c_{\bphi,n}^\sigma(t)-c_{\bphi}^\sigma(t)\big)\]
gleichm"a"sig in $\bar I$  gegen $0$ konvergieren. Ist $\bphi\in H^1_{0,\dv}(\Omega)$ und setzen
wir das zun"achst auf $\Omega_{\delta_\sigma}^I$ definierte Feld $\T_{\delta_\sigma}\bphi$ durch
$\boldsymbol{0}$ auf $I\times B_\alpha$ fort, so gilt
\begin{equation}\label{eqn:schnarch}
 \begin{aligned}
\norm{\T_{\delta_\sigma}\bphi}_{H^1(I,L^2(B_\alpha))\cap C(\bar I,H^1(B_\alpha))}\le
c\,\norm{\bphi}_{H^1_{0,\dv}(\Omega)}.
 \end{aligned}
\end{equation}
Das folgt unter Beachtung der verschwindenden Randwerte von $\T_{\delta_\sigma}\bphi$ aus den
entsprechenden Absch"atzungen auf $\Omega_{\delta_\sigma}^I$; siehe Bemerkung \ref{bem:tdelta}. Die
Identit"at \eqref{eqn:schwacht} mit den L"osungen $(\eta_n,\bu_n)$ und den
Testfunktionen $(0,\T_{\delta_\sigma}\bphi)\in T^I_{\eta_n}$ ($n$ hinreichend gro"s) zeigt unter
Beachtung von \eqref{ab:beschraenkt} und \eqref{eqn:schnarch}, dass die Funktionen
$c_{\bphi,n}^\sigma$ unabh"angig von $\norm{\bphi}_{H^1_{0,\dv}(\Omega)}\le 1$
und $n$ in $C^{0,1/12}(\bar I)$ beschr"ankt sind. Aus den Konvergenzen
\eqref{eqn:schwkonv} und Arzela-Ascoli erhalten wir
wiederum die gleichm"a"sige Konvergenz in $\bar I$ von $(c_{\bphi,n}^\sigma)_n$ gegen
$c_{\bphi}^\sigma$. Die
kompakte Einbettung \eqref{eqn:keinb} zeigt zudem, dass die Folge $(c_{\bphi,n}^\sigma(t))_n$ f"ur
$t\in I_0$ unabh"angig von $\norm{\bphi}_{H^1_{0,\dv}(\Omega)}\le 1$ konvergiert. Wie zuvor
schlie"sen wir daraus, dass $(h_n^\sigma)_n$ gleichm"a"sig gegen $0$ konvergiert.

Unter Verwendung von Lemma \ref{lemma:ehrling} folgt
\begin{equation}\label{eqn:gsigma}
 \begin{aligned}
\lim_n\int_I g^\sigma_n(t)\ dt=0,
 \end{aligned}
\end{equation}
wobei
\[g_n^\sigma(t):=\sup_{\norm{\bphi}_{H(\Omega)}\le
1} \big(c_{\bphi,n}^\sigma(t)-c_{\bphi}^\sigma(t)\big).\]
Die $L^2(\Omega_{\eta_n(t)})$-Normen der Felder $\bu_n(t,\cdot)-(\F_{\eta_n}\pa_t\eta_n)(t,\cdot)\in
H(\Omega_{\eta_n(t)})$ sind wegen \eqref{ab:beschraenkt} f"ur fast alle $t$ unabh"angig von $t$ und
$n$ beschr"ankt. Gem"a"s Lemma \ref{lemma:divdichtglm} existiert also f"ur jedes $\epsilon > 0$ ein
$\sigma>0$ derart, dass sich f"ur fast alle $t$ und alle hinreichend gro"sen $n$ Felder
$\bpsi_{t,n}\in H(\Omega_{\eta_n(t)})$ mit
unabh"angig von $t$ und $n$ beschr"ankten $L^2(\Omega_{\eta_n(t)})$-Normen, $\supp
\bpsi_{t,n}\subset
\Omega_{\delta_\sigma(t)}$ und
\begin{equation*}
 \begin{aligned}
  \norm{\bu_n(t,\cdot)-(\F_{\eta_n}\pa_t\eta_n)(t,\cdot) -
\bpsi_{t,n}}_{(H^{1/4}(\setR^3))'} <
\epsilon
 \end{aligned}
\end{equation*}
finden lassen; insbesondere gilt $\bpsi_{t,n}\in H(\Omega_{\delta_\sigma(t)})$ mit unabh"angig von
$t$ und $n$ beschr"ankter Norm. Durch Nulladdition erhalten wir die Identit"at
\begin{equation}\label{eqn:nulladd2}
 \begin{aligned}
&  \int_{\Omega_{\eta_n(t)}} \bu_n(t,\cdot)\cdot\big(\bu_n(t,\cdot) -
(\F_{\eta_n}\pa_t\eta_n)(t,\cdot)\big)\ dx -
\int_{\Omega_{\eta(t)}} \bu(t,\cdot)\cdot\big(\bu(t,\cdot) - (\F_{\eta}\pa_t\eta)(t,\cdot)\big)\
dx\\
&= \int_{\Omega_{\eta_n(t)}} \bu(t,\cdot)\cdot\big(\bu_n(t,\cdot) -
(\F_{\eta_n}\pa_t\eta_n)(t,\cdot)\big)\, dx -
\int_{\Omega_{\eta(t)}} \bu(t,\cdot)\cdot\big(\bu(t,\cdot) - (\F_{\eta}\pa_t\eta)(t,\cdot)\big)\,
dx\\
&\hspace{0.5cm} + \int_{\Omega_{\eta_n(t)}} \bu_n(t,\cdot)\cdot\bpsi_{t,n}\ dx -
\int_{\Omega_{\eta(t)}} \bu(t,\cdot)\cdot\bpsi_{t,n}\ dx\\
&\hspace{0.5cm}+ \int_{\Omega_{\eta_n(t)}} \bu_n(t,\cdot)\cdot\big(\bu_n(t,\cdot) -
(\F_{\eta_n}\pa_t\eta_n)(t,\cdot)-
\bpsi_{t,n}\big)\ dx\\
&\hspace{0.5cm} - \int_{\Omega_{\eta(t)}} \bu(t,\cdot)\cdot\big(\bu_n(t,\cdot) -
(\F_{\eta_n}\pa_t\eta_n)(t,\cdot)-
\bpsi_{t,n}\big)\ dx.
 \end{aligned}
\end{equation}
Der Betrag der zweiten Zeile der rechten Seite ist durch $c\,
g_n^\sigma(t)$ abgesch"atzt, w"ahrend die Betr"age der letzten beiden Zeilen durch
\[\epsilon\, ( \norm{\bu_n(t,\cdot)}_{H^{1/4}(\setR^3)} +
\norm{\bu(t,\cdot)}_{H^{1/4}(\setR^3)})\] dominiert sind; man beachte Proposition
\ref{lemma:nullfort}. Integrieren wir \eqref{eqn:nulladd2} "uber $I$ und verwenden
\eqref{eqn:gsigma}
sowie die aus \eqref{eqn:schwkonv} und Lemma \ref{lemma:konvergenzen} $(2.c)$ folgende schwache
Konvergenz der
Folge $(\chi_{\Omega^I_{\eta_n}}(\bu_n - \F_{\eta_n}\pa_t\eta_n))$ gegen $\chi_{\Omega^I_{\eta}}(\bu
-
F_{\eta}\pa_t\eta)$ in $L^2(I\times\setR^3)$, so erhalten wir \eqref{eqn:l2konv2}.
\qed\\

Beim Versuch der Konstruktion einer L"osung st"o"st man auf das Problem, dass der
Definitionsbereich der L"osung von der L"osung selber abh"angt. Das macht die Verwendung eines
Galerkin-Ansatzes zun"achst unm"oglich, da die Ansatzfunktionen ebenso von der L"osung
abh"angen w"urden. Diesen Zirkel wollen wir durch ein Fixpunktargument aufbrechen. Da wir keine
perturbative Aussage, sondern die Existenz zeitlich \emph{globaler} L"osungen zeigen m"ochten, wird
uns ein topologisches Resultat, eine Variante des Schauder'schen
Fixpunktsatzes, weiterhelfen. 

Naheliegend mag nun folgendes Vorgehen sein. Wir geben eine Bewegung $\delta$ des Randes
vor und l"osen die Fluidgleichungen auf dem resultierenden variablen Gebiet unter Beachtung der
no-slip-Bedingung \eqref{eqn:fluid}$_{3,4}$; vgl. zum Beispiel \cite{b30}. Anschlie"send
l"osen wir die Schalengleichung
\eqref{eqn:shell} mit der vom Fluid auf den Rand ausge"ubten Kraft auf der rechten Seite
und
erhalten eine Randbewegung $\eta$. Schlie"slich versuchen wir zu zeigen, dass die Abbildung
$\delta\mapsto\eta$ einen Fixpunkt besitzt. Der Haken an dieser Idee ist, dass wir globale, starke
L"osungen konstruieren m"ussen, damit die Auswertung des Spannungstensors auf dem Rand sinnvoll ist.
Das ist m"oglich, wenn wir den Wirbelterm regularisieren. Da wir aber letztlich an schwachen
L"osungen interessiert sind, wollen wir uns diesen Mehraufwand sparen. Zudem ist fraglich, ob die
auf diese Weise gewonnene Funktion $\delta\mapsto\eta$ eine geeignete beschr"ankte Menge in sich
abbildet, weil die Bewegung des Randes nat"urlich in die Energieabsch"atzungen der L"osungen
eingeht. Es ist mithin geschickter, eine schwache Formulierung zu verwenden, die die Kopplung an
die Wellengleichung bereits enth"alt.
% Der Haken an dieser Idee ist die Tatsache, dass die
% Fluidgleichungen auf dem Gebiet mit vorgegebener Bewegung des Randes keine Energieabsch"atzungen
% zulassen. Das liegt darin begr"undet, dass der Energietransfer vom Rand in das Fluid a-priori
% nicht
% beschr"ankt ist. Energie ist bekanntlich gleich Kraft mal Weg. W"ahrend der Weg des sich
% bewegenden Randes vorgegeben ist, h"angt die Kraft, die dieser aus"uben muss, vom Fluid, also der
% L"osung der Gleichungen, ab. Das Problem des Fluids in einem vorgegebenen zeitlich variablen
% Gebiet
% mit no-slip-Bedingung ist vermutlich nicht wohlgestellt. Das d"urfte der Grund daf"ur sein, dass
% beim Studium von Fluidgleichungen in Gebieten mit vorgegebener Bewegung fast ausschlie"slich
% die unphysikalische Bedingung verschwindender Randwerte anstelle der no-slip-Bedingung Verwendung
% findet; siehe zum Beispiel REF.\marginpar{STIMMT NICHT! SIEHE INSBES. MIYAKAWA,TERMAMOTO!} Es ist
% mithin essentiell,
% die Kopplung mit der Wellengleichung beizubehalten, um die "`energetische Konsistenz"' der
% Gleichungen zu bewahren.
Diese "Uberlegung f"uhrt uns auf
die Aufgabe, zu geeignetem gegebenem $\delta$ mit $\delta(0,\cdot)=\eta_0$ Funktionen $\eta\in Y^I$
mit $\norm{\eta}_{L^\infty(I\times M)}<\kappa$ und $\eta(0,\cdot)=\eta_0$ sowie $\bu\in X_\delta^I$
mit $\trd \bu =
\pa_t\eta\,\bnu$ zu finden, die der Gleichung
\begin{equation}\label{eqn:zwischen}
\begin{aligned}
% &\iot \bu(t,\cdot)\cdot\bphi(t,\cdot)\ dx 
&- \int_I\int_{\Omega_{\delta(t)}} \bu\cdot\pa_t\bphi\ dxdt -
\int_I\im\pa_t\eta\,\pa_t\delta\, b\, \gamma(\delta)\
dAdt+\int_I\int_{\Omega_{\delta(t)}}(\bu\cdot\nabla)\bu\cdot\bphi\ dxdt\\
&+\int_I\int_{\Omega_{\delta(t)}} \nabla\bu:\nabla\bphi\ dxdt
% +\im\pa_t\eta(t,\cdot)b(t,\cdot)\ dA
-\int_I\im\pa_t\eta\, \pa_tb\ dAdt + 2\int_I K(\eta,b)\ dt \\
&\hspace{1cm}=\int_I\int_{\Omega_{\delta(t)}} \ff\cdot\bphi\ dxdt + \int_I\im g\, b\
dAdt+\int_{\Omega_{\eta_0}}\bu_0\cdot\bphi(0,\cdot)\ dx + \im\eta_1\, b(0,\cdot)\ dA
 \end{aligned}
\end{equation}
f"ur alle Testfunktionen $(b,\bphi)\in T_\delta^I$ gen"ugen. Um die Energieabsch"atzungen zu
bewahren, muss hier allerdings noch der Wirbelterm abge"andert
werden. Wenn wir n"amlich obige Gleichung formal mit $(\pa_t\eta,\bu)$ testen,
sehen wir (vgl. \eqref{eqn:wirbelid}), dass das "`erste $\bu$"' im Term
\[\int_I\int_{\Omega_{\delta(t)}}(\bu\cdot\nabla)\bu\cdot\bu\ dxdt\] am
Rand mit der Geschwindigkeit des Randes \"ubereinstimmen m"usste.
Wir k"onnen dieses Problem beheben, indem wir in \eqref{eqn:schwach} den Wirbelterm durch partielle
Integration umformen: 
\begin{equation*}
 \begin{aligned}
\int_I\iot(\bu\cdot\nabla)\bu\cdot\bphi\ dxdt =
\frac{1}{2}\int_I\iot(\bu\cdot\nabla)\bu\cdot\bphi\ dxdt &-
\frac{1}{2}\int_I\iot(\bu\cdot\nabla)\bphi\cdot\bu\ dxdt\\
&+\frac{1}{2}\int_I\im(\pa_t\eta)^2\,b\,\gamma(\eta)\ dAdt.
 \end{aligned}
\end{equation*}
Wenn wir die gegebene Randbewegung $\delta$ \emph{nach} dieser Umformung
einf"uhren, erhalten wir anstelle von \eqref{eqn:zwischen} die Gleichung
\begin{equation}\label{eqn:zwischen2}
\begin{aligned}
% &\iot \bu(t,\cdot)\cdot\bphi(t,\cdot)\ dx 
&- \int_I\int_{\Omega_{\delta(t)}} \bu\cdot\pa_t\bphi\ dxdt -
\frac{1}{2}\int_I\im\pa_t\eta\,\pa_t\delta\, b\,\gamma(\delta)\ dAdt\\
&+\frac{1}{2}\int_I\int_{\Omega_{\delta(t)}}(\bu\cdot\nabla)\bu\cdot\bphi\ dxdt -
\frac{1}{2}\int_I\int_{\Omega_{\delta(t)}}(\bu\cdot\nabla)\bphi\cdot\bu\ dxdt\\
&+\int_I\int_{\Omega_{\delta(t)}} \nabla\bu:\nabla\bphi\ dxdt -\int_I\im\pa_t\eta\,\pa_t b\ dAdt +
2\int_I K(\eta,b)\ dt \\
% +\im\pa_t\eta(t,\cdot)b(t,\cdot)\ dA
&\hspace{4.5cm}=\int_I\int_{\Omega_{\delta(t)}} \ff\cdot\bphi\ dxdt + \int_I\im g\ b\
dAdt\\
&\hspace{5.0cm} +\int_{\Omega_{\eta_0}}\bu_0\cdot\bphi(0,\cdot)\ dx + \im\eta_1\, b(0,\cdot)\ dA.
 \end{aligned}
\end{equation}
Wenn wir diese formal mit $(\pa_t\eta,\bu)$ testen, heben sich die beiden "`Wirbelterme"'
gegenseitig weg, w"ahrend sich die ersten beiden Terme wie zuvor zu 
\[\int_I\frac{1}{2}\frac{d}{dt}\int_{\Omega_{\delta(t)}}|\bu|^2\ dxdt\]
zusammenfassen lassen.

Schlie"slich m"ussen wir noch eine Regularisierung von \eqref{eqn:zwischen2} vornehmen. Zun"achst
ben"otigen wir
eine Regularisierung der Randbewegung $\delta$, weil diese in einem Raum liegen muss,
in den $Y^I$ kompakt einbettet, und deshalb nicht hinreichend regul"ar sein kann. Zum Beispiel
w"aren
die Ansatzfunktionen, die wir in K"urze konstruieren werden, ohne eine Regularisierung von
$\delta$ unbrauchbar. Desweiteren ist eine Entsch"arfung des Wirbelterms (durch eine
Regularisierung des
"`ersten $\bu$"') unerl"asslich. Bekanntlich ist die Eindeutigkeit schwacher L"osungen
der Navier-Stokes-Gleichungen ein offenes Problem, sodass wir bez"uglich \eqref{eqn:zwischen2} erst
recht keine Eindeutigkeit erwarten k"onnen. Das ist ein Problem, weil wir dann
nicht wissen, auf welche der L"osungen die Funktion $\delta\mapsto\eta$, f"ur die wir die Existenz
eines Fixpunktes zeigen wollen, abbilden soll. Durch eine geeignete Entsch"arfung des Wirbelterms
wird dieses Problem gel"ost. Andererseits ist die Eindeutigkeit von L"osungen von
\eqref{eqn:zwischen2} auch mit entsch"arftem Wirbelterm nicht offensichtlich, und wir w"urden gerne
auf einen Beweis verzichten.\footnote{H"angt allerdings der Spannungstensor nichtlinear
vom Scherratentensor $D\bu$ ab, so ist eine Untersuchung in diese Richtung unerl"asslich.}
Stattdessen wollen wir eine
Variante des Schauder'schen Fixpunktsatzes f"ur mengenwertige Abbildungen, Theorem
\ref{theorem:kakutani}, verwenden.
Dazu muss allerdings die Menge der L"osungen zu gegebenen, festen Daten konvex sein. Diese
Eigenschaft
k"onnen wir erzwingen, indem wir den Wirbelterm "`linearisieren"'. Genauer: Wir f"uhren das
Fixpunktargument nicht nur in $\eta$, sondern zus"atzlich im ersten
Argument des Wirbelterms durch, d.h. wir ersetzen $(\bu\cdot\nabla)\bu$ durch $(\bv\cdot\nabla)\bu$,
wobei $\bv$ ein gegebenes Feld aus einem Funktionenraum ist, in dem wir Kompaktheit
der L"osungen der Fluidgleichungen bekommen k"onnen. Da also die Regularit"at von $\bv$ gering sein
wird, ist nach wie vor eine Regularisierung des ersten Arguments des Wirbelterms notwendig. Die
resultierende schwache Formulierung entspricht im Wesentlichen der, die in \cite{b27}
verwendet wird.

Wir wollen nun entsprechende Regularisierungsoperatoren $\R_\epsilon$, $\epsilon>0$, konstruieren.
% Da wir eine glatte
% Auslenkung des Randes w"unschen, m"ussen wir insbesondere die Anfangsauslenkung $\eta_0$ gl"atten.
% Wir nehmen zun"achst an, dass der Tr"ager von $\eta_0$ im Inneren vom $M$ enthalten sei. Dann
% k"onnen wir $\eta_0$ durch Funktionen $(\R_\epsilon\eta_0)\subset C_0^\infty(\inn M)$ in
% $H^2_0(M)$
% von
% oben approximieren, d.h. es gelte $\eta_0\le\R_\epsilon\eta_0<\kappa$ in $M$.\marginpar{WIE
% ALLGEMEIN?
% SIEHE S.104!} 
Wir w"ahlen dazu einen Gl"attungskern $\omega\in C_0^\infty(\setR^3)$
mit $\int_{\setR^3}\omega\ dx=1$ und $\supp\omega\subset\{(t,x)\in\setR^3\ |\ t>0\}$ und setzen
$\omega_\epsilon:=\epsilon^{-3}\omega(\epsilon^{-1}\cdot)$. 
% sowie $\widetilde\omega\in
% C_0^\infty((-1,0))$ mit $\int_{\setR}\widetilde\omega\ dx=1$
% und $\widetilde\omega_\epsilon:=\epsilon^{-1}\widetilde\omega(\epsilon^{-1}\cdot)$.
Zudem sei $(\phi_k,U_k)_k$ ein endlicher Atlas von $\pa\Omega$ mit untergeordneter
Zerlegung
der Eins $(\psi_k)_k$; siehe Anhang A.1. Funktionen $\delta\in C(\bar I\times\pa\Omega)$ mit
$\delta(0,\cdot)=\eta_0$ setzen wir durch $\eta_0$ auf $(-\infty,T]\times\pa\Omega$ fort, und wir
definieren
\begin{equation}\label{eqn:reg}
 \begin{aligned}
\R_\epsilon\delta:=\sum_k(\omega_\epsilon\ast((\psi_k\,
\delta)\circ\phi^{-1}_k))\circ\phi_k +\epsilon^{1/2}.
 \end{aligned}
\end{equation}
Der mit $k$ indizierte Summand sei dabei durch $0$ ausserhalb $U_{k}$ fortgesetzt.
Man beachte, dass $\R_\epsilon\eta_0:=(\R_\epsilon\delta)(0,\cdot)$ aufgrund der speziellen
Lokalisierung des Tr"agers von $\omega$ nur von $\eta_0$ (und $\epsilon$) abh"angt. Wegen
der
H"olderstetigkeit von $\eta_0$ gilt zudem $\R_\epsilon\eta_0\ge\eta_0$ f"ur alle $0<\epsilon\le
c(\eta_0)$ mit einer nur von $\eta_0$
abh"angenden (kleinen) Konstante $c(\eta_0)$. Das folgt aus der Absch"atzung
\begin{equation*}
 \begin{aligned}
  |\omega_\epsilon\ast((\psi_k\,
\delta)\circ\phi^{-1}_k)&-(\psi_k\,
\delta)\circ\phi^{-1}_k|(0,x)\\
&=|\int_{\setR^4}\omega_\epsilon(-s,x-y)\Big(((\psi_k\,
\eta_0)\circ\phi^{-1}_k)(y)-((\psi_k\,
\eta_0)\circ\phi^{-1}_k)(x)\Big)\ dyds|\\
&\le c\int_{\setR^4}|\omega_\epsilon(-s,x-y)||x-y|^{3/4}\ dyds=c\ \epsilon^{3/4}.
 \end{aligned}
\end{equation*}
F"ur die Ungleichung haben wir die $3/4$-H"olderstetigkeit von $\eta_0$ und die Regularit"at von
$\phi_k$, $\psi_k$ verwendet. F"ur hinreichend kleines
$\epsilon$ ist der Term \[(\omega_\epsilon\ast((\psi_k\,
\delta)\circ\phi^{-1}_k))\circ\phi_k\]
in $C^4(\bar I\times\pa\Omega)$, weil dann der Tr"ager in $\bar I\times U_k$ enthalten ist.
Aus Proposition \ref{theorem:faltung} folgt die Konvergenz von $(\omega_\epsilon\ast((\psi_k\,
\delta)\circ\phi^{-1}_k))$ gegen $(\psi_k\,
\delta)\circ\phi^{-1}_k$ in $L^\infty(I\times \setR^2)$, woraus durch Aufsummieren sofort
die
Konvergenz von $(\R_\epsilon\delta)$ gegen $\delta$ in $L^\infty(I\times \pa\Omega)$ folgt.
Falls
$\delta$ eine Zeitableitung $\pa_t\delta$ in $L^2(I\times\pa\Omega)$ besitzt, folgt wegen
$\pa_t\R_\epsilon\delta=\R_\epsilon\pa_t\delta$ die Konvergenz von
$(\pa_t\R_\epsilon\delta)$ gegen $\pa_t\delta$ in $L^2(I\times \pa\Omega)$ vollkommen
analog aus Proposition \ref{theorem:faltung}. Es seien schlie"slich $\widetilde\omega\in
C_0^\infty(\setR^4)$ mit $\int_{\setR^4}\omega\ dx=1$
und $\widetilde\omega_\epsilon:=\epsilon^{-4}\widetilde\omega(\epsilon^{-1}\cdot)$. F"ur
Funktionen $v\in L^2(I\times\setR^3)$ setzen wir $\R_\epsilon
v:=\widetilde\omega_\epsilon\ast v$, wobei $v$ durch $0$ auf $\setR^4$ fortgesetzt sei. Wegen
Proposition \ref{theorem:faltung} konvergiert $(\R_\epsilon v)$ gegen $v$ in $L^2(I\times\setR^3)$.

Der Definitionsbereich der Anfangsgeschwindigkeit $\bu_0$ des Fluids ist (im Allgemeinen)
von $\Omega_{\R_\epsilon\eta_0}$
verschieden. Wir m"ussen deshalb die Anfangswerte modifizieren. Die Tatsache, dass
$(\R_\epsilon\eta_0)$ die Funktion $\eta_0$ von oben approximiert, erleichtert uns dabei die
Arbeit. Wir setzen $\eta_1$ gem"a"s Bemerkung \ref{bem:nspur} zu einem divergenzfreien
$L^2$-Vektorfeld $\bpsi$ auf $B_\alpha$ fort, wobei $\alpha$ eine Zahl mit
$\norm{\eta_0}_{L^\infty(M)}<\alpha<\kappa$, und definieren 
\[\bu_0^\epsilon:=\left\{\begin{array}{cl} \bu_0 & \text{in }\Omega_{\eta_0}, \\
\bpsi & \text{in } \Omega_{\R_\epsilon\eta_0}\setminus\Omega_{\eta_0}\end{array}\right..\]
Offenbar gilt $\bu_0^\epsilon\in L^2(\Omega_{\R_\epsilon\eta_0})$, und aus \eqref{eqn:nspur} und 
$\trnormaln\bu_0=\eta_1\gamma(\eta_0)$ folgern wir mit Hilfe von Proposition
\ref{lemma:nspur} die Divergenzfreiheit von $\bu_0^\epsilon$. Mit der
Definition
\[\eta_1^\epsilon:= \exp\Big(\int_{\eta_0}^{\R_\epsilon\eta_0} \beta(\cdot+\tau\bnu)\ d\tau\Big)\
\eta_1\]
gilt zudem $\trnormale \bu_0^\epsilon=\eta_1^\epsilon\,\gamma(\R_\epsilon\eta_0)$, was sich anhand
der Definition von $\bpsi$ leicht
durch ein Approximationsargument wie im Beweis von Proposition \ref{lemma:FortVonRand} einsehen
l"asst.
Setzen wir $\bu_0$ und $\bu_0^\epsilon$ durch $\boldsymbol{0}$ auf $\setR^3$ fort, so folgen aus
der gleichm"a"sigen Konvergenz von $\R_\epsilon\eta_0\rightarrow\eta_0$ in $\pa\Omega$ die
Konvergenzen
\begin{equation}\label{eqn:konv7}
 \begin{aligned}
  \bu_0^\epsilon\rightarrow\bu_0 &\text{ in } L^2(\setR^3),\\
  \eta_1^\epsilon\rightarrow\eta_1&\text{ in } L^2(M).
 \end{aligned}
\end{equation}

Wir nehmen im Folgenden ein Intervall $I=(0,T)$, $T>0$, und beliebige, aber feste Funktionen $\bv\in
L^2(I\times \setR^3)$ und $\delta\in C(\bar I\times \pa\Omega)$ mit
$\norm{\delta}_{L^\infty(I\times\pa\Omega)}<\kappa$ und $\delta(0,\cdot)=\eta_0$ als
gegeben an. Unsere
Vor"uberlegungen f"uhren uns auf die folgende Definition. Der "Ubersichtlichkeit
halber werden wir vorerst den Parameter $\epsilon$ in der Notation unterdr"ucken. Insbesondere
stehen die Bezeichner $\bu_0$, $\eta_1$ vorerst f"ur die
modifizierten Anfangswerte $\bu_0^\epsilon$, $\eta_1^\epsilon$.

\begin{definition}\label{def:ent}
Ein Tupel $(\eta,\bu)$ hei"st schwache L"osung des entkoppelten, regularisierten
Systems zum Argument $(\delta,\bv)$ auf dem Intervall $I$, falls $\eta\in Y^I$
mit $\eta(0,\cdot)=\eta_0$, $\bu\in X_{\R\delta}^I$ mit $\trrd
\bu =\pa_t\eta\,\bnu$ und 
\begin{equation}\label{eqn:ent}
\begin{aligned}
% &\iot \bu(t,\cdot)\cdot\bphi(t,\cdot)\ dx 
&- \int_I\iortd \bu\cdot\pa_t\bphi\ dxdt - \frac{1}{2}\int_I\im\pa_t\eta\, 
\pa_t\mathcal{R}\delta\, b\, \gamma(\mathcal{R}\delta)\ dAdt\\
&+\frac{1}{2}\int_I\iortd(\mathcal{R}\bv\cdot\nabla)\bu\cdot\bphi\ dxdt -
\frac{1}{2}\int_I\iortd(\mathcal{R}\bv\cdot\nabla)\bphi\cdot\bu\ dxdt\\
& + \int_I\iortd \nabla\bu:\nabla\bphi\ dxdt -\int_I\im\pa_t\eta\, \pa_tb\ dAdt + 2\int_I K(\eta,b)\
dt\\
% +\im\pa_t\eta(t,\cdot)b(t,\cdot)\ dA
&\hspace{4.5cm}=\int_I\iortd \ff\cdot\bphi\ dxdt + \int_I\im g\, b\
dAdt\\
&\hspace{5cm} +\int_{\Omega_{\mathcal{R}\eta_0}}\bu_0\cdot\bphi(0,\cdot)\ dx +
\im \eta_1\, b(0,\cdot)\ dA
\end{aligned}
\end{equation}
f"ur alle Testfunktionen $(b,\bphi)\in T_{\R\delta}^I$.
\end{definition}
\vspace{0.2cm}

Die Existenz einer L"osung des entkoppelten, regularisierten Systems wird in \cite{b27} im Falle
einer einfachen Geometrie durch Transformation auf einen Raumzeitzylinder und einer anschlie"senden
Galerkin-Approximation gezeigt. Wie der Beweis der folgenden Proposition zeigt, ist die
Transformation
nicht sinnvoll, da sie f"ur die Galerkin-Approximation irrelevant ist und die nat"urliche Struktur
der Gleichungen zerst"ort. Die Existenz eines Diffeomorphismus des zeitlich variablen Gebiets auf
einen Raumzeitzylinder ist aber auch hier n"otig, und zwar zur Konstruktion der zeitlich variablen
Ansatzfunktionen.

\begin{proposition}\label{lemma:ent}
Es existiert eine schwache L"osung $(\eta,\bu)$ des entkoppelten, regularisierten
Systems zum Argument $(\delta,\bv)$ auf dem Intervall $I$, die die Absch"atzung
\begin{equation}\label{ab:ent}
 \begin{aligned}
 &\norm{\eta}_{Y^I}^2 + \norm{\bu}_{X_{\R\delta}^I}^2\le
c_0(T,\Omega_{\R\delta}^I,\ff,g,\eta_0,\eta_1)
  \end{aligned}
\end{equation}
erf"ullt. Insbesondere ist die linke Seite unabh"angig
vom Parameter $\epsilon$ und vom Argument $(\delta,\bv)$ beschr"ankt.
\end{proposition}

\noindent{\bf Vorarbeit.} Wir wollen zun"achst geeignete Ansatzfunktionen konstruieren. Dazu w"ahlen
wir eine Basis
$(\widehat\bX_k)_{k\in\setN}$ von $H^1_{0,\dv}(\Omega)$ und eine Basis
$(\widehat Y_k)_{k\in\setN}$ des Raumes \[\Big\{Y\in H^2_0(M)\big|\int_M Y\ dA =
0\Big\}.\] Durch L"osen des Stokes-Systems in $\Omega$ mit Randwerten $\widehat Y_k\,\bnu$
(wie immer durch 0 auf $\pa\Omega$ fortgesetzt), siehe Theorem \ref{theorem:stokes}, erhalten wir
divergenzfreie
Fortsetzungen $\widehat\bY_k$. F"ur $t\in\bar I$ setzen wir
\[\bX_k(t,\cdot):=\T_{\R\delta(t)}\widehat\bX_k,\quad
\bY_k(t,\cdot):=\T_{\R\delta(t)}\widehat\bY_k;\] siehe Bemerkung \ref{bem:tdelta}. Die Felder
$\bX_k(t,\cdot)$ bilden
offenbar eine Basis von $H^1_{0,\dv}(\Omega_{\R\delta(t)})$. Man beachte zudem, dass f"ur
$q\in\pa\Omega$ das Differential $d\Psi_{\R\delta(t)}(q)$ die Normale $\bnu(q)$ lediglich skaliert,
weshalb die Definition
\[Y_k(t,\cdot)\,\bnu:=\trrdt
\bY_k(t,\cdot)=d\Psi_{\R\delta(t)}\, (\det
d\Psi_{\R\delta(t)})^{-1}\,\widehat Y_k\,\bnu\]
sinnvoll ist. Aus der Identit"at $\trnormald
\bY_k=Y_k\,\gamma(\R\delta)$ und Proposition \ref{lemma:partInt} schlie"sen wir 
\[\int_M Y_k(t,\cdot)\,\gamma(\R\delta(t,\cdot))\ dA = 0.\]
Aus der Basiseigenschaft der $\widehat Y_k$ folgt, dass die Felder $Y_k(t,\cdot)$ eine
Basis des Raumes
\begin{equation}\label{eqn:raum}
 \begin{aligned}
\Big\{Y\in H^2_0(M)\ \big|\ \int_M
Y\, \gamma(\R\delta(t,\cdot))\ dA = 0\Big\}  
 \end{aligned}
\end{equation}
bilden. Um die Notation zu vereinfachen w"ahlen wir eine Aufz"ahlung\footnote{Wir k"onnen zum
Beispiel $\bW_{2k}:=\bX_k$ und $\bW_{2k-1}:=\bY_k$ w"ahlen.}
$(\bW_k)_{k\in\setN}$ der Felder $\bX_k,\bY_k$ und setzen $W_k\,\bnu:=\trrd
\bW_k$. 

Wir wollen nun zeigen, dass
\[\spann\{(\phi\, W_k,\phi\,\bW_k)\ |\ \phi\in
C_0^1([0,T)),\, k\in\setN\}\]
dicht liegt im Raum aller Tupel
\[(b,\bphi)\in \big(H^1(I,L^2(M))\cap L^2(I,H^2_0(M))\big)\times
H^1(\Omega_{\R\delta}^I)\]
mit $b(T,\cdot)=0$, $\bphi(T,\cdot)=0$, $\dv\bphi=0$ und $\trrd\bphi=b\,\bnu$. Offenbar bettet
$T_{\R\delta}^I$ in diesen Raum ein. Aufgrund der durch $\T_{\R\delta}$
induzierten Isomorphismen ist die Behauptung "aquivalent zur Dichtheitheit von
\[\spann\{(\phi\,\widehat W_k,\phi\,\widehat\bW_k)\ |\ \phi\in
C_0^1([0,T)),\, k\in\setN\}\] im Raum $T$ aller Tupel 
\[(b,\bphi)\in H^1(I,L^2(M))\cap L^2(I,H^2_0(M))\times H^1(I\times\Omega)\]
mit $b(T,\cdot)=0$, $\bphi(T,\cdot)=\boldsymbol{0}$, $\dv\bphi=0$ und $\bphi|_{I\times\pa\Omega} =
b\,\bnu$. Sei also $(b,\bphi)\in T$. Wir approximieren zun"achst $b$ durch Funktionen $\tilde
b\in C^\infty_0([0,T),H^2_0(M))$ in $H^1(I,L^2(M))\cap L^2(I,H^2_0(M))$ mit 
\[\int_M \tilde b(t,\cdot)\ dA=0\] f"ur $t\in \bar I$.\footnote{Wir k"onnen $\tilde b$ konstruieren,
indem wir die
vektorwertige Funktion $b$ durch $0$ auf $[0,\infty)$ fortsetzen und anschlie"send $b(\cdot+h)$,
$h>0$, mit einem Gl"attungskern falten. Die Mittelwertfreiheit von $\tilde b(t,\cdot)$ folgt dann
aus der Mittelwertfreiheit von $b(t,\cdot)$.}
Nun k"onnen wir $\pa_t\tilde b$ durch eine Folge $(\phi^k_n\,
\widehat Y_k)_n$ (Summation von $1$ bis $n$), $\phi^k_n\in
C_0^1([0,T))$, in $L^2(I,H^2_0(M))$ approximieren. Ist n"amlich $f\in L^2(I,H^2_0(M))$ mit
\[\int_M f(t,\cdot)\ dA=0\]
f"ur fast alle $t\in I$ und
\[\int_I\phi(t)\, (\widehat Y_k,f(t,\cdot))_{H^2(M)}\ dt=0\]
f"ur jedes $\phi\in C_0^1([0,T))$ und jedes $k\in\setN$, so verschwindet das
Skalarprodukt im Integranden fast "uberall, und aufgrund der Basiseigenschaft der Funktionen
$\widehat Y_k(t,\cdot)$ ist $f$ somit die Nullfunktion. Wegen
\begin{equation*}
 \begin{aligned}
 \Bignorm{\tilde b(t,\cdot) + \int_t^T\phi^k_n(s)\, ds\,
\widehat Y_k}_{H^2(M)}&\le \int_0^T \norm{\pa_s\tilde b(s,\cdot) - \phi^k_n(s)
\,\widehat Y_k}_{H^2(M)}\ ds
 \end{aligned}
\end{equation*}
konvergiert die Folge $(-\int_t^T\phi^k_n(s) ds\,
\widehat Y_k)_{n}$ gegen $\tilde b$ in $C(\bar I, H^2_0(M))$, insgesamt also in
$H^1(I,H^2_0(M))$. Wir haben somit Linearkombinationen der Felder $\widehat Y_k$ konstruiert, die
$b$ in $H^1(I,L^2(M))\cap L^2(I,H^2_0(M))$ approximieren. Aufgrund der
Stetigkeitseigenschaften des L"osungoperators des Stokes-Systems konvergieren die entsprechenden
Linearkombinationen der Felder $\widehat\bY_k$ gegen ein $\bY$ in $H^1(I\times\Omega)$.
Es gilt $(\bY-\bphi)|_{I\times\pa\Omega}=0$, sodass lediglich zu zeigen bleibt, dass wir jedes
$\bphi$
mit $(0,\bphi)\in T$ durch eine Folge $(\phi^k_n\,
\widehat \bX_k)_{n}$, $\phi^k_n\in C_0^1([0,T))$, in $H^1(I\times \Omega)$ approximieren
k"onnen. Dazu k"onnen wir aber genau wie bei der Approximation von $b$ vorgehen.\\

\beweis (von Proposition \ref{lemma:ent}) Wir verwenden die Galerkin-Methode,
d.h. wir approximieren das vorliegende
unendlichdimensionale dynamische System durch endlichdimensionale Systeme, also durch gew"ohnliche
Differentialgleichungen, indem wir die Gleichung auf endlichdimensionale Teilr"aume der jeweiligen
Funktionenr"aume "`projizieren"'.
Wir suchen dazu Funktionen $\alpha_n^k:[0,T]\rightarrow\setR$, $k,n\in\setN$, derart, dass
$\bu_n:=\alpha_n^k\,\bW_k$ und $\eta_n(t,\cdot):=\int_0^t\alpha^k_n\, W_k\ ds+\eta_0$ (jeweils
Summation
von $1$ bis $n$) die Gleichung\footnote{Wie "ublich unterdr"ucken wir die unabh"angigen Variablen
in der Notation, d.h. $\bu=\bu(t,\cdot)$, etc.}
\begin{equation}\label{eqn:gal}
\begin{aligned}
% &\iot \bu(t,\cdot)\cdot\bphi(t,\cdot)\ dx 
&\iortd \pa_t\bu_n\cdot\bW_j\ dx +\frac{1}{2} \im\pa_t\eta_n
\, \pa_t\mathcal{R}\delta\, W_j\, \gamma(\mathcal{R}\delta)\ dA\\
&+\frac{1}{2}\iortd(\mathcal{R}\bv\cdot\nabla)\bu_n\cdot\bW_j\ dx -
\frac{1}{2}\iortd(\mathcal{R}\bv\cdot\nabla)\bW_j\cdot\bu_n\ dx\\
&+\iortd \nabla\bu_n:\nabla\bW_j\ dx +\im\pa^2_t\eta_n\, W_j\ dA + 2K(\eta_n, W_j) \\
% +\im\pa_t\eta(t,\cdot)b(t,\cdot)\ dA 
&\hspace{5cm}=\iortd \ff_n\cdot\bW_j\ dx + \im g_n\, W_j\
dA
 \end{aligned}
\end{equation}
f"ur alle $1\le j \le n$ erf"ullen. Dabei sind $\ff_n$ und $g_n$ Funktionen hoher Regularit"at, die
gegen $\ff$
und $g$ in $L^2_\loc([0,\infty)\times \setR^3)$ bzw. $L^2_\loc([0,\infty)\times M)$ konvergieren.
Wir k"onnen zum
Beispiel $\ff_n:=\R_{1/n}\ff$ und
\[g_n:=\sum_\alpha(\omega_{1/n}\ast((\psi_\alpha
g)\circ\phi^{-1}_\alpha))\circ\phi_\alpha\]
setzen; vgl. Definition \eqref{eqn:reg}. Die Konvergenzen folgen
dann direkt aus Proposition \ref{theorem:faltung}.

Wir geben zudem Anfangswerte vor. Wir w"ahlen die $\alpha^k_n(0)$ derart, dass
\[\pa_t\eta_n(0,\cdot)\rightarrow\eta_1\text{ in } L^2(M)\]
und 
\[\bu_n(0,\cdot)\rightarrow \bu_0\text{ in }L^2(\Omega_{\R\eta_0})\]
gilt. Dazu w"ahlen
wir die Koeffizienten der $\bY_k$ bei $t=0$ so, dass die erste der beiden Konvergenzen
gilt. Das ist
m"oglich, weil die Felder $Y_k(0,\cdot)$ eine Basis von \eqref{eqn:raum} mit
$t=0$ bilden und wir aus Proposition \ref{lemma:nspur} die Identit"at                   
\[\int_M \eta_1\, \gamma(\R\eta_0)\ dA=\int_{\Omega_{\R\eta_0}}\dv\bu_0\ dx=0\]
folgern k"onnen. Theorem \ref{theorem:stokes} impliziert, dass der
L"osungsoperator des
Stokes-Systems stetig von 
\[\Big\{Y\,\bnu\in L^2(M)\ |\ \int_M Y\,\gamma(\R\eta_0)\ dA=0\Big\}\] nach $L^2(\Omega_{\R\eta_0})$
abbildet. Somit konvergieren nicht nur die Linearkombinationen der $Y_k(0,\cdot)$ gegen $\eta_1$ in
$L^2(M)$, sondern auch die entsprechenden Linearkombinationen der $\bY_k(0,\cdot)$ gegen ein $\bY$
in $L^2(\Omega_{\R\eta_0})$. Es gilt 
$\trnormaleo(\bu_0-\bY) = 0$. Da die $\bX_k(0,\cdot)$ inbesondere eine Basis von 
\[\big\{\bX\in L^2(\Omega_{\R\eta_0})\ |\ \dv \bX=0,\, \trnormaleo \bX=0\big\}\] 
bilden, k"onnen wir nun deren Koeffizienten bei $t=0$ so
w"ahlen, dass die Folge $(\alpha_n^k\,\bW_k)_{n}$ gegen $\bu_0$ in $L^2(\Omega_{\R\eta_0})$
konvergiert.

Es liegt somit ein Anfangswertproblem f"ur ein lineares System gew"ohnlicher
Integro-Differentialgleichungen der Form ($1\le j\le n$, Summation von $1$ bis $n$)
\begin{equation*}
   A_{jk}(t)\,\dot\alpha^k(t) = B_{jk}\,\alpha^k(t) + \int_0^t C_{jk}(t,s)\,\alpha^k(s)\ ds +
D_j(t),
% \begin{matrix}
% \begin{pmatrix}
%   1 & 0 \\
%  A & B+C
% \end{pmatrix} 
% &\hspace{-0.3cm} \begin{pmatrix}
%   \dot \alpha^k \\
%   \dot \beta^k
% \end{pmatrix}
% &=\hspace{1mm}\begin{pmatrix}
%   \beta^j \\
%   -D_{jk}\alpha^k + E_j
% \end{pmatrix}  
% \end{matrix}.
\end{equation*}
vor. Die Koeffizienten sind durch
\begin{equation*}
 \begin{aligned}
A_{jk}(t)&=\iortd \bW_k\cdot\bW_j\ dx + \im W_k\, W_j\ dA,\\
B_{jk}(t)&=\iortd \pa_t\bW_k\cdot\bW_j\ dx + \frac{1}{2}\im W_k\, W_j\
\pa_t\R\delta\, \gamma(\R\delta)\ dA\\
&\hspace{0.3cm} +\frac{1}{2} \iortd (\R\bv\cdot\nabla)\bW_k\cdot\bW_j\ dx - \frac{1}{2} \iortd
(\R\bv\cdot\nabla)\bW_j\cdot\bW_k\ dx \\
&\hspace{0.3cm} + \iortd \nabla\bW_k:\nabla\bW_j\ dx + \im \pa_t W_k\, W_j\ dA\\
C_{jk}(t,s)&=2 K(W_k(s),W_j(t)),\\
D_j(t)&=\iortd \ff_n\cdot \bW_j\ dx + \im g_n\, W_j\ dA
 \end{aligned}
\end{equation*}
gegeben und somit offenbar stetig. Die Matrix $A(t)$ ist symmetrisch, und aufgrund der linearen
Unabh"angigkeit unserer Ansatzfunktionen k"onnen wir aus der Identit"at ($\beta^j\in\setR$,
$j=1,\ldots,n$)
\[\beta^j\beta^k A_{jk}(t)=\iortd |\beta^j\bW_j|^2\ dx+\im |\beta^j W_j|^2\ dA\]
schlie"sen, dass $A(t)$ positiv definit ist. Insbesondere ist diese Matrix invertierbar,
womit Gleichung \eqref{eqn:gal} von der in Anhang A.3 behandelten Form ist. Sie besitzt somit f"ur
jedes
$n$ eine
L"osung auf dem Intervall $[0,T]$. Um Energieabsch"atzungen zu bekommen, testen wir die Gleichung
mit $(\pa_t\eta_n,\bu_n)$ und erhalten
\begin{equation*}
\begin{aligned}
% &\iot \bu(t,\cdot)\cdot\bphi(t,\cdot)\ dx 
&\iortd \pa_t\bu_n\cdot\bu_n\ dx +\frac{1}{2} \im(\pa_t\eta_n)^2
\, \pa_t\mathcal{R}\delta\, \gamma(\mathcal{R}\delta)\ dA\\
&+\iortd |\nabla\bu_n|^2\ dx
% +\im\pa_t\eta(t,\cdot)b(t,\cdot)\ dA
+\im\pa^2_t\eta_n\, \pa_t\eta_n\ dA + 2K(\eta_n,\pa_t\eta_n)\\
&\hspace{5cm}=\iortd \ff_n\cdot\bu_n\ dx + \im g_n\, \pa_t\eta_n\ dA.
\end{aligned}
\end{equation*}
Die ersten beiden Terme lassen sich wie gewohnt zusammenfassen, sodass sich
\begin{equation*}
\begin{aligned}
% &\iot \bu(t,\cdot)\cdot\bphi(t,\cdot)\ dx 
\frac{1}{2}\frac{d}{dt}\iortd |\bu_n|^2\ dx + \iortd |\nabla\bu_n|^2\ dx
% +\im\pa_t\eta(t,\cdot)b(t,\cdot)\ dA
+\frac{1}{2}\frac{d}{dt}&\im|\pa_t\eta_n|^2\ dA + \frac{d}{dt} K(\eta_n)\
dA\\ &=\iortd \ff_n\cdot\bu_n\ dx + \im g_n\, \pa_t\eta_n\ dA
\end{aligned}
\end{equation*}
ergibt. Nun k"onnen wir wie am Ende des letzten Abschnitts vorgehen und erhalten
\begin{equation*}
 \begin{aligned}
\norm{\eta_n}_{Y^I}^2 +
\norm{\bu_n}_{X_{\R\delta}^I}^2\le
c_0(T,\Omega_{\R\delta}^I,\ff_n,g_n,\bu_n(0,\cdot),\eta_0,\pa_t\eta_n(0,\cdot)).
\end{aligned}
\end{equation*}
Aus dieser Absch"atzung folgern wir f"ur eine Teilfolge die Konvergenzen\footnote{Man beachte, dass
\begin{equation*}
 \begin{aligned}
X_{\R\delta}^I&\simeq
L^\infty(I,L^2(\Omega))\cap L^2(I,H^1_{\dv}(\Omega))\simeq (L^1(I,L^2(\Omega))+
L^2(I,H^1_{\dv}(\Omega)))'\\ &\simeq (L^1(I,L^2(\Omega_{\R\delta(t)}))+
L^2(I,H^1_{\dv}(\Omega_{\R\delta(t)})))'
 \end{aligned}
\end{equation*}
gilt, wobei die erste und die dritte Isomorphie durch die Abbildung $\T_{\R\delta}$ induziert
werden.}
\begin{equation*}
 \begin{aligned}
  \eta_n&\rightarrow \eta \hspace{0.3cm}&&\text{ schwach$^*$ in } L^\infty(I,H^2_0(M)),\\
  \pa_t\eta_n&\rightarrow \pa_t\eta &&\text{ schwach$^*$ in } L^\infty(I,L^2(M)),\\
\bu_n&\rightarrow \bu &&\text{ schwach$^*$ in } X_{\R\delta}^I.
 \end{aligned}
\end{equation*}
Mit Hilfe der Unterhalbstetigkeit der Normen bez"uglich der schwach*-Konvergenz erhalten wir
\eqref{ab:ent}. Desweiteren folgt aus $\trrd \bu_n=\pa_t\eta_n\,\bnu$ und obigen Konvergenzen die
Identit"at $\trrd \bu=\pa_t\eta\,\bnu$.

Wir multiplizieren nun Gleichung \eqref{eqn:gal} mit $\phi(t)$, wobei $\phi\in C_0^1([0,T))$,
integrieren "uber $I$ und anschlie"send partiell in der Zeit und erhalten f"ur $1\le j\le n$
\begin{equation*}
\begin{aligned}
% &\iot \bu(t,\cdot)\cdot\bphi(t,\cdot)\ dx 
&-\int_I\iortd \bu_n\cdot\pa_t(\phi\,\bW_j)\ dxdt -\frac{1}{2}\int_I\im\pa_t\eta_n
\, \pa_t\mathcal{R}\delta\, \phi\, W_j\, \gamma(\mathcal{R}\delta)\
dAdt\\
& +\frac{1}{2}\int_I\iortd(\mathcal{R}\bv\cdot\nabla)\bu_n\cdot(\phi\,\bW_j)\ dxdt\\
&-\frac{1}{2}\int_I\iortd(\mathcal{R}\bv\cdot\nabla)(\phi\,\bW_j)\cdot\bu_n\ dxdt
+\int_I\iortd \nabla\bu_n:\nabla(\phi\,\bW_j)\ dxdt\\
% +\im\pa_t\eta(t,\cdot)b(t,\cdot)\ dA
&+\int_I\im\pa_t\eta_n\, \pa_t (\phi\, W_j)\ dAdt + 2\int_I K(\eta_n, \phi\, W_j)\
dAdt \\
&\hspace{2.5cm}=\int_I\iortd \ff_n\cdot(\phi\,\bW_j)\ dxdt + \int_I\im g_n\, \phi\, W_j\
dAdt\\
&\hspace{3cm} +\int_{\Omega_{\mathcal{R}\eta_0}}\bu_n(0)\cdot(\phi(0)\bW_j(0,\cdot))\ dx +
\im\pa_t\eta_n(0)\, \phi(0)W_j(0,\cdot)\ dA.
 \end{aligned}
\end{equation*}
Durch Grenz"ubergang in $n$ in obiger Gleichung sehen wir, dass $\eta$ und $\bu$ die
Identit"at \eqref{eqn:ent} f"ur alle Testfunktionen aus \[\spann\{(\phi\, W_j,\phi\,\bW_j)\
|\ \phi\in
C_0^1([0,T)),\, j\in\setN\}\] erf"ullen. Aufgrund der in der Vorarbeit gezeigten
Dichtheit dieser Funktionen folgt die G"ultigkeit der Identit"at f"ur alle Testfunktionen aus
$T^I_{\R\delta}$. Dabei ist zu beachten, dass f"ur den "`inneren
Anteil"' der Testfunktionen die Konvergenz in $H^1(\Omega_{\R\delta}^I)$ wegen der
Entsch"arfung des Wirbelterms gen"ugt.
\qed
\\

Man beachte, dass die Eindeutigkeit schwacher L"osungen des entkoppelten Systems nicht
offensichtlich ist. Diese Eigenschaft zeigt man im Falle einer parabolischen Gleichung
"ublicherweise, indem man die Gleichung mit der L"osung selber testet. Das ist m"oglich, wenn wir
wissen, dass schwache L"osungen Zeitableitungen im Dualraum ihrer eigenen Regularit"atsklasse
besitzen. Das ist bei dem vorliegenden parabolisch-dispersiven System selbst formal nicht der Fall.
Gleichung \eqref{eqn:ent} zeigt n"amlich, dass die formale Zeitableitung des Tupels
$(\pa_t\eta,\bu)$ nicht im Dualraum der Regularit"atsklasse von $(\pa_t\eta,\bu)$, sondern in dem
der Regularit"atsklasse von $(\eta,\bu)$ liegt. Andererseits zeigt man die Eindeutigkeit schwacher
L"osungen von Gleichungen vom Typ der Schalengleichung "ublicherweise, indem man mit einem Term der
Form $\int \eta(s,\cdot)\ ds$ testet, vgl. \cite{b6},
was im Falle unseres gekoppelten Systems auch nicht m"oglich ist. Dieses Problem kann durch eine
weitere Regularisierung des Systems umgangen werden. Da wir jedoch dem Begriff der Zeitableitung
bislang keine rigorose Bedeutung verliehen haben, ist die Eindeutigkeit auch dann nicht
offensichtlich, sodass wir, wie bereits erw"ahnt, gerne (vorerst) auf einen Beweis verzichten
w"urden. Aus der Linearit"at von Gleichung \eqref{eqn:ent} folgt, dass die L"osungsmenge konvex
ist. Dies erlaubt es uns, anstelle des Schauder'schen
Fixpunktsatzes eine Variante dieses Satzes f"ur mengenwertige Abbildungen,
den Fixpunktsatz von Kakutani-Glicksberg-Fan, Theorem \ref{theorem:kakutani}, zu verwenden.

Wir wollen zun"achst den nun naheliegenden L"osungsbegriff fixieren.

\begin{definition} \label{def:epsilon} Ein Tupel $(\eta,\bu)$ hei"st schwache L"osung des
regularisierten Systems zum
Parameter $\epsilon$ auf dem Intervall $I$, falls $\eta\in Y^I$
mit $\norm{\eta}_{L^\infty(I\times M)}<\kappa$ und $\eta(0,\cdot)=\eta_0$, $\bu\in
X_{\mathcal{R}_\epsilon\eta}^I$ mit
$\trre \bu =
\pa_t\eta\,\bnu$ und
\begin{equation*}
\begin{aligned}
% &\iot \bu(t,\cdot)\cdot\bphi(t,\cdot)\ dx 
&- \int_I\iorte \bu\cdot\pa_t\bphi\ dxdt - \frac{1}{2}\int_I\im\pa_t\eta
\, \pa_t\mathcal{R}_\epsilon\eta\ b\, \gamma(\mathcal{R}_\epsilon\eta)\ dAdt\\
&+\frac{1}{2}\int_I\iorte(\mathcal{R}_\epsilon\bu\cdot\nabla)\bu\cdot\bphi\ dxdt -
\frac{1}{2}\int_I\iorte(\mathcal{R}_\epsilon\bu\cdot\nabla)\bphi\cdot\bu\ dxdt\\
& + \int_I\iorte \nabla\bu:\nabla\bphi\ dxdt-\int_I\im\pa_t\eta\, \pa_tb\ dAdt + 2\int_I K(\eta,b)\
dt \\
% +\im\pa_t\eta(t,\cdot)b(t,\cdot)\ dA
&\hspace{4cm}=\int_I\iorte \ff\cdot\bphi\ dxdt + \int_I\im g\, b\
dAdt\\
&\hspace{4.5cm} +\int_{\Omega_{\R_\epsilon\eta_0}}\bu_0^\epsilon\cdot\bphi(0,\cdot)\ dx +
\im\eta_1^\epsilon\, b(0,\cdot)\ dA
 \end{aligned}
\end{equation*}
f"ur alle Testfunktionen $(b,\bphi)\in  T_{\R_\epsilon\eta}^I$.
\end{definition}
\vspace{0.2cm}

\begin{proposition}\label{lemma:approx}
Es existieren ein $T>0$ und f"ur jedes hinreichend kleine $\epsilon>0$ eine schwache L"osung $(\eta_\epsilon,\bu_\epsilon)$ des
regularisierten Systems zum Parameter $\epsilon$ auf dem Intervall $I=(0,T)$, die die
Absch"atzung
\begin{equation}\label{ab:approx}
 \begin{aligned}
 &\norm{\eta_\epsilon}_{Y^I}^2 + \norm{\bu_\epsilon}_{X_{\R_\epsilon\eta}^I}^2\le
 c_0(T,\Omega_{\R_\epsilon\eta}^I,\ff,g,\bu^\epsilon_0,\eta_0,\eta^\epsilon_1)
  \end{aligned}
\end{equation}
erf"ullt. Die Zeit $T$ h"angt nur von $\tau(\eta_0)$ und der durch \eqref{ab:approx} gegebenen
Schranke f"ur die $Y^I$-Norm von
$\eta_\epsilon$ ab.
\end{proposition}
\beweis Wir setzen $\alpha:=(\norm{\eta_0}_{L^\infty(M)}+\kappa)/2$ und fixieren ein
beliebiges, hinreichend kleines $\epsilon>0$, unterdr"ucken diesen
Parameter jedoch im Rahmen dieses Beweises in der Notation. Wir wollen Theorem
\ref{theorem:kakutani} anwenden. Dazu betrachten wir den Raum \[Z:=C(\bar I\times\pa\Omega)\times
L^2(I\times
\setR^3)\] 
und die konvexe Menge
\[D:=\{(\delta,\bv)\in Z\ |\  \delta(0,\cdot)=\eta_0,\,
\norm{\delta}_{L^\infty(I\times\pa\Omega)}\le\alpha,\,
\norm{\bv}_{L^2(I\times \setR^3)}\le c_1\}\] mit einer hinreichend gro"sen Konstante $c_1>0$.
Desweiteren betrachten wir
die Abbildung
\[F:D\subset Z\rightarrow 2^Z,\] die jedem Tupel $(\delta,\bv)$
die Menge der schwachen L"osungen $(\eta,\bu)$ des entkoppelten, regularisierten Systems zum
Argument $(\delta,\bv)$ zuordnet, die der Absch"atzung
\begin{equation}\label{ab:schranke}
\begin{aligned}
\norm{\eta}_{Y^I} +
\norm{\bu}_{X_{\R\delta}^I} \le c(\delta)  
\end{aligned}
\end{equation}
gen"ugen, wobei $c(\delta)$ die rechte Seite von $\eqref{ab:ent}$ bezeichne. Wegen Proposition
\ref{lemma:ent} ist $F(\delta,\bv)$ nichtleer. Es gelte
$(\eta,\bu)\in
F(\delta,\bv)$. Aufgrund von \eqref{ab:schranke} ist dann $\norm{\bu}_{L^2(I\times
\setR^3)}$  durch $c_1$ abgesch"atzt, und die Norm von $\eta$ in
\[Y^I\embedding C^{0,1-\theta}(\bar I, C^{0,2\theta -1}(\pa\Omega))\quad (1/2<\theta<1)\]
ist beschr"ankt. Da zudem $\eta(0,\cdot)=\eta_0$ gilt,
kann das Zeitintervall $I=(0,T)$ unabh"angig vom Regularisierungsparameter $\epsilon$ so klein
gew"ahlt werden, dass $\norm{\eta}_{L^\infty(I\times\pa\Omega)}\le\alpha$. Also bildet
$F$
die Menge $D$ in ihre Potenzmenge ab, $F(D)\subset
2^D$. Aufgrund der Linearit"at der Gleichung ist $F(\delta,\bv)$ konvex. Zudem ist diese Menge
abgeschlossen in $Z$. Falls n"amlich eine Folge $(\eta_n,\bu_n)\subset F(\delta,\bv)$
gegen ein $(\eta,\bu)$ in $Z$ konvergiert, so gilt wegen $\eqref{ab:schranke}$ f"ur eine Teilfolge
\begin{equation*}
 \begin{aligned}
  \eta_n&\rightarrow \eta \hspace{0.3cm}&&\text{ schwach$^*$ in } L^\infty(I,H^2_0(M)),\\
  \pa_t\eta_n&\rightarrow \pa_t\eta &&\text{ schwach$^*$ in } L^\infty(I,L^2(M)),\\
\bu_n&\rightarrow \bu &&\text{ schwach$^*$ in } X_{\R\delta}^I.
 \end{aligned}
\end{equation*}
Wir k"onnen also in Gleichung \eqref{eqn:ent} zum Grenzwert
"ubergehen
und erhalten $(\eta,\bu)\in F(\delta,\bv)$.

Als n"achstes wollen wir zeigen, dass $F(D)$ relativ kompakt ist in $Z$. Sei dazu
$(\delta_n,\bv_n)$ eine Folge aus $D$ und $(\eta_n,\bu_n)\in F(\delta_n,\bv_n)$. Es ist zu
zeigen, dass eine Teilfolge von $(\eta_n,\bu_n)$ in $Z$ konvergiert. Die gleichm"a"sige Konvergenz
einer Teilfolge von $(\eta_n)$ ist offensichtlich. Da der Regularisierungsoperator 
\[\mathcal{R}:\{\delta\in C(\bar I\times\pa\Omega)|\ \delta(0,\cdot)=\eta_0 \}\rightarrow
C^3(\bar I\times\pa\Omega)\compactembedding C^2(\bar I\times\pa\Omega)\]
Kompaktheit erzeugt, konvergiert eine Teilfolge von
$(\mathcal{R}\delta_n)$ gegen ein $\delta$ in $C^2(\bar I\times\pa\Omega)$. Der Beweis der
relativen Kompaktheit der Folge $(\bu_n)$ in $L^2(I\times\setR^3)$ l"asst sich nun fast w"ortlich
von
Proposition \ref{lemma:komp} "ubernehmen. Es ist lediglich die leicht abgewandelte Form
\eqref{eqn:ent}
des Systems zu beachten und insbesondere teilweise die Folge $(\eta_n)$ mit Grenzwert $\eta$ durch
die Folge $(\R\delta_n)$ mit Grenzwert $\R\delta$ zu ersetzen. Aufgrund der Regularisierungen
k"onnten nat"urlich einige Argumente vereinfacht werden. Dieses Argument zeigt zus"atzlich die
relative Kompaktkeit der Folge $(\pa_t\eta_n)$ in $L^2(I\times M)$.

Es bleibt zu zeigen, dass die Abbildung $F$ graphenabgeschlossen ist.\footnote{Wir f"uhren die
Argumente im Beweis der Graphenabgeschlossenheit so, dass sie auch f"ur R"ander von
asymptotisch geringer Regularit"at g"ultig sind. Das ist sinnvoll, weil diese Argumente im Beweis
von Theorem \ref{theorem:hs} wiederholt werden m"ussen.} Wir nehmen dazu eine Folge
$(\delta_n,\bv_n)\subset D$ mit $(\delta_n,\bv_n)\rightarrow(\delta,\bv)$ in $Z$ sowie eine Folge
$(\eta_n,\bu_n)\in F(\delta_n,\bv_n)$
mit $(\eta_n,\bu_n)\rightarrow(\eta,\bu)$ in $Z$ als gegeben an. Es ist zu zeigen, dass
$(\eta,\bu)\in F(\delta,\bv)$. Aus \eqref{ab:schranke} und der
soeben gezeigten relativen Kompaktheit folgern
wir f"ur eine Teilfolge die Konvergenzen 
\begin{equation}\label{eqn:wichtigekonv}
 \begin{aligned}
  \eta_n&\rightarrow\eta &&\text{gleichm"a"sig und schwach$^*$ in }L^\infty(I,H^2_0(M)),\\
  \pa_t\eta_n&\rightarrow\pa_t\eta &&\text{ in
}L^2(I\times M)\text{ und schwach$^*$ in }L^\infty(I,L^2(M)),\\
  \bu_n&\rightarrow\bu &&\text{ in
}L^2(I\times \setR^3)\text{ und schwach$^*$ in } L^\infty(I,L^2(\setR^3)),\\
  \nabla\bu_n&\rightarrow\nabla\bu &&\text{ schwach in }
L^2(I\times \setR^3).
\end{aligned}
\end{equation}
Die zun"achst auf $\Omega_{\R\delta_n}^I$ bzw. $\Omega_{\R\delta}^I$ definierten Felder
$\nabla\bu_n$ und $\nabla\bu$ werden dabei durch $0$ auf $I\times \setR^3$ fortgesetzt. Aus der
Unterhalbstetigkeit der Normen bez"uglich schwacher bzw. schwach*-Konvergenz folgern wir die
Absch"atzung
\eqref{ab:schranke} f"ur $\eta$ und $\bu$. Die Identit"at $\eta(0,\cdot)=\eta_0$ ist
aufgrund der gleichm"a"sigen Konvergenz von $(\eta_n)$ offensichtlich. Desweiteren gilt 
\begin{equation*}
 \begin{aligned}
   \pa_t\eta_n\,\bnu=\trrdn
\bu_n = \bw_n|_{\pa\Omega},
 \end{aligned}
\end{equation*}
wobei $\bw_n:=\bu_n\circ\Psi_{\R\delta_n}$. Wegen Lemma \ref{lemma:psi} konvergiert eine Teilfolge
von $(\bw_n)$ schwach gegen ein $\bw$ in $L^2(I,W^{1,1}(\Omega))$, sodass die Identit"at
$\pa_t\eta\,\bnu=\bw|_{\pa\Omega}$ folgt. Der erste Term auf der rechten Seite der
Absch"atzung 
\[\norm{\bw_n-\bu\circ\Psi_{\R\delta}}_{L^1(I\times \Omega)}\le
\norm{(\bu_n-\bu)\circ\Psi_{\R\delta_n}}_{L^1(I\times \Omega)}+
\norm{\bu\circ\Psi_{\R\delta_n}-\bu\circ\Psi_{\R\delta}}_{L^1(I\times \Omega)}\]
wird f"ur gro"se $n$ wegen Lemma \ref{lemma:psi} und \eqref{eqn:wichtigekonv}$_3$ klein, w"ahrend
der zweite Term wegen Bemerkung \ref{bem:konv} klein wird. Es folgt $\bw=\bu\circ\Psi_{\R\delta}$,
sodass wir
insgesamt die Identit"at $\pa_t\eta\,\bnu=\trrd\bu$ gezeigt haben. Es bleibt lediglich zu zeigen,
dass
Gleichung
\eqref{eqn:ent} erf"ullt ist. F"ur alle $n$ und alle Testfunktionen $(b_n,\bphi_n)\in
T_{\R\delta_n}^I$
gilt
\begin{equation}\label{eqn:graphab}
\begin{aligned}
% &\iot \bu(t,\cdot)\cdot\bphi(t,\cdot)\ dx 
&- \int_I\iortdn \bu_n\cdot\pa_t\bphi_n\ dxdt - \frac{1}{2}\int_I\im\pa_t\eta_n
\, \pa_t\mathcal{R}\delta_n\, b_n\, \gamma(\mathcal{R}\delta_n)\ dAdt\\
&+\frac{1}{2}\int_I\iortdn(\mathcal{R}\bv_n\cdot\nabla)\bu_n\cdot\bphi_n\ dxdt -
\frac{1}{2}\int_I\iortdn(\mathcal{R}\bv_n\cdot\nabla)\bphi_n\cdot\bu_n\ dxdt\\
& + \int_I\iortdn \nabla\bu_n:\nabla\bphi_n\ dxdt-\int_I\im\pa_t\eta_n\, \pa_tb_n\ dAdt + 2\int_I
K(\eta_n,b_n)\
dt \\
% +\im\pa_t\eta(t,\cdot)b(t,\cdot)\ dA 
&\hspace{4cm}=\int_I\iortdn \ff\cdot\bphi_n\ dxdt + \int_I\im g\, b_n\
dAdt\\
&\hspace{4.5cm} +\int_{\Omega_{\mathcal{R}\eta_0}}\mathcal{R}\bu_0\cdot\bphi_n(0,\cdot)\ dx +
\im \R\eta_1\, b_n(0,\cdot)\ dA.
\end{aligned}
\end{equation}
Wir k"onnen hier nicht ohne  weiteres zum Grenzwert "ubergehen, weil die Testfunktionen vom Index
$n$ abh"angen. F"ur $(b,\bphi)\in
T_{\R\delta}^I$ verwenden wir deshalb die speziellen Testfunktionen 
$(b_n,\bphi_n):=(\M_{\R\delta_n}b,\F_{\R\delta_n}\M_{\R\delta_n}b)\in T_{\R\delta_n}^I$ wie sie auch
im Beweis von
Proposition \ref{lemma:komp} zur Anwendung kamen. Die Zahl $\alpha$ in Propositon
\ref{lemma:FortVonRandZeit} w"ahlen wir dabei identisch $(\alpha+\kappa)/2$. Mit diesen
Testfunktionen k"onnen wir wegen der Konvergenzen \eqref{eqn:wichtigekonv} und wegen Lemma
\ref{lemma:konvergenzen} $(1.b)$, $(2.b)$ in \eqref{eqn:graphab} zum Grenzwert "ubergehen, sodass
wir die G"ultigkeit
von \eqref{eqn:ent} f"ur $(b,\bphi)=(b,\F_{\R\delta}b)\in
T_{\R\delta}^I$ erhalten. Per Definition von $T_{\R\delta}^I$ bleibt lediglich zu zeigen, dass
\eqref{eqn:ent} f"ur Testfunktionen $(0,\bphi)\in
T_{\R\delta}^I$ mit $\bphi(T,\cdot)=0$ und $\supp\bphi\subset\Omega_{\R\delta}^{\bar I}$ gilt.
Aufgrund der gleichm"a"sigen Konvergenz von $(\R\delta_n)$ gilt f"ur hinreichend gro"se $n$ in
diesem Fall aber $(0,\bphi)\in
T_{\R\delta_n}^I$, sodass wir wiederum in \eqref{eqn:graphab} zum Grenzwert "ubergehen k"onnen.

Theorem \ref{theorem:kakutani} garantiert uns nun die Existenz eines Fixpunktes von $F$, d.h. es
existiert ein Tupel $(\eta,\bu)\in D$ mit $(\eta,\bu)\in F(\eta,\bu)$. Dies zeigt die Behauptung der
Proposition.\qed\\

Wir k"onnen nun unsere zentrale Behauptung beweisen, indem wir den Regularisierungsparameter
$\epsilon$ in Definition \ref{def:epsilon} gegen $0$ gehen lassen.\\

\beweis (von Theorem \ref{theorem:hs}) Wir haben gezeigt, dass ein $T>0$ und
f"ur jedes $\epsilon=1/n$, $n\in\setN$ hinreichend gro"s, eine schwache L"osung
$(\eta_\epsilon,\bu_\epsilon)$ des
regularisierten Systems zum Parameter $\epsilon$ auf dem Intervall $I=(0,T)$ existieren. Aus
Absch"atzung \eqref{ab:approx} und der kompakten Einbettung $Y^I\compactembedding C(\bar
I\times\pa\Omega)$ erhalten wir f"ur eine Teilfolge die Konvergenzen
\begin{equation}\label{eqn:konv3}
 \begin{aligned}
   \eta_\epsilon,\,\R_\epsilon\eta_\epsilon&\rightarrow\eta &&\text{ gleichm"a"sig und
schwach$^*$ in }
L^\infty(I,H^2_0(M)),\\
   \pa_t\eta_\epsilon,\, \pa_t\R_\epsilon\eta_\epsilon&\rightarrow\pa_t\eta &&\text{ schwach$^*$
in }
L^\infty(I,L^2(M)),\\
  \bu_\epsilon&\rightarrow\bu &&\text{ schwach$^*$ in }
L^\infty(I,L^2(\setR^3)),\\
\nabla\bu_\epsilon&\rightarrow\nabla\bu &&\text{ schwach in }L^2(I\times \setR^3).
\end{aligned}
\end{equation}
Die zun"achst auf $\Omega_{\R_\epsilon\eta_\epsilon}^I$ bzw. $\Omega_{\eta}^I$ definierten Felder
$\nabla\bu_\epsilon$ und $\nabla\bu$ werden dabei durch $0$ auf $I\times \setR^3$
fortgesetzt. Die Konvergenzen der Folge $(\R_\epsilon\eta_\epsilon)$ schlie"sen wir unter
Verwendung der Faltungsungleichung aus der Absch"atzung
\[\norm{
\R_\epsilon\eta_\epsilon-\eta}_{L^\infty(I\times \pa\Omega)}
\le\norm{\R_\epsilon(\eta_\epsilon-\eta)}_{L^\infty(I\times
\pa\Omega)}+\norm{\R_\epsilon\eta-\eta}_{L^\infty(I\times \pa\Omega)}.\]
Wir k"onnen nun den Beweis von Proposition \ref{lemma:komp} fast w"ortlich
wiederholen, um die Konvergenzen
\begin{equation}\label{eqn:konv31}
 \begin{aligned}
   \pa_t\eta_\epsilon&\rightarrow\pa_t\eta &&\text{ in }
L^2(I\times M),\\
  \bu_\epsilon&\rightarrow\bu &&\text{ in }
L^2(I\times\setR^3)
\end{aligned}
\end{equation}
zu zeigen. Desweiteren erhalten wir wie im Beweis von Proposition
\ref{lemma:approx} die Identit"at $\tr\bu=\pa_t\eta\,\bnu$. Weiterhin zeigen wir unter Verwendung
von \eqref{eqn:konv31}, der Interpolationsungleichung
\[\norm{\bu_\epsilon-\bu}_{L^2(I,L^4(\setR^3))}\le
\norm{\bu_\epsilon-\bu}_{L^2(I\times
\setR^3)}^{1/6}\,\norm{\bu_\epsilon-\bu}_{L^2(I,L^5(\setR^3))}^{5/6}\]
und Korollar \ref{lemma:sobolev} die Konvergenz
\begin{equation}\label{eqn:konv20}
 \begin{aligned}
\bu_\epsilon\rightarrow\bu\text{ in }L^2(I,L^4(\setR^3)).  
 \end{aligned}
\end{equation}
Es folgt
\begin{equation}\label{eqn:konv8}
 \begin{aligned}
\pa_t\R_\epsilon\eta_\epsilon&\rightarrow\pa_t\eta &&\text{ in }L^2(I\times \pa\Omega),\\
\R_\epsilon\bu_\epsilon&\rightarrow\bu &&\text{ in }L^2(I,L^4(\setR^3)).
 \end{aligned}
\end{equation}
Aus der Unterhalbstetigkeit der Normen folgern wir die Absch"atzung \eqref{ab:hs}, und aufgrund der
gleichm"a"sigen Konvergenz von $(\eta_\epsilon)$ ist die Identit"at $\eta(0,\cdot)=\eta_0$
offensichtlich. F"ur alle $\epsilon$ und alle $(b_\epsilon,\bphi_\epsilon)\in
T_{\R_\epsilon\eta_\epsilon}^I$ gilt
\begin{equation}\label{eqn:sepp}
\begin{aligned}
% &\iot \bu(t,\cdot)\cdot\bphi(t,\cdot)\ dx 
&- \int_I\int_{\Omega_{R_\epsilon\eta_\epsilon(t)}} \bu_\epsilon\cdot\pa_t\bphi_\epsilon\ dxdt -
\frac{1}{2}\int_I\im\pa_t\eta_\epsilon\,
\pa_t\mathcal{R}_\epsilon\eta_\epsilon\, b_\epsilon\, \gamma(\mathcal{R}_\epsilon\eta_\epsilon)\
dAdt\\
&+\frac{1}{2}\int_I\int_{\Omega_{R_\epsilon\eta_\epsilon(t)}}(\mathcal{R}
_\epsilon\bu_\epsilon\cdot\nabla)\bu_\epsilon\cdot\bphi_\epsilon\ dxdt -
\frac{1}{2}\int_I\int_{\Omega_{R_\epsilon\eta_\epsilon(t)}}(\mathcal{R}
_\epsilon\bu_\epsilon\cdot\nabla)\bphi_\epsilon\cdot\bu_\epsilon\ dxdt\\
& + \int_I\int_{\Omega_{R_\epsilon\eta_\epsilon(t)}} \nabla\bu_\epsilon:\nabla\bphi_\epsilon\
dxdt-\int_I\im\pa_t\eta_\epsilon\, \pa_tb_\epsilon\ dAdt + 2\int_I K(\eta_\epsilon,b_\epsilon)\
dt \\
% +\im\pa_t\eta(t,\cdot)b(t,\cdot)\ dA
&\hspace{5cm}=\int_I\int_{\Omega_{R_\epsilon\eta_\epsilon(t)}} \ff\cdot\bphi_\epsilon\ dxdt +
\int_I\im g\, b_\epsilon\
dAdt\\
&\hspace{5.5cm} +\int_{\Omega_{\eta_0^\epsilon}}\bu_0^\epsilon\cdot\bphi_\epsilon(0,\cdot)\ dx +
\im\eta_1^\epsilon\, b_\epsilon(0,\cdot)\ dA.
 \end{aligned}
\end{equation}
F"ur $(b,\bphi)\in T_\eta^I$ verwenden wir wie im Beweis von Proposition \ref{lemma:approx} die
speziellen Testfunktionen
$(b_\epsilon,\bphi_\epsilon):=(\M_{\R_\epsilon\eta_\epsilon}b,\F_{\R_\epsilon\eta_\epsilon}\M_{
\R_\epsilon\eta_\epsilon}b)\in T_{\R_\epsilon\eta_\epsilon}^I$. Die Zahl $\alpha$ in Propositon
\ref{lemma:FortVonRandZeit} gen"uge dabei der Ungleichung
$\sup_{\epsilon}\norm{\eta_\epsilon}_{L^\infty(I\times M)}<\alpha<\kappa$. Da die Felder
$(\nabla\,\F_{\R_\epsilon\eta_\epsilon}\M_{\R_\epsilon\eta_\epsilon}b)$ in
$L^\infty(I,L^{2}(B_\alpha))$ beschr"ankt sind, gilt
\begin{equation}\label{eqn:konv5}
 \begin{aligned}
\nabla\,\F_{\R_\epsilon\eta_\epsilon}\M_{\R_\epsilon\eta_\epsilon}
b\rightarrow\nabla\,\F_{\eta}\M_{\eta}b\quad\text{ schwach$^*$ in
}L^\infty(I,L^2(B_\alpha)).  
 \end{aligned}
\end{equation}
Unter Ausnutzung von \eqref{eqn:konv3}, \eqref{eqn:konv31}, \eqref{eqn:konv20},
\eqref{eqn:konv8} und \eqref{eqn:konv7} sowie von Lemma \ref{lemma:konvergenzen} $(1.b)$, $(2.b)$
und \eqref{eqn:konv5} k"onnen wir nun in \eqref{eqn:sepp} den Grenz"ubergang vollziehen.
\eqref{eqn:konv31}$_1$ und \eqref{eqn:konv8}$_2$ werden dabei f"ur den zweiten Term ben"otigt,
w"ahrend
\eqref{eqn:konv20}, \eqref{eqn:konv8}$_3$ und \eqref{eqn:konv5} den Grenz"ubergang im dritten Term
und vierten Term
erm"oglichen. Wir k"onnen also auf die G"ultigkeit von \eqref{eqn:schwach} f"ur
$(b,\bphi)=(b,\F_{\eta}b)\in
T_{\eta}^I$ schlie"sen. Der Grenz"ubergang f"ur Testfunktionen
$(0,\bphi)\in T_\eta^I$ mit $\bphi(T,\cdot)=0$ und $\supp\bphi\subset\Omega^{\bar I}_{\eta}$ ist
aber jetzt klar.
% Wenn wir zun"achst annehmen, dass
% $b\in C^1(\overline{I\times M})$, so\footnote{Die Aussagen sind klar f"ur $W^{1,4}$ im Ort
% (sogar $W^{1,r}$ f"ur alle $r<\infty$), der nat"urlich in $L^\infty\cap W^{1,3}$ einbettet.} gilt
% zudem $b_n\rightarrow b$ in $L^\infty(I,L^\infty(M)\cap
% W^{1,3}(M))$ und somit auch
% $\widetilde\bphi_n\rightarrow\widetilde\bphi$ in $L^\infty(I,L^\infty(B)\cap
% W^{1,3}(B))$. Diese letzte Konvergenz wird f"ur den Grenz"ubergang in den beiden "`Wirbeltermen"'
% ben"otigt\erl. Wir k"onnen nun in \eqref{eqn:approx} zum Grenzwert "ubergehen und den Wirbelterm
% in seine urspr"ungliche Form bringen, sodass wir \eqref{eqn:schwach} erhalten. F"ur allgemeines
% $b$ k"onnen wir anschlie"send $b$ approximieren und in \eqref{eqn:schwach} zum Grenzwert
% "ubergehen.

Das Existenzintervall der konstruierten L"osung h"angt lediglich von der Supremumsnorm von $\eta$
zum Anfangszeitpunkt sowie der Schranke f"ur die
H"older-Norm von $\eta$ ab. Nach Konstruktion ist $\norm{\eta}_{L^\infty(I\times M)}<\kappa$. Die
Gr"o"sen
$\norm{\bu(t)}_{L^2(\oet)}$, $\norm{\eta(t)}_{H^2(M)}$ und $\norm{\pa_t\eta(t)}_{L^2(M)}$ sind
gleichm"a"sig f"ur fast alle $t\in I$ beschr"ankt. Wenn wir nun f"ur fast alle $t$ L"osungen zu
diesen Anfangsdaten konstruieren, sind die H"older-Normen der Auslenkungen des Randes wegen
\eqref{ab:hs} unabh"angig
vom Anfangszeitpunkt $t$ nach oben und somit die Lebensdauer der L"osungen nach unten beschr"ankt.
Liegt $t\in I=(0,T)$ hinreichend nahe bei $T$, erhalten wir also eine
schwache L"osung auf einem Intervall $(t,\widetilde
T)$, $\widetilde T>T$. Verketten dieser L"osung mit $(\eta,\bu)$ an der Stelle
$t$ ergibt wegen Bemerkung \ref{bem:grenzen} eine L"osung $(\widetilde\eta,\widetilde\bu)$ auf dem
Intervall $(0,\widetilde T)$. Zudem erf"ullt $(\widetilde\eta,\widetilde\bu)$ die Absch"atzung
\eqref{ab:hs} auf dem Intervall
$(0,\widetilde T)$, denn f"ur $t_0<t<\widetilde T$ gilt
\begin{equation*}
 \begin{aligned}
 &\onorm{\widetilde\bu(t,\cdot)}^2 +
\int_0^t\norm{\nabla\widetilde\bu(s,\cdot)}^p_{L^p(\Omega_{\widetilde\eta(s)})}\ ds +
\mnorm{\pa_t\widetilde\eta(t,\cdot)}^2 + \norm{\widetilde\eta(t,\cdot)}_{H^2(M)}^2 \\ 
&\hspace{1cm}\le \big(\norm{\bu(t_0)}_{L^2(\Omega_{\widetilde\eta(t_0)})}^2 +
\mnorm{\pa_t\eta(t_0)}^2 +
\norm{\eta(t_0)}_{H^2(M)}^2\big)e^{c(t-t_0)}\\ 
&\hspace{1.5cm} + \int_{t_0}^t
\big(\norm{\ff(s,\cdot)}_{L^2(\Omega_{\eta(s)})}^2 + \mnorm{g(s,\cdot)}^2 \big)e^{c(t-s)}\
ds + \int_0^{t_0}\norm{\nabla\bu(s,\cdot)}^p_{L^p(\Omega_{\widetilde\eta(s)})}\ ds\\
&\hspace{1cm}\le \big(\norm{\bu_0}_{L^2(\Omega_{\eta_0})}^2 + \mnorm{\eta_1}^2 +
\norm{\eta_0}_{H^2(M)}^2\big)e^{ct}\\ 
&\hspace{5cm} + \int_0^t
\big(\norm{\ff(s,\cdot)}_{L^2(\Omega_{\eta(s)})}^2 + \mnorm{g(s,\cdot)}^2\big)e^{c(t-s)}\ ds.
 \end{aligned}
\end{equation*}
Durch Wiederholen dieses Vorgangs erhalten wir eine maximale Zeit $T^*\in (0,\infty]$ und ein Tupel
$(\eta,\bu)$, das f"ur jedes $T<T^*$ L"osung auf dem Intervall $(0,T)$ ist und die Absch"atzung
\eqref{ab:hs} erf"ullt. Ist $T^*$ endlich, so ist aufgrund dieser Absch"atzung die H"older-Norm
von $\eta$ in $[0,T^*]\times M$ beschr"ankt, und somit muss
$\norm{\eta(t,\cdot)}_{L^\infty(M)}\nearrow\kappa$ f"ur $t\nearrow T^*$ gelten.
\qed\\
\chapter{Verallgemeinerte Newton'sche Fluide}

Wir wollen in diesem Abschnitt
erste Schritte in Richtung einer interessanten und in hohem Ma"se nichttrivialen
Verallgemeinerung der bisherigen Resultate skizzieren, die eng mit dem bereits
diskutierten Eindeutigkeitsproblem zusammenh"angen. Bisher war der viskose Spannungstensor eine
lineare Funktion des
Scherratentensors $D\bu$. F"ur eine wichtige Klasse von Fluiden, die verallgemeinerten
Newton'schen Fluide, die eine scherratenabh"angige
Viskosit"at besitzen, gilt dieser lineare Zusammenhang nicht mehr. Ein wichtiges Beispiel daf"ur
ist Blut. Blut zeigt bei kleinen Scherraten eine hohe Viskosit"at, was f"ur die schnelle
Gerinnung wichtig ist. Flie"st es aber durch d"unne Adern, so entstehen gr"o"sere Scherraten,
und die Viskosit"at wird kleiner, was den Durchfluss bei konstantem Druck erh"oht. F"ur eine
detaillierte Diskussion verallgemeinerter Newton'scher Fluide verweisen wir auf
\cite{b34} und die dort angef"uhrten Referenzen. Typische Beispiele f"ur die viskosen
Spannungstensoren solcher Fluide sind
\begin{equation*}
 \begin{aligned}
  S(D\bu)&=\mu_0(\delta+|D\bu|)^{p-2}D\bu,\\
  S(D\bu)&=\mu_0(\delta^2+|D\bu|^2)^{\frac{p-2}{2}}D\bu
 \end{aligned}
\end{equation*}
mit $\mu_0>0$, $\delta\ge 0$ und $1<p<\infty$. Wir sind deshalb an Abbildungen $S$ mit $p$-Struktur
interessiert, d.h. f"ur ein $1<p<\infty$ und ein $\delta\ge 0$ gelte
\begin{itemize}
 \item $S: M_{sym}\rightarrow M_{sym}$ stetig,
 \item Wachstum: $S(D)\le c_0(\delta+|D|)^{p-2}|D|$ f"ur alle $D\in M_{sym}$ mit $c_0>0$,
 \item Koerzivit"at: $S(D):D\ge c_1(\delta+|D|)^{p-2}|D|^2$ f"ur alle $D\in
M_{sym}$ mit $c_1>0$,
\item Strikte Monotonie: $(S(D)-S(E)):(D-E)>0$ f"ur alle $D,E\in M_{sym}, D\not=E$.
\end{itemize}
Dabei bezeichnet $M_{sym}$ den Raum der reellen, symmetrischen $3\times 3$-Matrizen. Wir geben Daten
genau
wie im Falle eines linearen Spannungstensors vor, wobei wir nun $\ff\in
L^{p'}_\loc([0,\infty)\times\setR^3)$ annehmen. Anstelle von
\eqref{eqn:fluid}$_1$ stellen wir die Gleichung 
\[\pa_t \bu + (\bu\cdot\nabla)\bu = \dv (S(D\bu) - \pi\id) + \ff \ \mbox{ in }\
\Omega_{\eta}^I\]
auf, w"ahrend die Kraftdichte ${\bf F}$ auf der rechten Seite von \eqref{eqn:shell}$_1$ die
Form
\[{\bf F}(t,\cdot) = \big(-S(D\bu(t,\cdot))\,\bnu_{\eta(t)} +
\pi(t,\cdot)\,\bnu_{\eta(t)}\big)\circ\Phi_{\eta(t)}\, |\det d\Phi_{\eta(t)}|\]
annimmt. Den Rest der Gleichungen "ubernehmen wir unver"andert. Was die Analysis der
Fluidgleichungen bei festem Rand (bzw. auf dem Torus) betrifft, verweisen
wir auf \cite{b32} und die dortigen Referenzen; siehe insbesondere auch \cite{b33}. Wir wollen uns
mit einigen knappen Anmerkungen begn"ugen. Abbildungen $S$ mit $p$-Struktur definieren in
geeigneten Sobolev-R"aumen durch die Zuordnung $v\mapsto\dv S(\nabla v)$ monotone Operatoren. Diese
Operatoren spielen als Erzeuger nichtlinearer Halbgruppen eine prominente Rolle; vgl. \cite{b16},
\cite{b7}, \cite{b31}. Allerdings ist die abstrakte Theorie auf die Fluidgleichungen nicht direkt
anwendbar, weil der Wirbelterm die Monotonie zerst"ort. Dennoch l"asst sich durch Kombination eines
Kompaktheitsarguments und eines Arguments aus der Theorie monotoner Operatoren f"ur hinreichend
gro"se $p$ die Existenz zeitlich globaler, schwacher L"osungen zeigen. Die Einschr"ankung an
den Exponenten kommt dadurch zustande, dass das Argument eine L"osung mit Zeitableitung im
Dualraum ihrer eigenen Regularit"atsklasse ben"otigt und der Wirbelterm dabei restriktiv wirkt. F"ur
kleinere Exponenten sind
deshalb andere Techniken vonn"oten.

Wir k"onnen genau wie zuvor, unter
Verwendung der Koerzivit"at von $S$, formale Energieabsch"atzungen herleiten. Dabei erhalten wir
anstelle von \eqref{ab:apriori} die Ungleichung
\begin{equation*}
 \begin{aligned}
&\onorm{\bu(t,\cdot)}^2 + \int_0^t\norm{D\bu(s,\cdot)}^p_{L^p(\Omega_{\eta(s)})}\ ds +
\mnorm{\pa_t\eta(t,\cdot)}^2 + \norm{\eta(t,\cdot)}_{H^2(M)}^2 \\ 
&\hspace{1cm}\le \big(\norm{\bu_0}_{L^2(\Omega_{\eta_0})}^2 + \mnorm{\eta_1}^2 +
\norm{\eta_0}_{H^2(M)}^2\big)\,e^{ct}\\ 
&\hspace{4cm} + \int_0^t
\big(\norm{\ff(s,\cdot)}_{L^{p'}(\Omega_{\eta(s)})}^{p'} + \mnorm{g(s,\cdot)}^2\big)\,e^{c(t-s)}\
ds.
 \end{aligned}
\end{equation*}
\"Ublicherweise wird der symmetrische Gradient in dieser Absch"atzung mit Hilfe einer
Korn-Ungleichung
\[\norm{\nabla\bu(s,\cdot)}^p_{L^p}\le c\,\norm{D\bu(s,\cdot)}^p_{L^p}\]
durch den vollen Gradienten ersetzt. Die Korn-Ungleichung ist allerdings in
Gebieten, deren Rand nicht Lipschitz-stetig ist, im Allgemeinen falsch; vgl. \cite{b26}. Deshalb
sehen wir davon
ab, begn"ugen uns mit der Absch"atzung
\begin{equation*}
 \begin{aligned}
\norm{\eta}_{W^{1,\infty}(I,L^2(M))\cap L^\infty(I,H_0^2(M))} +
\norm{\bu}_{L^\infty(I,L^2(\Omega_{\eta(t)}))} + \norm{D\bu}_{L^p(\Omega_{\eta}^I)}\le
c(T,\text{Daten})
 \end{aligned}
\end{equation*}
und setzen f"ur $\eta\in Y^I$ mit $\norm{\eta}_{L^\infty(I\times M)}<\kappa$
\begin{equation*}
 \begin{aligned}
L^p(I,W^{1,p}_{\dv,s}(\Omega_{\eta(t)}))&:=\{\bv\ |\ \bv,D\bv\in L^p(\Omega_{\eta}^I),\, \dv
\bv=0\},\\
X^I_{\eta,p}&:= L^\infty(I,L^2(\Omega_{\eta(t)}))\cap L^p(I,W^{1,p}_{\dv,s}(\Omega_{\eta(t)})).
\end{aligned}
\end{equation*}
Den Raum $T^I_{\eta,p}$ der Testfunktionen definieren wir so wie $T^I_\eta$, wobei wir
den Raum $L^2(I,H^2_0(M))$ durch $L^{\max(2,p)}(I,H^2_0(M))$ und den Raum $H^1(\Omega_{\eta}^I)$
durch \[H^1(I,L^2(\Omega_{\eta(t)}))\cap L^{\max(2,p)}(I,W^{1,\max(2,p)}(\Omega_{\eta(t)}))\]
ersetzen. Durch die Forderung $b\in L^{p}(I,H^2_0(M))$ stellen wir sicher, dass die
Fortsetzung $\F_\eta b$ in $L^{p}(I,W^{1,p}(\Omega_{\eta(t)}))$ liegt.

Wir geben nun das Analogon von Definition \ref{def:ent} an. Wir konstruieren genau wie zuvor
modifizierte Anfangswerte $\R_\epsilon\eta_0$,
$\eta_1^\epsilon$ und $\bu^\epsilon$ und geben ein Zeitintervall $I=(0,T)$, $T>0$,
und Funktionen $\bv\in L^2(I\times B)$ und $\delta\in C(\bar I\times \pa\Omega)$ mit
$\norm{\delta}_{L^\infty(I\times\pa\Omega)}<\kappa$ und $\delta(0,\cdot)=\eta_0$
vor. Zudem setzen wir
\[\widetilde Y^I:=\{b\in Y^I\ |\ b\in H^1(I,H^2_0(M))\}.\]
Wie zuvor unterdr"ucken wir zun"achst den Parameter $\epsilon$ in der Notation.

\begin{definition}\label{def:pS}
Ein Tupel $(\eta,\bu)$ hei"st schwache L"osung des entkoppelten, regularisierten
$p$-Systems zum Argument $(\delta,\bv)$ auf dem Intervall $I$, falls $\eta\in\widetilde  Y^I$
mit $\eta(0,\cdot)=\eta_0$, $\bu\in X_{\R\delta,p}^I$ mit $\trrd
\bu =\pa_t\eta\,\bnu$ und 
\begin{equation}\label{eqn:entp}
\begin{aligned}
% &\iot \bu(t,\cdot)\cdot\bphi(t,\cdot)\ dx 
&- \int_I\iortd \bu\cdot\pa_t\bphi\ dxdt - \frac{1}{2}\int_I\im\pa_t\eta\ 
\pa_t\mathcal{R}\delta\, b\, \gamma(\mathcal{R}\delta)\ dAdt\\
&+\frac{1}{2}\int_I\iortd(\mathcal{R}\bv\cdot\nabla)\R\bu\cdot\bphi\ dxdt -
\frac{1}{2}\int_I\iortd(\mathcal{R}\bv\cdot\nabla)\R\bphi\cdot\bu\ dxdt\\
& + \int_I\iortd S(D\bu):D\bphi\ dxdt -\int_I\im\pa_t\eta\, \pa_tb\ dAdt + 2\int_I
K(\eta+\pa_t\eta,b)\
dt\\
% +\im\pa_t\eta(t,\cdot)b(t,\cdot)\ dA
&\hspace{5.0cm}=\int_I\iortd \ff\cdot\bphi\ dxdt + \int_I\im g\, b\
dAdt\\
&\hspace{5.5cm} +\int_{\Omega_{\mathcal{R}\eta_0}}\bu_0\cdot\bphi(0,\cdot)\ dx +
\im \eta_1\, b(0,\cdot)\ dA
\end{aligned}
\end{equation}
f"ur alle Testfunktionen $(b,\bphi)\in T_{\R\delta,p}^I$.
\end{definition}
\vspace{0.2cm}

Aufgrund der Wachstumsbedingung an $S$ ist
das Integral
\[\int_I\iortd S(D\bu):D\bphi\ dxdt\]
endlich. Eine wichtige Modifikation in dieser Definition gegen"uber Definition \ref{def:ent} ist der
Zusatzterm
\[2\int_I K(\pa_t\eta,b)\ dt,\]
den man sich mit einem Faktor $\epsilon$ versehen denke. Dieser entspricht einem zus"atzlichen
Ausdruck $\grad_{L^2}K(\pa_t\eta)$ in der Schalengleichung. Ein "ahnlicher Term wird in \cite{b27}
zur D"ampfung der Plattengleichung verwendet. Wie bereits angedeutet werden
wir
die Eindeutigkeit obiger L"osungen ben"otigen. Im Anschluss an den Beweis von Proposition
\ref{lemma:ent} wurde aber klar, dass ein Eindeutigkeitsbeweis aufgrund des gemischten Charakters
der Gleichungen schwierig sein kann. Der Zusatzterm stellt nun gewisserma"sen eine
"`Parabolisierung"'
der Schalengleichung dar. Insbesondere hat er die h"ohere Regularit"at $\eta\in
\widetilde Y^I$ zur Folge, wodurch, zusammen mit der zus"atzlichen Regularisierung der Wirbelterme,
nun zumindest formal die Zeitableitung von $(\pa_t\eta,\bu)$ im Dualraum von $L^2(I,H^2_0(M)))\times
L^p(I,W^{1,p}_{\dv,s}(\Omega_{\R\delta(t)}))$, der Regularit"atsklasse von $(\pa_t\eta,\bu)$, liegt.
Das Ziel ist, die Eindeutigkeit durch Aufstellen einer Energiegleichung, des Analogons von
\eqref{eqn:energiesatz}, zu beweisen. Diese Gleichung spielt auch beim
Beweis der Existenz obiger L"osungen eine wichtige Rolle. Augenscheinlich besitzen die schwachen
L"osungen unseres
parabolisch-dispersiven Systems nicht gen"ugend Regularit"at, um der Energieidentit"at Sinn zu
verleihen. Das legt die Vermutung nahe, dass die Energiegleichung f"ur derartige L"osungen nicht
gilt, was die N"utzlichkeit der Parabolisierung des Systems unterstreicht.

Sei nun $(\eta,\bu)$ ein schwache L"osung gem"a"s Definition \ref{def:pS}, wobei allerdings der Term
$S(D\bu)$ durch eine beliebige Matrix $\xi\in L^{p'}(\Omega_{\R\delta}^I)$ ersetzt sei. Wir m"ussen
die Energiegleichung ohne R"uckgriff auf die Theorie der Bochner-R"aume und das
Konzept der distributionellen Zeitableitung im Dualraum beweisen. Wie im Beweis von Proposition
\ref{lemma:komp} wollen wir stattdessen direkt mit der Gleichung arbeiten. Wir setzen dazu
\begin{equation*}
 \begin{aligned}
V:=L^2(I,H^2_0(M))) \times L^p(I,W^{1,p}_{\dv,s}(\Omega_{\R\delta(t)}))
 \end{aligned}
\end{equation*} 
und stellen die folgende Behauptung auf.

\begin{behauptung}\label{beh}
Es existieren hinreichend glatte Regularisierungen $(\eta_k,\bu_k)_{k\in\setN}$ mit
$\trrd\bu_k=\pa_t\eta_k\,\bnu$ und
\begin{equation}\label{eqn:konvp1}
 \begin{aligned}
  \bu_k&\rightarrow\bu&&\text{ in } L^2(\Omega_{\R\delta}^I),\\
  (\pa_t\eta_k,\bu_k)&\rightarrow(\pa_t\eta,\bu)&&\text{ in } V,\\
  (\pa_t^2\eta_k,\pa_t\bu_k)&\rightarrow \Sigma&&\text{ in } V'.
 \end{aligned}
\end{equation}
\end{behauptung}
\vspace{0.2cm}

\noindent Dabei ist die "`Zeitableitung"' $(\pa_t^2\eta_k,\pa_t\bu_k)\in V'$ definiert durch
\begin{equation*}
 \begin{aligned}
\langle(\pa_t^2\eta_k,\pa_t\bu_k),(b,\bphi)\rangle_V&:=\int_I\iortd \pa_t\bu_k\cdot\bphi\ dxdt +
\frac{1}{2}\int_I\im\pa_t\eta_k\, 
\pa_t\mathcal{R}\delta\, b\, \gamma(\mathcal{R}\delta)\ dAdt\\
&\hspace{0.5cm} + \int_I\im\pa_t^2\eta_k\, b\ dAdt,\\
 \end{aligned}
\end{equation*}
w"ahrend $\Sigma\in V'$ durch
\begin{equation*}
 \begin{aligned}
\langle\Sigma,(b,\bphi)\rangle_V&:=-\frac{1}{2}\int_I\iortd(\mathcal{R}
\bv\cdot\nabla)\R\bu\cdot\bphi\ dxdt +
\frac{1}{2}\int_I\iortd(\mathcal{R}\bv\cdot\nabla)\R\bphi\cdot\bu\ dxdt\\
&\hspace{0.55cm} - \int_I\iortd
\xi:D\bphi\ dxdt - 2\int_I K(\eta+\pa_t\eta,b)\ dt\\
&\hspace{0.55cm} + \int_I\iortd \ff\cdot\bphi\ dxdt +
\int_I\im g\, b\
dAdt\\
 \end{aligned}
\end{equation*}
gegeben ist. Die Konstruktion derartiger Approximationen sollte zwar keine gr"o"seren Probleme
bereiten, d"urfte aber aufw"andig und technisch sein. Wir wollen sie in
dieser Arbeit
nicht durchf"uhren. Offenbar gilt f"ur $s,t\in\bar I$, $s<t$
\begin{equation}\label{eqn:zeitabl}
 \begin{aligned}
\langle(\pa_t^2\eta_k,\pa_t\bu_k),(\pa_t\eta_k,\bu_k)\chi_
{(s,t)}\rangle_V&=\frac{1}{2} \int_{\Omega_{\R\delta(t)}} |\bu_k(t,\cdot)|^2\ dx +
\frac{1}{2} \im |\pa_t\eta_k(t,\cdot)|^2\
dA\\
&\hspace{0.5cm}-\frac{1}{2} \int_{\Omega_{\R\delta(s)}} |\bu_k(s,\cdot)|^2\ dx +
\frac{1}{2} \im |\pa_t\eta_k(s,\cdot)|^2\ dA.
 \end{aligned}
\end{equation}
Wenn wir in dieser Gleichung $\eta_k$ durch $\eta_k-\eta_l$ und $\bu_k$ durch $\bu_k-\bu_l$
ersetzen und die resultierende Identit"at bez"uglich $s$
integrieren, so folgt
\begin{equation*}
 \begin{aligned}
&\int_{\Omega_{\R\delta(t)}} |(\bu_k-\bu_l)(t,\cdot)|^2\ dx + \im
|\pa_t(\eta_k-\eta_l)(t,\cdot)|^2\ dA \\
&\hspace{0.5cm}\le
c\,\big(\norm{(\pa_t^2(\eta_k-\eta_l),\pa_t(\bu_k-\bu_l))}_{V'}^2+\norm{(\pa_t(\eta_k-\eta_l),
\bu_k-\bu_l)
}_V^2 +\norm{\bu_k-\bu_l}_{L^2(\Omega_{\R\delta}^I)}^2\big).
 \end{aligned}
\end{equation*}
Setzen wir die Felder $\bu_k$ und $\bu$ durch $\boldsymbol{0}$ auf $I\times \setR^3$ fort, so folgt
aus dieser Absch"atzung zusammen mit den Konvergenzen \eqref{eqn:konvp1}, dass die Folgen $(\bu_k)$
und $(\pa_t\eta_k)$ gegen $\bu$ in $C(\bar I, L^2(\setR^3))$ bzw. gegen $(\pa_t\eta)$ in $C(\bar I,
L^2(M))$ konvergieren. Setzen wir in \eqref{eqn:zeitabl} $s=0$, so konvergiert die linke Seite gegen
\begin{equation*}
 \begin{aligned}
  \langle\Sigma,(\pa_t\eta,\bu)\chi_{(0,t)}\rangle_V&=- \int_0^t\iortd
\xi:D\bu\ dxds - \int_0^t\frac{d}{ds} K(\eta) + 2 K(\pa_t\eta)\ ds \\
&\hspace{0.5cm} + \int_0^t\iortd \ff\cdot\bu\ dxds +
\int_0^t\im g\, \pa_t\eta\ dAds,\\
 \end{aligned}
\end{equation*}
w"ahrend die rechte Seite gegen die Funktion
\begin{equation*}
 \begin{aligned}
 & \frac{1}{2} \int_{\Omega_{\R\delta(t)}} |\bu(t,\cdot)|^2\ dx +
\frac{1}{2} \im |\pa_t\eta(t,\cdot)|^2\
dA\\
&-\frac{1}{2} \int_{\Omega_{\R\delta(0)}} |\bu(0,\cdot)|^2\ dx +
\frac{1}{2} \im |\pa_t\eta(0,\cdot)|^2\
dA
 \end{aligned}
\end{equation*}
konvergiert. Mit der Definition
\[E_{\eta,\bu}(t):= \frac{1}{2} \int_{\Omega_{\R\delta(t)}} |\bu(t,\cdot)|^2\ dx + \frac{1}{2} \im
|\pa_t\eta(t,\cdot)|^2\
dA + K(\eta(t,\cdot))\]
gilt also f"ur alle $t\in\bar I$ die Energiegleichung
\begin{equation*}
 \begin{aligned}
E_{\eta,\bu}(t)-E_{\eta,\bu}(0)&=- \int_0^t\iortd
\xi:D\bu\ dxds - 2\int_0^t K(\pa_t\eta)\ ds \\
&\hspace{0.5cm} + \int_0^t\iortd \ff\cdot\bu\ dxds +
\int_0^t\im g\, \pa_t\eta\
dAds.
 \end{aligned}
\end{equation*}
Im Unterschied zur Energiebilanz \eqref{eqn:energiesatz} tritt hier durch die Parabolisierung des
Systems ein zweiter dissipativer Term auf. Mithin ist der korrespondierende Ausdruck in der
Schalengleichung als D"ampfungsterm zu interpretieren. 
Wir wollen noch zeigen, dass
\begin{equation}\label{eqn:id}
 \begin{aligned}
E_{\eta,\bu}(0)=\frac{1}{2} \int_{\Omega_{\R\eta_0}} |\bu_0|^2\ dx + \frac{1}{2} \im
|\eta_1|^2\
dA + K(\eta_0)  
 \end{aligned}
\end{equation}
gilt. Wir wissen bereits, dass $\eta\in C(\bar I,H^2_0(M))$ und $\eta(0,\cdot)=\eta_0$ gilt, sodass
nur die ersten beiden Terme zu identifizieren sind. Wenn wir analog zum Beweis von
\eqref{eqn:schwacht} vorgehen, k"onnen wir mit Hilfe der Stetigkeit von $\bu$ und $\pa_t\eta$
zeigen, dass Gleichung \eqref{eqn:entp} mit
$\bu(0,\cdot)$ und $\pa_t\eta(0,\cdot)$ anstelle von $\bu_0$ und $\eta_1$ gilt, woraus wir
\eqref{eqn:id} folgern.
% Zudem k"onnen wir genau wie im Beweis von
% Theorem \ref{theorem:hs} zeigen, dass f"ur fast alle $t\in \bar I$, nach Obigem also sogar f"ur
% alle
% $t\in\bar I$, gilt
% \begin{equation}\label{eqn:entpt}
% \begin{aligned}
% % &\iot \bu(t,\cdot)\cdot\bphi(t,\cdot)\ dx 
% &- \int_t^T\iortd \bu\cdot\pa_t\bphi\ dxds - \frac{1}{2}\int_t^T\im\pa_t\eta\ 
% \pa_t\mathcal{R}\delta\ b\ \gamma(\mathcal{R}\delta)\ dAds\\
% &+\frac{1}{2}\int_t^T\iortd(\mathcal{R}\bv\cdot\nabla)\R\bu\cdot\bphi\ dxds -
% \frac{1}{2}\int_t^T\iortd(\mathcal{R}\bv\cdot\nabla)\R\bphi\cdot\bu\ dxds\\
% & + \int_t^T\iortd \chi:D\bphi\ dxds -\int_t^T\im\pa_t\eta\ \pa_tb\ dAds +
% 2\int_t^T K(\eta+\pa_t\eta,b)\
% ds\\
% % +\im\pa_t\eta(t,\cdot)b(t,\cdot)\ dA
% &\hspace{3.0cm}=\int_t^T\iortd \ff\cdot\bphi\ dxdt + \int_t^T\im g\ b\
% dAdt\\
% &\hspace{3.5cm} +\int_{\Omega_{\R\delta(t)}}\bu(t,\cdot)\cdot\bphi(t,\cdot)\ dx +
% \im \pa_t\eta(t,\cdot)\ b(t,\cdot)\ dA.
% \end{aligned}
% \end{equation}
% Zusammen mit \eqref{eqn:entp} folgt \eqref{eqn:id}. 

Um nun die Eindeutigkeit schwacher L"osungen im Sinne von Definition \ref{def:pS} einzusehen,
bemerken wir, dass die Differenz zweier L"osungen $(\eta,\bu)$ und $(\tilde\eta,\tilde\bu)$
Gleichung \eqref{eqn:entp} mit verschwindenden Daten $(\ff,g,\bu_0,\eta_0,\eta_1)$ erf"ullt, wenn
wir den Term $S(D\bu)$ durch $S(D\bu)-S(D\tilde\bu)$ ersetzen. Es gilt somit f"ur
alle $t\in\bar I$ die Energiegleichung
\begin{equation*}
 \begin{aligned}
E_{\eta-\tilde\eta,\bu-\tilde\bu}(t)&=- \int_0^t\iortd
S(D\bu)-S(D\tilde\bu):(D\bu-D\tilde\bu)\ dxds - 2\int_0^t K(\pa_t\eta-\pa_t\tilde\eta)\ ds
 \end{aligned}
\end{equation*}
Die linke Seite ist nichtnegativ, w"ahrend die rechte Seite aufgrund der Monotonie von $S$
nichtpositiv ist. Mithin verschwinden beide Seiten und somit auch die Differenz der L"osungen.

Kommen wir nun zur Existenz schwacher L"osungen. Der Wirbelterm spielt infolge der
Regularisierungen praktisch keine Rolle. Aufgrund der komplexen Struktur des Systems ist der
abstrakte Begriff des monotonen Operators hier dennoch
unbrauchbar. Wir k"onnen aber den klassischen Beweis ohne R"uckgriff auf die abstrakte Theorie
gewisserma"sen elementar nachzeichnen. Somit wird die Energiegleichung auch hier eine wichtige
Rolle spielen. Wir gehen zun"achst genau wie im
Falle des Navier-Stokes-Fluids vor, wobei die Felder
$(\widehat\bX_k)_{k\in\setN}$ nun eine Basis von $W^{1,p}_{0,\dv}(\Omega)$ bilden. Der
Galerkin-Ansatz f"uhrt uns auf nichtlineare Systeme gew"ohnlicher
Integro-Differentialgleichungen der in Anhang A.3 behandelten Form, zu denen wir lokale L"osungen
$(\eta_n,\bu_n)$ auf
Intervallen $I_n$ erhalten. Wir k"onnen allerdings wie zuvor Energieabsch"atzungen
\begin{equation*}
 \begin{aligned}
 \norm{\eta_n}_{\widetilde Y^{I_n}} +
\norm{\bu_n}_{X_{\R\delta,p}^{I_n}}\le c(T,\text{Daten})
  \end{aligned}
\end{equation*}
herleiten, die zeigen, dass die L"osungen auf dem ganzen Intervall $I$ existieren. Wegen
der Wachstumsbedingung an $S$ schlie"sen wir zudem, dass die Folge $(S(D\bu_n))$ in
$L^{p'}(\Omega_{\R\delta}^I)$ beschr"ankt ist. Somit erhalten
wir f"ur eine Teilfolge die Konvergenzen
\begin{equation*}
 \begin{aligned}
  \eta_n&\rightarrow \eta \hspace{0.3cm}&&\text{ schwach in } H^1(I,H^2_0(M)),\\
  \pa_t\eta_n&\rightarrow \pa_t\eta &&\text{ schwach$^*$ in } L^\infty(I,L^2(M)),\\
\bu_n&\rightarrow \bu &&\text{ schwach in } L^p(I,W^{1,p}_{\dv}(\Omega_{\R\delta(t)}))\\
& &&\text{ und
schwach$^*$ in } L^\infty(I,L^2(\Omega_{\R\delta(t)})),\\
S(D\bu_n)&\rightarrow \xi &&\text{ schwach in } L^{p'}(\Omega_{\R\delta}^I).
 \end{aligned}
\end{equation*}
Bez"uglich der dritten Konvergenz ist zu beachten, dass wegen des regul"aren Randes die Korn'sche
Ungleichung, Proposition \ref{theorem:korn}, und somit 
\[L^p(I,W^{1,p}_{\dv,s}(\Omega_{\R\delta(t)}^I))=
L^p(I,W^{1,p}_{\dv}(\Omega_{\R\delta(t)}^I))\]
gilt. Wie zuvor folgt $\trrd \bu=\pa_t\eta\,\bnu$. 
% Die
% Galerkin-L"osungen erf"ullen f"ur $\phi\in C^1_0([0,T))$ und $1\le j\le n$ die Gleichung
% \begin{equation*}
% \begin{aligned}
% % &\iot \bu(t,\cdot)\cdot\bphi(t,\cdot)\ dx 
% &-\int_I\iortd \bu_n\cdot\pa_t(\phi\bW_j)\ dxdt -\frac{1}{2}\int_I\im\pa_t\eta_n\
% \pa_t\mathcal{R}\delta\ \phi W_j\ \gamma(\mathcal{R}\delta)\
% dAdt\\
% & +\frac{1}{2}\int_I\iortd(\mathcal{R}\bv\cdot\nabla)\R\bu_n\cdot(\phi\bW_j)\ dxdt\\
% &-\frac{1}{2}\int_I\iortd(\mathcal{R}\bv\cdot\nabla)\R(\phi\bW_j)\cdot\bu_n\ dxdt
% +\int_I\iortd S(D\bu_n):D(\phi\bW_j)\ dxdt\\
% % +\im\pa_t\eta(t,\cdot)b(t,\cdot)\ dA
% &+\int_I\im\pa_t\eta_n\ \pa_t (\phi W_j)\ dAdt + 2\int_I\im K(\eta_n+\pa_t\eta_n, \phi W_j)\
% dAdt \\
% &\hspace{0.5cm}=\int_I\iortd \ff_n\cdot(\phi\bW_j)\ dxdt + \int_I\im g_n\ \phi W_j\
% dAdt\\
% &\hspace{1.5cm} +\int_{\Omega_{\mathcal{R}\eta_0}}\bu_n(0)\cdot(\phi(0)\bW_j(0,\cdot))\ dx +
% \im\pa_t\eta_n(0)\ \phi(0)W_j(0,\cdot)\ dA.
%  \end{aligned}
% \end{equation*}
Wenn wir in der Galerkin-Gleichung den Grenz"ubergang vollziehen, sehen wir, dass das
Tupel $(\eta,\bu)$ die Identit"at \eqref{eqn:entp} mit $\xi$ anstelle von $S(D\bu)$ f"ur alle
Testfunktionen erf"ullt. Es bleibt also lediglich,
das Feld $\xi$ zu identifizieren. Teilfolgen von $(\bu_n(T,\cdot))$ und $(\pa_t\eta_n(T,\cdot))$
konvergieren schwach gegen Felder $\bu^*$ in $L^2(\Omega_{\R\delta(T)})$ bzw. $\eta^*$ in $L^2(M)$.
Wenn wir die Funktionen $\bW_j$ und $W_j$ wie zuvor definieren, folgt daraus die Identit"at
\begin{equation*}
\begin{aligned}
% &\iot \bu(t,\cdot)\cdot\bphi(t,\cdot)\ dx 
&-\int_I\iortd \bu\cdot\pa_t(\phi\,\bW_j)\ dxdt -\frac{1}{2}\int_I\im\pa_t\eta\,
\pa_t\mathcal{R}\delta\ \phi\, W_j\, \gamma(\mathcal{R}\delta)\
dAdt\\
& +\frac{1}{2}\int_I\iortd(\mathcal{R}\bv\cdot\nabla)\R\bu\cdot(\phi\,\bW_j)\ dxdt\\
&-\frac{1}{2}\int_I\iortd(\mathcal{R}\bv\cdot\nabla)\R(\phi\,\bW_j)\cdot\bu\ dxdt
+\int_I\iortd \xi:D(\phi\,\bW_j)\ dxdt\\
% +\im\pa_t\eta(t,\cdot)b(t,\cdot)\ dA
&+\int_I\im\pa_t\eta\, \pa_t (\phi\, W_j)\ dAdt + 2\int_I\im K(\eta+\pa_t\eta, \phi\, W_j)\
dAdt \\
&\hspace{3.0cm}=\int_I\iortd \ff\cdot(\phi\,\bW_j)\ dxdt + \int_I\im g\, \phi\, W_j\
dAdt\\
&\hspace{3.5cm} -\int_{\Omega_{\mathcal{R}\delta(T)}}\bu^*\cdot(\phi(T)\,\bW_j(T,\cdot))\ dx -
\im\eta^*\, \phi(T)\,W_j(T,\cdot)\ dA
 \end{aligned}
\end{equation*}
f"ur alle $j\ge 1$ und $\phi\in C_0^1((0,T])$. Wir k"onnen analog zum Beweis von
\eqref{eqn:schwacht} vorgehen und unter Verwendung der Stetigkeit von $\bu$ und $\pa_t\eta$ zeigen,
dass diese Identit"at auch mit $\bu(T,\cdot)$ und
$\pa_t\eta(T,\cdot)$ anstelle von $\bu^*$ und $\eta^*$ gilt, woraus wir $\bu^*=\bu(T,\cdot)$ und
$\pa_t\eta(T,\cdot)=\eta^*$ folgern. Zudem konvergiert eine Teilfolge von $(\eta_n(T,\cdot))$
schwach gegen $\eta(T,\cdot)$ in $H^2_0(M)$. F"ur die Galerkin-L"osungen gilt der Energiesatz
\begin{equation*}
 \begin{aligned}
&\int_I\int_{\Omega_{\R\delta(t)}}S(D\bu_n):D\bu_n\ dxdt + 2\int_I K(\pa_t\eta_n)\ dt\\ 
&\hspace{1cm}=-E_{\eta_n,\bu_n}(T)+E_{\eta_n,\bu_n}(0) + \int_I\iortd \ff_n\cdot\bu_n\ dxdt
+\int_I\im g_n\,\pa_t\eta_n\ dAdt.
 \end{aligned}
\end{equation*}
Wenn wir den $\limsup$ dieser Gleichung nehmen und die schwache Unterhalbstetigkeit der Energie $E$
verwenden,\footnote{Man beachte, dass jede stetige, nichtnegative quadratische Form, insbesondere
also $K$, schwach unterhalbstetig ist. Das folgt, wenn man den $\liminf$ der Ungleichung
\[0\le K(\eta_n-\eta,\eta_n-\eta)=K(\eta_n)-2K(\eta,\eta_n)+K(\eta)\] nimmt.} so folgt
\begin{equation*}
 \begin{aligned}
&\limsup_n \int_I\int_{\Omega_{\R\delta(t)}}S(D\bu_n):D\bu_n\ dxdt + 2\int_I K(\pa_t\eta_n)\ dt\\ 
&\hspace{1cm}\le-E_{\eta,\bu}(T)+E_{\eta,\bu}(0) + \int_I\iortd \ff\cdot\bu\ dxdt +\int_I\im g\,
\pa_t\eta\ dAdt.
 \end{aligned}
\end{equation*}
Man beachte dabei, dass $\eta_n(0,\cdot)=\eta_0$ f"ur alle $n\in\setN$ gilt. Aus der
Energiegleichung f"ur die schwache L"osung $(\eta,\bu)$ (mit $\xi$ anstelle von $S(D\bu)$) folgt
\begin{equation*}
 \begin{aligned}
&\limsup_n \int_I\int_{\Omega_{\R\delta(t)}}S(D\bu_n):D\bu_n\ dxdt + 2\int_I K(\pa_t\eta_n)\ dt\\ 
&\hspace{4cm}\le\int_I\int_{\Omega_{\R\delta(t)}}\xi:D\bu\ dxdt + 2\int_I K(\pa_t\eta)\ dt.
 \end{aligned}
\end{equation*} 
Unter Ausnutzung dieser Absch"atzung und der schwachen Konvergenzen schlie"sen wir
\begin{equation*}
 \begin{aligned}
  0&\le\limsup_n\Big(\int_I\int_{\Omega_{\R\delta(t)}}(S(D\bu_n)-S(D\bu)):(D\bu_n-D\bu)\ dxdt\\
&\hspace{4.3cm}+2\int_I K(\pa_t\eta_n-\pa_t\eta,\pa_t\eta_n-\pa_t\eta)\ dt\Big)\\
&=\limsup_n\Big(\int_I\int_{\Omega_{\R\delta(t)}}S(D\bu_n):D\bu_n + S(D\bu):D\bu\\
&\hspace{4.75cm} - S(D\bu_n):D\bu -
S(D\bu):D\bu_n\ dxdt\\
&\hspace{2cm}+2\int_I K(\pa_t\eta_n,\pa_t\eta_n) + K(\pa_t\eta,\pa_t\eta)-
2K(\pa_t\eta_n,\pa_t\eta)\ dt\Big)\\
&\le \int_I\int_{\Omega_{\R\delta(t)}}\xi:D\bu + S(D\bu):D\bu - \xi:D\bu - S(D\bu):D\bu\ dxdt\\
&\hspace{2.2cm}+2\int_I K(\pa_t\eta,\pa_t\eta) + K(\pa_t\eta,\pa_t\eta)-
2K(\pa_t\eta,\pa_t\eta)\ dt\\
&=0.
 \end{aligned}
\end{equation*}
Somit gilt f"ur eine Teilfolge die Konvergenz
\[(S(D\bu_n)-S(D\bu)):(D\bu_n-D\bu)\rightarrow 0\]
fast "uberall in $\Omega_{\R\delta}^I$. Mit Hilfe von Proposition \ref{lemma:dmm}
schlie"sen wir
$D\bu_n\rightarrow D\bu$  und somit $S(D\bu_n)\rightarrow S(D\bu)$ fast "uberall in
$\Omega_{\R\delta}^I$. Proposition \ref{lemma:vitali} liefert uns nun die Identifizierung
$\xi=S(D\bu)$. Unter der Annahme der Existenz obiger Approximationen ist damit die Existenz und
Eindeutigkeit schwacher L"osungen des entkoppelten, regularisierten p-Systems gezeigt.

Auch bei der Durchf"uhrung des Fixpunktarguments gehen wir zun"achst genau wie im Falle des
Navier-Stokes-Fluids vor. Wir ersetzen lediglich die Absch"atzung
\eqref{ab:schranke} durch 
\begin{equation*}
\begin{aligned}
\norm{\eta}_{\widetilde Y^I} +
\norm{\bu}_{X_{\R\delta,p}^I} \le c.  
\end{aligned}
\end{equation*}
F"ur $(\delta,\bv)\in D$ ist die Menge
$F(\delta,\bv)$ aufgrund der Eindeutigkeit der schwachen L"osungen einelementig, insbesondere
konvex und abgeschlossen. Desweiteren k"onnen wir den Beweis von Proposition
\ref{lemma:komp} fast w"ortlich "ubernehmen, um die Kompaktheit von $F$ zu zeigen. Zus"atzlich zu
\eqref{ab:b} verwenden wir dabei die Absch"atzung
\begin{equation*}
 \begin{aligned}
\norm{\F_{\eta_n}\M_{\eta_n}b}_{L^\infty(I,W^{1,p}(B_\alpha))}
\le c\,\norm{b}_{H^2_0(M)},
 \end{aligned}
\end{equation*}
die aus Lemma \ref{lemma:mittelwert} und Proposition \ref{lemma:FortVonRand} folgt. Zudem ersetzen
wir in der Definition von $h_n^\sigma$ den Raum $H^1_0(\Omega)$ durch $W^{1,\max(2,p)}_0(\Omega)$,
sodass wir zus"atzlich zu \eqref{eqn:schnarch} "uber die Absch"atzung
\[\norm{\T_{\delta_\sigma}\bphi}_{L^\infty(I,W^{1,p}(B_\alpha))}\le
c\,\norm{\bphi}_{W^{1,p}_0(\Omega)}\]
verf"ugen. Bei der Anwendung von Lemma \ref{lemma:ehrling} tritt die Einschr"ankung $p>3/2$
auf.\footnote{Im Prinzip l"asst sich die Einschr"ankung an dieser Stelle auf $p>6/5$ absenken,
da die Kompaktheit von $(\pa_t\eta_n)$ in $L^2(I\times M)$ an dieser Stelle noch nicht ben"otigt
wird. Die Einschr"ankung kommt n"amlich dadurch zustande, dass wir die r"aumliche Regularit"at der
Folge $(\pa_t\eta_n)$ durch Spurbildung aus der r"aumlichen Regularit"at von $(\bu_n)$ gewinnen und
somit der Raum $W^{1-1/r,r}(M)$ f"ur ein beliebiges $r<p$ kompakt nach $L^2(M)$ einbetten muss. Beim
finalen Grenz"ubergang wird die Einschr"ankung (mit dieser Beweismethode) jedoch vermutlich nicht zu
umgehen
sein.} Den Exponenten $p$ in Proposition
\ref{lemma:nullfort} w"ahlen wir identisch $\min(2,p)$. Auch der Beweis der Graphenabgeschlossenheit
von $F$ funktioniert fast genauso wie zuvor. Die Konvergenz der Folge $(\F_{\eta_n}\M_{\eta_n}b)$ in
$L^p(I,W^{1,p}(B_\alpha))$ erh"alt man dabei wie im Beweis von Lemma \ref{lemma:konvergenzen}. F"ur
die Identitfizierung des Grenzwerts der viskosen Spannungstensoren $(S(D\bu_n))\subset
L^{p'}(I\times \setR^3)$ k"onnen wir die Argumentation in obigem Existenzbeweis fast w"ortlich
"ubernehmen. Im Wesentlichen ersetzen wir die dortigen
Integrale "uber $\Omega_{\R\delta}^I$ durch Integrale "uber $I\times \setR^3$, wobei wir die
entsprechenden Felder durch $0$ auf $I\times \setR^3$ fortsetzen.

Unter der Annahme der Existenz der Approximationen und $3/2<p<\infty$ haben wir somit gezeigt, dass
sich ein Intervall $I$ und zu jedem hinreichend kleinen Parameter $\epsilon>0$ eine schwache
L"osungen des regularisierten p-Systems zum Parameter $\epsilon$ auf dem Intervall $I$ gem"a"s der
nachfolgenden Definition finden l"asst.

\begin{definition}
Ein Tupel $(\eta,\bu)$ hei"st schwache L"osung des
regularisierten $p$-Systems zum
Parameter $\epsilon$ auf dem Intervall $I$, falls $\eta\in \widetilde Y^I$
mit $\norm{\eta}_{L^\infty(I\times M)}<\kappa$ und $\eta(0,\cdot)=\eta_0$, $\bu\in
X_{\mathcal{R}_\epsilon\eta,p}^I$ mit
$\trre \bu =\pa_t\eta\,\bnu$ und
\begin{equation*}
\begin{aligned}
% &\iot \bu(t,\cdot)\cdot\bphi(t,\cdot)\ dx 
&- \int_I\iorte \bu\cdot\pa_t\bphi\ dxdt - \frac{1}{2}\int_I\im\pa_t\eta
\, \pa_t\mathcal{R}_\epsilon\eta\ b\, \gamma(\mathcal{R}_\epsilon\eta)\ dAdt\\
&+\frac{1}{2}\int_I\iorte(\mathcal{R}_\epsilon\bu\cdot\nabla)\mathcal{R}_\epsilon\bu\cdot\bphi\ dxdt
-\frac{1}{2}\int_I\iorte(\mathcal{R}_\epsilon\bu\cdot\nabla)\mathcal{R}_\epsilon\bphi\cdot\bu\
dxdt\\
& + \int_I\iorte S(D\bu):D\bphi\ dxdt-\int_I\im\pa_t\eta\, \pa_tb\ dAdt + 2\int_I
K(\eta+\epsilon\,\pa_t\eta,b)\ dt\\ 
% +\im\pa_t\eta(t,\cdot)b(t,\cdot)\ dA
&\hspace{4cm}=\int_I\iorte \ff\cdot\bphi\ dxdt + \int_I\im g\, b\
dAdt\\
&\hspace{4.5cm} +\int_{\Omega_{\R_\epsilon\eta_0}}\bu_0^\epsilon\cdot\bphi(0,\cdot)\ dx +
\im\eta_1^\epsilon\, b(0,\cdot)\ dA
 \end{aligned}
\end{equation*}
f"ur alle Testfunktionen $(b,\bphi)\in  T_{\R_\epsilon\eta,p}^I$.
\end{definition}
\vspace{0.2cm}

Um den Grenz"ubergang $\epsilon\rightarrow0$ vollziehen zu k"onnen, werden nun allerdings andere
Techniken ben"otigt. F"ur $p<11/5$ ist die Energieidentit"at selbst f"ur die
Fluidgleichungen in einem Raumzeitzylinder ohne jegliche Kopplung vermutlich falsch. Dementsprechend
sind f"ur kleine $p$ schon in diesem einfacheren Fall andere
Techniken vonn"oten. Da mit $\epsilon$
auch die Parabolisierung verschwindet, ist die Energieidentit"at f"ur das gekoppelte,
nichtregularisierte System selbst f"ur gro"se $p$ vermutlich nicht richtig. Auch die asymptotisch
geringe Randregularit"at und das resultierende Versagen der Korn'schen Gleichung k"onnten Probleme
bereiten.

Zur Behandlung der Fluidgleichungen f"ur kleine $p$ in einem Raumzeitzylinder ohne Kopplung haben
sich die Methoden
der lokalen Druckfelder und der parabolischen $L^\infty$- bzw. $W^{1,\infty}$-Abschneidungen als
besonders leistungsf"ahig und flexibel erwiesen; siehe \cite{b32}. Mit Hilfe dieser Techniken l"asst
sich die Existenz schwacher L"osungen f"ur $p>6/5$ zeigen. Das Ziel ist, auch hier diese Techniken
erfolgreich einzubringen. Dies wird Gegenstand weiterer Forschung sein.
\chapter{Ausblick}

Der n"achste Schritt wird der vollst"andige Existenzbeweis im Falle der verallgemeinerten
Newton'schen Fluide sein. Hier ist noch einige Arbeit zu leisten. 
% Zudem soll gepr"uft
% werden, ob sich die Einschr"ankung der Bewegung der Schale in Normalenrichtung f"ur kleine
% Schalendicken und kleine Auslenkungen physikalisch rechtfertigen l"asst. 
Zudem k"onnte man
unter Beibehaltung der Einschr"ankung der Auslenkungen in Richtung der Normale an $\pa\Omega$ zur
Koiter-Energie f"ur nichtlinear elastische Schalen
"ubergehen. Trotz der
Einschr"ankung wird dabei allerdings ein Elliptizit"atsverlust auftreten, da die
Entartungsrichtungen im nichtlinearen Fall ja
mit der L"osung variieren und nicht l"anger tangential an $\pa\Omega$ liegen werden. Interessant
k"onnte auch die Konstruktion starker L"osungen f"ur kurze
Zeiten bei Newton'schen wie bei verallgemeinerten
Newton'schen Fluiden sein. Insbesondere in letzterem Fall scheint
sich die Technik der maximalen $L^p$-Regularit"at anzubieten, wie sie in \cite{b33} zur
Konstruktion starker Kurzzeit-L"osungen f"ur verallgemeinerte
Newton'sche Fluide in Raumzeitzylindern ohne zus"atzliche Kopplung verwendet wird. Auch bei
der Konstruktion von Kurzzeitl"osungen sollte man die M"oglichkeit, die Koiter-Energie f"ur
nichtlinear elastische Schalen zu verwenden, in Betracht ziehen. Eine weitere m"ogliche
Sto"srichtung besteht darin, die Einschr"ankung der Schalenauslenkung auf die Normalenrichtung
aufzuheben. Hier tritt, wie bereits angemerkt, das
Problem auf, dass der Gradient der Koiter-Energie in tangentiale Richtungen degeneriert
ist. Zudem ist der Rand in diesem Fall im Allgemeinen kein Graph "uber $\pa\Omega$ mehr.
\begin{appendix}
\addchap{Anhang}
\addtocounter{chapter}{1}
%\setcounter{Def}{0} % Falls Equations und Rest unterschiedlich nummeriert werden sollen.
%\setcounter{equation}{0} 
%\markboth{Anhang}{Anhang}

\section{Differentialgeometrie} Details und gegebenenfalls Beweise zu den folgenden
Ausf"uhrungen finden sich in \cite{b1}, \cite{b2}, \cite{b3}. Wir nehmen im Folgenden an, dass alle
auftretenden mathematischen Objekte so regul"ar sind, dass die Definitionen sinnvoll und die
durchgef"uhrten Operationen zul"assig
sind. Sei $M$ eine kompakte Riemann'sche Mannigfaltigkeit (berandet oder nicht) endlicher Dimension
mit Riemann'scher Metrik $g=\langle\cdot,\cdot\rangle$. Bez"uglich beliebiger Koordinaten
bezeichnen wir die Komponenten von
$g$ mit $g_{\alpha\beta}$, die Koordinatenvektorfelder mit $\pa_\alpha$ und die
Koordinaten-1-Formen mit $dx^\alpha$. Das Skalarprodukt auf dem Tangentialb"undel induziert ein
Skalarprodukt $\langle\cdot,\cdot\rangle$ auf s"amtlichen Tensorb"undeln; zum Beispiel f"ur
$(0,2)$-Tensorfelder $T,S$ in Koordinaten $\langle T,S \rangle =
g^{\alpha\gamma} g^{\beta\delta} T_{\alpha\beta} S_{\gamma\delta}$. Dabei ist $(g^{\alpha\gamma})$
die inverse Matrix von $(g_{\alpha\gamma})$. F"ur beliebige Tensorfelder $T$ setzen wir
$|T|^2:=\langle
T,T\rangle$. Desweiteren l"asst sich ein Tensorprodukt definieren, das ein $(k,l)$-Tensorfeld und
ein $(r,s)$-Tensorfeld auf ein $(k+r,l+s)$-Tensorfeld abbildet; zum Beispiel f"ur zwei
$(1,1)$-Tensorfelder
$T,S$ in Koordinaten $(T\otimes S)_{\alpha\gamma}^{\beta\delta}=T_\alpha^\beta S_\gamma^\delta$. 
Die Riemann'sche Metrik induziert einen kanonischen Isomorphismus zwischen Tangential- und
Kotangentialb"undel und damit auch zwischen ko- und kontravarianten Tensorfeldern. Zum Beispiel
l"asst sich ein $(0,2)$-Tensorfeld $T$ in ein $(1,1)$-Tensorfeld umwandeln; in
Koordinaten $T_\alpha^\beta=g^{\beta\delta}T_{\alpha\delta}$. Dieser Vorgang nennt sich Indexziehen.
Desweiteren ist ein $(1,1)$-Tensorfeld $T$ zu einem skalaren Feld kontrahierbar. Ein solches
Feld l"asst sich n"amlich auch so interpretieren, dass punktweise Vektoren linear auf lineare
Funktionale abgebildet werden, die Kovektoren aufnehmen. Ein lineares Funktional, das Kovektoren
aufnimmt, ist (da jeder endlichdimensionale, normierte Raum reflexiv ist) aber ein Vektor. D.h. ein
solches Tensorfeld definiert punktweise einen Endomorphismus des Tangentialraums. Von diesem
Endomorphismenfeld nehmen wir punktweise die Spur, was ein skalares Feld ergibt; in Koordinaten
$\trace T=T^\alpha_\alpha$. Allgemein l"asst sich die Spur beliebiger Tensorfelder bez"uglich
eines ko-
und eines kontravarianten Index nehmen, was den Rang des Tensorfelds um $2$ vermindert; zum Beispiel
f"ur ein $(2,1)$-Tensorfeld $T$ in Koordinaten $T^{\alpha\beta}_\beta$. Durch
Verkn"upfen von Indexziehen und Spurbildung l"asst sich jedes Tensorfeld "uber zwei
verschiedene Indizes kontrahieren; zum Beispiel f"ur ein $(0,2)$-Tensorfeld $T$ in
Koodinaten $\trace_g T=g^{\alpha\beta}T_{\alpha\beta}$.

Es existiert ein kanonischer linearer Zusammenhang $\nabla$ auf dem
Tangentialb"undel, der Levi-Civita-Zusammenhang, der dadurch charakterisiert ist, dass er
symmetrisch und Riemann'sch ist. Sei $X$ ein Vektorfeld auf $M$; in Koordinaten
$X=X^\alpha\pa_\alpha$. Dann bezeichnen wir mit $\nabla X$ die totale kovariante Ableitung von $X$
bez"uglich $\nabla$. Dies ist ein $(1,1)$-Tensorfeld, d.h. es nimmt (punktweise) einen Vektor (die
Richtung, in die abgeleitet wird) und einen Kovektor (der das Ergebnis der Ableitung aufnimmt) als
Argument auf. In Koordinaten gilt $(\nabla X)(\pa_\alpha,dx^\beta)=
dx^\beta(\nabla_{\pa_\alpha} X)=\pa_\alpha X^\beta + \Gamma_{\alpha\delta}^\beta X^\delta$, wobei
die Christoffel-Symbole $\Gamma_{\alpha\beta}^\gamma$ durch
$\Gamma_{\alpha\beta}^\gamma\pa_\gamma:=\nabla_{\pa_\alpha}\pa_\beta$ definiert sind. F"ur
skalare Funktionen $f$ setzen wir zudem $\nabla f:=df$ mit dem Differential $df$ von $f$,
ein Kovektorfeld; in Koordinaten $df(\pa_\alpha)=\pa_\alpha f$. Durch Indexziehen erhalten wir aus
$\nabla f$ das Gradientenfeld $\grad f$; in Koordinaten $(\grad
f)^\alpha=g^{\alpha\beta}\pa_\beta f$. Die kovariante Ableitung von
Vektorfeldern induziert eine kovariante Ableitung beliebiger Tensorfelder; f"ur ein
Kovektorfeld $\omega$ und Vektorfelder $X,Y$ zum Beispiel gilt $(\nabla\omega)(X,Y)=\nabla_X
(\omega(Y))-\omega(\nabla_X Y)$, in Koordinaten 
$(\nabla\omega)(\pa_\alpha,\pa_\beta)=\pa_\alpha\omega_\beta - \Gamma_{\alpha\beta}^\delta
\omega_\delta$. Der $L^2$-adjungierte Operator der totalen kovarianten Ableitung ist die Divergenz
$\nabla^*=\trace\nabla$. Dabei wird "uber den Index der Ableitung und einen Index des
Tensorfeldes, von dem wir die Divergenz nehmen wollen, kontrahiert; zum Beispiel f"ur ein
Vektorfeld $X$ in Koordinaten 
\begin{equation}\label{eqn:div}
 \begin{aligned}
\nabla^* X=\pa_\alpha X^\alpha + \Gamma_{\alpha\delta}^\alpha
X^\delta=(\det(g_{\beta\gamma}))^{-\frac{1}{2}}\pa_\alpha((\det(g_{\beta\gamma}))^{\frac{1}{2}}
X^\alpha).
 \end{aligned}
\end{equation}
Zudem setzen wir $\nabla^2=\nabla\nabla$. F"ur ein skalares Feld $f$ zum Beispiel ist $\nabla^2 f$
die totale kovariante Ableitung des Differentials, ein symmetrisches $(0,2)$-Tensorfeld, und $\Delta
f=\nabla^*\nabla f=\trace_g \nabla^2 f$ der
Laplace(-Beltrami)-Operator auf $f$
angewendet. Dieser ist f"ur beliebige Tensorfelder durch $\Delta=\nabla^*\nabla=\trace_g \nabla^2$
definiert.

Der Integralsatz von Stokes (f"ur Differentialformen) beinhaltet den Spezialfall
\begin{equation}\label{eqn:divsatz}
 \begin{aligned}
\int_M \dv X\ dV=0,  
 \end{aligned}
\end{equation}
falls das Vektorfeld $X$ auf dem (m"oglicherweise leeren) Rand von $M$ verschwindet. Dabei
ist $dV$ das durch die Metrik induzierte Ma"s. Sind nun $T$
ein $(0,k)$-Tensorfeld und $S$ ein $(0,k+1)$-Tensorfeld, beide mit verschwindenden Randwerten,
so gilt
\begin{equation}\label{eqn:pisatz}
 \begin{aligned}
  \int_M \langle \nabla T,S\rangle\ dV = -\int_M\langle T,\trace_g\nabla S\rangle\ dV,
 \end{aligned}
\end{equation}
wobei "uber die ersten beiden Indizes von $\nabla S$ kontrahiert wird. Diese Aussage folgt, wenn wir
\eqref{eqn:divsatz} auf das Vektorfeld $X:=(\trace_g)^k\ T\otimes S$ anwenden, wobei wir den ersten
Index von $T$ mit dem zweiten Index von $S$ kontrahieren, den zweiten mit dem dritten, etc. Dabei
ist lediglich zu beachten, dass die Operatoren $\nabla$ und $\trace_g$ ebenso wie Kontraktionen
bez"uglich verschiedener Indizes kommutieren und dass
\[\nabla X=(\trace_g)^k\ (\nabla T)\otimes S + (\trace_g)^k\ T\otimes (\nabla S)\]
gilt. 

Wir definieren f"ur Vektorfelder $X,Y,Z$ den Riemann'schen
Kr"ummungstensor $R$ durch 
\[R(X,Y)Z:=\nabla^2_{X,Y}Z-\nabla^2_{Y,X}Z\]
und den Ricci-Tensor $Rc$ durch Kontraktion von $R$ "uber den ersten und den letzten
Index; in Koordinaten $Rc_{\alpha\beta}=R_{\gamma\alpha\beta}^{\gamma}$. F"ur ein skalares Feld $f$
haben wir
\begin{equation}\label{eqn:kommutator}
 [\Delta,\nabla] f=Rc(\grad f,\,\cdot\,),
\end{equation}
denn in lokalen Koordinaten gilt
\begin{equation*}
 \begin{aligned}
  (\Delta\nabla f) (\pa_\gamma)&=g^{\alpha\beta}\,(\nabla^3 f)(\pa_\alpha,\pa_\beta,\pa_\gamma)=
g^{\alpha\beta}\,(\nabla^3 f)(\pa_\alpha,\pa_\gamma,\pa_\beta)\\
&=g^{\alpha\beta}\,(\nabla^3 f)(\pa_\gamma,\pa_\alpha,\pa_\beta) + g^{\alpha\beta}\,\langle
R(\pa_\alpha,\pa_\gamma)\grad f,\pa_\beta\rangle\\
&=(\nabla\Delta f)(\pa_\gamma) + Rc(\grad f,\pa_\gamma).
\end{aligned}
\end{equation*}

F"ur eine weitere Riemann'sche Mannigfaltigkeit $N$ und eine Abbildung
$\Phi:M\rightarrow N$ bezeichnen wir mit $d\Phi: TM \rightarrow TN$ das Differential von $\Phi$.
Dabei sind $TM$ und $TN$ die
Tangentialb"undel von $M$ bzw. $N$. Ist $v$ ein Tangentialvektor an $M$ und $c:t\mapsto c(t)$ ein
Kurve in $M$ mit $\frac{d}{dt}\big|_{t=0}\ c(t)=v$, so gilt
$d\Phi\, v=\frac{d}{dt}\big|_{t=0}\ \Phi\circ c(t)$. F"ur Funktionen
$f:N\rightarrow\setR$ gilt der Transformationssatz
\begin{equation}\label{eqn:trafo}
 \begin{aligned}
\int_N f\ dA = \int_M f\circ \Phi\ |\det d\Phi|\ dA.  
 \end{aligned}
\end{equation}
Dabei wird die Determinante von $d\Phi$ im Punkt $q\in M$ bez"uglich zweier
beliebiger Orthonormalbasen von $T_qM$ und $T_{\Phi(q)}N$ gebildet. Dieser Kunstgriff ist
notwendig, weil die Determinante von Homomorphismen zwischen \emph{verschiedenen} Vektorr"aumen
keine koordinateninvariante Bedeutung besitzt. Die Determinanten bez"uglich beliebiger
Orthonormalbasen hingegen k"onnen h"ochstens um den Faktor $-1$ differieren. Ist $T$ ein Tensorfeld
auf $N$, so bezeichnen wir mit
$\Phi^*T$ den Pullback von $T$ oder den zur"uckgeholten Tensor. Ist $T$ ein skalares Feld, so ist
dieser durch $\Phi^*T=T\circ\Phi$ definiert; f"ur ein ein $(0,2)$-Tensorfeld $T$ und Vektorfelder
$X,Y$ auf $M$ gilt $(\Phi^*T)(X,Y)=T(d\Phi\, X, d\Phi\, Y)$. Die
kovariante
Ableitung verh"alt sich nat"urlich unter Isometrien. Ist $\Phi$ eine Isometrie, so bedeutet das
speziell f"ur die Divergenz eines Vektorfeldes $X$ auf $M$ die Identit"at $\nabla^*
X=(\widetilde\nabla^* Y)\circ\Phi$, wobei $\widetilde\nabla$ der Levi-Civita-Zusammenhang von $N$
und $Y:=(d\Phi\, X)\circ\Phi^{-1}$ der
Pushforward von $X$ unter $\Phi$ ist. 

Ist $(\phi_k,U_k)_{k}$ ein endlicher Atlas (insbesondere $U_k\subset M$ offen), so existiert eine
untergeordnete
Zerlegung der Eins, d.h. es existieren differenzierbare Funktionen $\psi_k: M\rightarrow
\setR$ mit $0\le\psi_k\le 1$, $\supp\psi_k\subset U_k$ und
$\sum_k\psi_k(q)=1$ f"ur alle $q\in M$.

Ist $M$ speziell eine kompakte, orientierte Fl"ache in $\setR^3$, so wird durch Einschr"anken des
euklidischen Skalarprodukts eine Riemann'sche Metrik auf $M$ definiert. Der Levi-Civita-Zusammenhang
ist in diesem Fall durch Differenzieren von Vektorfeldern (l"angs Kurven auf $M$) im $\setR^3$ und
anschlie"sende Orthogonalprojektion auf den Tangentialraum gegeben. Die Kr"ummung der
Fl"ache wird durch die zweite Fundamentalform $h$, ein
symmetrisches $(0,2)$-Tensorfeld, beschrieben. Diese ist f"ur Vektorfelder $X,Y$ durch
$h(X,Y):=\frac{\pa Y}{\pa
X}\cdot\nu$ definiert, wobei $\frac{\pa Y}{\pa X}$ die Ableitung von $Y$ in Richtung $X$ im
$\setR^3$ und $\nu$ die Normale an $M$ bezeichnet. Durch Indexziehen erhalten wir ein
Endomorphismenfeld, die Weingarten-Abbildung $W$. Fassen wir die Normale $\nu$ als
Abbildung von
$M$ in die $2$-Sph"are auf, so l"asst sich das Differential $d\nu$ als Endomorphismus des
Tangentialb"undels von $M$ interpretieren. Dann gilt $W=-d\nu$. Die Weingarten-Abbildung
ist symmetrisch, und ihre  Eigenwerte hei"sen Hauptkr"ummungen. Die
Determinante der Weingarten-Abbildung ist die Gau"s-Kr"ummung $G$; in Koordinaten $G=\det
(h^\alpha_\beta)$. Die gemittelte Spur ist die
mittlere Kr"ummung $H$; in Koordinaten
$H=\frac{1}{2}g^{\alpha\beta}h_{\alpha\beta}=\frac{1}{2}h^\alpha_\alpha$. Der Riemann'sche
Kr"ummungstensor l"asst sich durch die zweite Fundamentalform ausdr"ucken; in lokalen Koordinaten
$R_{\alpha\beta\gamma}^\delta= h_{\beta\gamma}\, h_\alpha^\delta - h_{\alpha\gamma}\,
h_\beta^\delta$. Entsprechend hat der Ricci-Tensor in lokalen
Koordinaten die Form 
\begin{equation}\label{eqn:ricci}
 \begin{aligned}
Rc_{\beta\gamma}=h_{\beta\gamma}\, h_\alpha^\alpha - h_{\alpha\gamma}\,h_\beta^\alpha.  
 \end{aligned}
\end{equation}
\\

\section{Sobolev-R"aume}

Eine Einf"uhrung in die Sobolev-R"aume auf euklidischen Gebieten findet sich in
\cite{b4}; siehe auch \cite{b48}, \cite{b55}, \cite{b46}. F"ur eine Einf"uhrung in die
Bochner-R"aume sei auf \cite{b52}, \cite{b7} verwiesen. F"ur $d,k\in\setN$,
$\Omega\subset\setR^d$ offen und $1\le p\le\infty$
bezeichnen wir mit $W^{k,p}(\Omega)$ den Sobolev-Raum der
reellwertigen $L^p(\Omega)$-Funktionen, deren distributionelle Ableitungen bis zur Ordnung $k$ in
$L^p(\Omega)$
liegen. F"ur $s\in (k-1,k)$ und $1\le p<\infty$ setzen wir zudem
$W^{s,p}(\Omega):=B^s_{pp}(\Omega):=(L^p(\Omega),W^{k,p}(\Omega))_{\frac{s}{k},p}$. Ist $\Omega$
der $\setR^d$ oder ein
beschr"anktes Gebiet mit Lipschitz-Rand, so ist der Raum $W^{s,p}(\Omega)$ isomorph zum (kanonisch
normierten) Raum der $W^{k-1,p}(\Omega)$-Funktionen $f$, f"ur die die Gr"o"se
\begin{equation*}
 \begin{aligned}
|\pa^\alpha f|_{\sigma,p;\Omega}^p:=\int_{\Omega}\int_{\Omega}
\frac{|\pa^\alpha f(x)-\pa^\alpha f(y)|^p}{|x-y|^{d+\sigma p}}\ dxdy  
 \end{aligned}
\end{equation*}
f"ur alle $|\alpha|=k-1$ endlich ist; siehe \cite{b4}, \cite{b48}, \cite{b55}, \cite{b46}.

\begin{theorem}\label{theorem:emb}
Es sei $\Omega\subset\setR^d$, $d\in\setN$, ein beschr"anktes Gebiet mit Lipschitz-Rand. F"ur $1\le
p<\infty$ und $s>0$ mit $0<\alpha:=s-d/p< 1$ gilt
\[W^{s,p}(\Omega)\embedding C^{0,\alpha}(\overline\Omega).\]
F"ur $1\le p, \tilde p<\infty$ und $0<\tilde s<s$ mit $s-d/p>\tilde s-d/\tilde p$ gilt
\[W^{s,p}(\Omega)\embedding W^{\tilde s,\tilde p}(\Omega).\] 
F"ur $1\le p<\infty$ und $s>0$ haben wir die kompakte Einbettung
\[W^{s,p}(\Omega)\compactembedding L^p(\Omega).\]
\end{theorem}
\beweis Die ersten beiden Behauptungen folgen mit Hilfe der in \cite{b4} konstruierten
Fortsetzungsoperatoren aus Theorem 2.8.1 in \cite{b55}. Die kompakte Einbettung ist wegen Theorem
3.8.1 in \cite{b48} eine Konsequenz der "ublichen Sobolev-Einbettungen.  
\qed\\

Ist $N$ eine geschlossene $C^k$-Mannigfaltigkeit und $M\subset N$, so bestehe der Raum $W^{s,p}(M)$,
$0 \le s \le k$, $1\le p<\infty$ oder $s=0$, $p=\infty$  aus den Funktionen(klassen) $f:
M\rightarrow\setR$, f"ur die in lokalen Koordinaten $(\phi,U)$ von $N$ f"ur alle $\psi\in C^k_0(U)$
gilt \[(f\,\psi)\circ
\phi^{-1}\in W^{s,p}(\phi(U\cap \inn M)).\]
Wir setzen $L^p(M):=W^{0,p}(M)$. Ist $(\phi_k,U_k)_{k}$ ein endlicher Atlas mit
untergeordneter Zerlegung der Eins $(\psi_k)$, so ist die Norm von $W^{s,p}(M)$
durch
\[\norm{f}_{W^{s,p}(M)}:=\sum_k \norm{(f\,\psi_k)\circ
\phi_k^{-1}}_{W^{s,p}(\phi_k(U_k\cap \inn M))}\]
gegeben. Ein $C^k$-Diffeomorphismus, $k\in\setN$, zwischen Gebieten $\Omega$ und
$\widetilde\Omega$ des
$\setR^d$ induziert Isomorphismen zwischen den R"aumen $W^{s,p}(\Omega)$ und
$W^{s,p}(\widetilde\Omega)$, falls $0\le s \le k$ und $1\le p<\infty$. Durch Wahl unterschiedlicher
Atlanten erhalten wir also "aquivalente Normen. Ist $s$
ganzzahlig und $N$ Riemann'sch mit Levi-Civita-Zusammenhang $\nabla$ und
induziertem Ma"s $dV$, so definiert auch die Gr"o"se
\[\sum_{j=0}^s\Big(\int_{\inn M} |\nabla^j v|^p\ dV\Big)^{1/p}\] 
eine "aquivalente Norm. Das ist eine simple Konsequenz der Kompaktheit von $N$. Durch Verwenden
eines endlichen Atlas von $N$ mit
untergeordneter Zerlegung der Eins sieht man leicht, dass
\[W^{s,p}(M)=(L^p(M),W^{k,p}(M))_{\frac{s}{k},p}\]
gilt, wenn $p$, $k$ und $s$ wie zu Beginn des Abschnitts gew"ahlt sind. Desweiteren sei
$W^{s,p}_0(M)$ der Abschluss von $C^k_0(\inn M)$ in $W^{s,p}(M)$. Offenbar gilt
\[W^{s,p}_0(M)\embedding W^{s,p}(N),\]
wenn wir die Funktionen durch $0$ auf $N$ fortsetzen. Eigenschaften wie Vollst"andigkeit,
Reflexivit"at, Dichtheit regul"arer Funktionen, etc. folgen ebenso wie Sobolev-Einbettungen sofort
aus
dem euklidischen Fall, wobei der Rand von $M$ gegebenenfalls hinreichend regul"ar oder leer
vorausgesetzt
werden muss. Aus der $L^2$-Theorie des
Laplace-Operators auf Mannigfaltigkeiten, folgt, dass $\norm{\Delta\, \cdot\,
}_{L^2(M)}$ eine "aquivalente Norm auf $H^2_0(M):=W^{2,2}_0(M)$ definiert, falls $M$ einen
nichtleeren $C^{1,1}$-Rand besitzt; siehe zum Beispiel \cite{b50}.\\

% F"ur
% $1\le p\le \tilde p<\infty$ und $s,\tilde s\ge 0$ mit $s-d/p=\tilde s-d/{\tilde p}$ gilt dann die
% Einbettung
% \[W^{s,p}(\Omega)\embedding W^{\tilde s,\tilde p}(\Omega).\]
% Die resultierende Einbettung
% \[W^{1,p}(N)\embedding L^{r}(N)\]
% f"ur $r<p^*$ ist kompakt. Desweiteren gilt f"ur $d<p<\infty$ und $\alpha=1-d/p$ 
% \[W^{1,p}(N)\embedding C^{0,\alpha}(N).\]

% F"ur eine Einf"uhrung in die Sobolev-R"aume auf Teilmengen des $\setR^d$ sowie die Bochner- und
% Bochner-Sobolev-R"aume sei auf \cite{b4}, \cite{b5} bzw. auf \cite{b6}, \cite{b16}, \cite{b7}
% verwiesen. Dort finden sich auch, falls nicht explizit anders vermerkt, die Beweise der folgenden
% ben"otigten Fakten dieser Theorien.

\begin{theorem}\label{theorem:spur}
Es seien $\Omega\subset\setR^d$, $d\in\setN$, ein beschr"anktes Gebiet mit Lipschitz-Rand und
$1<p<\infty$. Dann besitzt die Abbildung $v\mapsto v|_{\pa\Omega}$, die f"ur $v\in
C^1(\overline\Omega)$ wohldefiniert
ist, eine stetige Fortsetzung von $W^{1,p}(\Omega)$ nach $W^{1-1/p,p}(\pa\Omega)$.
\end{theorem}         
\noindent Ein Beweis findet sich in \cite{b5}.\\

\begin{proposition}\label{theorem:faltung}
Es sei $\omega\in C^\infty_0(\setR^d)$, $d\in\setN$, mit $\int_{\setR^d}\omega\ dx=1$. F"ur
$\epsilon>0$ setzen wir $\omega_\epsilon:=\epsilon^{-d}\omega(\epsilon^{-1}\cdot)$. Falls
$1\le p <\infty$ und $f\in L^p(\setR^d)$, so gilt $\omega_\epsilon\ast f\rightarrow f$
in $L^p(\setR^d)$ und fast "uberall f"ur $\epsilon\rightarrow 0$.
Ist $f$ stetig in Umgebung einer kompakten Menge $K\subset\setR^d$, so gilt 
$\omega_\epsilon\ast f\rightarrow f$ in $L^\infty(K)$ f"ur $\epsilon\rightarrow 0$.
F"ur $d\ge 2$, $1\le p,r <\infty$ und $f\in L^p(\setR,L^r(\setR^{d-1})$ gilt
\[\norm{\omega_\epsilon\ast f}_{L^p(\setR,L^r(\setR^{d-1})}\le
c\,\norm{f}_{L^p(\setR,L^r(\setR^{d-1})}\]
und $\omega_\epsilon\ast f\rightarrow f$ in $L^p(\setR,L^r(\setR^{d-1})$.
\end{proposition}
\noindent Ein Beweis der Aussagen "uber skalarwertige Funktionen findet sich in \cite{b14},
Theorem 1.2.19 und Korollar 2.1.17. F"ur den Beweis der Aussagen "uber vektorwertige Funktionen
verweisen wir auf \cite{b43}.\\

\begin{proposition}\label{theorem:dicht}
 Es seien $\Omega\subset\setR^d$, $d\in\setN$, ein beschr"anktes Gebiet mit
$C^0$-Rand und $1\le p<\infty$. Dann ist
$C^\infty_0(\setR^d)$ dicht in $W^{1,p}(\Omega)$ und in
$E^p(\Omega)$.
\end{proposition}
\noindent {\bf Beweisskizze:}\  Es seien $B\subset\setR^{d-1}$ ein offener Ball, $g:\overline
B\rightarrow\setR$ eine stetige Funktion und $\widetilde\Omega:=\{(x',x_d)\in\setR^d\ |\ x'\in B,\,
g(x')< x_d\}$. Zudem
sei $v\in W^{1,p}(\widetilde\Omega)$ eine Funktion, deren Tr"ager eine beschr"ankte Teilmenge von
$\{(x',x_d)\in\setR^d\ |\ x'\in B,\, g(x')\le x_d\}$ sei und die wir durch $0$ auf $\setR^d$
fortsetzen. Wir setzen $v_t(x):=v(x',x_d+t)$ mit $t>0$, sodass in $\widetilde\Omega$ f"ur
hinreichend kleine $\epsilon>$ die Identit"at $\nabla(\omega_\epsilon\ast
v_t)=\omega_\epsilon\ast\nabla v_t$, $\omega_\epsilon$ wie in Proposition \ref{theorem:faltung},
gilt. Aus Proposition \ref{theorem:faltung} folgt dann f"ur
$\epsilon\rightarrow 0$ die
Konvergenz von $\omega_\epsilon\ast v_t\in C^\infty_0(\setR^d)$ gegen $v_t$ in
$W^{1,p}(\widetilde\Omega)$. Die Konvergenz von $v_t$ gegen $v$ in $W^{1,p}(\widetilde\Omega)$ ist
eine direkte Konsequenz der Stetigkeit der Translation in $L^p(\setR^d)$. Die erste Behauptung 
folgt nun aus Obigem mittels Lokalisierung. Details finden sich in \cite{b5}, \cite{b49}.
F"ur die zweite Behauptung ersetzen wir einfach $\nabla$ durch $\dv$.
\qed\\

\begin{theorem}{\bf (de Rham)}\label{theorem:derham}
Es seien $\Omega\subset\setR^d$, $d\in\setN$, ein beschr"anktes Gebiet mit Lipschitz-Rand und
$f\in (H^1_0(\Omega,\setR^d))'$. Gilt $\langle f,\bphi\rangle =0$ f"ur alle $\bphi\in
C^\infty_0(\Omega)$ mit $\dv\bphi=0$, so existiert genau eine Funktion $p\in L^2(\Omega)$ mit
$\int_\Omega p\ dx=0$ und $f=\nabla p$.
\end{theorem}
\noindent F"ur einen Beweis siehe \cite{b13}.\\

\begin{proposition}{\bf (Korn'sche Ungleichung)}\label{theorem:korn}
Es seien $\Omega\subset\setR^d$, $d\in\setN$, ein beschr"anktes Gebiet mit Lipschitz-Rand und $1\le
p<\infty$. Dann definiert $\norm{\,\cdot\,}_{L^p(\Omega)}+\norm{D\,\cdot\,}_{L^p(\Omega)}$
eine "aquivalente Norm auf $W^{1,p}(\Omega,\setR^d)$.
\end{proposition}
\noindent Diese Aussage wird in \cite{b56} bewiesen.\\

\begin{proposition}\label{theorem:bochneremb}
Es seien $I\subset\setR$ ein offenes, beschr"anktes Intervall und die komplexen Hilbert-R"aume
$H_0$, $H_1$
ein Interpolationspaar.
Dann gilt die Einbettung
\[\Big\{v\in L^2(I,H_0)\ \big|\ \frac{dv}{dt}\in L^2(I,H_1)\Big\}\embedding C(\bar
I,[H_0,H_1]_{\frac{1}{2}}).\]
\end{proposition}
\noindent F"ur einen Beweis siehe \cite{b51}.\\

\begin{proposition}\label{theorem:bochnerinterpol}
Es seien $I\subset\setR$ ein offenes Intervall und die Banach-R"aume $B_0$, $B_1$
ein Interpolationspaar. F"ur $1\le p,r,q \le \infty$, $0<\theta<1$ und
$\frac{1}{s}=\frac{1-\theta}{p}+\frac{\theta}{r}$ gilt die Einbettung
\[L^p(I,B_0)\cap L^r(I,B_1)\embedding L^s(I,(B_0,B_1)_{\theta,q}).\]
\end{proposition}
\beweis Unter Verwendung H"older-Ungleichung mit Exponenten
$\tilde p=\frac{p}{(1-\theta)s}$ und $\tilde p'=\frac{r}{\theta s}$ erhalten wir
\begin{equation*}
 \begin{aligned}
  \int_I\norm{u}_{(B_0,B_1)_{\theta,q}}^s\ dt \le \int_I
\norm{u}_{B_0}^{(1-\theta)s}\norm{u}_{B_1}^{\theta s}\ dt\le
\norm{u}_{L^p(I,B_0)}^{(1-\theta)s}\norm{u}_{L^r(I,B_1)}^{\theta s}.
 \end{aligned}
\end{equation*}
\qed\\

\section{Gew"ohnliche Integro-Differentialgleichungen}
Es seien $d\in\setN$, $\bA\in C([0,\infty)\times\setR^d,\setR^d)$ und $\bB\in
C([0,\infty)^2\times\setR^d,\setR^d)$. Wir suchen L"osungen $\balpha\in
C^1([0,T^*),\setR^d)$ der Gleichung
\begin{equation}\label{eqn:gidgl}
 \begin{aligned}
\dot \balpha(t)=\bA(t,\balpha(t)) + \int_0^t \bB(t,s,\balpha(s))\ ds  
 \end{aligned}
\end{equation}
f"ur $t\in [0,T^*)$. Wir zeigen nun, dass zu jedem $\balpha_0\in\setR^d$ ein
$T^*\in (0,\infty]$ und eine L"osung auf dem Intervall $[0,T^*)$ mit $\balpha(0)=\balpha_0$
existieren derart, dass aus $T^*<\infty$ die Divergenz $\lim_{t\nearrow T^*}|\balpha(t)|=\infty$
folgt. Sind $\bA$ und $\bB$ affin-linear in $\balpha$, so ist $T^*=\infty$. 

Die Konstruktion der zeitlich lokalen L"osung folgt
de facto w"ortlich dem Beweis des Existenzsatzes von Peano. Der Vollst"andigkeit halber skizzieren
wir sie dennoch. Wir k"onnen zun"achst ohne Einschr"ankung annehmen, dass die rechte Seite von
\eqref{eqn:gidgl} unabh"angig von $t$ und $\balpha$ beschr"ankt ist, da wir sonst die beiden
Summanden mit Abschneidefunktionen multiplizieren
k"onnen, die in einer Umgebung von $(t=0$, $\balpha=\balpha_0)$ identisch $1$ sind. Die so gewonnene
L"osung ist f"ur hinreichend kleine Zeiten eine L"osung der urspr"unglichen Gleichung. Wir
definieren N"aherungsl"osungen $\balpha_k$, $k\in\setN$, indem wir $\balpha_k(t):=\balpha_0$
f"ur $t<0$ und 
\[\balpha_k(t):=\balpha_0+\int_0^t\bA(s,\balpha_k(s-1/k)) + \int_0^s
\bB(s,\tau,\balpha_k(\tau-1/k))\ d\tau\ ds\]
f"ur $0\le t\le 1$ setzen. Aus der Beschr"anktheit der rechten Seite von \eqref{eqn:gidgl} folgt
die Beschr"anktheit von $(\dot\balpha_k)$ in $C([0,1],\setR^d)$ unabh"angig von $k$. Der
Satz von Arzela-Ascoli liefert uns einen gleichm"a"sigen Grenzwert $\balpha$ einer Teilfolge von
$(\balpha_k)$. Aufgrund der Absch"atzung $|\balpha_k(t-1/k)-\balpha_k(t)|\le c/k$
konvergiert auch die entsprechende Teilfolge von $(\balpha_k(\cdot-1/k))$ gleichm"a"sig gegen
$\balpha$. Nun k"onnen wir in der Definition von $\balpha_k$ den Grenz"ubergang vollziehen und
erhalten eine L"osung von \eqref{eqn:gidgl} f"ur kleine Zeiten. Jede beschr"ankte L"osung $\balpha$
auf einem Intervall $[0,T)$, $T<\infty$, l"asst sich auf ein Intervall $[0,\widetilde T)$,
$\widetilde T>T$, fortsetzen. Mit Hilfe von \eqref{eqn:gidgl} folgt n"amlich die Beschr"ankheit von
$\dot\balpha$ und mithin die gleichm"a"sige Stetigkeit von $\balpha$ auf $[0,T)$. Somit existiert
der Grenzwert $\balpha(T)=\lim_{t\nearrow T}\balpha(t)$, und wir k"onnen zum Anfangswert
$\balpha(T)$ eine L"osung konstruieren, die, mit der urspr"unglichen L"osung verkn"upft, die
Fortsetzung liefert. Die Behauptung $T^*=\infty$ im affin-linearen Fall folgt mit Hilfe einer
Variation des zum Beweis des Gronwall-Lemmas verwendeten Arguments. Wir wollen uns hier
die Details sparen, weil wir die Langzeitexistenz der L"osungen $\alpha_n^k$ im Beweis von
Proposition \ref{lemma:ent} ebenso gut aus den dortigen Energieabsch"atzungen folgern k"onnen.

\section{Weitere verwendete Fakten}

\begin{proposition}{\bf (Gronwall'sches Lemma)}\label{lemma:gronwall}
Es seien $T>0$, $\alpha\in C([0,T))$ und $c\in\setR$.
Erf"ullt $\phi\in C^1([0,T))$ f"ur alle $t\in [0,T)$ die Ungleichung
\[\phi'(t)\le\alpha(t)+c\,\phi(t),\]
so folgt f"ur alle $t\in [0,T)$
\[\phi(t)\le\phi(0)\,e^{ct}+\int_0^t\alpha(s)\,e^{c(t-s)}\ ds.\]
\end{proposition}
\vspace{0.1cm}
% \begin{theorem}{\bf (Arzela-Ascoli)}\label{theorem:arzasc}
%  Es seien $K$ ein kompakter metrischer Raum und $F$ ein Banachraum. Eine Menge $M\subset C(K,F)$
% ist
% genau dann relativ kompakt, wenn $M$ gleichgradig stetig und f"ur
% alle $y\in K$ die Menge $\{f(y)\ |\ f\in M\}$
% relativ kompakt in $F$ ist.
% \end{theorem}
% \vspace{0.3cm}
\begin{proposition}{\bf (Reynolds'sches Transporttheorem)}\label{theorem:reynolds}
Es seien $\Omega\subset\setR^3$ ein beschr"anktes Gebiet mit $C^1$-Rand, $I\subset\setR$ ein
Intervall und $\Psi\in C^1(I\times\overline\Omega,\setR^3)$ derart, dass
\[\Psi_t:=\Psi(t,\cdot):\overline\Omega\rightarrow \Psi_t(\overline\Omega)\] f"ur alle $t\in I$ ein
Diffeomorphismus ist. Wir setzen $\Omega_t:=\Psi_t(\Omega)$ und
$\bv:=(\pa_t\Psi)\circ\Psi_t^{-1}$. Dann gilt f"ur $\xi\in C^1(\bigcup_{t\in I}\{t\}\times
\overline\Omega_t)$ und $t\in I$ die Identit"at
\[\frac{d}{dt}\int_{\Omega_t} \xi(t,x)\ dx= \int_{\Omega_t} \pa_t
\xi(t,x)\ dx + \int_{\pa\Omega_t} \bv\cdot \bnu_t\ \xi(t,\cdot)\ dA_t.\]
Dabei bezeichnen $dA_t$ das Fl"achenma"s und $\bnu_t$ die "au"sere
Einheitsnormale von $\pa\Omega_t$.
\end{proposition}
\noindent Ein Beweis dieser Aussage findet sich in \cite{b13}.\\

Wir betrachten nun f"ur ein beschr"anktes Gebiet $\Omega\subset\setR^3$ das Stokes-System
\begin{equation*}
 \begin{aligned}
  -\Delta\bu+\nabla\pi&=0&&\text{ in }\Omega,\\
\dv \bu&=0&&\text{ in }\Omega,\\
\bu&=\bg&&\text{ auf }\pa\Omega.
 \end{aligned}
\end{equation*}

\begin{theorem}\label{theorem:stokes}
 Es seien $\Omega\subset\setR^3$ ein beschr"anktes Gebiet mit $C^{\max(2,k)}$-Rand, $k\in\setN$,
und $1<p<\infty$. Desweiteren sei $\bg\in W^{k-1/p,p}(\pa\Omega)$ mit
\[\int_{\pa\Omega}\bg\cdot\bnu\ dA=0,\] 
wobei
$\bnu$ und $dA$ Einheitsnormale bzw. Fl"achenma"s von $\pa\Omega$ bezeichnen. Dann existiert
genau eine (sehr schwache) L"osung $\bu\in W^{k,p}(\Omega)$ des Stokes-Systems, d.h. es gilt
$\dv\bu=0$ und
\[-\int_\Omega\bu\cdot\Delta\bphi\ dx+\int_{\pa\Omega}\bg\cdot\pa_{\bnu}\bphi\ dA=0\]
f"ur alle $\bphi\in C^2(\overline\Omega)$ mit $\dv\bphi=0$ und $\bphi=0$ auf $\pa\Omega$. Die
Zuordnung $\bg\mapsto\bu$ definiert einen stetigen, linearen Operator von $W^{k-1/p,p}(\pa\Omega)$
nach $W^{k,p}(\Omega)$.

Ist $\pa\Omega$ von der Klasse $C^{2,1}$, so gilt die Aussage auch dann, wenn wir
$W^{k-1/p,p}(\pa\Omega)$ durch $L^p(\pa\Omega)$ und $W^{k,p}(\Omega)$ durch $L^p(\Omega)$ ersetzen.
\end{theorem}
\noindent Die erste Aussage wird in \cite{b18} bewiesen. Die zweite Aussage ist eine direkte
Folgerung aus
Theorem 3 in \cite{b19}.\\

Ein Beweis der folgenden Variante des Schauder'schen Fixpunktsatzes f"ur mengenwertige Abbildungen
ist in \cite{b29} nachzulesen. Die endlichdimensionale Version dieses Theorems, der
Kakutani-Fixpunktsatz, wurde von John Nash bei der Beschreibung des Nash-Gleichgewichts verwendet.
\begin{theorem}{\bf (Kakutani-Glicksberg-Fan)}\label{theorem:kakutani}
Es seien $C$ eine konvexe Teilmenge eines normierten Raumes $Z$ und $F: C\rightarrow 2^C$ eine
oberhalbstetige mengenwertige Abbildung, d.h. f"ur jede offene Teilmenge $W\subset C$ sei die Menge
$\{c\in C\ |\ F(c)\subset W\}\subset C$ offen. Zudem sei $F(C)$ in einer kompakten
Teilmenge von $C$ enthalten, und f"ur alle $z\in C$ sei die Menge $F(z)$ nichtleer, konvex und
kompakt. Dann besitzt $F$ einen Fixpunkt, d.h. es existiert ein $c_0\in C$ mit $c_0\in F(c_0)$.
\end{theorem}
Es ist nicht schwierig zu sehen, dass die Forderung der Oberhalbstetigkeit "aquivalent durch die
Bedingung der Graphenabgeschlossenheit ersetzt werden kann. Letztere besagt, dass aus den
Konvergenzen $c_n\rightarrow c$ in $C$ und mit $F(c_n)\ni z_n\rightarrow z$ in $Z$
folgt, dass $z\in F(c)$.\\

Wir bezeichnen mit $M_{sym}$ den Raum der reellen, symmetrischen $3\times 3$-Matrizen.

\begin{proposition}{\bf (Dal Maso-Murat)}\label{lemma:dmm}
Sei $S:M_{sym}\rightarrow M_{sym}$ stetig und strikt monoton, d.h. es gelte $(S(A)-S(B)):(A-B)>0$
f"ur $A,B\in M_{sym}$, $A\not= B$. Desweiteren sei $(A_n)_{n\in\setN}\subset M_{sym}$ eine Folge
mit
\[\lim_{n}\ (S(A_n)-S(A)):(A_n-A)=0\]
f"ur ein $A\in M_{sym}$. Dann gilt $\lim_{n}A_n=A$.
\end{proposition}
\noindent Dies wird in \cite{b15} gezeigt.\\

\begin{proposition}{\bf (Vitali)}\label{lemma:vitali}
Es seien $\Omega\subset\setR^d$, $d\in\setN$, ein beschr"anktes Gebiet und
$(f_n)_{n\in\setN}\subset L^1(\Omega)$.
Zudem konvergiere die Folge $(f_n)$ fast "uberall in $\Omega$. Dann ist
Konvergenz der Folge in $L^1(\Omega)$ "aquivalent zu der Aussage, dass zu
jedem $\epsilon>0$ ein $\delta>0$ existiert derart, dass 
\[\sup_n\int_K |f_n|\ dx<\epsilon\]
f"ur alle messbaren Mengen $K\subset\Omega$ mit $|K|<\delta$ gilt.
\end{proposition}
\noindent F"ur einen Beweis sei auf \cite{b17} verwiesen.\\

\section{Ausgelagertes}

\begin{bemerkung}\label{bem:korn}
Im Kontext von Abschnitt $4$ gilt formal, d.h. unter Vernachl"assigung von Regularit"atsfragen,
\[\int_{\Omega_{\eta(t)}} (\nabla\bu)^T:\nabla\bphi\ dx =0,\]
falls $\bphi$ ein Feld mit $\tr \bphi=b\,\bnu$ f"ur eine skalare Funktion $b$ ist, insbesondere also
falls
$\bphi=\bu$. Partielle
Integration liefert n"amlich aufgrund der Divergenzfreiheit von $\bu$
\[\int_{\Omega_{\eta(t)}} (\nabla\bu)^T:\nabla\bphi\ dx = \int_{\pa\Omega_{\eta(t)}}
(\bphi\cdot\nabla)\bu\cdot\bnu_{\eta(t)}\ dA.\]
Es gen"ugt also zu zeigen, dass auf $\pa\Omega_{\eta(t)}$ gilt
\[((\bnu\circ q)\cdot\nabla)\bu\cdot\bnu_{\eta(t)}=0. \]
Um diese Identit"at zu beweisen, setzen wir $\be_1:=\bnu\circ q$. Zudem w"ahlen wir auf
$\pa\Omega_{\eta(t)}$ zwei linear
unabh"angige, tangentiale Vektorfelder\footnote{Zur Erinnerung: Diese
existieren lokal.} und setzen diese konstant l"angs $\be_1$ fort. Die so konstruierten Felder nennen
wir $\be_2$ und $\be_3$. Per Konstruktion gilt $\Gamma_{1i}^k\,\be_k:=\nabla_{\be_1}\be_i=0$ und
somit
$\Gamma_{1i}^k=0$
f"ur alle $i,k$.
Schreiben wir $\bu=u^i\,\be_i$, so gilt
\[((\bnu\circ
q)\cdot\nabla)\bu=\nabla_{\be_1}\bu=du^i\,\be_1\ \be_i+u^i\,\nabla_{\be_1}\be_i=du^i\,\be_1\
\be_i,\]
also auf $\pa\Omega_{\eta(t)}$
\begin{equation}\label{eqn:zw}
 \begin{aligned}
((\bnu\circ
q)\cdot\nabla)\bu\cdot\bnu_{\eta(t)}=du^1\,\be_1\ \be_1\cdot\bnu_{\eta(t)}.  
 \end{aligned}
\end{equation}
Andererseits gilt
\[0=\dv\bu=du^i\,\be_i+ u^k\,\Gamma_{ik}^i.\]
Die Komponenten $u^2$ und $u^3$ und ihre
tangentialen Ableitungen $du^2\,\be_2$ bzw. $du^3\,\be_3$ verschwinden auf $\pa\Omega_{\eta(t)}$.
Wir
folgern, dass auf $\pa\Omega_{\eta(t)}$
\[0=du^1\,\be_1+u^1\,\Gamma_{1i}^i=du^1\,\be_1\]
gilt, was zusammen mit \eqref{eqn:zw} die Behauptung zeigt.
\end{bemerkung}
\vspace{0.2cm}

\begin{bemerkung}\label{bem:tdelta}
F"ur $\delta\in C^2(\pa\Omega)$ mit $\norm{\delta}_{L^\infty(\pa\Omega)}<\kappa$ und
$\bphi:\Omega\rightarrow\setR^3$ sei $\T_\delta\bphi$ der \emph{Pushforward} von $(\det
d\Psi_\delta)^{-1} \bphi$ unter $\Psi_\delta$, d.h. $\T_\delta\bphi:=\big(d\Psi_\delta\, (\det
d\Psi_\delta)^{-1}\bphi\big)\circ\Psi^{-1}_\delta.$ Die Abbildung $\T_\delta$ mit der Inversen
$\T_\delta^{-1}\bphi:=\big(d\Psi_\delta^{-1}\, (\det
d\Psi_\delta)\,\bphi\big)\circ\Psi_\delta$
definiert offenbar Isomorphismen zwischen den Lebesgue- und Sobolev-R"aumen auf $\Omega$ bzw.
$\Omega_\delta$, solange die Differenzierbarkeitsordnung kleiner oder gleich $1$ bleibt. Zudem
erh"alt die Abbildung Nullrandwerte. 

Der Diffeomorphismus $\Psi_\delta$ ist eine Isometrie von
$\overline\Omega$ mit der Riemann'schen
Metrik $h:=(d\Psi_\delta)^T d\Psi_\delta$ nach $\overline{
\Omega_\delta}$ mit der euklidischen Metrik. Aus der Nat"urlichkeit des Levi-Civita-Zusammenhangs
unter Isometrien sowie den Identit"aten $\sqrt{\det h_{\beta\gamma}}=\det d\Psi_t$ und
\eqref{eqn:div} folgt
\[(\dv \T_\delta\bphi)\circ\Psi_\delta=\dv_h((\det d\Psi_\delta)^{-1}\bphi)=(\det
d\Psi_\delta)^{-1}\dv\bphi.\]
$\T_\delta$ erh"alt also auch die Divergenzfreiheit und definiert somit Isomorphismen zwischen den
entsprechenden Funktionenr"aumen auf $\Omega$ bzw. $\Omega_\delta$.

Ist $\delta\in C^2(I\times\pa\Omega)$ mit $\norm{\delta}_{L^\infty(I\times\pa\Omega)}<\kappa$, so
definiert die Anwendung von $\T_{\delta(t)}$ zu jedem Zeitpunkt $t\in I$ Isomorphismen zwischen
geeigneten Funktionenr"aumen auf $I\times\Omega$ bzw. $\Omega_\delta^I$, solange die
Differenzierbarkeitsordnung wiederum kleiner oder gleich $1$ bleibt.
\end{bemerkung}
\vspace{0.2cm}

\begin{lemma}\label{lemma:mittelwert}
Es sei $\eta\in Y^I$ mit $\norm{\eta}_{L^\infty(I\times M)}<\kappa$. Dann existiert ein linearer
Operator $\M_\eta$ mit
\begin{equation*}
 \begin{aligned}
  \norm{\M_\eta b}_{L^r(I\times M)}&\le c\, \norm{b}_{L^r(I\times M)},\\
  \norm{\M_\eta b}_{C(\bar I,L^r(M)}&\le c\, \norm{b}_{C(\bar I,L^r(M)},\\
  \norm{\M_\eta b}_{L^r(I,H^2_0(M))}&\le c\, \norm{b}_{L^r(I,H^2_0(M))},\\
  \norm{\M_\eta b}_{H^1(I,L^2(M))}&\le c\, \norm{b}_{H^1(I,L^2(M))}\\
 \end{aligned}
\end{equation*}
f"ur $1\le r\le\infty$ und
\[\int_M (\M_\eta b)(t,\cdot)\,\gamma(\eta(t,\cdot))\ dA = 0\]
f"ur fast alle $t\in I$. Die Konstante $c$ h"angt nur von $\Omega$, $\norm{\eta}_{Y^I}$ und
$\tau(\eta)$ ab; sie bleibt beschr"ankt, wenn
$\norm{\eta}_{Y^I}$ und $\tau(\eta)$ beschr"ankt bleiben.
\end{lemma}
\beweis Wir fixieren eine beliebige Funktion $\psi\in C_0^\infty(\inn M)$ mit $\psi\ge
0$, $\psi\not\equiv 0$ und setzen
\begin{equation*}
 \begin{aligned}
  (M_\eta b)(t,\cdot):= b- \psi\, \frac{a(b(t,\cdot),\eta(t,\cdot))}{a(\psi,\eta(t,\cdot))}
 \end{aligned}
\end{equation*}
mit der Abk"urzung 
\[a(b(t,\cdot),\eta(t,\cdot)):=\int_M b\,\gamma(\eta(t,\cdot))\ dA.\]
Der Beweis der behaupteten Eigenschaften ist dann sehr einfach zu f"uhren, wenn wir beachten, dass
die Ungleichung $a(\psi,\eta)\ge c$ mit einer Konstante $c>0$, die nur von $\tau(\eta)$
abh"angt, gilt. Das folgt aber aus Bemerkung \ref{bem:groessernull}, da $\gamma$ eine stetige
Funktion von $\eta$ ist. 
\qed\\

\begin{lemma}\label{lemma:konvergenzen}
F"ur die Folge $(\eta_n)\subset Y^I$ mit $\sup_n\norm{\eta_n}_{L^\infty(I\times
M)}<\alpha<\kappa$ gelten die Konvergenzen \eqref{eqn:schwkonv}$_{(1,2)}$.
\begin{itemize}
 \item[$(1.a)$] Ist $b\in
C(\bar I,L^2(M))$, so konvergiert $(\M_{\eta_n}b)$
gegen $\M_{\eta}b$ in $C(\bar I,L^2(M))$ unabh"angig von $\norm{b}_{C(\bar I,L^2(M))}\le
1$.
\item[$(1.b)$] Konvergiert zus"atzlich $(\pa_t\eta_n)$ in $L^2(I\times M)$ und ist $b\in
H^1(I,L^2(M))\cap L^2(I,H^2_0(M))$, so konvergiert $(\M_{\eta_n}b)$ gegen $\M_\eta b$ in
$H^1(I,L^2(M))\cap L^2(I,H^2_0(M))$.
\item[$(2.a)$] Ist $b\in
C(\bar I,L^2(M))$, so konvergiert $(\F_{\eta_n}\M_{\eta_n}b)$ gegen
$\F_{\eta}\M_{\eta}b$ in $C(\bar I,L^2(B_\alpha))$ unabh"angig von $\norm{b}_{C(\bar I,L^2(M))}\le
1$. 
 \item[$(2.b)$] Unter den Voraussetzungen von $(1.b)$ konvergiert $(\F_{\eta_n}\M_{\eta_n}b)$
gegen $\F_{\eta}\M_\eta b$ in
$H^1(I\times B_\alpha)\cap L^\infty(I,L^4(B_\alpha))$.
\item[$(2.c)$] Konvergiert $(b_n)$ gegen $b$ schwach in $L^2(I\times M)$, so konvergiert
die Folge $(\F_{\eta_n}b_n)$ gegen $\F_\eta b$ schwach in $L^2(I\times B_\alpha)$.
\end{itemize}
\end{lemma}
\beweis Behauptung $(1.a)$ folgt aus
\begin{equation*}
 \begin{aligned}
  \norm{a(b,\eta_n)-a(b,\eta)}_{L^\infty(I)}&=\Bignorm{\int_M b\,
(\gamma(\eta_n)-\gamma(\eta))\ dA}_{L^\infty(I)}\\
&\le c\,\norm{b}_{C(\bar I,L^2(M))}\,\norm{\gamma(\eta_n)-\gamma(\eta)}_{L^\infty(I\times M)},
 \end{aligned}
\end{equation*}
der analogen Absch"atzung mit $\psi$ anstelle von $b$ und der Ungleichung $a(\psi,\eta_n)\ge c>0$,
wenn wir beachten, dass $(\gamma(\eta_n))$ gleichm"a"sig gegen $\gamma(\eta)$ konvergiert.
Behauptung $(1.b)$ ist klar, wenn wir
wissen, dass die Folgen $(a(b,\eta_n))$ und $(a(\psi,\eta_n))$ gegen $a(b,\eta)$ bzw. $a(\psi,\eta)$
in $H^1(I)$ konvergieren. Der Beweis
dieser Konvergenzen ist aber sehr einfach zu f"uhren. 

Kommen wir zu $(2.a)$. Wir schlie"sen aus
\eqref{eqn:fort} und $(1.a)$, dass $(\F_{\eta_n}\M_{\eta_n}b)$ in $C(\bar I,L^2(S_\alpha))$
unabh"angig von $\norm{b}_{C(\bar I,L^2(M))}\le 1$ konvergiert. Ebenso konvergiert die formale Spur
\begin{equation*}
 \begin{aligned}
\exp\Big(\int_{\eta_n(t,q)}^{-\alpha}
\beta(q+\tau\,\bnu\circ q))\ d\tau\Big)\,
(\M_{\eta_n}b)(t,q)\,\bnu\circ q  
 \end{aligned}
\end{equation*}
von $\F_{\eta_n}\M_{\eta_n}b$ auf $I\times\pa(\Omega\setminus\overline{S_\alpha})$ in
$C(\bar I,L^2(\Omega\setminus\overline{S_\alpha}))$ unabh"angig von $\norm{b}_{C(\bar I,L^2(M))}\le
1$. Aus den Stetigkeitseigenschaften des L"osungsoperators des Stokes-Systems folgt die Behauptung;
vgl.
Bemerkung \ref{bem:nspur}. Vollkommen analog erhalten wir $(2.c)$. 

Wir zeigen nun $(2.b)$. Aus \eqref{eqn:fort}, $(1.b)$ und der Einbettung
\[H^1(I,L^2(M))\cap L^2(I,H^2_0(M))\embedding C(\bar I,H^1(M))\embedding L^\infty(I,L^4(M))\]
schlie"sen wir die Konvergenz von
$(\F_{\eta_n}\M_{\eta_n}b)$ in $L^\infty(I,L^4(S_\alpha))$. Um die Konvergenz in
$H^1(I\times S_\alpha)$ einzusehen, folgern wir zun"achst aus Proposition
\ref{theorem:bochnerinterpol} mit $\theta=2/3$ und
Theorem \ref{theorem:emb} die Einbettung
\begin{equation}\label{eqn:lurchi}
 \begin{aligned}
L^\infty(I,L^2(M))\cap L^2(I,H^2_0(M)) \embedding L^{3}(I,H^{4/3}(M))\embedding
L^3(I,L^\infty(M)).
 \end{aligned}
\end{equation}
Die Konvergenz in $L^2(I,H^1(S_\alpha))$ schlie"sen wir nun aus \eqref{eqn:fortnabla}, $(1.b)$,
\eqref{eqn:lurchi} und der
kompakten Einbettung
\[Y^I\compactembedding L^6(I,H^1_0(M)),\]
die aus Proposition \ref{theorem:aubinlions} folgt. Die Konvergenz in
$H^1(I,L^2(S_\alpha))$ ist eine Konsequenz von \eqref{eqn:klopps}, $(1.b)$, \eqref{eqn:lurchi} und
der Konvergenz von
$(\pa_t\eta_n)$ in $L^6(I,L^2(M))$. Letztere folgt aus der Interpolationsungleichung
\[\norm{\pa_t\eta_n-\pa_t\eta}_{L^6(I,L^2(M))}\le
\norm{\pa_t\eta_n-\pa_t\eta}_{L^\infty(I,L^2(M))}^{2/3}\, \norm{\pa_t\eta_n-\pa_t\eta}_{L^2(I\times
M)}^{1/3}.\]
Ebenso zeigen wir die
Konvergenz der Spur in $H^1(I\times\pa(\Omega\setminus\overline{S_\alpha}))\cap
L^\infty(I,L^4(\pa(\Omega\setminus\overline{S_\alpha})))$. Wie zuvor gen"ugt
es nun, die Stetigkeitseigenschaften des L"osungsoperators des Stokes-Systems zu beachten, um die
Behauptung zu erhalten.
\qed\\

Wir verwenden im Folgenden die zeitunabh"angige Variante des Operators $\M$, die auf
naheliegende
Weise definiert ist.
\begin{lemma}\label{lemma:ehrling} F"ur alle $N\in\setN$, $3/2<p\le\infty$ und $\epsilon>0$ 
existiert eine Konstante $c$ derart, dass f"ur alle $\eta,\,\tilde\eta\in H^2_0(M)$ mit
$\norm{\eta}_{H^2_0(M)}+\norm{\tilde\eta}_{H^2_0(M)}+\tau(\eta)+\tau(\tilde\eta)\le N$
und alle $\bv\in W^{1,p}(\Omega_\eta)$, $\tilde\bv\in W^{1,p}(\Omega_{\tilde\eta})$ die Absch"atzung
\begin{equation*}
 \begin{aligned}
&\sup_{\norm{b}_{L^2(M)}\le
1}\bigg(\int_{\Omega_{\eta}} \bv\cdot\F_{\eta}\M_{\eta}b\ dx - \int_{\Omega_{\tilde\eta}}
\tilde\bv\cdot\F_{\tilde\eta}\M_{\tilde\eta}b\ dx\\
&\hspace{2cm} + \int_M \tr\bv\cdot\bnu\ \M_{\eta}b - \trt\tilde\bv\cdot\bnu\ \M_{\tilde\eta}b\
dA\bigg)\\
&\le \epsilon\, \big(\norm{\bv}_{W^{1,p}(\Omega_{\eta})} +
\norm{\tilde\bv}_{W^{1,p}(\Omega_{\tilde\eta})}\big) + c\sup_{\norm{b}_{H^2_0(M)}\le
1}\bigg(\int_{\Omega_{\eta}} \bv\cdot\F_{\eta}\M_{\eta}b\ dx - \int_{\Omega_{\tilde\eta}}
\tilde\bv\cdot\F_{\tilde\eta}\M_{\tilde\eta}b\ dx\\
&\hspace{7.8cm}  + \int_M \tr\bv\cdot\bnu\ \M_{\eta}b -
\trt\tilde\bv\cdot\bnu\ \M_{\tilde\eta}b\ dA\bigg)
 \end{aligned}
\end{equation*}
gilt.

Ebenso existiert f"ur alle $N\in\setN$, $6/5<p,r\le\infty$ und $\epsilon>0$ eine Konstante $c$
derart, dass f"ur alle $\eta,\,\tilde\eta\in
H^2_0(M)$ und $\delta\in C^4(\pa\Omega)$
mit
\[\norm{\eta}_{H^2_0(M)}+\norm{\tilde\eta}_{H^2_0(M)}+
\norm{\delta}_{C^4(\pa\Omega)}+\tau(\eta)+\tau(\tilde\eta)+\tau(\delta)\le N\] 
und alle $\bv\in W^{1,p}(\Omega_\eta)$, $\tilde\bv\in W^{1,p}(\Omega_{\tilde\eta})$ die
Absch"atzung
\begin{equation*}
 \begin{aligned}
&\sup_{\norm{\bphi}_{H(\Omega)}\le
1}\bigg(\int_{\Omega_{\eta}} \bv\cdot\T_\delta\bphi\ dx - \int_{\Omega_{\tilde\eta}}
\tilde\bv\cdot\T_\delta\bphi\ dx\bigg)\\
&\le \epsilon\, \big(\norm{\bv}_{W^{1,p}(\Omega_{\eta})} +
\norm{\tilde\bv}_{W^{1,p}(\Omega_{\tilde\eta})}\big) +
c\sup_{\norm{\bphi}_{W^{1,r}_{0,\dv}(\Omega)}\le
1}\bigg(\int_{\Omega_{\eta}} \bv\cdot\T_\delta\bphi\ dx - \int_{\Omega_{\tilde\eta}}
\tilde\bv\cdot\T_\delta\bphi\ dx\bigg)
 \end{aligned}
\end{equation*}
gilt.
\end{lemma}
\beweis Wir beweisen diese Aussagen vom Typ des Ehrling-Lemmas mittels des "ublichen
Widerspruchsarguments. Wir zeigen zun"achst die erste Behauptung. W"are diese falsch, so g"abe es
ein
$3/2<p\le\infty$, ein $\epsilon>0$, beschr"ankte Folgen
$(\eta_n),\,(\tilde\eta_n)\subset H^2_0(M)$ mit
$\sup_n\big(\tau(\eta_n)+\tau(\tilde\eta_n)\big)<\infty$ sowie Folgen
$(\bv_n),\,(\tilde\bv_n)$ mit
\[\norm{\bv_n}_{W^{1,p}(\Omega_{\eta_n})}
+ \norm{\tilde\bv_n}_{W^{1,p}(\Omega_{\tilde\eta_n})}=1\]
und 
\begin{equation}\label{eqn:wid}
 \begin{aligned}
&\sup_{\norm{b}_{L^2(M)}\le
1}\bigg(\int_{\Omega_{\eta_n}} \bv_n\cdot\F_{\eta_n}\M_{\eta_n}b\ dx - \int_{\Omega_{\tilde\eta_n}}
\tilde\bv_n\cdot\F_{\tilde\eta_n}\M_{\tilde\eta_n}b\ dx\\
&\hspace{2cm} + \int_M \tren\bv_n\cdot\bnu\ \M_{\eta_n}b - \trent\tilde\bv_n\cdot\bnu\
\M_{\tilde\eta_n}b\
dA\bigg)\\
&\hspace{1cm}> \epsilon + n\sup_{\norm{b}_{H^2_0(M)}\le
1}\bigg(\int_{\Omega_{\eta_n}} \bv_n\cdot\F_{\eta_n}\M_{\eta_n}b\ dx - \int_{\Omega_{\tilde\eta_n}}
\tilde\bv_n\cdot\F_{\tilde\eta_n}\M_{\tilde\eta_n}b\ dx\\
&\hspace{4.5cm}  + \int_M \tren\bv_n\cdot\bnu\ \M_{\eta_n}b -
\trent\tilde\bv_n\cdot\bnu\ \M_{\tilde\eta_n}b\ dA\bigg) .
\end{aligned}
\end{equation}
Wegen Korollar \ref{lemma:spur} sind die Folgen $(\tren\bv_n)$, $(\trent\tilde\bv_n)$ in
$W^{1-1/r,r}(M)$ f"ur ein $r>3/2$ und wegen Theorem \ref{theorem:emb} insbesondere in $H^s(M)$ f"ur
ein $s>0$ beschr"ankt. Es existieren also Teilfolgen mit
\begin{equation*}
 \begin{aligned}
\tren\bv_n\cdot\bnu\rightarrow d,\ \trent\tilde\bv_n\cdot\bnu\rightarrow \tilde d&\text{\quad in
}L^2(M),\\
\eta_n\rightarrow\eta,\ \tilde\eta_n\rightarrow\tilde\eta &\text{\quad schwach in }
H^2_0(M) \text{, insbesondere gleichm"a"sig}. 
 \end{aligned}
\end{equation*}
Lemma \ref{lemma:psi} und die "ublichen
Sobolev-Einbettungen zeigen, dass eine Teilfolge von
$(\bw_n:=\bv_n\circ\Psi_{\eta_n})$ gegen eine Funktion $\bw$ in $L^3(\Omega)$
konvergiert. Setzen wir alle beteiligten Funktionen durch $\boldsymbol{0}$ auf $\setR^3$ fort, so
folgt aus der Absch"atzung 
\begin{equation*}
 \begin{aligned}
  \norm{\bv_n-\bw\circ \Psi^{-1}_{\eta}}_{L^2(\setR^3)} \le \norm{(\bw_n-\bw)\circ
\Psi^{-1}_{\eta_n}}_{L^2(\setR^3)} + \norm{\bw\circ\Psi^{-1}_{\eta_n} -\bw\circ
\Psi^{-1}_{\eta}}_{L^2(\setR^3)},
 \end{aligned}
\end{equation*}
Lemma \ref{lemma:psi} und Bemerkung \ref{bem:konv}, dass
die Folge $(\bv_n)$ gegen $\bv:=\bw\circ \Psi^{-1}_{\tilde\eta}$ in $L^2(\setR^3)$ konvergiert.
Ebenso konvergiert die Folge $(\tilde\bv_n)$ gegen ein $\tilde \bv$ in $L^2(\setR^3)$. Aus Lemma
\ref{lemma:konvergenzen} $(1.a)$, $(2.a)$ schlie"sen wir zudem, dass die Folgen
$(\M_{\eta_n}b)$, $(\M_{\tilde\eta_n}b)$ in $L^2(M)$ und die Folgen $(\F_{\eta_n}\M_{\eta_n}b)$,
$(\F_{\tilde\eta_n}\M_{\tilde\eta_n}b)$ in $L^2(\setR^3)$ unabh"angig von $\norm{b}_{L^2(M)}\le 1$
konvergieren. Somit strebt das Supremum auf der rechten Seite
von \eqref{eqn:wid} gegen
\begin{equation*}
 \begin{aligned}
\sup_{\norm{b}_{H^2_0(M)}\le
1}\bigg(\int_{\Omega_{\eta}} \bv\cdot\F_{\eta}\M_{\eta}b\ dx - \int_{\Omega_{\tilde\eta}}
\tilde\bv\cdot\F_{\tilde\eta}\M_{\tilde\eta}b\ dx + \int_M d\ \M_{\eta}b - \tilde d\
\M_{\tilde\eta}b\ dA\bigg).
\end{aligned}
\end{equation*}
Da die linke Seite von \eqref{eqn:wid} beschr"ankt ist, muss dieser Grenzwert jedoch verschwinden.
Aufgrund der Dichtheit von $H^2_0(M)$ in $L^2(M)$ muss dann aber auch der Grenzwert 
\begin{equation*}
 \begin{aligned}
\sup_{\norm{b}_{L^2(M)}\le
1}\bigg(\int_{\Omega_{\eta}} \bv\cdot\F_{\eta}\M_{\eta}b\ dx - \int_{\Omega_{\tilde\eta}}
\tilde\bv\cdot\F_{\tilde\eta}\M_{\tilde\eta}b\ dx + \int_M d\ \M_{\eta}b - \tilde d\
\M_{\tilde\eta}b\ dA\bigg) 
\end{aligned}
\end{equation*}
der linken Seite von \eqref{eqn:wid} identisch $0$ sein; im Widerspruch zu $\epsilon>0$.

Der Beweis der zweiten Behauptung geht vollkommen analog. Wir wollen hier deshalb lediglich zeigen,
dass f"ur beschr"ankte Folgen $(\eta_n)\subset H^2_0(M)$, $(\delta_n)\subset C^4(\pa\Omega)$ und
eine Folge $(\bv_n)$
mit $\sup_n\big(\norm{\bv_n}_{W^{1,p}(\Omega_{\eta_n})}+\tau(\eta_n)+\tau(\delta_n)\big)<\infty$
f"ur eine Teilfolge die Konvergenz
\begin{equation*}
 \begin{aligned}
\int_{\Omega_{\eta_n}} \bv_n\cdot\T_{\delta_n}\bphi\ dx\quad\rightarrow\quad 
\int_{\Omega_{\eta}} \bv\cdot\T_{\delta}\bphi\ dx
 \end{aligned}
\end{equation*}
f"ur $n\rightarrow\infty$ unabh"angig von $\norm{\bphi}_{H(\Omega)}\le 1$ gilt. Wie zuvor finden
wir Teilfolgen mit
\begin{equation*}
 \begin{aligned}
\bv_n\rightarrow \bv&\text{\quad in }L^2(\setR^3),\\
\eta_n\rightarrow\eta &\text{\quad schwach in }
H^2_0(M) \text{, insbesondere gleichm"a"sig},\\
\delta_n\rightarrow\delta &\text{\quad in } C^3(\pa\Omega).
 \end{aligned}
\end{equation*}
Durch Nulladdition erhalten wir die Absch"atzung
\begin{equation*}
 \begin{aligned}
  \big|\int_{\Omega_{\eta_n}}\bv_n\cdot\T_{\delta_n}\bphi\
dx-\int_{\Omega_{\eta}}\bv\cdot\T_{\delta}\bphi\ dx\big|&\le
\norm{\bv_n-\bv}_{L^2(\setR^3)}\,\norm{\T_{\delta_n}\bphi}_{L^2(\setR^3)}\\
&\hspace{0.5cm}+\norm{\bv}_{W^{1,p}(\Omega_\eta)}\,\norm{\T_{\delta_n}\bphi-\T_\delta\bphi}_{
(W^{1,p}(\Omega_\eta))'}.
 \end{aligned}
\end{equation*}
Wir m"ussen also lediglich zeigen, dass die Folge $(\T_{\delta_n}\bphi)$ gegen $\T_{\delta}\bphi$
in $(W^{1,p}(\Omega_\eta))'$ unabh"angig von $\norm{\bphi}_{H(\Omega)}\le 1$ konvergiert. W"are
dies nicht der Fall, so g"abe es aber ein $\epsilon>0$ und eine Folge $(\bphi_n)\subset H(\Omega)$,
die schwach gegen ein $\bphi$ in $H(\Omega)$ konvergiert, und f"ur die
\[\norm{\T_{\delta_n}\bphi_n-\T_{\delta}\bphi_n}_{(W^{1,p}(\Omega_\eta))'}>\epsilon\]
gilt. Das steht aber im Wiederspruch zur Absch"atzung
\begin{equation}\label{ab:fin}
 \begin{aligned}
\norm{\T_{\delta_n}\bphi_n-\T_{\delta}\bphi_n}_{(W^{1,p}(\Omega_\eta))'}\le
\norm{\T_{\delta_n}(\bphi_n-\bphi)}_{(W^{1,p}(\Omega_\eta))'} +
\norm{\T_{\delta_n}\bphi-\T_{\delta}\bphi_n}_{(W^{1,p}(\Omega_\eta))'}.  
 \end{aligned}
\end{equation}
Die Folgen $(\T_{\delta}\bphi_n)$ und $(\T_{\delta_n}\bphi)$ konvergieren n"amlich schwach
bzw. stark in $L^2(\setR^3)$ gegen $\T_\delta\bphi$. Die starke Konvergenz von
$(\T_{\delta_n}\bphi)$ l"asst sich leicht durch Approximation von $\bphi$ durch glatte Funktionen
einsehen; vgl. Bemerkung \ref{bem:konv}. Zudem zeigt die Identit"at
\[\int_{\setR^3} \T_{\delta_n}(\bphi_n-\bphi)\ h\ dx = \int_\Omega
d\Psi_{\delta_n}(\bphi_n-\bphi)\ h\circ\Psi_{\delta_n}\ dx\]
f"ur $h\in L^2(\setR^3)$, dass $(\T_{\delta_n}(\bphi_n-\bphi))$ schwach in $L^2(\setR^3)$ gegen
$\boldsymbol{0}$ konvergiert. Aus der kompakten Einbettung
\[L^2(\Omega_\eta)\compactembedding (W^{1,p}(\Omega_\eta))',\]
siehe Korollar \ref{lemma:sobolev}, folgt nun, dass die rechte Seite von \eqref{ab:fin} f"ur gro"se
$n$
klein wird.
\qed\\

\begin{lemma}\label{lemma:divdichtglm}
F"ur alle $N\in\setN$, $s>0$ und $\epsilon>0$ existiert eine kleine Zahl $\sigma>0$ derart, dass
f"ur alle $\eta\in H^2_0(M)$ mit $\norm{\eta}_{H^2_0(M)}+\tau(\eta)\le N$ und alle $\bphi\in
H(\Omega_\eta)$ mit $\norm{\bphi}_{L^2(\Omega_\eta)}\le
1$ ein $\bpsi\in H(\Omega_\eta)$ existiert mit
$\supp \bpsi\subset\Omega_{\eta-\sigma}$, $\norm{\bpsi}_{L^2(\Omega_\eta)}\le
2$ und $\norm{\bphi-\bpsi}_{(H^{s}(\setR^3))'}<\epsilon$.\footnote{Bei der Konstruktion von
$\Omega_{\eta-\sigma}$ wird \emph{zuerst} $\eta$ durch $0$ auf $\pa\Omega$ fortgesetzt und
\emph{anschlie"send} $\sigma$ subtrahiert. Zudem sei nochmals explizit darauf hingewiesen, dass wir
die
Felder $\bphi,\bpsi$ wie immer durch $\boldsymbol{0}$ auf $\setR^3$ fortsetzen.}
\end{lemma}
\beweis W"are die Behauptung falsch, so g"abe es $s,\,\epsilon>0$, eine
Nullfolge $(\sigma_n)_{n\in\setN}$ positiver Zahlen, eine beschr"ankte
Folge $(\eta_n)\subset H^2_0(M)$ mit $\sup_n\tau(\eta_n)<\infty$ und eine Folge
$(\bphi_n)\subset H(\Omega_{\eta_n})$ mit
$\norm{\bphi_n}_{L^2(\Omega_{\eta_n})}\le 1$, $\norm{\bphi_n-\bpsi}_{(H^{s}(\setR^3))'}\ge
\epsilon$ f"ur alle $\bpsi\in H(\Omega_{\eta_n})$ mit $\norm{\bpsi}_{L^2(\Omega_{\eta_n})}\le
2$ und $\supp \bpsi\subset\Omega_{\eta_n-\sigma_n}$ sowie
\begin{equation*}
 \begin{aligned}
  \eta_n&\rightarrow\eta &&\text{ schwach in }H^2_0(M), \text{insbesondere gleichm"a"sig},\\
  \bphi_n&\rightarrow \bphi &&\text{ schwach in }L^2(\setR^3).
\end{aligned}
\end{equation*}
Aus der kompakten Einbettung
\[L^2(B)\compactembedding  (H^{s}(B))',\]
siehe Theorem \ref{theorem:emb}, f"ur einen geeigneten Ball $B\subset\setR^3$ folgt die
Konvergenz von $(\bphi_n)$ in $(H^{s}(\setR^3))'$ unter Verwendung der in \cite{b4}
konstruierten Fortsetzungsoperatoren. Ist $\psi$ eine $C^1$-Funktion in einer
Umgebung von $\Omega_\eta$, so folgt aus Proposition
\ref{lemma:nspur}
\begin{equation*}
 \begin{aligned}
  \int_{\Omega_\eta} \bphi\cdot\nabla\psi\ dx = \lim_n \int_{\Omega_{\eta_n}}
\bphi_n\cdot\nabla\psi\ dx =0,
 \end{aligned}
\end{equation*}
da jedes Folgenglied verschwindet. Somit gilt $\trnormal \bphi = 0$, sodass gem"a"s Proposition
\ref{lemma:divdicht} ein $\bpsi\in H(\Omega_\eta)$ mit $\supp \bpsi\subset\Omega_\eta$ und
$\norm{\bphi-\bpsi}_{L^2(\Omega_\eta)}<\epsilon/2$ existiert. Es folgt
\begin{equation*}
 \begin{aligned}
  \norm{\bphi_n-\bpsi}_{(H^{s}(\setR^3))'}\le
\norm{\bphi_n-\bphi}_{(H^{s}(\setR^3))'} +
\norm{\bphi-\bpsi}_{L^2(\Omega_{\eta})}< \epsilon,
 \end{aligned}
\end{equation*}
falls $n$ hinreichend gro"s ist. Das ist ein Widerspruch,
sofern $\norm{\bpsi}_{L^2(\Omega_{\eta_n})}\le
2$ gilt. Diese Absch"atzung ist aber eine Konsequenz von
$\norm{\bphi-\bpsi}_{L^2(\Omega_\eta)}<\epsilon/2$,
falls $\epsilon$ hinreichend klein ist. Letzteres kann ohne Einschr"ankung angenommen
werden.
\qed\\

\begin{bemerkung}\label{bem:ausw}
Wir zeigen nun, dass f"ur $\eta\in Y^I$, $\norm{\eta}_{L^\infty(I\times
M)}<\alpha<\kappa$ und $(b,\bphi)\in T^I_\eta$ die Fortsetzung von
$\bphi$ durch
$(b\,\bnu)\circ q$ auf $I\times B_\alpha$ in $H^1(I,L^2(B_\alpha))$ liegt. Wir approximieren
dazu
$\bphi$ durch Funktionen
$(\bphi_k)\subset C_0^\infty(\setR^4)$ in $H^1(\Omega_\eta^I)$ und $\eta$ durch $(\eta_n)\subset
C^4(\bar I\times\pa\Omega)$ derart, dass die Folgen $(\eta_n)$ gegen $\eta$ in
$L^\infty(I\times\pa\Omega)$ und $(\pa_t\eta_n)$ gegen $\pa_t\eta$ in $L^2(I\times\pa\Omega)$
konvergieren; vgl. Definition \eqref{eqn:reg}. Unter Verwendung des Reynolds'schen Transporttheorems
mit $\bv=(\pa_t\eta_n\bnu)\circ
\Phi_{\eta_n(t)}^{-1}$ und
$\boldsymbol{\xi}=\bphi_k\,\psi$, $\psi\in C_0^\infty(I\times\setR^3)$, erhalten wir die Identit"at
\begin{equation*}
 \begin{aligned}
0 = \int_I\frac{d}{dt}\int_{\Omega_{\eta_n(t)}}
\bphi_k\,\psi\ dxdt&=\int_I\int_{\Omega_{\eta_n(t)}} \pa_t\bphi_k\, \psi +
\bphi_k\, \pa_t\psi\ dxdt\\
&\hspace{0.5cm}+ \int_I\int_{\pa\Omega} \tren(\bphi_k\,\psi)\, \pa_t\eta_n\, \gamma(\eta_n)\ dAdt.
% +\int_I\int_{\pa\Omega_{\eta_n(t)}}\bphi_k\psi\
% ((\pa_t\eta_n\bnu)\circ
% \Phi_{\eta_n(t)}^{-1})\cdot\bnu_{\eta_n(t)}\
% dA_{\eta_n(t)}dt.
 \end{aligned}
\end{equation*}
Lassen wir zun"achst $n$ und anschlie"send $k$ gegen unendlich gehen, so ergibt sich
\begin{equation*}
 \begin{aligned}
\int_I\int_{\Omega_{\eta(t)}} \bphi\, \pa_t\psi\ dxdt=-
\int_I\int_{\Omega_{\eta(t)}} \pa_t\bphi\, \psi\ dxdt-\int_I\int_{\pa\Omega}b\,\bnu\,
\tr\psi\ \pa_t\eta\, \gamma(\eta)\ dAdt.
 \end{aligned}
\end{equation*}
Analog zeigen wir
\begin{equation*}
 \begin{aligned}
\int_I\int_{B_\alpha\setminus\overline{\Omega_{\eta(t)}}} (b\,\bnu)\circ q\ \pa_t\psi\ dxdt&=-
\int_I\int_{B_\alpha\setminus\overline{\Omega_{\eta(t)}}} (\pa_tb\ \bnu)\circ q\ \psi\ dxdt\\
&\hspace{0.5cm}+\int_I\int_{\pa\Omega}b\,\bnu\, \tr\psi\ \pa_t\eta\, \gamma(\eta)\ dAdt.
 \end{aligned}
\end{equation*}
Die Addition der letzten beiden Gleichungen zeigt die Behauptung.
\end{bemerkung}
\end{appendix}
\bibliographystyle{alphadin}
\bibliography{biblio}

\begin{thebibliography}{DEGLT01}

% this bibliography is generated by alphadin.bst [8.2] from 2005-12-21

\providecommand{\url}[1]{\texttt{#1}}
\expandafter\ifx\csname urlstyle\endcsname\relax
  \providecommand{\doi}[1]{doi: #1}\else
  \providecommand{\doi}{doi: \begingroup \urlstyle{rm}\Url}\fi

\bibitem[ADLG]{b26}
\textsc{Acosta}, Gabriel ; \textsc{Dur{\'a}n}, Ricardo~G.  ;
  \textsc{L{\'o}pez~Garc{\'i}a}, Fernando:
\newblock Korn inequality and divergence operator: counterexamples and
  optimality of weighted estimates.
\newblock  \url{http:/\/mate.dm.uba.ar/~rduran/papers/adlg.pdf}. --
\newblock Preprint

\bibitem[AF03]{b4}
\textsc{Adams}, Robert~A. ; \textsc{Fournier}, John J.~F.:
\newblock \emph{Pure and Applied Mathematics (Amsterdam)}. Bd. 140:
  {\emph{Sobolev spaces}}.
\newblock Second.
\newblock Elsevier/Academic Press, Amsterdam, 2003. --
\newblock  xiv+305 S. --
\newblock ISBN 0--12--044143--8

\bibitem[Alt06]{b17}
\textsc{Alt}, H.~W.:
\newblock \emph{Lineare Funktionalanalysis}.
\newblock F"unfte, "uberarbeitete Auflage.
\newblock Berlin : Springer-Verlag, 2006. --
\newblock  xiv+431 S. --
\newblock ISBN 3--540--34186--2

\bibitem[Aub63]{b58}
\textsc{Aubin}, Jean-Pierre:
\newblock Un th\'eor\`eme de compacit\'e.
\newblock {In: }\emph{C. R. Acad. Sci. Paris} 256 (1963), S. 5042--5044

\bibitem[B{\"a}r01]{b3}
\textsc{B{\"a}r}, Christian:
\newblock \emph{Elementare {D}ifferentialgeometrie}.
\newblock Berlin : Walter de Gruyter \& Co., 2001 (de Gruyter Lehrbuch. [de
  Gruyter Textbook]). --
\newblock  xii+281 S. --
\newblock ISBN 3--11--015519--2

\bibitem[Bar10]{b31}
\textsc{Barbu}, Viorel:
\newblock \emph{Nonlinear differential equations of monotone types in {B}anach
  spaces}.
\newblock New York : Springer, 2010 (Springer Monographs in Mathematics). --
\newblock  x+272 S.
\newblock \url{http://dx.doi.org/10.1007/978-1-4419-5542-5}.
\newblock \url{http://dx.doi.org/10.1007/978-1-4419-5542-5}. --
\newblock ISBN 978--1--4419--5541--8

\bibitem[BF06]{b13}
\textsc{Boyer}, Franck ; \textsc{Fabrie}, Pierre:
\newblock \emph{Math\'ematiques \& Applications (Berlin) [Mathematics \&
  Applications]}. Bd.~52: {\emph{\'{E}l\'ements d'analyse pour l'\'etude de
  quelques mod\`eles d'\'ecoulements de fluides visqueux incompressibles}}.
\newblock Berlin : Springer-Verlag, 2006. --
\newblock  xii+398 S. --
\newblock ISBN 978--3--540--29818--2; 3--540--29818--5

\bibitem[BL76]{b48}
\textsc{Bergh}, J{\"o}ran ; \textsc{L{\"o}fstr{\"o}m}, J{\"o}rgen:
\newblock \emph{Interpolation spaces. {A}n introduction}.
\newblock Berlin : Springer-Verlag, 1976. --
\newblock  x+207 S. --
\newblock Grundlehren der Mathematischen Wissenschaften, No. 223

\bibitem[Bou05]{b38}
\textsc{Boulakia}, Muriel:
\newblock Existence of weak solutions for an interaction problem between an
  elastic structure and a compressible viscous fluid.
\newblock {In: }\emph{J. Math. Pures Appl. (9)} 84 (2005), Nr. 11, 1515--1554.
\newblock \url{http://dx.doi.org/10.1016/j.matpur.2005.08.004}. --
\newblock DOI 10.1016/j.matpur.2005.08.004. --
\newblock ISSN 0021--7824

\bibitem[BP07]{b33}
\textsc{Bothe}, Dieter ; \textsc{Pr{\"u}ss}, Jan:
\newblock {$L_P$}-theory for a class of non-{N}ewtonian fluids.
\newblock {In: }\emph{SIAM J. Math. Anal.} 39 (2007), Nr. 2, 379--421
  (electronic).
\newblock \url{http://dx.doi.org/10.1137/060663635}. --
\newblock DOI 10.1137/060663635. --
\newblock ISSN 0036--1410

\bibitem[CDEG05]{b27}
\textsc{Chambolle}, Antonin ; \textsc{Desjardins}, Beno{\^{\i}}t ;
  \textsc{Esteban}, Maria~J.  ; \textsc{Grandmont}, C{\'e}line:
\newblock Existence of weak solutions for the unsteady interaction of a viscous
  fluid with an elastic plate.
\newblock {In: }\emph{J. Math. Fluid Mech.} 7 (2005), Nr. 3, 368--404.
\newblock \url{http://dx.doi.org/10.1007/s00021-004-0121-y}. --
\newblock DOI 10.1007/s00021--004--0121--y. --
\newblock ISSN 1422--6928

\bibitem[Cia97]{b44}
\textsc{Ciarlet}, Philippe~G.:
\newblock \emph{Studies in Mathematics and its Applications}. Bd.~27:
  {\emph{Mathematical elasticity. {V}ol. {II}}}.
\newblock Amsterdam : North-Holland Publishing Co., 1997. --
\newblock  lxiv+497 S. --
\newblock ISBN 0--444--82570--3. --
\newblock Theory of plates

\bibitem[Cia00]{b45}
\textsc{Ciarlet}, Philippe~G.:
\newblock \emph{Studies in Mathematics and its Applications}. Bd.~29:
  {\emph{Mathematical elasticity. {V}ol. {III}}}.
\newblock Amsterdam : North-Holland Publishing Co., 2000. --
\newblock  lxii+599 S. --
\newblock ISBN 0--444--82891--5. --
\newblock Theory of shells

\bibitem[Cia05]{b23}
\textsc{Ciarlet}, Philippe~G.:
\newblock \emph{An introduction to differential geometry with applications to
  elasticity}.
\newblock Dordrecht : Springer, 2005. --
\newblock  iv+209 S. --
\newblock ISBN 978--1--4020--4247--8; 1--4020--4247--7. --
\newblock Reprinted from J. Elasticity {{\bf{7}}8/79} (2005), no. 1-3
  [MR2196098]

\bibitem[CS05]{b40}
\textsc{Coutand}, Daniel ; \textsc{Shkoller}, Steve:
\newblock Motion of an elastic solid inside an incompressible viscous fluid.
\newblock {In: }\emph{Arch. Ration. Mech. Anal.} 176 (2005), Nr. 1, 25--102.
\newblock \url{http://dx.doi.org/10.1007/s00205-004-0340-7}. --
\newblock DOI 10.1007/s00205--004--0340--7. --
\newblock ISSN 0003--9527

\bibitem[CS06]{b41}
\textsc{Coutand}, Daniel ; \textsc{Shkoller}, Steve:
\newblock The interaction between quasilinear elastodynamics and the
  {N}avier-{S}tokes equations.
\newblock {In: }\emph{Arch. Ration. Mech. Anal.} 179 (2006), Nr. 3, 303--352.
\newblock \url{http://dx.doi.org/10.1007/s00205-005-0385-2}. --
\newblock DOI 10.1007/s00205--005--0385--2. --
\newblock ISSN 0003--9527

\bibitem[CS10]{b42}
\textsc{Cheng}, C. H.~A. ; \textsc{Shkoller}, Steve:
\newblock The interaction of the 3{D} {N}avier-{S}tokes equations with a moving
  nonlinear {K}oiter elastic shell.
\newblock {In: }\emph{SIAM J. Math. Anal.} 42 (2010), Nr. 3, 1094--1155.
\newblock \url{http://dx.doi.org/10.1137/080741628}. --
\newblock DOI 10.1137/080741628. --
\newblock ISSN 0036--1410

\bibitem[DE99]{b35}
\textsc{Desjardins}, B. ; \textsc{Esteban}, M.~J.:
\newblock Existence of weak solutions for the motion of rigid bodies in a
  viscous fluid.
\newblock {In: }\emph{Arch. Ration. Mech. Anal.} 146 (1999), Nr. 1, 59--71.
\newblock \url{http://dx.doi.org/10.1007/s002050050136}. --
\newblock DOI 10.1007/s002050050136. --
\newblock ISSN 0003--9527

\bibitem[DEGLT01]{b39}
\textsc{Desjardins}, B. ; \textsc{Esteban}, M.~J. ; \textsc{Grandmont}, C.  ;
  \textsc{Le~Tallec}, P.:
\newblock Weak solutions for a fluid-elastic structure interaction model.
\newblock {In: }\emph{Rev. Mat. Complut.} 14 (2001), Nr. 2, S. 523--538. --
\newblock ISSN 1139--1138

\bibitem[DMM98]{b15}
\textsc{Dal~Maso}, Gianni ; \textsc{Murat}, Fran{\c{c}}ois:
\newblock Almost everywhere convergence of gradients of solutions to nonlinear
  elliptic systems.
\newblock {In: }\emph{Nonlinear Anal.} 31 (1998), Nr. 3-4, 405--412.
\newblock \url{http://dx.doi.org/10.1016/S0362-546X(96)00317-3}. --
\newblock DOI 10.1016/S0362--546X(96)00317--3. --
\newblock ISSN 0362--546X

\bibitem[Dob06]{b49}
\textsc{Dobrowolski}, M.:
\newblock \emph{Angewandte Funktionalanalysis: Funktionalanalysis,
  Sobolev-Raume Und Elliptische Differentialgleichungen}.
\newblock Springer, 2006 (Springer-Lehrbuch Masterclass).
\newblock \url{http://books.google.com/books?id=o1X\_Ma-MMvMC}. --
\newblock ISBN 9783540253952

\bibitem[DR05]{b59}
\textsc{Diening}, Lars ; \textsc{R{\r{u}}{\v{z}}i{\v{c}}ka}, Michael:
\newblock Strong solutions for generalized {N}ewtonian fluids.
\newblock {In: }\emph{J. Math. Fluid Mech.} 7 (2005), Nr. 3, 413--450.
\newblock \url{http://dx.doi.org/10.1007/s00021-004-0124-8}. --
\newblock DOI 10.1007/s00021--004--0124--8. --
\newblock ISSN 1422--6928

\bibitem[Dro]{b43}
\textsc{Droniou}, J{\'e}r{\^o}me:
\newblock Int{\'e}gration et Espaces de Sobolev {\`a} Valeurs Vectorielles.
\newblock  \url{http://www-gm3.univ-mrs.fr/polys/gm3-02/gm3-02.pdf}. --
\newblock Polycopi{\'e} de l'Ecole Doctorale de Maths-Info de Marseille

\bibitem[DRW10]{b32}
\textsc{Diening}, Lars ; \textsc{R{\r{u}}{\v{z}}i{\v{c}}ka}, Michael  ;
  \textsc{Wolf}, J{\"o}rg:
\newblock Existence of weak solutions for unsteady motions of generalized
  {N}ewtonian fluids.
\newblock {In: }\emph{Ann. Sc. Norm. Super. Pisa Cl. Sci. (5)} 9 (2010), Nr. 1,
  S. 1--46. --
\newblock ISSN 0391--173X

\bibitem[Eva10]{b6}
\textsc{Evans}, Lawrence~C.:
\newblock \emph{Graduate Studies in Mathematics}. Bd.~19: {\emph{Partial
  differential equations}}.
\newblock Second.
\newblock Providence, RI : American Mathematical Society, 2010. --
\newblock  xxii+749 S. --
\newblock ISBN 978--0--8218--4974--3

\bibitem[FJMM03]{b22}
\textsc{Friesecke}, Gero ; \textsc{James}, Richard~D. ; \textsc{Mora}, Maria~G.
   ; \textsc{M{\"u}ller}, Stefan:
\newblock Derivation of nonlinear bending theory for shells from
  three-dimensional nonlinear elasticity by {G}amma-convergence.
\newblock {In: }\emph{C. R. Math. Acad. Sci. Paris} 336 (2003), Nr. 8,
  697--702.
\newblock \url{http://dx.doi.org/10.1016/S1631-073X(03)00028-1}. --
\newblock DOI 10.1016/S1631--073X(03)00028--1. --
\newblock ISSN 1631--073X

\bibitem[FS69]{b30}
\textsc{Fujita}, Hiroshi ; \textsc{Sauer}, Niko:
\newblock Construction of weak solutions of the {N}avier-{S}tokes equation in a
  noncylindrical domain.
\newblock {In: }\emph{Bull. Amer. Math. Soc.} 75 (1969), S. 465--468. --
\newblock ISSN 0002--9904

\bibitem[Gal94]{b18}
\textsc{Galdi}, Giovanni~P.:
\newblock \emph{Springer Tracts in Natural Philosophy}. Bd.~38: {\emph{An
  introduction to the mathematical theory of the {N}avier-{S}tokes equations.
  {V}ol. {I}}}.
\newblock New York : Springer-Verlag, 1994. --
\newblock  xii+450 S. --
\newblock ISBN 0--387--94172--X. --
\newblock Linearized steady problems

\bibitem[GD03]{b29}
\textsc{Granas}, Andrzej ; \textsc{Dugundji}, James:
\newblock \emph{Fixed point theory}.
\newblock New York : Springer-Verlag, 2003 (Springer Monographs in
  Mathematics). --
\newblock  xvi+690 S. --
\newblock ISBN 0--387--00173--5

\bibitem[GM00]{b36}
\textsc{Grandmont}, C{\'e}line ; \textsc{Maday}, Yvon:
\newblock Existence for an unsteady fluid-structure interaction problem.
\newblock {In: }\emph{M2AN Math. Model. Numer. Anal.} 34 (2000), Nr. 3,
  609--636.
\newblock \url{http://dx.doi.org/10.1051/m2an:2000159}. --
\newblock DOI 10.1051/m2an:2000159. --
\newblock ISSN 0764--583X

\bibitem[Gra08a]{b14}
\textsc{Grafakos}, Loukas:
\newblock \emph{Graduate Texts in Mathematics}. Bd. 249: {\emph{Classical
  {F}ourier analysis}}.
\newblock Second.
\newblock New York : Springer, 2008. --
\newblock  xvi+489 S. --
\newblock ISBN 978--0--387--09431--1

\bibitem[Gra08b]{b28}
\textsc{Grandmont}, C{\'e}line:
\newblock Existence of weak solutions for the unsteady interaction of a viscous
  fluid with an elastic plate.
\newblock {In: }\emph{SIAM J. Math. Anal.} 40 (2008), Nr. 2, 716--737.
\newblock \url{http://dx.doi.org/10.1137/070699196}. --
\newblock DOI 10.1137/070699196. --
\newblock ISSN 0036--1410

\bibitem[GSS05]{b19}
\textsc{Galdi}, G.~P. ; \textsc{Simader}, C.~G.  ; \textsc{Sohr}, H.:
\newblock A class of solutions to stationary {S}tokes and {N}avier-{S}tokes
  equations with boundary data in {$W^{-1/q,q}$}.
\newblock {In: }\emph{Math. Ann.} 331 (2005), Nr. 1, 41--74.
\newblock \url{http://dx.doi.org/10.1007/s00208-004-0573-7}. --
\newblock DOI 10.1007/s00208--004--0573--7. --
\newblock ISSN 0025--5831

\bibitem[KJF77]{b5}
\textsc{Kufner}, Alois ; \textsc{John}, Old{\v{r}}ich  ;
  \textsc{Fu{\v{c}}{\'{\i}}k}, Svatopluk:
\newblock \emph{Function spaces}.
\newblock Leyden : Noordhoff International Publishing, 1977. --
\newblock  xv+454 S. --
\newblock ISBN 90--286--0015--9. --
\newblock Monographs and Textbooks on Mechanics of Solids and Fluids;
  Mechanics: Analysis

\bibitem[Koi60]{b54}
\textsc{Koiter}, W.~T.:
\newblock A consistent first approximation in the general theory of thin
  elastic shells.
\newblock {In: }\emph{Proc. {S}ympos. {T}hin {E}lastic {S}hells ({D}elft,
  1959)}.
\newblock Amsterdam : North-Holland, 1960, S. 12--33

\bibitem[Koi66]{b20}
\textsc{Koiter}, W.~T.:
\newblock On the nonlinear theory of thin elastic shells. {I}, {II}, {III}.
\newblock {In: }\emph{Nederl. Akad. Wetensch. Proc. Ser. B} 69 (1966), S.
  1--17, 18--32, 33--54. --
\newblock ISSN 0023--3366

\bibitem[LDR00]{b21}
\textsc{Le~Dret}, H. ; \textsc{Raoult}, A.:
\newblock The membrane shell model in nonlinear elasticity: a variational
  asymptotic derivation [ {MR}1375820 (97b:73028)].
\newblock {In: }\emph{Mechanics: from theory to computation}.
\newblock New York : Springer, 2000, S. 59--84

\bibitem[Lee97]{b2}
\textsc{Lee}, John~M.:
\newblock \emph{Graduate Texts in Mathematics}. Bd. 176: {\emph{Riemannian
  manifolds}}.
\newblock New York : Springer-Verlag, 1997. --
\newblock  xvi+224 S. --
\newblock ISBN 0--387--98271--X. --
\newblock An introduction to curvature

\bibitem[Lee03]{b1}
\textsc{Lee}, John~M.:
\newblock \emph{Graduate Texts in Mathematics}. Bd. 218: {\emph{Introduction to
  smooth manifolds}}.
\newblock New York : Springer-Verlag, 2003. --
\newblock  xviii+628 S. --
\newblock ISBN 0--387--95495--3

\bibitem[LM72]{b51}
\textsc{Lions}, J.-L. ; \textsc{Magenes}, E.:
\newblock \emph{Non-homogeneous boundary value problems and applications.
  {V}ol. {I}}.
\newblock New York : Springer-Verlag, 1972. --
\newblock  xvi+357 S. --
\newblock Translated from the French by P. Kenneth, Die Grundlehren der
  mathematischen Wissenschaften, Band 181

\bibitem[MRR95]{b34}
\textsc{M{\'a}lek}, J. ; \textsc{Rajagopal}, K.~R.  ;
  \textsc{R{\r{u}}{\v{z}}i{\v{c}}ka}, M.:
\newblock Existence and regularity of solutions and the stability of the rest
  state for fluids with shear dependent viscosity.
\newblock {In: }\emph{Math. Models Methods Appl. Sci.} 5 (1995), Nr. 6,
  789--812.
\newblock \url{http://dx.doi.org/10.1142/S0218202595000449}. --
\newblock DOI 10.1142/S0218202595000449. --
\newblock ISSN 0218--2025

\bibitem[Ne{\v{c}}66]{b56}
\textsc{Ne{\v{c}}as}, J.:
\newblock Sur les normes {\'e}quvalentes dans $W^k_p$ et sur la coerivit{\'e}
  des formes formellement positives.
\newblock Montr{\'e}al : Les Presses de l'Universit{\'e} de Montr{\'e}al, 1966,
  S. 102--128

\bibitem[R{\r{u}}{\v{z}}04]{b7}
\textsc{R{\r{u}}{\v{z}}i{\v{c}}ka}, M.:
\newblock \emph{Nichtlineare Funktionalanalysis: Eine Einf{\"u}hrung}.
\newblock Berlin : Springer, 2004

\bibitem[Tar07]{b46}
\textsc{Tartar}, Luc:
\newblock \emph{Lecture Notes of the Unione Matematica Italiana}. Bd.~3:
  {\emph{An introduction to {S}obolev spaces and interpolation spaces}}.
\newblock Berlin : Springer, 2007. --
\newblock  xxvi+218 S. --
\newblock ISBN 978--3--540--71482--8; 3--540--71482--0

\bibitem[Tay11]{b50}
\textsc{Taylor}, Michael~E.:
\newblock \emph{Applied Mathematical Sciences}. Bd. 115: {\emph{Partial
  differential equations {I}. {B}asic theory}}.
\newblock Second.
\newblock New York : Springer, 2011. --
\newblock  xxii+654 S. --
\newblock ISBN 978--1--4419--7054--1

\bibitem[Tem01]{b57}
\textsc{Temam}, Roger:
\newblock \emph{Navier-{S}tokes equations}.
\newblock AMS Chelsea Publishing, Providence, RI, 2001. --
\newblock  xiv+408 S. --
\newblock ISBN 0--8218--2737--5. --
\newblock Theory and numerical analysis, Reprint of the 1984 edition

\bibitem[Tri78]{b55}
\textsc{Triebel}, Hans:
\newblock \emph{North-Holland Mathematical Library}. Bd.~18:
  {\emph{Interpolation theory, function spaces, differential operators}}.
\newblock Amsterdam : North-Holland Publishing Co., 1978. --
\newblock  528 S. --
\newblock ISBN 0--7204--0710--9

\bibitem[Vei04]{b37}
\textsc{Veiga}, H. Beirao~d.:
\newblock On the Existence of Strong Solutions to a Coupled Fluid-Structure
  Evolution Problem.
\newblock {In: }\emph{Journal of Mathematical Fluid Mechanics} 6 (2004), 21-52.
\newblock \url{http://dx.doi.org/10.1007/s00021-003-0082-5}. --
\newblock ISSN 1422--6928. --
\newblock 10.1007/s00021-003-0082-5

\bibitem[Zei90a]{b52}
\textsc{Zeidler}, Eberhard:
\newblock \emph{Nonlinear functional analysis and its applications. {II}/{A}}.
\newblock New York : Springer-Verlag, 1990. --
\newblock  xviii+467 S. --
\newblock ISBN 0--387--96802--4. --
\newblock Linear monotone operators, Translated from the German by the author
  and Leo F. Boron

\bibitem[Zei90b]{b16}
\textsc{Zeidler}, Eberhard:
\newblock \emph{Nonlinear functional analysis and its applications. {II}/{B}}.
\newblock New York : Springer-Verlag, 1990. --
\newblock  i--xvi and 469--1202 S. --
\newblock ISBN 0--387--97167--X. --
\newblock Nonlinear monotone operators, Translated from the German by the
  author and Leo F. Boron

\end{thebibliography}

%%%%%%%%%%%%%%%%%%%%%%%%%%%%%%%%%%%%%%%%%%%%%%%%%%%%%%%%%%%%%%%%%
\end{document}